\documentclass[11pt]{amsart}
\usepackage[T1]{fontenc}
\usepackage[utf8]{inputenc}
\usepackage{geometry}
\usepackage{amssymb, amsmath, amsthm, mathtools}
\usepackage[mathscr]{euscript}
\usepackage{stmaryrd}
\usepackage{xcolor}
\usepackage[inline]{enumitem}

\usepackage{hyperref}
\hypersetup{
	colorlinks=true, linktocpage=true, hidelinks, hypertexnames=false
}

\newtheorem{theorem}{Theorem}[section]
\newtheorem{lemma}[theorem]{Lemma}
\newtheorem{corollary}[theorem]{Corollary}
\newtheorem{question}[theorem]{Question}
\newtheorem*{theoremA}{Theorem A}
\newtheorem*{theoremB}{Theorem B}
\newtheorem*{theoremC}{Theorem C}
\newtheorem{definition}[theorem]{Definition}
\theoremstyle{definition}
\newtheorem{example}[theorem]{Example}
\newtheorem{remark}[theorem]{Remark}
\numberwithin{equation}{section}
\newcommand{\rea}{\mathbb{R}}
\newcommand{\ste}{\mathbb{V}}
\newcommand{\std}{\mathbb{V}^{\ast}}
\newcommand{\projp}{\mathbb{P}^{+}}
\newcommand{\Ga}{\Gamma}
\newcommand{\Om}{\Omega}
\newcommand{\ka}{\kappa}
\newcommand{\la}{\lambda}
\newcommand{\ep}{\epsilon}
\newcommand{\pr}{\partial}
\newcommand{\lb}{\langle}
\newcommand{\ra}{\rangle}
\newcommand{\eul}{\mathbb{X}}
\newcommand{\cA}{\mathscr{A}}
\newcommand{\cK}{\mathscr{K}}
\newcommand{\cAs}{\mathscr{A}^{\infty}}
\newcommand{\mahler}{\mathsf{s}}
\newcommand{\cP}{\mathcal{P}}
\newcommand{\cS}{\mathcal{S}}
\newcommand{\vol}{\operatorname{vol}}
\newcommand{\MA}{\operatorname{MA}}
\newcommand{\tr}{\operatorname{tr}}
\newcommand{\gb}{\bar{g}}
\newcommand{\cinf}{C^{\infty}}
\newcommand{\be}{\beta}
\newcommand{\si}{\sigma}
\newcommand{\lie}{\mathcal{L}}
\newcommand{\monge}{\mathcal{M}}
\newcommand{\Bp}{\mathcal{B}}
\newcommand{\al}{\alpha}
\newcommand{\Det}{\operatorname{Det}}
\newcommand{\sign}{\operatorname{sgn}}
\providecommand{\en}{[\nabla]}
\providecommand{\adj}{\operatorname{adj}}
\providecommand{\vr}{\delta}
\newcommand{\fD}{\mathbb{D}}
\newcommand{\Vrel}{V_{\mathrm{rel}}}
\newcommand{\Vrelpolar}{V^{\ast}_{\mathrm{rel}}}
\newcommand{\Srel}[1]{\cS^{\mathrm{rel}}_{#1}}
\newcommand{\cVrel}{\mathcal{V}_{\mathrm{rel}}}
\newcommand{\dift}{\mathcal{K}}
\newcommand{\Lop}{\mathsf{L}}
\newcommand{\Efun}{\mathsf{E}}
\newcommand{\nrm}{\mathsf{n}}
\newcommand{\Fset}{\mathsf{F}}
\newcommand{\cV}{\mathcal{V}}

\makeatletter
\renewcommand{\l@section}{\@tocline{1}{3pt}{0pt}{1.8pc}{}}
\renewcommand{\l@subsection}{\@tocline{2}{0pt}{1.8pc}{2.4pc}{}}
\makeatother

\begin{document}
\title[Covolume products and hyperbolic affine spheres]{Covolume products of invariant convex bodies and hyperbolic affine spheres}
\author{Daniel J. F. Fox}
\address{Departamento de Matem\'atica Aplicada\\ ETSAM\\ Universidad Polit\'ecnica de Madrid\\Av. Juan de Herrera 4 \\ 28040 Madrid España}
\email{daniel.fox@upm.es}

\date{26 November 2026}

\begin{abstract}
Let $C$ be an $N$-dimensional properly convex open cone divided by a discrete group $\Gamma$ of unimodular linear transformations, and let $M = \mathbb{P}^{+}(C)/\Gamma$. The hyperbolic affine spheres of Cheng and Yau asymptotic to the boundary of $C$ determine a centroaffine metric $h$ on $M$. The covolume product of a $\Gamma$-convex body, a generalization of a Fuchsian convex body, is the product of its covolume and that of its polar. It is bounded below by $\mathrm{vol}(M,h)^{2}/N^{2}$, with equality exactly for the bodies bounded by these affine spheres. This is a hyperbolic counterpart of the Blaschke--Santaló inequality for Santaló-centered Mahler products, with the inequality reversed. The approach gives a direct variational proof of the existence of invariant affine spheres and a sharp Harnack inequality for invariant concave functions. The minimum characterizes the projectively invariant equiaffine volume $\mathrm{vol}(M,h)$, admits a stability estimate, and satisfies an explicit product formula for suitably divided product cones. The maximum of the covolume product is a new projective invariant, bounded in terms of the minimum and the sharp Harnack constant. For convex projective structures on closed surfaces, the minimum is expressed through the Euler characteristic and the squared $L^{2}$ norm of the cubic form of the affine sphere, recovering known bounds, and attains its least value exactly on the Fuchsian locus.
\end{abstract}

\maketitle

\enlargethispage{2\baselineskip}
\tableofcontents

\section{Introduction}
By the Blaschke--Santaló inequality, the Mahler product of a convex body, the product of the volumes of the body and of its polar with respect to its Santaló point, is bounded from above, and its maximizers are exactly the ellipsoids, which are the elliptic affine spheres. This article describes a hyperbolic counterpart of the Blaschke--Santaló inequality, in which a covolume product defined on the invariant convex bodies in a properly convex cone divided by a discrete group is bounded from below, and its minimizers are exactly the invariant convex bodies whose boundaries are hyperbolic affine spheres. The minimum is $\cV(M)^{2}/N^{2}$, where $\cV(M)$ is the volume of the quotient manifold $M$ in the metric induced on $M$ by the centroaffine metric of the affine sphere, which is a projective invariant of the flat projective structure of $M$. This gives a novel characterization of this invariant as an extremal value, with the affine spheres as the only minimizers among all invariant convex bodies, and it leads to a direct variational proof of the existence of affine spheres foliating a divisible properly convex cone that is independent of the classical proof of Cheng and Yau.

Throughout the article $\ste$ is a real vector space of dimension $N = n+1 \geq 2$ equipped with a fixed parallel volume form $\Psi$, $C \subset \ste$ is an open properly convex cone, and $\Ga \subset SL^{\pm}(\ste)$ is a discrete group of unimodular linear transformations preserving $C$ and acting freely, properly, and cocompactly on $\Om = \projp(C)$, so that $M = \Om/\Ga$ is a closed properly convex flat projective manifold. The cone $C$ is said to be \emph{divisible} by $\Ga$ which is said to \emph{divide} $C$.

A \emph{$C$-convex body} is a nonempty convex set $L \subset C$, closed in $\ste$, such that the Minkowski sum $L + C$ satisfies $L + C \subset L$. In the terminology of Schneider \cite{Schneider-pseudocones}, a $C$-convex body $L$ is a \emph{$C$-pseudo-cone}, and, when $C \setminus L$ has finite volume, $L$ is \emph{$C$-close} and $C \setminus L$ is a \emph{$C$-coconvex set} \cite{Schneider-coconvex}.

A \emph{$\Ga$-convex body} is defined to be a $C$-convex body $L$ such that $\Ga L = L$. The set of $\Ga$-convex bodies in $C$ is written $\cK(C, \Ga)$. As consequences of the group action, the conditions $L + C \subset L$ and closedness in $\ste$ of the definition are in fact automatic; by Lemma \ref{gammaconvexbodylemma}, a $\Ga$-convex body is the same thing as a nonempty convex $\Ga$-invariant subset of $C$ that is relatively closed in $C$ and different from $C$. Thus the definition of $\Ga$-convex body generalizes Fillastre's definition \cite{Fillastre-Fuchsian} of Fuchsian convex bodies, which are the $\Ga$-convex bodies in the case that $C$ is the Minkowski future cone and $\Ga$ is Fuchsian. 
Fillastre developed the Brunn--Minkowski theory of Fuchsian convex bodies and introduced their polars. The present article extends that setting from Fuchsian groups to arbitrary groups dividing a properly convex cone, and studies the product of the covolume of a body with that of its polar, which was not considered in \cite{Fillastre-Fuchsian}.

The subset of the dual vector space $\std$ defined by $C^{\ast} = \{\zeta \in \std: \zeta(x) > 0 \text{ for all } x \in \bar{C}\setminus\{0\}\}$ is the open dual cone of $C$.
The \emph{polar} $L^{\circ} = \{\zeta \in C^{\ast}: \zeta(x) \geq 1 \text{ for all } x \in L\}$ of $L$ is a convex subset of $C^{\ast}$, relatively closed in $C^{\ast}$ and satisfying $L^{\circ} + C^{\ast} \subset L^{\circ}$. For a $\Ga$-convex body $L$, the polar $L^{\circ}$ is invariant under the action of the dual group $\Ga^{\ast} = \{(g^{-1})^{\ast}: g \in \Ga\}$ and is a $\Ga^{\ast}$-convex body (Lemmas \ref{gammaconvexbodylemma} and \ref{polartransformlemma}). Polarity is a bijection $\cK(C, \Ga) \to \cK(C^{\ast}, \Ga^{\ast})$, and, under the identification of the dual of $\std$ with $\ste$, it is involutive, $L^{\circ\circ} = L$ \cite{Schneider-copolarity}.

The \emph{covolume} of a $\Ga$-convex body $L$ is defined to be $|(C \setminus L)/\Ga|$, the volume of a fundamental domain for the action of $\Ga$ on the coconvex set $C \setminus L$. The covolume is finite and positive (Lemmas \ref{conevolumelemma}, \ref{harnacklemma}, and \ref{gammaconvexbodylemma}). The covolume of the polar $L^{\circ}$ is defined with respect to the volume form on $\std$ that gives $1$ when paired with $\Psi$. The \emph{covolume product} $\cP(L)$ of $L$ is the product of the covolumes of $L$ and $L^{\circ}$,
\begin{align}\label{volprodintro}
\cP(L) = |(C\setminus L)/\Ga|\,|(C^{\ast}\setminus L^{\circ})/\Ga^{\ast}|.
\end{align}
The covolume product $\cP(L)$ is unchanged under either rescaling the volume form on $\ste$ or homotheties of $L$. Without a group action the covolume of a $C$-convex body is in general infinite, and no inequality for the product of the covolumes of a $C$-close pseudo-cone and its polar seems to be known. For general $C$-pseudo-cones, Schneider obtains finite substitutes for the covolume, the surface area measure, and the cone-volume measure by introducing weights \cite{Schneider-weightedminkowski, Schneider-conevolume}; here finiteness is obtained instead by passing to the quotient by $\Ga$.

For a convex body $K \subset \ste$ with the origin in its interior, the gauge $\gamma_{K}(x) = \inf\{t > 0: x \in tK\}$ is convex, the gauge of the polar body is the support function $\gamma_{K^{\circ}}(\zeta) = \sup_{x \neq 0}\zeta(x)/\gamma_{K}(x)$, and the Mahler product of $K$ is
\begin{align}\label{mahlerproduct}
\mahler(K) = |K||K^{\circ}| = |\{x \in \ste: \gamma_{K}(x) \leq 1\}|\,|\{\zeta \in \std: \gamma_{K^{\circ}}(\zeta) \leq 1\}|.
\end{align}
For a $\Ga$-convex body $L$, the \emph{cogauge} $\psi_{L}(x) = \sup\{t > 0: x \in tL\}$ is a $\Ga$-invariant concave function on $C$, homogeneous of degree one, satisfying $L = \{x \in C: \psi_{L}(x) \geq 1\}$ and $\pr L = \{x \in C: \psi_{L}(x) = 1\}$, and the cogauge of $L^{\circ}$ equals $\psi_{L}^{\ast}(\zeta) = \inf_{x \in C}\zeta(x)/\psi_{L}(x)$ (Lemma \ref{polartransformlemma}). Since $C \setminus L = \{x \in C: \psi_{L}(x) < 1\}$, $C^{\ast}\setminus L^{\circ} = \{\zeta \in C^{\ast}: \psi_{L}^{\ast}(\zeta) < 1\}$, and the level sets $\{x \in C: \psi_{L}(x) = 1\}$ and $\{\zeta \in C^{\ast}: \psi_{L}^{\ast}(\zeta) = 1\}$ are null sets,
\begin{align}\label{volprodcogauge}
\cP(L) = |\{x \in C: \psi_{L}(x) \leq 1\}/\Ga|\,|\{\zeta \in C^{\ast}: \psi_{L}^{\ast}(\zeta) \leq 1\}/\Ga^{\ast}|.
\end{align}
Thus the supremum defining the polar gauge becomes an infimum, the convex body is replaced by the coconvex region between the origin and $\pr L$, and both factors are computed modulo the group.
The name \emph{covolume product} for $\cP(L)$ is chosen because its definition is parallel to that of the Mahler volume product of a convex body, with the gauge replaced by the cogauge and the volume replaced by covolume. 

As for the Mahler product, the covolume product is linearly invariant. For $g \in GL(\ste)$, the image $gL$ is a $g\Ga g^{-1}$-convex body in $gC$, and its covolume product, computed with respect to $gC$ and $g\Ga g^{-1}$, equals $\cP(L)$, because the two covolumes rescale by the reciprocal factors $|\det g|$ and $|\det g|^{-1}$, and it is invariant under polarity, $\cP(L^{\circ}) = \cP(L)$, because $L^{\circ\circ} = L$. Consequently the minimum and the maximum of $\cP$ over $\cK(C, \Ga)$ are invariants of the flat projective structure of $M$, unchanged when $(C, \Ga)$ is replaced by $(C^{\ast}, \Ga^{\ast})$.

For a $\Ga$-convex body $L$ with smooth, locally strongly convex boundary $\Sigma = \pr L$, the \emph{centroaffine metric} of $\Sigma$ is its second fundamental form with respect to the position vector as transversal field, equivalently the restriction to $T\Sigma$ of $-\fD^{2}\psi_{L}$, where $\psi_{L}$ is the cogauge of $L$ (Section \ref{functionalssection}). Its definition uses only the vector space structure, so it is invariant under $GL(\ste)$ and unchanged when $L$ is replaced by its image under a homothety, and it descends to a Riemannian metric on $M$.

The conormal map $x \mapsto d\psi_{L}(x)$ is a diffeomorphism from $\Sigma$ onto the dual hypersurface $\pr L^{\circ}$ (Lemma \ref{polarmaplemma}) and an isometry between their centroaffine metrics (Lemma \ref{centroaffinemetriclemma}), so polarity preserves the centroaffine metric. The conormal Gauss map $\Sigma \to \projp(\std)$ associating with a point of $\Sigma$ the oriented annihilator of the tangent space to $\Sigma$ at the point is the composition of the conormal map with the projection to $\projp(\std)$. Because $\Sigma$ is locally strongly convex, the conormal Gauss map is a local diffeomorphism, so the pullback via the conormal Gauss map of the flat projective structure on $\projp(\std)$ yields a flat projective structure on $\Sigma$ that, like its centroaffine metric, depends only on $\Sigma$ and the underlying vector space structure \cite{Fox-coupled}.

The equiaffine metric of $\Sigma$, also called the Blaschke metric, is the representative of its second fundamental form determined by the equiaffine normal, whose definition depends instead on the volume form $\Psi$, and $\Sigma$ is a hyperbolic affine sphere centered at the origin exactly when its equiaffine metric is a constant multiple of its centroaffine metric (Lemma \ref{affinespherelemma}). 

By theorems of Cheng and Yau \cite{Cheng-Yau-dirichlet, Cheng-Yau-real, Cheng-Yau-affinespheres}, resolving a conjecture of Calabi \cite{Calabi-completeaffine} (see also \cite{Gigena, Li-calabi, Li-calabi2, Sasaki-hyperbolic} for refinements), $C$ is foliated by hyperbolic affine spheres centered at the vertex of $C$ and asymptotic to $\pr C$, and every hyperbolic affine sphere with these two properties is a leaf of this foliation; see \cite{Loftin-affinespheres} for the statement in this form. (Definition \ref{asymptoticdefinition} makes the meaning of \emph{asymptotic to $\pr C$} precise.) The leaves are the images of any one of them under homotheties preserving $C$, and the affine mean curvature, which is constant on each leaf, takes each negative value on exactly one leaf. Because an element of $\Ga \subset SL^{\pm}(\ste)$ preserves $C$, it permutes the leaves of the foliation, and because it preserves the affine mean curvature, it maps each leaf to itself. That a hyperbolic affine sphere centered at the vertex of $C$ and asymptotic to $\pr C$ is unique up to homothety, and so $\Ga$-invariant, is elementary (Lemma \ref{generaluniquenesslemma}). Let $\Sigma_{0}$ be any one of these affine spheres. Nothing below depends on which leaf is used, because the covolume product and the metric induced on $M$ by the centroaffine metric of a leaf are unchanged by homotheties.

\begin{theoremA}
Let $\ste$ be a real vector space of dimension $N \geq 2$
equipped with a fixed parallel volume form, let $C \subset \ste$ be an open properly convex cone, and let $\Ga \subset SL^{\pm}(\ste)$ be a discrete group of unimodular linear transformations preserving $C$ and acting freely, properly, and cocompactly on $\Om = \projp(C)$, with quotient $M = \Om/\Ga$. Let $\Sigma_{0} \subset C$ be a hyperbolic affine sphere centered at the vertex of $C$ and asymptotic to $\pr C$, 
and let $h$ be the Riemannian metric induced on $M$ by its centroaffine metric. For every $\Ga$-convex body $L \subset C$,
\begin{align}\label{theoremAineq}
\cP(L) \geq \frac{\vol(M, h)^{2}}{N^{2}},
\end{align}
with equality if and only if $\pr L$ is a homothetic image of $\Sigma_{0}$.
\end{theoremA}
The metric $h$ on $M$, which does not depend on the choice of leaf $\Sigma_{0}$, equals the metric induced on $M$ by the equiaffine (Blaschke) metric of the leaf of affine mean curvature $-1$, on which the equiaffine and centroaffine metrics coincide. Its volume $\cV(M) = \vol(M, h)$ is the \emph{equiaffine volume} of $M$, and it depends only on the flat projective structure of $M$. Theorem A gives
\begin{align}\label{extremalintro}
\cV(M) = N\min_{L \in \cK(C, \Ga)}\cP(L)^{1/2},
\end{align}
and the minimum is attained exactly at the bodies bounded by the affine spheres.

An affine sphere $\Sigma_{0}$ as in Theorem A exists by the theorems of Cheng and Yau and, independently of them, by Theorem B below. It is $\Ga$-invariant by Lemma \ref{generaluniquenesslemma}, so its centroaffine metric induces a metric on $M$. A third proof of existence, by means of the theorem of Cheng and Yau on compact affine manifolds, is described in Remark \ref{compactquotientremark}. The proofs of Theorems A and B do not use the theorems of Cheng and Yau. 

When $\Ga$ is Fuchsian, $h$ is a hyperbolic metric, and \eqref{theoremAineq} becomes a lower bound by the square of the hyperbolic volume of $M$ (Corollary \ref{fuchsiancorollary}). In this case, \eqref{theoremAineq} and its equality case also follow from the first reversed Minkowski inequality of \cite[Theorem $5.5$]{Fillastre-Fuchsian}, applied to $L$ and the body bounded by the hyperboloid, together with Jensen's inequality for the covolume of the polar (Remark \ref{fillastreremark}). 
Granted the existence of the affine sphere, Theorem A has a short proof, given in Section \ref{reversesantalosection}, which extends that argument to arbitrary divisible cones, with the convexity of the polar covolume (Lemma \ref{covolumeconvexitylemma}) in place of the reversed Minkowski inequality. That proof gives the inequality \eqref{theoremAineq}, the characterization of equality among all $\Ga$-convex bodies, including those with nonsmooth boundary, and the uniqueness of the affine sphere up to homothety (Theorem \ref{stabilitytheorem}). It does not show that an affine sphere, or a minimizer of $\cP$, exists; that is the content of Theorem B, or of the theorems of Cheng and Yau.

The manifolds to which Theorem A applies include the quotients of symmetric cones, in particular all closed hyperbolic manifolds, and their deformations by bending \cite{Cooper-Long-Thistlethwaite, Johnson-Millson}; convex projective structures on surfaces of negative Euler characteristic \cite{Choi-Goldman, Goldman-convex}; structures on Gromov--Thurston manifolds constructed by Kapovich \cite{Kapovich-GT}; examples obtained from Coxeter groups \cite{Benoist-divisibles3, Choi-Lee-Marquis, Choi-Lee-Marquis-truncation, Vinberg-reflections}; three-manifolds of Benoist \cite{Benoist-divisibles4} and of Ballas, Danciger, and Lee \cite{Ballas-Danciger-Lee}, whose universal covers $\Om$ have non-strictly convex boundary; and, in every dimension at least $3$, examples of Blayac and Viaggi \cite{Blayac-Viaggi}, which are divided by Zariski dense groups and for which $\pr\Om$ is likewise not strictly convex. See the surveys \cite{Benoist-survey, Choi-Lee-Marquis-survey, Marquis-survey} and \cite[Section $1.1$]{Blayac-Viaggi}. Note that nothing in Theorem A, or in Theorems B and C later in the introduction, needs the convexity of $\pr\Om$ to be strict.

By the Blaschke--Santaló inequality, for a convex body $K$ whose Santal\'o point is the origin, in particular for $K$ symmetric about the origin, the Mahler product satisfies $\mahler(K) \leq \mahler(B)$, where $B$ is the Euclidean unit ball \cite{Santalo-affineinvariant}, with equality if and only if $K$ is a centered ellipsoid \cite{Lutwak-BS, Meyer-Pajor, Petty-affineisoperimetric}. For the functional form of the Blaschke--Santal\'o inequality, see \cite{AKM-santalo}. The inequality of Theorem A is a hyperbolic counterpart of the Blaschke--Santaló inequality, in which covolumes of $\Ga$-convex bodies replace volumes of convex bodies, and for this reason it is called here the \emph{Blaschke--Santal\'o inequality for the covolume product}. Its direction is opposite to that of the Blaschke--Santal\'o inequality, in the same way that the Brunn--Minkowski and Aleksandrov--Fenchel inequalities have opposite directions for coconvex bodies \cite{Khovanskii-Timorin, Schneider-coconvex} and for Fuchsian convex bodies \cite{Fillastre-Fuchsian}. 

This change of direction is of the kind that accompanies the passage from the elliptic to the hyperbolic setting, and that pervades the passage from Euclidean to affine notions, as when the locally strongly convex surfaces critical for the affine area are maximal rather than minimal, the second variation of the affine area at such a surface being negative definite, by a theorem of Calabi \cite[Theorem $1.3$]{Calabi-maximal}. It is not a reversal in the sense in which an inequality that bounds a quantity on one side is said to be reversed by a bound for the same quantity on the other side. In that sense, the name \emph{reverse Santal\'o inequality} is commonly used for the lower bound $\mahler(K) \geq c^{N}\mahler(B)$ of Bourgain and Milman \cite{Bourgain-Milman}, valid for $K$ symmetric about the origin and an absolute constant $c > 0$, and its counterpart for the covolume product is the upper bound \eqref{maxminintro} discussed below, not Theorem A. 

For a general properly convex cone, divisible or not, the existence of the foliation by hyperbolic affine spheres needed in Theorem A follows from theorems of Cheng and Yau. For a cone divided by a cocompact group, Theorem B gives an alternative, variational existence proof that does not use those theorems. What is new in Theorem B is that the extremal property is strong enough to produce the affine sphere: minimizing $\cP$ over all $\Ga$-convex bodies yields a body whose polar solves the affine sphere equation weakly, and regularity theory shows that its boundary is the affine sphere. 

\begin{theoremB}
Let $C$ and $\Ga$ be as in Theorem A. On the set of $\Ga$-convex bodies in $C$ the covolume product $\cP$ attains its minimum. Every minimizer has smooth, locally strongly convex boundary, meaning that the second fundamental form of its boundary is positive definite, and its boundary is a $\Ga$-invariant hyperbolic affine sphere centered at the origin and asymptotic to $\pr C$. Any two minimizers are homothetic.
\end{theoremB}

Cheng and Yau remark that the real Monge--Amp\`ere equation does not admit a natural variational characterization \cite[Section $0$]{Cheng-Yau-real}. For the affine sphere equation on a divisible cone, Theorem B provides one. Obtaining a solution of a Monge--Amp\`ere equation as the minimizer of a functional of convex functions is also the strategy of Berman and Berndtsson \cite{Berman-Berndtsson} and of Cordero-Erausquin and Klartag \cite{Cordero-Klartag} on $\rea^{n}$. Here the functional is the covolume product $\cP$. In the analogy with K\"ahler geometry described in Section \ref{kahlereinsteinsection}, $\cP$ corresponds formally, at the level of normalized first variations, to the energy functional whose minimizers are the K\"ahler--Einstein metrics, so the extremal characterization \eqref{extremalintro} is the counterpart of a familiar one. The characterization of the Cheng--Yau solution as the supremum of subsolutions in \cite[Theorem $4.1$]{Fox-schwarz} is extremal in a different sense, resting on a comparison principle rather than on a functional, and does not give an independent proof of existence.

The proof of Theorem B has three steps. A minimizing sequence, normalized to remove the scale invariance of $\cP$, has a locally uniformly convergent subsequence, because of a Harnack inequality valid for all positive concave $\Ga$-invariant functions of degree one (Lemma \ref{harnacklemma}), whose sharp form is described below. By the variational lemma of Aleksandrov, a minimizer satisfies the Euler--Lagrange equation of $\cP$ in a weak sense; this is the centroaffine Minkowski equation \cite{Chou-Wang} for the polar set, written as a Monge--Amp\`ere equation on the dual manifold $M^{\ast} = \projp(C^{\ast})/\Ga^{\ast}$ (Lemma \ref{weakeulerlagrangeproposition}). The regularity theory of Caffarelli \cite{Caffarelli-strictconvexity, Caffarelli-interior} then applies, because the Harnack inequality forces the solution to vanish at the boundary, which excludes degenerate contact sets (Lemma \ref{regularityproposition}). Theorem A follows by combining Theorem B with the elementary uniqueness of affine spheres (Lemma \ref{uniquenessproposition}).

Theorem C, which follows from Theorem \ref{secondvariationtheorem} proved in Section \ref{variationsection}, calculates the second variation of $\cP$ at a $\Ga$-convex body bounded by an affine sphere, in the parametrization of $\Ga$-convex bodies by their cogauges, and shows that the minimum in Theorem A is nondegenerate in directions transverse to homotheties. It is the infinitesimal form of the extremal characterization \eqref{extremalintro}; a global counterpart, valid for all $\Ga$-convex bodies, is the stability estimate of Theorem \ref{stabilitytheorem}.

\begin{theoremC}
Let $C$, $\Ga$, $M$, $\Sigma_{0}$, and $h$ be as in Theorem A, let $\psi_{0}$ be the cogauge of the $\Ga$-convex body $L_{0}$ bounded by $\Sigma_{0}$, and let $V = \vol(M, h)$. For a smooth function $f$ on $M$ with $h$-mean $\bar{f}$, and $\psi_{t} = (1 + tf)\psi_{0}$,
\begin{align}\label{theoremCformula}
\tfrac{d^{2}}{dt^{2}}_{|t = 0}\cP(\{x \in C: \psi_{t}(x) \geq 1\}) = \frac{V}{N}\int_{M}\left(|df|^{2}_{h} + 2N(f - \bar{f})^{2}\right)d\vol_{h}.
\end{align}
In particular the second variation of $\cP$ at $L_{0}$ is a positive semidefinite quadratic form whose radical consists of the constant functions, which correspond to homotheties preserving $C$.
\end{theoremC}

For the Mahler product of a centered ellipsoid in $\ste$, the corresponding second variation in the direction of a function $f$ on $S^{n}$ is a positive multiple of $\int_{S^{n}}(-|df|^{2} + 2N(f - \bar{f})^{2})\,d\vol$, which is positive on the spherical harmonics of degree $1$, which describe translations of the body, zero on those of degree $2$, which describe its images under linear maps, and negative on those of higher degree (Remark \ref{ellipticremark}). The ellipsoid is thus a minimum in the translation directions and a maximum in the directions transverse to the translations and the linear maps, whereas by Theorem C all directions transverse to homotheties are directions of increase of $\cP$.

For convex bodies $K_{i} \subset \rea^{m_{i}}$, $i = 1, 2$, containing the origin in their interiors, the Mahler product satisfies the rescaled multiplicativity
\begin{align}\label{mahlerproductintro}
(m_{1} + m_{2})!\,\mahler(K_{1}\times K_{2}) = m_{1}!\,m_{2}!\,\mahler(K_{1})\mahler(K_{2}).
\end{align}
The covolume product is multiplicative in a similar sense under products of cones. If the $N_{i}$-dimensional convex cones $C_{i} \subset \ste_{i}$ are divided by $\Ga_{i}$ for $i = 1, 2$, then, for every $\tau > 0$, the cone $C_{1}\times C_{2}$ is divided by $\Ga = \Ga_{1}\times\Ga_{2}\times\lb g_{\tau}\ra$, where $g_{\tau}$ acts on $\ste_{1}$ by $e^{N_{2}\tau}$ and on $\ste_{2}$ by $e^{-N_{1}\tau}$ (Lemma \ref{productdivisiblelemma}). The quotient is homeomorphic to $M_{1}\times M_{2}\times S^{1}$, and $\tau$ is proportional to the length of the circle factor. The product $L_{1}\boxtimes L_{2}$ of $\Ga_{i}$-convex bodies is the $\Ga$-convex body whose cogauge is the weighted geometric mean $\psi_{1}^{N_{1}/N}\psi_{2}^{N_{2}/N}$ of their cogauges. By Theorem \ref{producttheorem},
\begin{align}\label{productintro}
N^{N}\cP(L_{1}\boxtimes L_{2}) = \tau^{2} N_{1}^{N_{1}+2}N_{2}^{N_{2}+2}\,\cP(L_{1})\,\cP(L_{2}).
\end{align}
The constants are explained after Theorem \ref{producttheorem}; the factor $\tau^{2}$ has no counterpart for convex bodies.

The minimizers of $\cP$ for a product of divisible cones are exactly the products of minimizers for the factors, so the minimal covolume product factors as in \eqref{productintro} (Corollary \ref{productinvariantcorollary}). This reformulates the Calabi product of hyperbolic affine spheres \cite{Calabi-completeaffine} (see also \cite[Section $3$]{Dillen-Vrancken}) in terms of cogauges. 

The extremal problem for the covolume product has two ends. The covolume product also has a maximum $\mathfrak{M}(M)$ on $\cK(C, \Ga)$, again an invariant of the flat projective structure of $M$ and of that of its dual. Unlike the minimum, it is not expressed by any previously known quantity, and it is not attained at any body with smooth, locally strongly convex boundary (Corollary \ref{nosmoothmaximizercorollary}), because a smooth maximizer would be critical, so an affine sphere by Lemma \ref{firstvariationproposition}, and hence a minimizer.

The maximum is bounded in terms of the minimum by means of a sharp form of the Harnack inequality used to prove Theorem B. Let $\psi_{0}$ be the cogauge of the $\Ga$-convex body bounded by the affine sphere $\Sigma_{0}$ of Theorem A, which is determined by $C$ and $\Ga$ up to a positive constant factor. For a positive concave $\Ga$-invariant function $\psi$ on $C$, homogeneous of degree one, the ratio $\psi/\psi_{0}$ descends to the compact manifold $M$, so the ratio $(\inf_{C}\psi/\psi_{0})/(\sup_{C}\psi/\psi_{0})$ is defined, and it follows from Lemma \ref{harnacklemma} that it is bounded from below by a positive constant depending only on $C$ and $\Ga$. The infimum $\mathfrak{c}(M)$ of this ratio over all such $\psi$ is the \emph{sharp Harnack constant} of $M$ (Definition \ref{maxmindefinition}). It depends only on the flat projective structure of $M$, and $\mathfrak{c}(M^{\ast}) = \mathfrak{c}(M)$, because polarity maps $\psi_{0}$ to the cogauge of the affine sphere of $C^{\ast}$ (Corollary \ref{maximumcorollary}). To compute it, it suffices to take the infimum over the simplest $\Ga$-convex bodies, those whose vertices form a single $\Ga$-orbit. Precisely, for every such $\psi$ there is $x \in C$ such that
\begin{align}\label{sharpharnackintro}
\frac{\inf_{C}\psi/\psi_{0}}{\sup_{C}\psi/\psi_{0}} \geq \frac{\inf_{C}\psi_{x}/\psi_{0}}{\sup_{C}\psi_{x}/\psi_{0}} = \frac{\psi_{0}(x)}{\sup_{\pr L_{x}}\psi_{0}},
\end{align}
where $L_{x}$ is the closure in $C$ of $\operatorname{conv}(\Ga x) + C$, with $\operatorname{conv}$ denoting the convex hull, and $\psi_{x}$ is the cogauge of $L_{x}$ (Theorem \ref{sharpharnacktheorem}, which treats the ratio relative to an arbitrary positive concave $\Ga$-invariant function homogeneous of degree one in place of $\psi_{0}$).

By Corollary \ref{maximumcorollary}, $\cP$ attains its supremum $\mathfrak{M}(M)$ on $\cK(C, \Ga)$, this supremum satisfies $\mathfrak{M}(M^{\ast}) = \mathfrak{M}(M)$, and 
\begin{align}\label{maxminintro}
\mathfrak{m}(M) < \mathfrak{M}(M) \leq \mathfrak{c}(M)^{-N}\mathfrak{m}(M),
\end{align}
where $\mathfrak{m}(M) = \min_{\cK(C, \Ga)}\cP = \cV(M)^{2}/N^{2}$.
The use of $\psi_{0}$ as the reference function makes $\mathfrak{c}(M)$ invariant under duality and gives the exponent $N$ in \eqref{maxminintro}. The same constant governs the stability of the minimum (Theorem \ref{stabilitytheorem}).

The upper bound in \eqref{maxminintro} is an analogue of the Bourgain--Milman inequality $\mahler(K) \geq c^{N}\mahler(B)$ \cite{Bourgain-Milman} mentioned above, which, with the Blaschke--Santal\'o inequality, bounds the ratio of the maximum to the minimum of the Mahler product of bodies symmetric about the origin by $c^{-N}$. The bound \eqref{maxminintro} has the same exponential dependence on $N$, but its base $\mathfrak{c}(M)^{-1}$ depends on the geometry of $M$, and, unlike the constant $c$ of Bourgain and Milman, it cannot be replaced by an absolute constant. Indeed, already for $N = 2$ the ratio $\mathfrak{M}(M)/\mathfrak{m}(M)$ is unbounded. For a simplicial cone, $\mathfrak{c}(M)$ is given by the explicit formula \eqref{simplicialharnack}. For $N = 2$, with $\Ga$ generated by an element with eigenvalues $e^{s}$ and $e^{-s}$, 
this formula gives $\mathfrak{c}(M) = 1/\cosh(s/2)$, and the body $L_{x}$, bounded by the polygonal line through the orbit $\Ga x$, gives the lower bound in
\begin{align}\label{noneintro}
\frac{\sinh^{2}(s/2)}{(s/2)^{2}} \leq \frac{\mathfrak{M}(M)}{\mathfrak{m}(M)} \leq \cosh^{2}(s/2)
\end{align}
(Example \ref{simplicialharnackexample}), which tends to infinity with $s$. The minimizers have smooth boundary, whereas no maximizer has smooth, locally strongly convex boundary (Corollary \ref{nosmoothmaximizercorollary}); for the Mahler product it is the maximizers that are smooth, namely the ellipsoids, and the minimizers are conjecturally polytopes.

This leads to Question \ref{maximizerquestion}: is there always a maximizer of $\cP$ whose vertices form a single orbit, such as some $L_{x}$, or equivalently, by polarity, one whose facets do? The behavior of $\mathfrak{M}/\mathfrak{m}$ under products of cones gives some evidence for an affirmative answer, described after Question \ref{maximizerquestion}. Already for $N = 2$, whether $L_{x}$ is a maximizer is not known.

Section \ref{reversesantalosection} gives a second proof of Theorem A, independent of the regularity theory used to prove Theorem B, first for bodies with smooth, locally strongly convex boundary and then, by approximation (Lemma \ref{densitylemma}), for all $\Ga$-convex bodies. If $\psi_{0}$ is the cogauge of the affine sphere and $\psi$ is the cogauge of a $\Ga$-convex body $L$, then
\begin{align}\label{linearintro}
\frac{|K_{\psi}/\Ga|}{|K_{\psi_{0}}/\Ga|} + \frac{|K^{\ast}_{\psi}/\Ga^{\ast}|}{|K^{\ast}_{\psi_{0}}/\Ga^{\ast}|} \geq 2,
\end{align}
where $K_{\psi} = \{x \in C: \psi(x) \leq 1\}$ and $K^{\ast}_{\psi} = \{\zeta \in C^{\ast}: \psi^{\ast}(\zeta) \leq 1\}$, so that, by \eqref{volprodcogauge}, the numerators are the covolumes of $L$ and $L^{\circ}$ (Theorem \ref{suminequalityproposition}). Replacing $\psi_{0}$ by $t\psi_{0}$, whose level set is again an affine sphere, and minimizing over $t > 0$ gives \eqref{theoremAineq}, while, conversely, \eqref{theoremAineq} implies \eqref{linearintro} for every $\Ga$-convex body by the inequality between arithmetic and geometric means. Retaining the remainder in the convexity of $x \mapsto x^{-N}$ in this proof gives the equality case of \eqref{theoremAineq} for all $\Ga$-convex bodies and a quantitative form of Theorem A: the difference $\cP(L)^{1/2} - \vol(M, h)/N$ is bounded below by a positive multiple, depending only on $N$ and the sharp Harnack constant $\mathfrak{c}(M)$, of the squared $L^{2}$ distance, with respect to $h$, between the cogauge of $L$, suitably rescaled, and that of the body bounded by the affine sphere (Theorem \ref{stabilitytheorem}).

By Corollary \ref{weightedcorollary}, for arbitrary constants $a, b > 0$ the weighted sum $a|K_{\psi}/\Ga| + b|K^{\ast}_{\psi}/\Ga^{\ast}|$ of the two covolumes in \eqref{linearintro} is bounded below, with equality exactly for the affine sphere of the foliation determined by the ratio $a/b$. This is the hyperbolic counterpart of the construction of affine hemispheres by Klartag \cite{Klartag-hemispheres}, in which a hemisphere is obtained as the minimizer of a sum of a volume term and a polar volume term (Remark \ref{klartagremark}). 

The inequality \eqref{linearintro} is a consequence of the strict convexity of the polar covolume as a function of the cogauge $\psi$, which can be seen as a Brunn--Minkowski inequality for coconvex sets (Lemma \ref{covolumeconvexitylemma}). This generalizes the convexity of the covolume of Fuchsian convex bodies proved by Fillastre \cite{Fillastre-Fuchsian}. The same argument applies to a properly convex cone that is not divided by any group, although the covolumes are then infinite. For a hypersurface that coincides with an affine sphere of the foliation of $C$ outside the cone over a compact subset of $\projp(C)$, the differences between its covolumes and those of the affine sphere are finite, and among such hypersurfaces the affine sphere is the unique minimizer of a weighted sum of these relative covolumes (Theorem \ref{relativetheorem}).

The two covolumes in \eqref{volprodintro}, multiplied by $N$, are the members $\cS_{0}$ and $\cS_{1}$ of a family of functionals depending on a parameter $\be \in \rea$, defined for $\Ga$-convex bodies $L$ with smooth, locally strongly convex boundary by
\begin{align}\label{sbetaintro}
\cS_{\be}(L) = \int_{\pr L/\Ga}\lb x, \nrm_{L}(x)\ra^{1-(n+2)\be}|\varkappa_{L}(x)|^{\be}\,d\si(x), \qquad \be \in \rea,
\end{align}
where $\nrm_{L}$ is the unit normal pointing into $L$, so that $\lb x, \nrm_{L}(x)\ra$ is the support function of $\pr L$ at $x$ in the inward-normal convention, in which it is positive, $\varkappa_{L}$ is the Gauss--Kronecker curvature of $\pr L$, and $d\si$ is the induced hypersurface measure, all with respect to an inner product on $\ste$ whose volume form is $\pm\Psi$. The integrand does not depend on the choice of such an inner product (Lemma \ref{centroaffineproposition}).

The functional $\cS_{\be}$ is an equivariant version of the $L_{p}$ affine surface area of Lutwak \cite{Lutwak-Lp}, with $p = N\be/(1-\be)$ for $\be \neq 1$, while $\be = 1$ corresponds to $p = \pm\infty$. Polarity exchanges $\cS_{\be}$ and $\cS_{1-\be}$ (Lemma \ref{dualityproposition}). The functional $\cS_{1/2}(L) = \vol(M, \gb_{L})$ is the volume of the centroaffine metric $\gb_{L}$ of $\pr L$, and, by the Cauchy--Schwarz inequality, $N^{2}\cP(L) = \cS_{1}(L)\cS_{0}(L)\geq \vol(M, \gb_{L})^{2}$, with equality exactly for affine spheres (Lemma \ref{cauchyschwarzproposition}). 
This differs from Theorem A, in which the metric is that of the affine sphere and not that of $\pr L$, and comparing the two leads to the question of whether the affine sphere maximizes the centroaffine volume (Question \ref{centroaffineconjecture}), which is answered affirmatively here for $n \leq 2$ (Theorem \ref{lowdimcentroaffinetheorem}). An affirmative answer would characterize $\cV(M)$ also as the maximum of $\vol(M, \gb_{L})$ over the $\Ga$-convex bodies with smooth, locally strongly convex boundary, so that the same invariant would be the minimum of one functional and the maximum of another, both attained exactly at the affine spheres. For $n \leq 2$ this holds, and the deficit is computed exactly (Remark \ref{energyremark}).

Section \ref{familysection} studies scale invariant ratios $\mathcal{T}_{\be}$ of these functionals. For every $\be$ some such ratio has only the affine spheres as critical points (Theorem \ref{criticalpointstheorem}). The second variation at the affine sphere is computed for all $\be$ (Corollary \ref{tbetaspectralcorollary}). For $\be \geq 1$ the affine sphere minimizes $\mathcal{T}_{\be}$, by Theorem A and the H\"older inequality, and for $\be \leq 1/2$ the corresponding global inequalities reduce, by the H\"older inequality, to a single conjectural upper bound for the centroaffine volume (Theorem \ref{reductionproposition}).

Section \ref{invariantsection} draws consequences of the identity $\mathfrak{m}(M) = \cV(M)^{2}/N^{2}$ for the minimum $\mathfrak{m}(M) = \min_{L \in \cK(C, \Ga)}\cP(L)$, both sides of which are unchanged under projective duality. For closed surfaces it reformulates known identities and bounds for the area of the equiaffine metric in terms of covolumes. For closed orientable surfaces Theorem \ref{surfacetheorem} gives
\begin{align}\label{mintro}
3\mathfrak{m}(M)^{1/2} = 2\pi|\chi(M)| + \tfrac{1}{2}\int_{M}|\dift|^{2}_{h}\,d\vol_{h} \geq 2^{1/3}\int_{M}|U|^{2/3},
\end{align}
where $\dift$ is the difference tensor of the equiaffine connection of the affine sphere and the Levi-Civita connection of $h$, whose norm is that of the cubic form of the affine sphere, and $U$ is the holomorphic cubic differential whose real part is that cubic form. In particular, $\mathfrak{m}$ is minimized exactly on the Fuchsian locus. Theorem \ref{surfacetheorem} also bounds $3\mathfrak{m}(M)^{1/2}$ above by $2\pi|\chi(M)| + 2^{1/3}\int_{M}|U|^{2/3}$, recovering the two-sided bound of \cite{Nie-entropy, Ouyang-Tamburelli}. Further consequences for surfaces are described in Remark \ref{surfaceremark}. In dimension $n \geq 3$, if $M$ is homeomorphic to a closed hyperbolic manifold $M_{0}$, the theorem of Besson, Courtois, and Gallot \cite{BCG}, applied to the equiaffine metric by means of the lower bound on its Ricci curvature, gives the sharp bound $\cV(M) \geq \vol_{\mathrm{hyp}}(M_{0})$, that is, $\mathfrak{m}(M) \geq \vol_{\mathrm{hyp}}(M_{0})^{2}/N^{2}$, with equality exactly for the hyperbolic structure (Theorem \ref{hyperbolicvolumetheorem}).

The affine sphere equation is the Einstein equation, with negative Einstein constant, in a real analogue of K\"ahler geometry described in Section \ref{kahlereinsteinsection}, in which a flat projective structure replaces the complex structure and the projective Hessian replaces $\partial\bar{\partial}$ (Lemma \ref{kelemma}). Through the tube domain over the cone (Lemma \ref{tubelemma}), the covolume product corresponds formally, at the level of normalized first variations, to the energy functional whose minimizers are the K\"ahler--Einstein metrics of a compact K\"ahler manifold with negative first Chern class, so that Theorems A, B, and C are the counterparts of the variational theory of that functional. Under this correspondence $-N\log\psi_{0}$ is, up to an additive constant, the canonical potential of $C$, which is a self-concordant barrier for $C$ with parameter $N$ (Remark \ref{barrierremark}).

\section{Preliminaries}\label{settingsection}
This section collects elementary facts about integration on $M$, the characteristic function of $C$, positive concave functions of degree one, and polarity. 

The following conventions are used throughout. The space $\ste$ is a real vector space of dimension $N = n+1$ with a fixed volume form $\Psi \in \Lambda^{N}\std$, and $dx = |\Psi|$ is the associated Lebesgue measure. The dual volume form $\Psi^{\ast} \in \Lambda^{N}\ste$ is defined by $\Psi^{\ast}(\zeta_{1}, \dots, \zeta_{N})\Psi(e_{1}, \dots, e_{N}) = 1$ for dual bases $e_{1}, \dots, e_{N}$ of $\ste$ and $\zeta_{1}, \dots, \zeta_{N}$ of $\std$, and $d\zeta = |\Psi^{\ast}|$ is the Lebesgue measure on $\std$.

Everywhere, $C \subset \ste$ denotes an open properly convex cone, and $\Ga \subset SL^{\pm}(\ste)$ is a subgroup that preserves $C$ and acts on $\Om = \projp(C)$ freely, properly, and cocompactly. Here an open convex cone $C$ is \emph{properly convex} if its closure $\bar{C}$ contains no line. The space $\projp(\ste)$ is the sphere of rays of $\ste\setminus\{0\}$ with its flat projective structure, and $\Om = \projp(C)$ is the set of rays of $C$, an open subset of $\projp(\ste)$. The action of $\Ga$ on $\Om$ is \emph{proper} if for every compact $Q \subset \Om$ only finitely many $\gamma \in \Ga$ satisfy $\gamma Q \cap Q \neq \emptyset$, so that $\Ga$ is discrete, and it is \emph{cocompact} if $\Om/\Ga$ is compact. A group $\Ga$ with these properties is said to \emph{divide} $\Om$. Then $M = \Om/\Ga$ is a closed manifold on which the flat projective structure of $\projp(\ste)$ descends, and $M$ is a \emph{closed properly convex flat projective manifold}. A compact set $\Fset \subset C$ and a constant $c_{\Fset}$ with the properties of Lemma \ref{fundamentalsetlemma} below are also fixed throughout.

For $g \in GL(\ste)$ write $g^{-\ast} = (g^{-1})^{\ast}$. The dual group $\Ga^{\ast} = \{g^{-\ast}: g \in \Ga\} \subset SL^{\pm}(\std)$ preserves the open dual cone $C^{\ast} = \{\zeta \in \std: \zeta(x) > 0 \text{ for all } x \in \bar{C}\setminus\{0\}\}$. Since elements of $\Ga$ have determinant $\pm 1$, they preserve $\Psi$ up to sign, hence they preserve $dx$, and, likewise, elements of $\Ga^{\ast}$ preserve $d\zeta$. The canonical isomorphism from $\ste$ to the dual of $\std$, sending $x$ to evaluation at $x$, maps $C$ onto the dual cone of $C^{\ast}$ and conjugates $\Ga$ onto the dual group of $\Ga^{\ast}$. Throughout, objects obtained by applying a duality construction twice are transported to $\ste$ by its inverse, so that for $L \subset C$ the double polar $L^{\circ\circ}$ is a subset of $C$, and for a function $\psi$ on $C$ the double conjugate $\psi^{\ast\ast}$ is a function on $C$. 

A function $F$ on $C$ is homogeneous of degree $k$ if $F(sx) = s^{k}F(x)$ for $s > 0$. Functions on $M = \Om/\Ga$ are identified with $\Ga$-invariant functions on $C$ homogeneous of degree $0$.
The Euler vector field generating dilations of $\ste$ is $\eul_{x} = x$. The flat affine connection of $\ste$ is written $\fD$. For a function $F$, $\fD F = dF$ is its differential and $\fD^{2}F$ its Hessian, and $\fD_{X}Y$ is the derivative of a vector field $Y$ along $X$. The ordinary script letter $D$ is reserved for the Levi-Civita connection of an equiaffine metric (Lemma \ref{gausslemma} and Section \ref{invariantsection}). Where convenient, in particular in the next paragraph and in Sections \ref{familysection} and \ref{kahlereinsteinsection}, tensors are written in abstract index notation, in which $\fD_{i}F$ and $\fD_{i}\fD_{j}F$ are $dF$ and $\fD^{2}F$, a repeated index indicates a contraction, and no coordinates are implied. Enclosure of indices in round (respectively square) brackets indicates complete symmetrization (respectively skew-symmetrization) over the enclosed indices, so that, for instance, $2\nabla_{[i}h_{j]k} = \nabla_{i}h_{jk} - \nabla_{j}h_{ik}$.

For a symmetric bilinear form $B$ on $\ste$, $\det B$ denotes the determinant of the matrix of $B$ in a basis $e_{1}, \dots, e_{N}$ with $\Psi(e_{1}, \dots, e_{N}) = 1$, which does not depend on the choice of such a basis. Intrinsically the determinant of $B$ is an element of $(\Lambda^{N}\std)^{\otimes 2}$, a $2$-density, and $\det B$ is its quotient by $\Psi\otimes\Psi$. The \emph{adjugate tensor} $\adj(B)$ is the symmetric contravariant tensor such that the contraction of $\adj(B)$ with $B$ in one slot is $\det B$ times the identity endomorphism of $\ste$. Intrinsically the adjugate takes values in $S^{2}\ste\otimes(\Lambda^{N}\std)^{\otimes 2}$, and $\adj(B)$ is its contraction with $\Psi^{\ast}\otimes\Psi^{\ast}$. In abstract index notation, $(N-1)!\,\adj(B)^{ij} = \Psi^{\ast\,i i_{2}\dots i_{N}}\Psi^{\ast\,j j_{2}\dots j_{N}}B_{i_{2}j_{2}}\cdots B_{i_{N}j_{N}}$. Consequently $\det B$ and $\adj(B)$ are multiplied by $c^{-2}$ when $\Psi$ is replaced by $c\Psi$. When $B$ is nondegenerate, $(\det B)^{-1}\adj(B)$ is the bilinear form on $\std$ inverse to $B$, meaning that the linear maps $\ste \to \std$ and $\std \to \ste$ they induce are mutually inverse.

For a symmetric contravariant two-tensor $U$ and a symmetric bilinear form $Q$ on $\ste$, $\lb U, Q\ra$ denotes their complete contraction, the natural pairing of $S^{2}\ste$ with $S^{2}\std$. In this notation the Jacobi formula for a smooth one-parameter family of symmetric bilinear forms $B_{s}$ reads $\tfrac{d}{ds}\det B_{s} = \lb\adj(B_{s}), \tfrac{d}{ds}B_{s}\ra$. If $B = \fD^{2}F$ is the Hessian of a function $F$, then $\adj(B)$ is divergence free, $\fD_{i}\adj(\fD^{2}F)^{ij} = 0$, because $\fD_{i}\fD_{i_{2}}\fD_{j_{2}}F$ is symmetric in $i$ and $i_{2}$ while $\Psi^{\ast\,i i_{2}\dots i_{N}}$ is antisymmetric. Equivalently, for every one-form $\al$, the vector field $\adj(\fD^{2}F)(\al, \cdot)$ has divergence $\lb\adj(\fD^{2}F), \fD\al\ra$ with respect to $\Psi$.

If $F$ is a continuous $\Ga$-invariant function on $C$ homogeneous of degree $-N$, the signed density $F\,|\iota_{\eul}\Psi|$, where $\iota_{\eul}\Psi = \Psi(\eul, \cdot, \dots, \cdot)$ is the interior multiplication of $\Psi$ by $\eul$, is annihilated by $\eul$ and is invariant under the flow of $\eul$, because $\lie_{\eul}\Psi = N\Psi$ and $F$ is homogeneous of degree $-N$, so it descends to a $\Ga$-invariant signed density on $\Om$, hence to a signed density on $M$. Its integral is written $\int_{M}F$. The same definition and notation apply when $F$ is only bounded and measurable.

\begin{lemma}\label{conevolumelemma}
Let $\psi$ be a positive continuous $\Ga$-invariant function on $C$ homogeneous of degree $1$, let $K_{\psi} = \{x \in C: \psi(x) \leq 1\}$, and let $f$ be a continuous function on $M$. Then
\begin{align}\label{conevolume}
\int_{M}f\psi^{-N} = N\int_{K_{\psi}/\Ga}f\,dx.
\end{align}
\end{lemma}

\begin{proof}
Let $S \subset C$ be a smooth hypersurface transverse to $\eul$ that meets each ray of $C$ at most once, let $E \subset S$ be a Borel set, and let $R_{E} = \{sy: s > 0,\ y \in E\}$ be the cone over $E$. In the coordinates $(s, y) \in (0, \infty)\times E$ on $R_{E}$ given by $x = sy$, $\Psi_{sy}(y, sv_{1}, \dots, sv_{n}) = s^{n}(\iota_{\eul}\Psi)_{y}(v_{1}, \dots, v_{n})$ for $v_{i} \in T_{y}S$, so $dx = s^{n}ds\,\big|(\iota_{\eul}\Psi)|_{S}\big|$. Since $\psi$ is homogeneous of degree $1$, $sy \in K_{\psi}$ if and only if $s \leq \psi(y)^{-1}$, and since $f$ is homogeneous of degree $0$, $f(sy) = f(y)$. Hence
\begin{align}\label{conevolumelocal}
\int_{K_{\psi}\cap R_{E}}f\,dx = \int_{E}f(y)\int_{0}^{\psi(y)^{-1}}s^{n}\,ds\,|\iota_{\eul}\Psi| = \frac{1}{N}\int_{E}f\psi^{-N}|\iota_{\eul}\Psi|.
\end{align}
Since $M$ is compact, there are finitely many such hypersurfaces $S_{1}, \dots, S_{m}$ and Borel sets $E_{j} \subset S_{j}$ such that the cones $R_{E_{1}}, \dots, R_{E_{m}}$ are pairwise disjoint and their union $R$ is a fundamental domain for the action of $\Ga$ on $C$. Equivalently, radial projection maps the $E_{j}$ bijectively onto disjoint Borel subsets of $\Om$ whose images in $M$ are pairwise disjoint and cover $M$. Since $K_{\psi}$ is $\Ga$-invariant, $K_{\psi}\cap R$ is a fundamental domain for the action of $\Ga$ on $K_{\psi}$, so $\int_{K_{\psi}/\Ga}f\,dx = \sum_{j=1}^{m}\int_{K_{\psi}\cap R_{E_{j}}}f\,dx$, while $\int_{M}f\psi^{-N} = \sum_{j=1}^{m}\int_{E_{j}}f\psi^{-N}|\iota_{\eul}\Psi|$ by the definition of $\int_{M}$ preceding the lemma. Summing \eqref{conevolumelocal} with $E_{j}$ in place of $E$ over $j$ gives \eqref{conevolume}.
\end{proof}

\begin{lemma}\label{divergencelemma}
Let $Y$ be a $\Ga$-invariant $C^{1}$ vector field on $C$ whose coefficients are homogeneous of degree $1-N$. Then $\int_{M}\operatorname{div}Y = 0$, where $\operatorname{div}$ is the divergence with respect to $dx$.
\end{lemma}
\begin{proof}
Since $[\eul, Y] = -NY$ and $\lie_{\eul}\Psi = N\Psi$, $\lie_{\eul}(\iota_{Y}\Psi) = 0$, so $\iota_{\eul}((\operatorname{div}Y)\Psi) = \iota_{\eul}d(\iota_{Y}\Psi) = -d(\iota_{\eul}\iota_{Y}\Psi)$. The form $\iota_{\eul}\iota_{Y}\Psi$ is annihilated by $\eul$ and $\eul$-invariant, so it descends to $M$ (with values in the orientation bundle), and the claim follows from the Stokes theorem on the closed manifold $M$.
\end{proof}

\begin{definition}
Let $\cA(C, \Ga)$ be the set of functions $\psi: C \to (0, \infty)$ that are concave, homogeneous of degree $1$, and $\Ga$-invariant. Let $\cAs(C, \Ga) \subset\cA(C, \Ga)$ be the subset comprising smooth $\psi\in \cA(C, \Ga)$ such that $\fD^{2}\psi$ is negative semidefinite with radical spanned by $\eul$. 
\end{definition}

\begin{remark}
For $\psi \in \cAs(C, \Ga)$, $-\psi$ can be identified with a negative section of the dual of the tautological bundle of $M$ with positive definite Hessian. By a theorem of Loftin \cite{Loftin-met}, a closed flat projective manifold admitting such a section is properly convex.
\end{remark}

The \emph{characteristic function} of $C$ is defined by
\begin{align}
\phi_{C}(x) = \int_{C^{\ast}}e^{-\zeta(x)}d\zeta,
\end{align}
and $\rho_{C} = \phi_{C}^{-1/N}$. By Lemma \ref{vinberglemma}, $\rho_{C} \in \cAs(C, \Ga)$. The properties of $\phi_{C}$ used in Lemma \ref{vinberglemma} are classical and due to Vinberg \cite{Vinberg}, and the concavity of $\rho_{C}$ is a short consequence of them.

\begin{lemma}\label{vinberglemma}
The function $\rho_{C} = \phi_{C}^{-1/N}$ is smooth, positive, homogeneous of degree $1$, invariant under $\{g \in SL^{\pm}(\ste): gC = C\}$, and concave, and the Hessian $\fD^{2}\rho_{C}$ is negative semidefinite with radical spanned by $\eul$. Moreover $\rho_{C}(x) \to 0$ as $x \to \pr C\setminus\{0\}$. If a linear function $\ell \in \std$ and a constant $c$ satisfy $|\ell| \leq c\rho_{C}$ on $C$, then $\ell = 0$.
\end{lemma}

\begin{proof}
By \cite[\S $1$, n$^{\circ}$ $3$]{Vey} (see also \cite[Chapter I, \S 2]{Vinberg} and \cite[Propositions $4.3$--$4.5$]{Shima}), $\phi_{C}$ is analytic and positive on $C$, $\phi_{C}(gx) = |\det g|^{-1}\phi_{C}(x)$ for every linear automorphism $g$ of $C$ and $x \in C$, the Hessian of $\Lambda = \log\phi_{C}$ is positive definite, and $\phi_{C}(x) \to \infty$ as $x \to \pr C\setminus\{0\}$. Taking $g$ to be a homothety shows that $\phi_{C}$ is homogeneous of degree $-N$. Hence $\rho_{C} = \phi_{C}^{-1/N}$ is smooth, positive, homogeneous of degree $1$, and invariant under $\{g \in SL^{\pm}(\ste): gC = C\}$, and $\rho_{C}(x) \to 0$ as $x \to \pr C\setminus\{0\}$.

Let $x \in C$. Since $\Lambda(tx) = \Lambda(x) - N\log t$ for $t > 0$, $d\Lambda(x)(x) = -N$, and differentiating this identity gives $\fD^{2}\Lambda(x)(x, \cdot) = -d\Lambda(x)$, so that $\fD^{2}\Lambda(x)(x, x) = N$. By the Cauchy--Schwarz inequality for the inner product $\fD^{2}\Lambda(x)$ on $\ste$, $d\Lambda(x)(v)^{2} = \fD^{2}\Lambda(x)(x, v)^{2} \leq N\fD^{2}\Lambda(x)(v, v)$ for $v \in \ste$, with equality if and only if $v$ is a multiple of $x$. Since $\rho_{C} = e^{-\Lambda/N}$, it follows that
\begin{align}\label{rhochessian}
\fD^{2}\rho_{C}(x)(v, v) = -\tfrac{1}{N}\rho_{C}(x)\Big(\fD^{2}\Lambda(x)(v, v) - \tfrac{1}{N}d\Lambda(x)(v)^{2}\Big) \leq 0, \qquad v \in \ste,
\end{align}
with equality if and only if $v$ is a multiple of $x$. Hence $\fD^{2}\rho_{C}$ is negative semidefinite with radical spanned by $\eul$, and in particular $\rho_{C}$ is concave.

Finally, let $\ell \in \std$ and a constant $c$ satisfy $|\ell| \leq c\rho_{C}$ on $C$. Since $\ell$ is continuous and $\rho_{C}(x) \to 0$ as $x \to \pr C\setminus\{0\}$, $\ell$ vanishes on $\pr C$. Let $x \in C$ and let $P \subset \ste$ be a two-dimensional subspace containing $x$, which exists because $N \geq 2$. The set $\bar{C}\cap P$ is a closed convex cone in $P$ that contains no line and has nonempty interior in $P$, so it is the sector bounded by the rays spanned by two linearly independent vectors $y_{1}, y_{2} \in P$. Since every point of $C \cap P$ is an interior point of $\bar{C}\cap P$ in $P$, $y_{1}, y_{2} \in \pr C$. Since $\ell(y_{1}) = \ell(y_{2}) = 0$ and $y_{1}$ and $y_{2}$ span $P$, $\ell(x) = 0$. So $\ell$ vanishes on the open set $C$, and therefore on $\ste$.
\end{proof}

By the homogeneity of degree $1$ and the concavity of $\psi \in \cA(C, \Ga)$, for $a, b \in C$,
\begin{align}\label{superadditive}
\psi(a+b) = 2\psi(\tfrac{a+b}{2}) \geq \psi(a) + \psi(b),
\end{align}
so every $\psi \in \cA(C, \Ga)$ is superadditive.

\begin{lemma}\label{fundamentalsetlemma}
There are a compact set $\Fset \subset C$ and a constant $c_{\Fset} \in (0, 1]$ such that
\begin{align}\label{fundamentalset}
C = \{s\gamma x : s > 0,\ \gamma \in \Ga,\ x \in \Fset\}
\end{align}
and
\begin{align}\label{fsetproperty}
y - tx \in C \qquad \text{for all } x, y \in \Fset \text{ and } 0 < t < c_{\Fset}.
\end{align}
Precisely, $c_{\Fset} = \operatorname{dist}(\Fset, \partial C)/\max_{x \in \Fset}|x|$.
\end{lemma}
\begin{proof}
Since $\Ga$ acts cocompactly on $\Om$, there is a compact set $\Fset \subset C$ satisfying \eqref{fundamentalset}. Fix a norm $|\cdot|$ on $\ste$. Since $\Fset$ is a compact subset of the open cone $C$, the numbers $r = \operatorname{dist}(\Fset, \partial C)$ and $R = \max_{x \in \Fset}|x|$ are positive and finite, and $r \leq R$ because $0 \in \pr C$. Set $c_{\Fset} = r/R$. For $x, y \in \Fset$ and $0 < t < c_{\Fset}$, $|tx| < c_{\Fset}R = r$, so $y - tx$ lies in the open ball of radius $r$ about $y$, which is contained in $C$.
\end{proof}

\begin{lemma}[Harnack inequality]\label{harnacklemma}
There is $c \in (0, 1]$, depending only on $C$ and $\Ga$, such that every $\psi \in \cA(C, \Ga)$ satisfies
\begin{align}\label{harnack}
\sup_{C}\frac{\psi}{\rho_{C}} \leq c^{-1}\inf_{C}\frac{\psi}{\rho_{C}}.
\end{align}
In particular every $\psi \in \cA(C, \Ga)$ is bounded above by a multiple of $\rho_{C}$, and so tends to $0$ at $\pr C\setminus\{0\}$.
\end{lemma}
\begin{proof}
Let $x, y \in \Fset$ and $0 < t < c_{\Fset}$, with $\Fset$ and $c_{\Fset}$ as in Lemma \ref{fundamentalsetlemma}. By \eqref{fsetproperty}, $y - tx \in C$, so, by superadditivity, homogeneity, and nonnegativity of $\psi$,
\begin{align}\label{superadditivestep}
\psi(y)  = \psi(tx + (y -tx)) \geq \psi(tx) + \psi(y - tx) \geq t\psi(x).
\end{align}
Taking the supremum of \eqref{superadditivestep} over $t \in (0, c_{\Fset})$ yields $\psi(y) \geq c_{\Fset}\psi(x)$. Consequently, for $x, y \in \Fset$,
\begin{align}\label{harnackonF}
\frac{\psi(y)}{\rho_{C}(y)} \geq c_{\Fset}\frac{\rho_{C}(x)}{\rho_{C}(y)}\frac{\psi(x)}{\rho_{C}(x)} \geq c\,\frac{\psi(x)}{\rho_{C}(x)}, \qquad c = c_{\Fset}\frac{\min_{\Fset}\rho_{C}}{\max_{\Fset}\rho_{C}}.
\end{align}
Since $\psi/\rho_{C}$ is $\Ga$-invariant and homogeneous of degree $0$, it follows from \eqref{fundamentalset} of Lemma \ref{fundamentalsetlemma} that each value of $\psi/\rho_{C}$ on $C$ is attained at some point of $\Fset$. Hence \eqref{harnackonF} holds for all $x, y \in C$, which is \eqref{harnack}. The last claim in the statement of the lemma follows from Lemma \ref{vinberglemma}.
\end{proof}

For $\psi \in \cA(C, \Ga)$ write $L_{\psi} = \{x \in C: \psi(x) \geq 1\}$, $K_{\psi} = \{x \in C: \psi(x) \leq 1\}$, and $\Sigma_{\psi} = \{x \in C: \psi(x) = 1\}$, so that $C\setminus L_{\psi} = \{x \in C: \psi(x) < 1\}$ and $K_{\psi}$ differ by $\Sigma_{\psi}$, which is the boundary of the convex set $L_{\psi}$, hence a null set, and which meets each ray of $C$ exactly once, because $\psi$ is positive, finite, and homogeneous of degree $1$. For $\psi \in \cAs(C, \Ga)$, $\Sigma_{\psi}$ is a smooth hypersurface, since $d\psi(\eul) = \psi \neq 0$. Lemma \ref{smoothboundarylemma} shows that it is locally strongly convex, in the sense defined in Section \ref{functionalssection}, and that these two properties of $\Sigma_{\psi}$ characterize $\cAs(C, \Ga)$ within $\cA(C, \Ga)$.

Lemma \ref{orbitlemma} shows that every $\Ga$-orbit in $\Om$ accumulates at every extreme point of $\pr\Om$. This is due to Vey \cite[Section $3$, Lemme $5$]{Vey}, who proves it using the Hilbert metric and deduces from it that the convex hull of every $\Ga$-orbit in $\Om$ is $\Om$ \cite[Section $3$, Proposition $3$]{Vey}, a result that he attributes to Koszul. The proof given here uses only Lemma \ref{fundamentalsetlemma}. Together with the Harnack inequality, Lemma \ref{harnacklemma}, it yields Lemma \ref{gammaconvexbodylemma}, which describes the $\Ga$-convex bodies of the introduction.

\begin{lemma}\label{orbitlemma}
Fix $\zeta_{0} \in C^{\ast}$, let $\mathsf{B} = \{x \in \bar{C}: \zeta_{0}(x) = 1\}$, and let $\pi(x) = x/\zeta_{0}(x)$ for $x \in \bar{C}\setminus\{0\}$. Then $\mathsf{B}$ is a compact convex set, and for every $x \in C$ and every extreme point $p$ of $\mathsf{B}$ there is a sequence $\gamma_{k} \in \Ga$ such that $\pi(\gamma_{k}x) \to p$.
\end{lemma}
\begin{proof}
Since $\zeta_{0} > 0$ on $\bar{C}\setminus\{0\}$, for some $c > 0$ there holds $\zeta_{0} \geq c|\cdot|$ on $\bar{C}$, so $\mathsf{B}$ is closed and bounded, and it is convex as the intersection of the convex set $\bar{C}$ with a hyperplane. Let $\Fset$ and $c_{\Fset}$ be as in Lemma \ref{fundamentalsetlemma}, and fix $\ep$ with $0 < \ep < \min\{1, c_{\Fset}\}$, so that $f - \ep f' \in C$ for all $f, f' \in \Fset$ by \eqref{fsetproperty}.

First it is shown that the following property holds.
\begin{itemize}
\item[(P)] If $\gamma_{k} \in \Ga$ and $f_{k} \in \Fset$ are such that $\pi(\gamma_{k}f_{k})$ converges to an extreme point $p$ of $\mathsf{B}$, then $\pi(\gamma_{k}f'_{k}) \to p$ for every sequence $f'_{k} \in \Fset$.
\end{itemize}
By the choice of $\ep$, $c_{k} = f_{k} - \ep f'_{k}$ and $c'_{k} = f'_{k} - \ep f_{k}$ are contained in $C$. Since $\gamma_{k}$ preserves $C$ and $\zeta_{0}$ is positive on $C$, applying $\gamma_{k}$ and then $\zeta_{0}$ to the identity $f_{k} = \ep f'_{k} + c_{k}$ gives
\begin{align}
\zeta_{0}(\gamma_{k}f_{k}) = \ep\zeta_{0}(\gamma_{k}f'_{k}) + \zeta_{0}(\gamma_{k}c_{k}) > \ep\zeta_{0}(\gamma_{k}f'_{k}),
\end{align}
and the same argument applied to $f'_{k} = \ep f_{k} + c'_{k}$ gives $\zeta_{0}(\gamma_{k}f'_{k}) > \ep\zeta_{0}(\gamma_{k}f_{k})$. By these inequalities, $w_{k} = \ep\zeta_{0}(\gamma_{k}f'_{k})/\zeta_{0}(\gamma_{k}f_{k})$ is contained in $(\ep^{2}, 1)$. Dividing $\gamma_{k}f_{k} = \ep\gamma_{k}f'_{k} + \gamma_{k}c_{k}$ by $\zeta_{0}(\gamma_{k}f_{k})$, and noting that $\zeta_{0}(\gamma_{k}c_{k})/\zeta_{0}(\gamma_{k}f_{k}) = 1 - w_{k}$, gives
\begin{align}\label{extremestep}
\pi(\gamma_{k}f_{k}) = w_{k}\pi(\gamma_{k}f'_{k}) + (1 - w_{k})\pi(\gamma_{k}c_{k}).
\end{align}
Suppose that $\pi(\gamma_{k}f'_{k})$ does not converge to $p$. Then there are an open neighborhood $U$ of $p$ in $\mathsf{B}$ and a subsequence $k_{j}$ such that $\pi(\gamma_{k_{j}}f'_{k_{j}}) \notin U$ for all $j$. The points $\pi(\gamma_{k_{j}}f'_{k_{j}})$ and $\pi(\gamma_{k_{j}}c_{k_{j}})$ are contained in the compact set $\mathsf{B}$, and $w_{k_{j}}$ is contained in the compact interval $[\ep^{2}, 1]$, so, after passing to a further subsequence, $\pi(\gamma_{k_{j}}f'_{k_{j}}) \to q$, $\pi(\gamma_{k_{j}}c_{k_{j}}) \to r$, and $w_{k_{j}} \to w$ as $j \to \infty$, for some $q, r \in \mathsf{B}$ and $w \in [\ep^{2}, 1]$, with $q \notin U$ because $\mathsf{B}\setminus U$ is closed. Letting $j \to \infty$ in \eqref{extremestep} gives $p = wq + (1 - w)r$. If $w = 1$ this says $q = p$, and if $0 < w < 1$ then $q = p$ because $p$ is an extreme point of $\mathsf{B}$. Either way $q = p \in U$, a contradiction. Hence $\pi(\gamma_{k}f'_{k}) \to p$, which proves (P).

Now let $x \in C$ and let $p$ be an extreme point of $\mathsf{B}$. Since $p \in \bar{C}\setminus\{0\}$, there are $z_{k} \in C$ with $\pi(z_{k}) \to p$. By \eqref{fundamentalset} of Lemma \ref{fundamentalsetlemma}, write $z_{k} = s_{k}\gamma_{k}f_{k}$ with $s_{k} > 0$, $\gamma_{k} \in \Ga$, and $f_{k} \in \Fset$, and $x = s\gamma f$ with $f \in \Fset$. Then $\pi(\gamma_{k}f_{k}) = \pi(z_{k}) \to p$, so property (P) with $f'_{k} = f$ for all $k$ gives $\pi(\gamma_{k}f) \to p$, and $\pi(\gamma_{k}\gamma^{-1}x) = \pi(s\gamma_{k}f) = \pi(\gamma_{k}f) \to p$.
\end{proof}

Lemma \ref{gammaconvexbodylemma} shows that the superlevel sets $L_{\psi}$, $\psi \in \cA(C, \Ga)$, are exactly the $\Ga$-convex bodies of the introduction, and that in the definition of a $\Ga$-convex body the conditions that $L + C \subset L$ and that $L$ be closed in $\ste$ are consequences of the other conditions.

\begin{lemma}\label{gammaconvexbodylemma}
Let $L \subset C$ be a nonempty convex $\Ga$-invariant set that is relatively closed in $C$ and different from $C$. Then $L + C \subset L$, $L$ is closed in $\ste$, and $0 \notin L$, that is, $L$ is a $\Ga$-convex body. Every $\zeta \in \std$ that is bounded from below on $L$ lies in $\bar{C}^{\ast}$. Its cogauge $\psi_{L}(x) = \sup\{s > 0: x \in sL\}$ is contained in $\cA(C, \Ga)$, and $L = L_{\psi_{L}}$. Conversely, for every $\psi \in \cA(C, \Ga)$, $L_{\psi}$ is a $\Ga$-convex body with cogauge $\psi$.
\end{lemma}

\begin{proof}
Let $\zeta_{0}$, $\mathsf{B}$, and $\pi$ be as in Lemma \ref{orbitlemma}. Let $x \in L$ and let $\zeta \in \std$ be bounded from below on $L$. Suppose $\zeta(u) < 0$ for some $u \in C$. Since $\mathsf{B}$ is the convex hull of its extreme points and $\zeta$ is affine on $\mathsf{B}$, there is an extreme point $p$ of $\mathsf{B}$ with $\zeta(p) < 0$, and $p \in \pr C$ because an extreme point of $\mathsf{B}$ is contained in the relative boundary of $\mathsf{B}$, which is $\mathsf{B} \cap \pr C$. By Lemma \ref{orbitlemma} there are $\gamma_{k} \in \Ga$ such that $x_{k} = \gamma_{k}x$ satisfies $\pi(x_{k}) \to p$, and $x_{k} \in L$ by the $\Ga$-invariance of $L$.

Write $x_{k} = t_{k}\pi(x_{k})$ with $t_{k} = \zeta_{0}(x_{k}) > 0$. If $t_{k}$ were bounded along a subsequence, then along a further subsequence $x_{k}$ would converge to a point $tp$ with $t \geq 0$, which is contained in $\pr C \cup \{0\}$. This is impossible, because $\rho_{C}(x_{k}) = \rho_{C}(x) > 0$ for all $k$ by the invariance in Lemma \ref{vinberglemma}, whereas $\rho_{C}$ extends continuously by $0$ to $\bar{C}$. Indeed, $\rho_{C}(y) \to 0$ as $y \to \pr C\setminus\{0\}$ by Lemma \ref{vinberglemma}, so the extension of $\rho_{C}$ by $0$ is continuous on $\bar{C}\setminus\{0\}$, hence bounded on the compact set of unit vectors of $\bar{C}$ for a norm $|\cdot|$ on $\ste$, and so, by homogeneity, there is a constant $c_{0}$ such that $\rho_{C}(y) \leq c_{0}|y|$ for all $y \in C$, which shows that the extension is continuous at $0$. Hence $t_{k} \to \infty$, so $\zeta(x_{k}) = t_{k}\zeta(\pi(x_{k})) \to -\infty$, contradicting the boundedness of $\zeta$ from below on $L$.

Therefore every $\zeta \in \std$ bounded from below on $L$ is nonnegative on $C$, that is, lies in $\bar{C}^{\ast}$. The closure $\bar{L}$ of $L$ in $\ste$ is the intersection of the closed half-spaces $\{y \in \ste: \zeta(y) \geq a\}$, with $\zeta \in \std\setminus\{0\}$ and $a \in \rea$, that contain $L$; each such $\zeta$ is bounded below on $L$, so lies in $\bar{C}^{\ast}$, and therefore $\{y \in \ste: \zeta(y) \geq a\} + C \subset \{y \in \ste: \zeta(y) \geq a\}$, so $\bar{L} + C \subset \bar{L}$, and $L + C \subset \bar{L}\cap C = L$ because $L$ is relatively closed in $C$.

Since $x + C \subset L$ and, for $c \in C$, $tc - x \in C$ for $t$ large, every ray of $C$ meets $L$, so $\psi_{L} > 0$ on $C$. For $0 < s' < s$, convexity of $L$ gives $sL \subset s'L + (s - s')L \subset s'L + C \subset s'L$, so $\{x \in C: \psi_{L}(x) \geq s\} = \bigcap_{s' < s}s'L = sL$, the last equality because $L$ is relatively closed. In particular $L = \{x \in C: \psi_{L}(x) \geq 1\}$. Suppose $\psi_{L}(y) = \infty$ for some $y \in C$. Then $y/s \in L$ for all $s > 0$, so $0 \in \bar{L}$. For every $y' \in L$ the segments from $y/s$ to $y'$ are contained in $L$, so $(0, 1]y' \subset \bar{L}\cap C = L$, and with $L + C \subset L$ and the fact that every ray meets $L$ this gives $L = C$, which is excluded. Hence $\psi_{L}$ is finite. It is homogeneous of degree $1$ and $\Ga$-invariant by construction, and superadditive, hence concave, because $y/\psi_{L}(y), y'/\psi_{L}(y') \in L$ and the convexity of $L$ give $(y + y')/(\psi_{L}(y) + \psi_{L}(y')) \in L$. Thus $\psi_{L} \in \cA(C, \Ga)$. By Lemma \ref{harnacklemma}, $\psi_{L} \leq c_{1}\rho_{C}$ for some $c_{1} > 0$, so $L \subset \{x \in C: \rho_{C}(x) \geq c_{1}^{-1}\}$. Since $\rho_{C}$ extends continuously by $0$ to $\bar{C}$, as shown above, this set is closed in $\ste$ and disjoint from $\pr C \cup \{0\}$. Hence $\bar{L} \subset C$, so $L = \bar{L}\cap C$ is closed in $\ste$ and $0 \notin L$.

Conversely, if $\psi \in \cA(C, \Ga)$, then $L_{\psi}$ is convex, $\Ga$-invariant, relatively closed in $C$, nonempty, and different from $C$, since $\psi(sx) = s\psi(x) < 1$ for $x \in C$ and $s > 0$ sufficiently small, so it is a $\Ga$-convex body by what has been shown, and $\{x \in C: \psi(x) \geq s\} = sL_{\psi}$ by homogeneity gives $\psi_{L_{\psi}} = \psi$.
\end{proof}

The best constant in Lemma \ref{harnacklemma} is the infimum of the corresponding ratio over the $\Ga$-convex bodies $L_{x}$, $x \in C$, defined below, each of which is generated by a single $\Ga$-orbit.
\begin{definition}\label{sharpharnackdefinition}
For $\phi \in \cA(C, \Ga)$, the \emph{sharp Harnack constant of $M$ relative to $\phi$} is
\begin{align}\label{sharpharnackdefined}
\mathfrak{c}_{\phi}(M) = \inf_{\psi \in \cA(C, \Ga)}\frac{\inf_{C}\psi/\phi}{\sup_{C}\psi/\phi}.
\end{align}
\end{definition}
By Lemma \ref{harnacklemma}, applied to $\psi$ and to $\phi$, $\mathfrak{c}_{\phi}(M) > 0$, and $\mathfrak{c}_{\rho_{C}}(M)$ is the best constant in \eqref{harnack}. Since the ratio in \eqref{sharpharnackdefined} is unchanged when $\phi$ is multiplied by a positive constant, $\mathfrak{c}_{\phi}(M)$ depends on $\phi$ only up to such a factor. For $x \in C$ let $L_{x}$ be the closure in $C$ of $\operatorname{conv}(\Ga x) + C$, where $\operatorname{conv}$ denotes the convex hull, and let $\psi_{x}$ be its cogauge.
\begin{theorem}\label{sharpharnacktheorem}
Let $\phi \in \cA(C, \Ga)$. For every $x \in C$, $L_{x}$ is a $\Ga$-convex body, and for every $\psi \in \cA(C, \Ga)$ there is $x \in C$ such that
\begin{align}\label{sharpharnack}
\frac{\inf_{C}\psi/\phi}{\sup_{C}\psi/\phi} \geq \frac{\inf_{C}\psi_{x}/\phi}{\sup_{C}\psi_{x}/\phi} = \frac{\phi(x)}{\sup_{\pr L_{x}}\phi}.
\end{align}
Consequently $\mathfrak{c}_{\phi}(M) = \inf_{x \in C}\phi(x)/\sup_{\pr L_{x}}\phi$.
\end{theorem}
\begin{proof}
The set $\{\phi \geq \phi(x)\}$ is convex, relatively closed in $C$, $\Ga$-invariant, and contains $\Ga x$, and by \eqref{superadditive} it is preserved by addition of elements of $C$, so it contains $L_{x}$. In particular $L_{x} \neq C$, so $L_{x}$ is a $\Ga$-convex body, and $\psi_{x} \leq \phi/\phi(x)$, with equality at $x$. Let $\psi \in \cA(C, \Ga)$. The function $\psi/\phi$ is identified with a continuous function on $M$, so it attains its maximum at some $x \in C$, and, since \eqref{sharpharnack} is unchanged by rescaling $\psi$, it can be supposed that $\psi(x) = 1$. Then $\psi \geq 1$ on $\Ga x$, hence on $\operatorname{conv}(\Ga x)$ by concavity, on $\operatorname{conv}(\Ga x) + C$ by \eqref{superadditive}, and on $L_{x}$ by continuity, so $\psi \geq \psi_{x}$. Hence $\sup_{C}\psi/\phi = 1/\phi(x) = \sup_{C}\psi_{x}/\phi$ and $\inf_{C}\psi/\phi \geq \inf_{C}\psi_{x}/\phi$, which gives the inequality in \eqref{sharpharnack}. Since $\psi_{x} = 1$ on $\pr L_{x}$, which meets every ray of $C$, $\inf_{C}\psi_{x}/\phi = 1/\sup_{\pr L_{x}}\phi$.
\end{proof}

For $\psi \in \cA(C, \Ga)$ define
\begin{align}\label{polardefined}
\psi^{\ast}(\zeta) = \inf_{x \in C}\frac{\zeta(x)}{\psi(x)}, \qquad \zeta \in C^{\ast}.
\end{align}
As in the introduction, the \emph{polar} of a $\Ga$-convex body $L \subset C$ is $L^{\circ} = \{\zeta \in C^{\ast}: \zeta(x) \geq 1 \text{ for all } x \in L\}$. For $\zeta \in C^{\ast}$ and $x \in C$, $\zeta(x)/\psi(x) = \zeta(x/\psi(x))$ with $x/\psi(x) \in \Sigma_{\psi} \subset L_{\psi}$, while for $x \in L_{\psi}$, $\zeta(x) \geq \zeta(x)/\psi(x)$ because $\psi(x) \geq 1$. Consequently $L^{\circ}_{\psi} = (L_{\psi})^{\circ}$ satisfies $L^{\circ}_{\psi} = \{\zeta \in C^{\ast}: \psi^{\ast}(\zeta) \geq 1\}$. This shows that $\psi^{\ast}$ is the support function $s_{L}(\zeta) = \inf_{x \in L}\zeta(x)$ of $L = L_{\psi}$, restricted to $C^{\ast}$. 

In general $L^{\circ}$ need not be closed in $\std$. Up to the reflection $\zeta \mapsto -\zeta$, $L^{\circ}$ is the copolar set of Rashkovskii \cite{Rashkovskii-copolar} and Schneider \cite{Schneider-copolarity}, and polarity is the order reversing transform of $L$ studied by Artstein-Avidan, Sadovsky, and Wyczesany \cite{AASW-zoo}; see also \cite{Xu-Li-Leng}, where this duality is characterized among the dualities of pseudo-cones.

For $\psi \in \cAs(C, \Ga)$ the following notation is used throughout the rest of the paper:
\begin{align}\label{phiqdefined}
\Phi = \tfrac{1}{2}\psi^{2}, \qquad \fD^{2}\Phi = \psi\,\fD^{2}\psi + d\psi\otimes d\psi, \qquad G_{\psi} = \fD\Phi = \psi\,d\psi: C \to \std.
\end{align}
When more than one function $\psi \in \cAs(C, \Ga)$ is under consideration, $\Phi$ is written $\Phi_{\psi}$, to indicate the dependence on $\psi$. Since $\Phi$ is homogeneous of degree $2$, $\fD^{2}\Phi(\eul, \cdot) = G_{\psi}$ and $G_{\psi}(\eul) = 2\Phi = \psi^{2}$.

\begin{lemma}\label{polarmaplemma}
For $\psi \in \cAs(C, \Ga)$:
\begin{enumerate}[label=$(\arabic*)$, ref=(\arabic*)]
\item\label{pml1} $\fD^{2}\Phi$ is nondegenerate of signature $(1, n)$, with $\fD^{2}\Phi(\eul, \eul) = \psi^{2}$ and $\fD^{2}\Phi$ negative definite on $\ker d\psi$. Conversely, if $\psi: C \to (0, \infty)$ is smooth, $\Ga$-invariant, and homogeneous of degree $1$, and $\fD^{2}(\tfrac{1}{2}\psi^{2})$ has signature $(1, n)$ at every point, then $\psi \in \cAs(C, \Ga)$.
\item\label{pml2} $G_{\psi}$ is homogeneous of degree $1$, satisfies $G_{\psi}(gx) = g^{-\ast}G_{\psi}(x)$ for $g \in \Ga$, has Jacobian $\fD G_{\psi} = \fD^{2}\Phi$, and satisfies $G_{\psi}(x)(z) \geq \psi(x)\psi(z)$ for $x, z \in C$, with $G_{\psi}(x)(x) = \psi(x)^{2}$.
\item\label{pml3} $G_{\psi}$ is a diffeomorphism of $C$ onto $C^{\ast}$, and $\Ga^{\ast}$ acts freely, properly, and cocompactly on $\Om^{\ast} = \projp(C^{\ast})$.
\item \label{pml4}$\psi^{\ast}\circ G_{\psi} = \psi$; $\Phi^{\ast} = \tfrac{1}{2}(\psi^{\ast})^{2}$ satisfies $\fD\Phi^{\ast} = G_{\psi}^{-1}$ and $\fD^{2}\Phi^{\ast}\circ G_{\psi} = (\det \fD^{2}\Phi)^{-1}\adj(\fD^{2}\Phi)$, the bilinear form on $\std$ inverse to $\fD^{2}\Phi$; and $\psi^{\ast} \in \cAs(C^{\ast}, \Ga^{\ast})$.
\end{enumerate}
\end{lemma}
\begin{proof}
Write $G = G_{\psi}$. First \ref{pml1} is proved. Since $\fD^{2}\psi(\eul, \cdot) = 0$ and $d\psi(\eul) = \psi$, $\fD^{2}\Phi(\eul, \eul) = \psi^{2} > 0$ and $\fD^{2}\Phi(\eul, v) = 0$ for $v \in \ker d\psi$. On $\ker d\psi$, which is complementary to $\eul$, $\fD^{2}\Phi = \psi \fD^{2}\psi$ is negative definite. Conversely, let $\psi$ be smooth, positive, $\Ga$-invariant, and homogeneous of degree $1$, with $\fD^{2}(\tfrac{1}{2}\psi^{2})$ of signature $(1, n)$. The same identities give $\fD^{2}\Phi(\eul, \eul) > 0$ and $\fD^{2}\Phi(\eul, v) = 0$ for $v \in \ker d\psi$, so, by the signature, $\fD^{2}\Phi = \psi\fD^{2}\psi$ is negative definite on the complement $\ker d\psi$ of $\eul$. Since $\fD^{2}\psi(\eul, \cdot) = 0$ by homogeneity, $\fD^{2}\psi$ is negative semidefinite with radical spanned by $\eul$. In particular $\psi$ is concave, so $\psi \in \cAs(C, \Ga)$.

Next \ref{pml2} is proved. Homogeneity and equivariance are immediate. By concavity and the Euler identity, $\psi(z) \leq \psi(x) + d\psi(x)(z - x) = d\psi(x)(z)$, so $G(x)(z) = \psi(x)d\psi(x)(z) \geq \psi(x)\psi(z)$, and $G(x)(x) = \psi(x)^{2}$.

Next \ref{pml3} is proved. By \ref{pml1}, $G$ is a local diffeomorphism, so $G(C)$ is open. By \ref{pml2}, $G(x)$ is positive on $C$, so $G(C)$ is contained in the closed dual cone $\bar{C}^{\ast}$, and, being open, in its interior $C^{\ast}$. For $x \in \Sigma_{\psi}$ let $\zeta_{x} = G(x)$. Then $\zeta_{x}(x) = 1$ and $\zeta_{x} \geq 1$ on $L_{\psi}$, so the hyperplane $\{y \in \ste: \zeta_{x}(y) = 1\}$ supports $L_{\psi}$ exactly at $x$, because, by the negative definiteness of $\fD^{2}\Phi$ on $\ker d\psi$ in \ref{pml1}, $\pr L_{\psi} = \Sigma_{\psi}$ contains no segments. If $G(x) = sG(y)$ with $x, y \in \Sigma_{\psi}$ and $s > 0$, then $s = \zeta_{x}(y) \geq 1$ and $1 = s\zeta_{y}(x) \geq s$, so $s = 1$ and $x = y$. By homogeneity $G$ is injective.

To show that the rays of $G(C)$ form a closed subset of $\Om^{\ast}$, let $x_{k} \in \Sigma_{\psi}$ and $s_{k} > 0$ with $s_{k}\zeta_{x_{k}} \to \zeta \in C^{\ast}$. For fixed $z \in C$, $\zeta_{x_{k}}(z) \geq \psi(z)$, so $s_{k}$ is bounded. Fix a norm. Since $\zeta \in C^{\ast}$, $\zeta(v) \geq c_{1}|v|$ on $\bar{C}$. If $|x_{k}| \to \infty$ along a subsequence, then $s_{k} = s_{k}\zeta_{x_{k}}(x_{k}) = |x_{k}|\,s_{k}\zeta_{x_{k}}(x_{k}/|x_{k}|)$, and $s_{k}\zeta_{x_{k}}(x_{k}/|x_{k}|) \geq c_{1}/2$ for large $k$ by uniform convergence on the unit sphere, contradicting the boundedness of $s_{k}$. So $|x_{k}|$ is bounded. By Lemma \ref{harnacklemma}, $1 = \psi(x_{k}) \leq c_{2}\rho_{C}(x_{k}) \leq c_{3}|x_{k}|$, so $|x_{k}|$ is bounded below, and a subsequence converges to $x \in \bar{C}\setminus\{0\}$. By Lemma \ref{harnacklemma} again, $x \notin \pr C$, so $x \in \Sigma_{\psi}$ and $\zeta$ is a multiple of $G(x)$.

The rays of $G(C)$ also form an open subset of $\Om^{\ast}$, because $G(C)$ is open, so, since $\Om^{\ast}$ is connected, $G$ maps the rays of $C$ onto those of $C^{\ast}$, and, by homogeneity, $G(C) = C^{\ast}$. Since $G$ is equivariant, it induces a homeomorphism $\Om \to \Om^{\ast}$ that conjugates the action of $\Ga$ on $\Om$ to the action of $\Ga^{\ast}$ on $\Om^{\ast}$, so the latter action is free, proper, and cocompact.

Finally \ref{pml4} is proved. By \ref{pml2}, $G(x)(z)/\psi(z) \geq \psi(x)$ with equality at $z = x$, so $\psi^{\ast}(G(x)) = \psi(x)$. Since $G$ is a diffeomorphism of $C$ onto $C^{\ast}$ by \ref{pml3}, $\psi^{\ast} = \psi\circ G^{-1}$ is smooth on $C^{\ast}$, and $\Phi^{\ast}\circ G = \Phi$. Differentiating, $\fD^{2}\Phi(\fD\Phi^{\ast}\circ G, \cdot) = \fD\Phi = \fD^{2}\Phi(x, \cdot)$, so $\fD\Phi^{\ast}\circ G = \mathrm{id}$ by nondegeneracy of $\fD^{2}\Phi$, and differentiating again, since $\fD G$ is the linear map $\ste \to \std$, $v \mapsto \fD^{2}\Phi(v, \cdot)$, gives that $\fD^{2}\Phi^{\ast}\circ G$ is the inverse bilinear form $(\det \fD^{2}\Phi)^{-1}\adj(\fD^{2}\Phi)$. The function $\psi^{\ast}$ is concave, as an infimum of linear functions, homogeneous of degree $1$, and $\Ga^{\ast}$-invariant. If $\fD^{2}\psi^{\ast}(v, \cdot) = 0$ at $\zeta$, then $\fD^{2}\Phi^{\ast}(v, \cdot) = d\psi^{\ast}(v)\,d\psi^{\ast}$, and $\fD^{2}\Phi^{\ast}(\zeta, \cdot) = \psi^{\ast}d\psi^{\ast}$ by the Euler identity, so $v$ is a multiple of $\zeta$ by nondegeneracy.
\end{proof}

\begin{lemma}\label{polartransformlemma}
For $\psi \in \cA(C, \Ga)$, $\psi^{\ast} \in \cA(C^{\ast}, \Ga^{\ast})$ and $\psi^{\ast\ast} = \psi$.
\end{lemma}

\begin{proof}
The function $\psi^{\ast}$ is concave, homogeneous of degree $1$, $\Ga^{\ast}$-invariant, and nonnegative. By Lemma \ref{harnacklemma}, $\psi \leq c_{1}\rho_{C}$, so $\psi^{\ast} \geq \rho_{C}^{\ast}/c_{1}$, and $\rho^{\ast}_{C} > 0$ by Lemma \ref{polarmaplemma}\ref{pml4} applied to $\rho_{C} \in \cAs(C, \Ga)$. By Lemma \ref{gammaconvexbodylemma}, $L = L_{\psi}$ is a closed convex subset of $\ste$ and every $\zeta \in \std$ bounded below on $L$ lies in $\bar{C}^{\ast}$, so $L$ is the intersection of the half-spaces $\{y \in \ste: \zeta(y) \geq s_{L}(\zeta)\}$ with $\zeta \in \bar{C}^{\ast}\setminus\{0\}$. Those with $\zeta \in C^{\ast}$ suffice. If $\zeta \in \pr C^{\ast}\setminus\{0\}$ and $\eta \in C^{\ast}$, then $\zeta_{t} = (1-t)\zeta + t\eta \in C^{\ast}$ for $0 < t \leq 1$, and an $x$ satisfying $\zeta_{t}(x) \geq s_{L}(\zeta_{t})$ for all such $t$ satisfies $\zeta_{t}(x) \geq (1-t)s_{L}(\zeta) + ts_{L}(\eta)$, by the concavity of $s_{L}$, hence $\zeta(x) \geq s_{L}(\zeta)$ in the limit $t \to 0$. Since $s_{L} = \psi^{\ast}$ on $C^{\ast}$, $x \in L$ if and only if $\psi^{\ast\ast}(x) = \inf_{\zeta \in C^{\ast}}\zeta(x)/\psi^{\ast}(\zeta) \geq 1$, that is, $L_{\psi^{\ast\ast}} = L_{\psi}$.

A nonnegative function $\varphi$ on $C$ that is homogeneous of degree $1$ is determined by its superlevel set $L_{\varphi} = \{x \in C: \varphi(x) \geq 1\}$, because $\varphi(x) = \sup\{s > 0: x/s \in L_{\varphi}\}$ for $x \in C$, with the convention that the supremum of the empty set is $0$. Applied to $\psi$ and $\psi^{\ast\ast}$, this gives $\psi^{\ast\ast} = \psi$.
\end{proof}

By Lemma \ref{polarmaplemma}\ref{pml3}, $M^{\ast} = \Om^{\ast}/\Ga^{\ast}$ is a closed properly convex flat projective manifold, and all results of this section applicable to $M$ apply equally to $M^{\ast}$. That $\Ga^{\ast}$ divides $\Om^{\ast}$ when $\Ga$ divides $\Om$ is \cite[Section $1$, Lemme $1$]{Vey}, where it is deduced, as here, from an equivariant diffeomorphism of $C$ onto $C^{\ast}$, there the one determined by the characteristic function. See also \cite[Lemme $2.8$]{Benoist-divisibles1}.

Every $\psi \in \cA(C, \Ga)$ is approximated, in the sense that $\psi_{k}/\psi \to 1$ uniformly, by cogauges of $\Ga$-convex bodies cut out by finitely many $\Ga$-orbits of half-spaces, and by elements of $\cAs(C, \Ga)$. For a nonempty finite set $S \subset C^{\ast}$ and $\psi \in \cA(C, \Ga)$, let
\begin{align}\label{finiteorbitcogauge}
\psi_{S}(x) = \inf_{\gamma \in \Ga,\, \zeta \in S}\frac{\zeta(\gamma x)}{\psi^{\ast}(\zeta)}.
\end{align}
\begin{lemma}\label{densitylemma}
Let $\psi \in \cA(C, \Ga)$ and $\ep \in (0, 1)$.
\begin{enumerate}[label=$(\arabic*)$, ref=(\arabic*)]
\item\label{dl1} There is a nonempty finite set $S \subset C^{\ast}$ such that $\psi_{S} \in \cA(C, \Ga)$ and $\psi \leq \psi_{S} \leq (1+\ep)\psi$.
\item\label{dl2} There is $\phi \in \cAs(C, \Ga)$ such that $(1-\ep)\psi \leq \phi \leq (1+\ep)\psi$.
\end{enumerate}
\end{lemma}
\begin{proof}
Since $\psi^{\ast}$ is $\Ga^{\ast}$-invariant, $\zeta(\gamma x)/\psi^{\ast}(\zeta) = (\zeta\circ\gamma)(x)/\psi^{\ast}(\zeta\circ\gamma)$, so the infimum in \eqref{finiteorbitcogauge} is taken over a subset of the $\eta \in C^{\ast}$ over which $\psi(x) = \psi^{\ast\ast}(x) = \inf_{\eta \in C^{\ast}}\eta(x)/\psi^{\ast}(\eta)$ (Lemma \ref{polartransformlemma}) is the infimum. Hence $\psi_{S} \geq \psi > 0$. As an infimum of linear functions $\psi_{S}$ is concave, and it is $\Ga$-invariant and homogeneous of degree $1$, so $\psi_{S} \in \cA(C, \Ga)$.

Fix a norm on $\std$, let $\Fset$ be as in Lemma \ref{fundamentalsetlemma}, and let $A = \max_{\Fset}\psi$. There is $c > 0$ with $\eta(x) \geq c|\eta|$ for $x \in \Fset$ and $\eta \in \bar{C}^{\ast}$, and, by Lemma \ref{harnacklemma} applied to $\psi^{\ast} \in \cA(C^{\ast}, \Ga^{\ast})$, there is $c' > 0$ with $\psi^{\ast} \leq c'\rho_{C^{\ast}}$. For $x \in \Fset$ and $|\eta| = 1$ it follows that $\eta(x)/\psi^{\ast}(\eta) \geq c/(c'\rho_{C^{\ast}}(\eta))$, which exceeds $A \geq \psi(x)$ unless $\eta$ lies in $K^{\ast} = \{\eta \in C^{\ast}: |\eta| = 1,\ \rho_{C^{\ast}}(\eta) \geq c/(c'A)\}$. Because $\rho_{C^{\ast}}$ is continuous on $\bar{C}^{\ast}$ and vanishes on $\pr C^{\ast}$ (Lemma \ref{vinberglemma}), $K^{\ast}$ is a compact subset of $C^{\ast}$, and $\psi(x) = \min_{\eta \in K^{\ast}}\eta(x)/\psi^{\ast}(\eta)$ for $x \in \Fset$. The function $(x, \eta) \mapsto \eta(x)/\psi^{\ast}(\eta)$ is uniformly continuous on $\Fset\times K^{\ast}$, so, if $S$ is a sufficiently fine finite net in $K^{\ast}$, then $\psi_{S} \leq \psi + \ep\min_{\Fset}\psi \leq (1+\ep)\psi$ on $\Fset$. Since $\psi_{S}/\psi$ is $\Ga$-invariant and homogeneous of degree $0$, $\psi_{S} \leq (1+\ep)\psi$ on $C$. This proves \ref{dl1}.

By \ref{dl1}, to prove \ref{dl2} it suffices to approximate $\psi_{S}$. Write $\zeta_{s} = s/\psi^{\ast}(s)$ for $s \in S$, so that $\psi_{S}(x) = \inf_{\gamma \in \Ga,\, s \in S}\zeta_{s}(\gamma x)$, and, for $p > N$, let
\begin{align}\label{softmin}
&\Phi_{p}(x) = \sum_{\gamma \in \Ga}\sum_{s \in S}\zeta_{s}(\gamma x)^{-p},& &\psi_{p} = \Phi_{p}^{-1/p}.&
\end{align}
First, $\Phi_{p}$ converges, with all its derivatives, uniformly on compact subsets of $C$. Let $x_{0} \in C$ and $r > 0$ be such that the closed ball of radius $2r$ about $x_{0}$ lies in $C$, and let $W$ be the closed ball of radius $r$ about $x_{0}$. For $x, w \in W$, $3x - w = 2(x + \tfrac{1}{2}(x - w)) \in C$, so $\eta(w) \leq 3\eta(x)$ for $\eta \in \bar{C}^{\ast}$, and, since $\gamma$ preserves Lebesgue measure,
\begin{align}\label{softminbound}
\zeta_{s}(\gamma x)^{-p} \leq \frac{3^{p}}{|W|}\int_{W}\zeta_{s}(\gamma w)^{-p}\,dw = \frac{3^{p}}{|W|}\int_{\gamma W}\zeta_{s}(y)^{-p}\,dy, \qquad x \in W.
\end{align}
Because $\Ga$ acts properly on $\Om$ and preserves $\rho_{C}$, it acts properly on $C$, so each $y \in C$ lies in at most $m$ of the sets $\gamma W$, for some $m$, and, because $\rho_{C}$ is $\Ga$-invariant and bounded between positive constants $a$ and $b$ on $W$, all $\gamma W$ lie in $\{a \leq \rho_{C} \leq b\}$. Let $\Sigma = \{\rho_{C} = 1\}$ and $\nu = |\iota_{\eul}\Psi|$ on $\Sigma$. Summing \eqref{softminbound} and computing as in \eqref{conevolumelocal} gives, for $x \in W$,
\begin{align}\label{softminintegral}
\Phi_{p}(x) \leq \frac{3^{p}m}{|W|}\sum_{s \in S}\int_{\{a \leq \rho_{C} \leq b\}}\zeta_{s}(y)^{-p}\,dy = \frac{3^{p}m}{|W|}\int_{a}^{b}t^{n-p}\,dt\,\sum_{s \in S}\int_{\Sigma}\zeta_{s}^{-p}\,d\nu.
\end{align}
There is $c_{1} > 0$ with $\zeta_{s}(u) \geq c_{1}|u|$ for $u \in \bar{C}$ and $s \in S$, $\Sigma$ is bounded away from $0$, and, by \eqref{conevolumelocal}, $\nu(\Sigma \cap \{|u| \leq R\})$ is $N$ times the volume of the cone over $\Sigma\cap\{|u| \leq R\}$ truncated by $\Sigma$, which lies in the ball of radius $R$, so $\nu(\Sigma \cap \{|u| \leq R\}) \leq c_{2}R^{N}$. Summing over the sets where $2^{k} \leq |u| < 2^{k+1}$ shows that $\int_{\Sigma}\zeta_{s}^{-p}\,d\nu < \infty$ because $p > N$. Hence $\Phi_{p}$ converges uniformly on $W$. There is $c_{W} > 0$ with $\eta(x) \geq c_{W}|\eta|$ for $x \in W$ and $\eta \in \bar{C}^{\ast}$, so the $k$th derivative of $\zeta_{s}(\gamma x)^{-p}$, which is a multiple of $\zeta_{s}(\gamma x)^{-p-k}(\zeta_{s}\circ\gamma)^{\otimes k}$, is bounded in norm on $W$ by a constant depending only on $k$, $p$, and $c_{W}$ times $\zeta_{s}(\gamma x)^{-p}$, and the differentiated series also converge uniformly on $W$. Thus $\Phi_{p}$ is smooth and positive, and $\psi_{p}$ is smooth, positive, homogeneous of degree $1$, and $\Ga$-invariant.

Each function $x \mapsto \zeta_{s}(\gamma x)^{-p}$ is convex on $C$, so $\{\psi_{p} \geq 1\} = \{\Phi_{p} \leq 1\}$ is convex. A positive function $f$ on $C$ homogeneous of degree $1$ with $\{f \geq 1\}$ convex is concave, because, for $x, y \in C$ and $t \in [0, 1]$, with $\al = tf(x)$ and $\be = (1-t)f(y)$, the point $(tx + (1-t)y)/(\al + \be) = \tfrac{\al}{\al + \be}\tfrac{x}{f(x)} + \tfrac{\be}{\al + \be}\tfrac{y}{f(y)}$ lies in $\{f \geq 1\}$. Hence $\psi_{p} \in \cA(C, \Ga)$. The supremum of the terms of $\Phi_{p}(x)$ is $\psi_{S}(x)^{-p}$, so $\psi_{p} \leq \psi_{S}$. Fix $p_{0} > N$. For $p \geq p_{0}$, since $\zeta_{s}(\gamma x)/\psi_{S}(x) \geq 1$, $\psi_{S}^{p}\Phi_{p} \leq \psi_{S}^{p_{0}}\Phi_{p_{0}} \leq B$, where $B$ is the maximum on $\Fset$ of the continuous, $\Ga$-invariant function $\psi_{S}^{p_{0}}\Phi_{p_{0}}$, homogeneous of degree $0$. Hence $B^{-1/p}\psi_{S} \leq \psi_{p} \leq \psi_{S}$.

Finally, for $\delta > 0$, $\phi = \psi_{p} + \delta\rho_{C}$ is smooth and lies in $\cA(C, \Ga)$, and $\fD^{2}\phi \leq \delta\fD^{2}\rho_{C}$, whose radical is spanned by $\eul$ (Lemma \ref{vinberglemma}), while $\eul$ lies in the radical of $\fD^{2}\phi$ by homogeneity, so $\phi \in \cAs(C, \Ga)$. Since $\rho_{C}/\psi_{p}$ is bounded, $\psi_{p} \leq \phi \leq (1 + \delta\sup_{C}\rho_{C}/\psi_{p})\psi_{p}$. Choosing $S$ by \ref{dl1}, then $p$, then $\delta$ gives \ref{dl2}.
\end{proof}

\section{The functionals \texorpdfstring{$\cS_{\be}$}{S beta}, volumes, centroaffine curvature, and affine spheres}\label{functionalssection}
Recall from \eqref{phiqdefined} that, for $\psi \in \cAs(C, \Ga)$, $\Phi = \tfrac{1}{2}\psi^{2}$, $\fD^{2}\Phi = \psi\,\fD^{2}\psi + d\psi\otimes d\psi$, and $G_{\psi} = \fD\Phi$.
Because $\fD^{2}\Phi$ is nondegenerate by Lemma \ref{polarmaplemma}\ref{pml1}, both $\det\fD^{2}\Phi$ and $\Psi\otimes\Psi$ are nowhere vanishing sections over $C$ of $(\Lambda^{N}\std)^{\otimes 2}$, so there is a positive function $\mu(\psi)$ on $C$ satisfying
\begin{align}\label{mudefined}
\det \fD^{2}\Phi = \pm\, \mu(\psi)\, \Psi \otimes \Psi.
\end{align}
That is, $\mu(\psi) = |\det \fD^{2}\Phi|$ in the notation of the beginning of Section \ref{settingsection}. Concretely, $\mu(\psi) = |\det(\fD^{2}\Phi(e_{i}, e_{j}))|$ for any basis $e_{1}, \dots, e_{N}$ of $\ste$ such that $\Psi(e_{1}, \dots, e_{N}) = 1$.

Since $\fD^{2}\Phi$ is homogeneous of degree $0$ and $\Ga \subset SL^{\pm}(\ste)$, $\mu(\psi)$ is the lift of a positive function on $M$.
Lemma \ref{centroaffineproposition} below expresses $\mu(\psi)$, the centroaffine curvature of $\Sigma_{\psi}$, in Euclidean terms.

\begin{definition}\label{sbetadefinition}
For $\psi \in \cAs(C, \Ga)$ and $\be \in \rea$, let
\begin{align}\label{sbetadefined}
\cS_{\be}(\psi) = \int_{M}\mu(\psi)^{\be}\psi^{-N} = N\int_{K_{\psi}/\Ga}\mu(\psi)^{\be}\,dx,
\end{align}
where the second equality is Lemma \ref{conevolumelemma}.
\end{definition}
Since $\mu(s\psi) = s^{2N}\mu(\psi)$, 
\begin{align}\label{csscaling}
\cS_{\be}(s\psi) = s^{N(2\be - 1)}\cS_{\be}(\psi). 
\end{align}
In particular $\cS_{1/2}$ and $\cS_{\be}\cS_{1-\be}$ are invariant under rescaling.

\begin{lemma}\label{chartlemma}
Let $\psi \in \cAs(C, \Ga)$. Choose linear coordinates $x = (z, t) \in \rea^{n}\times\rea$ on $\ste$ such that $\Psi = dz^{1}\wedge\dots\wedge dz^{n}\wedge dt$, let $\Om_{0} = \{z \in \rea^{n}: (z, 1) \in C\}$, an open convex set, and let $p(z) = \psi(z, 1)$ for $z \in \Om_{0}$. The map $z \mapsto [(z, 1)]$ identifies $\Om_{0}$ with the open subset $\{[x] \in \Om: t(x) > 0\}$ of $\Om$, and under this identification, at $(z, 1)$,
\begin{align}\label{chartformula}
&\det \fD^{2}\Phi = p^{n+2}\det \fD^{2}p,& &\mu(\psi)\,\psi^{-N} = p\,\det(-\fD^{2}p)\,dz,& &\psi^{-N} = p^{-N}dz,&
\end{align}
where the last two equalities are between densities on $\Om_{0}$.
\end{lemma}
\begin{proof}
Write $\pr_{i} = \pr/\pr z^{i}$, $p_{i} = \pr_{i}p$, and $p_{ij} = \pr_{i}\pr_{j}p$. At $x = (z, 1)$ the vectors $\pr_{1}, \dots, \pr_{n}, \eul_{x}$ form a parallel frame on $\ste$ such that $\Psi(\pr_{1}, \dots, \pr_{n}, \eul_{x}) = dt(\eul_{x}) = 1$. Since $\psi(z, 1) = p(z)$, $d\psi(\pr_{i}) = p_{i}$ and $\fD^{2}\psi(\pr_{i}, \pr_{j}) = p_{ij}$, and by homogeneity $d\psi(\eul) = \psi = p$ and $\fD^{2}\psi(\eul, \cdot) = 0$. Hence $\fD^{2}\Phi = \psi\fD^{2}\psi + d\psi\otimes d\psi$ has, in this basis, the matrix
\begin{align}
\begin{pmatrix} p\,p_{ij} + p_{i}p_{j} & p\,p_{i}\\ p\,p_{j} & p^{2}\end{pmatrix}, \qquad 1 \leq i, j \leq n.
\end{align}
Subtracting $p_{i}/p$ times the last row from the $i$th row for $1 \leq i \leq n$ leaves the determinant unchanged and gives a block lower triangular matrix with diagonal blocks $p\,p_{ij}$ and $p^{2}$, so the determinant is $p^{2}\det(p\fD^{2}p) = p^{n+2}\det \fD^{2}p$. Since $\fD^{2}p$ is negative definite, $|\det\fD^{2}p| = \det(-\fD^{2}p)$. The density $|\iota_{\eul}(F\Psi)|$ restricted to the slice $\{(z, t) \in \ste: t = 1\}$, which is identified with $\Om_{0}$, is $F(z, 1)\,dz$, and applying this with $F = \mu(\psi)\psi^{-N} = p^{n+2}\det(-\fD^{2}p)\,p^{-N}$ and with $F = \psi^{-N}$ gives the other two identities.
\end{proof}

Let $\Sigma \subset \ste$ be a smooth hypersurface and let $\xi$ be a \emph{transversal} along $\Sigma$, meaning a smooth vector field along $\Sigma$ such that $T_{x}\ste = T_{x}\Sigma\oplus\rea\xi_{x}$ for every $x \in \Sigma$. For vector fields $X$ and $Y$ tangent to $\Sigma$, the decompositions
\begin{align}\label{gaussformula}
\fD_{X}Y = \nabla^{\xi}_{X}Y + h_{\xi}(X, Y)\xi, \qquad \fD_{X}\xi = -S_{\xi}(X) + \tau_{\xi}(X)\xi,
\end{align}
with $\nabla^{\xi}_{X}Y$ and $S_{\xi}(X)$ tangent to $\Sigma$, define the \emph{induced connection} $\nabla^{\xi}$, a torsion-free affine connection on $\Sigma$, the \emph{second fundamental form} $h_{\xi}$, a symmetric bilinear form on $T\Sigma$, the \emph{shape operator} $S_{\xi}$, and the one-form $\tau_{\xi}$. The \emph{induced volume form} $\theta_{\xi}$ is the restriction to $T\Sigma$ of the interior multiplication $\iota_{\xi}\Psi$. Differentiating $\theta_{\xi}$ with \eqref{gaussformula} and the parallelism of $\Psi$ gives $\nabla^{\xi}_{X}\theta_{\xi} = \tau_{\xi}(X)\theta_{\xi}$. The transversal $\xi$ is \emph{equiaffine} if $\tau_{\xi} = 0$, that is, if $\fD_{X}\xi$ is tangent to $\Sigma$ for every tangent $X$, equivalently if $\nabla^{\xi}$ preserves $\theta_{\xi}$. For a positive function $\la$ on $\Sigma$, $h_{\la\xi} = \la^{-1}h_{\xi}$ and $\theta_{\la\xi} = \la\theta_{\xi}$, and for a tangent vector field $Z$, $h_{\xi + Z} = h_{\xi}$. The hypersurface $\Sigma$ is said to be \emph{locally strongly convex} if $h_{\xi}$ is definite for some, equivalently for every, transversal $\xi$. Replacing $\xi$ by $-\xi$ if necessary, $h_{\xi}$ is then positive definite. 

For $\psi \in \cAs(C, \Ga)$ and $\Sigma = \Sigma_{\psi}$, since $d\psi(\eul) = \psi = 1$ on $\Sigma$, the Euler field $\eul$ is transverse to $\Sigma$ and $T\Sigma = \ker d\psi|_{\Sigma}$. For vector fields $X$ and $Y$ tangent to $\Sigma$, applying $X$ to $d\psi(Y) = 0$ gives $\fD^{2}\psi(X, Y) + d\psi(\fD_{X}Y) = 0$, and applying $d\psi$ to \eqref{gaussformula} with $\xi = \eul$ gives $d\psi(\fD_{X}Y) = h_{\eul}(X, Y)$. Hence the second fundamental form of $\Sigma$ relative to the transversal $\eul$ is
\begin{align}\label{hxdefined}
h_{\eul} = -\fD^{2}\psi|_{T\Sigma},
\end{align}
which is positive definite because $\fD^{2}\psi$ is negative semidefinite with radical spanned by $\eul$. Thus $\Sigma_{\psi}$ is locally strongly convex, and $h_{\eul}$ is a Riemannian metric on $\Sigma_{\psi}$, the \emph{centroaffine metric}, denoted $\gb_{\psi}$. Via the identification $\Sigma_{\psi}/\Ga \simeq M$ by radial projection, $\gb_{\psi}$ is a Riemannian metric on $M$. The restriction of $G_{\psi}$ to $\Sigma_{\psi}$ is the \emph{conormal map} $x \mapsto d\psi(x)$, which associates with $x \in \Sigma_{\psi}$ the covector that annihilates $T_{x}\Sigma_{\psi}$ and equals $1$ on $\eul_{x}$. Its composition with the projection of $C^{\ast}$ onto $\Om^{\ast}$ is the \emph{conormal Gauss map} of $\Sigma_{\psi}$, which, by Lemma \ref{polarmaplemma}\ref{pml3}, is a diffeomorphism of $\Sigma_{\psi}$ onto $\Om^{\ast}$. The pullback by the conormal Gauss map of the flat projective structure of $\Om^{\ast}$ is a flat projective structure on $\Sigma_{\psi}$ that, like $\gb_{\psi}$, depends only on $\Sigma_{\psi}$ and the vector space structure of $\ste$, and not on the volume form $\Psi$ \cite{Fox-coupled}.

The next lemma shows that the class $\cAs(C, \Ga)$ consists exactly of the cogauges of the $\Ga$-convex bodies with smooth, locally strongly convex boundary.

\begin{lemma}\label{smoothboundarylemma}
Let $L$ be a $\Ga$-convex body with cogauge $\psi = \psi_{L} \in \cA(C, \Ga)$, and suppose that $\pr L = \Sigma_{\psi}$ is a smooth hypersurface. Then $\eul$ is transverse to $\Sigma_{\psi}$, $\psi$ is smooth on $C$, $\fD^{2}\psi$ is negative semidefinite with radical containing $\eul$, and the second fundamental form of $\Sigma_{\psi}$ relative to $\eul$ is $h_{\eul} = -\fD^{2}\psi|_{T\Sigma_{\psi}}$. Consequently $\psi \in \cAs(C, \Ga)$ if and only if $\Sigma_{\psi}$ is locally strongly convex.
\end{lemma}
\begin{proof}
Let $x \in \Sigma_{\psi}$. Since $\Sigma_{\psi}$ is smooth, the tangent hyperplane $x + T_{x}\Sigma_{\psi}$ is the unique supporting hyperplane of the convex set $L$ at $x$. Let $\zeta \in \std$ be a nonzero linear functional with kernel $T_{x}\Sigma_{\psi}$, with the sign chosen so that $\zeta(y) \geq \zeta(x)$ for all $y \in L$, so that the tangent hyperplane is $\{y \in \ste: \zeta(y) = \zeta(x)\}$. Then $\zeta$ is bounded from below on $L$, so $\zeta \in \bar{C}^{\ast}$ by Lemma \ref{gammaconvexbodylemma}, and $\zeta(x) > 0$ because $x \in C$ and $\zeta \neq 0$. Hence $\eul_{x} = x \notin \ker\zeta = T_{x}\Sigma_{\psi}$, which shows the claimed transversality.

The map $(s, y) \mapsto sy$ from $(0, \infty)\times\Sigma_{\psi}$ to $C$ is a bijection, because $\Sigma_{\psi}$ meets each ray of $C$ exactly once, and it is a local diffeomorphism by transversality, so it is a diffeomorphism, and so $\psi$ is smooth because $\psi(sy) = s$. By the concavity of $\psi$, $\fD^{2}\psi \leq 0$, and by its homogeneity, $\fD^{2}\psi(\eul, \cdot) = 0$. The formula for $h_{\eul}$ follows as in the derivation of \eqref{hxdefined}, which used only that $\psi$ is smooth, homogeneous of degree $1$, and equal to $1$ on $\Sigma_{\psi}$.

Finally, at $x \in \Sigma_{\psi}$ the radical of the negative semidefinite form $\fD^{2}\psi(x)$ is exactly the span of $\eul_{x}$ if and only if $\fD^{2}\psi(x)$ is negative definite on the complementary subspace $T_{x}\Sigma_{\psi} = \ker d\psi(x)$, that is, if and only if $h_{\eul}$ is positive definite at $x$. By homogeneity the same holds at every point of the ray through $x$. Since $h_{\eul} \geq 0$, definiteness of $h_{\eul}$ means positive definiteness, which gives the last claim.
\end{proof}

\begin{lemma}\label{centroaffinemetriclemma}
Let $\psi \in \cAs(C, \Ga)$ and write $\mu = \mu(\psi)$. Then:
\begin{enumerate}[label=$(\arabic*)$, ref=(\arabic*)]
\item\label{cml1} $\mu^{1/2}\psi^{-N}$ is the Riemannian density of $\gb_{\psi}$. In particular $\cS_{1/2}(\psi) = \vol(M, \gb_{\psi})$.
\item\label{cml2} For a $C^{1}$ function $f$ on $M$, both sides of the following identity being $\Ga$-invariant functions homogeneous of degree $0$,
\begin{align}\label{dualmetricadjugate}
|df|^{2}_{\gb_{\psi}} = -\psi^{2}(\det\fD^{2}\Phi)^{-1}\adj(\fD^{2}\Phi)(df, df) = (-1)^{n+1}\mu^{-1}\psi^{2}\adj(\fD^{2}\Phi)(df, df).
\end{align}
\item\label{cml3} The restriction of $G_{\psi}$ to $\Sigma_{\psi}$, which is the conormal map $x \mapsto d\psi(x)$, is an isometry from $(\Sigma_{\psi}, \gb_{\psi})$ onto $(\Sigma_{\psi^{\ast}}, \gb_{\psi^{\ast}})$.
\item\label{cml4} With the determinant computed with respect to $\Psi$ and with $\gb_{\psi}$ pulled back to $C$ by radial projection onto $\Sigma_{\psi}$,
\begin{align}\label{mulogdet}
&\mu(\psi) = \psi^{2N}\det(-\fD^{2}\log\psi),& &-\fD^{2}\log\psi = (d\log\psi)^{2} + \gb_{\psi}.&
\end{align}
\end{enumerate}
\end{lemma}
\begin{proof}
At $x \in \Sigma_{\psi}$ choose a basis $\eul_{x}, e_{1}, \dots, e_{n}$ with $e_{i} \in T_{x}\Sigma_{\psi} = \ker d\psi(x)$. In this basis $\fD^{2}\Phi$ is block diagonal with entries $\psi^{2} = 1$ and $\psi \fD^{2}\psi(e_{i}, e_{j}) = -\gb_{\psi}(e_{i}, e_{j})$, so
\begin{align}\label{mucentroaffine}
&\mu = \det(\gb_{\psi}(e_{i}, e_{j}))\,\Psi(\eul_{x}, e_{1}, \dots, e_{n})^{-2},& &\text{that is,}& &\mu^{1/2}\big|(\iota_{\eul}\Psi)|_{T\Sigma_{\psi}}\big| = \vol_{\gb_{\psi}}.&
\end{align} Since $\psi = 1$ on $\Sigma_{\psi}$ and the density $\psi^{-N}$ is represented on $\Sigma_{\psi}$ by $|\iota_{\eul}\Psi|$, \ref{cml1} follows. Next \ref{cml2} is proved. Since $f$ is homogeneous of degree $0$, $\al = df$ annihilates $\eul$. Write $\al = \fD^{2}\Phi(v, \cdot)$. Then $\fD^{2}\Phi(v, \eul) = \al(\eul) = 0$, and $\fD^{2}\Phi(\cdot, \eul) = \psi d\psi$, so $d\psi(v) = 0$ and $\al = \psi \fD^{2}\psi(v, \cdot) = -\gb_{\psi}(v, \cdot)$ on $T\Sigma_{\psi}$ at points where $\psi = 1$. Since $\al = \fD^{2}\Phi(v, \cdot)$, $\adj(\fD^{2}\Phi)(\al, \cdot) = (\det\fD^{2}\Phi)\,v$, so $(\det\fD^{2}\Phi)^{-1}\adj(\fD^{2}\Phi)(\al, \al) = \al(v) = -\gb_{\psi}(v, v) = -|\al|^{2}_{\gb_{\psi}}$. Both sides of \eqref{dualmetricadjugate} are homogeneous of degree $0$, which gives \eqref{dualmetricadjugate} everywhere, and the second expression follows from $\sign\det\fD^{2}\Phi = (-1)^{n}$ (Lemma \ref{polarmaplemma}\ref{pml1}).
Next \ref{cml3} is proved. On $\Sigma_{\psi}$, $G_{\psi} = \psi\,d\psi = d\psi$, and by Lemma \ref{polarmaplemma}\ref{pml3}--\ref{pml4}, $G_{\psi}$ maps $\Sigma_{\psi} = \{\psi = 1\}$ diffeomorphically onto $\Sigma_{\psi^{\ast}} = \{\psi^{\ast} = 1\}$. Its differential at $x \in \Sigma_{\psi}$ is the linear map $\ste \to \std$, $v \mapsto \fD^{2}\Phi(x)(v, \cdot)$, which in the basis above sends $e_{i}$ to $\fD^{2}\Phi(e_{i}, \cdot) = -\gb_{\psi}(e_{i}, \cdot)$, a covector annihilating $\eul_{x}$. By Lemma \ref{polarmaplemma}\ref{pml4}, $\fD^{2}\Phi^{\ast}$ at $G_{\psi}(x)$ is the bilinear form on $\std$ inverse to $\fD^{2}\Phi(x)$, so in the dual basis it is block diagonal with blocks $1$ and $-\gb_{\psi}^{-1}$, and since $\psi^{\ast} = 1$ at $G_{\psi}(x)$, by \eqref{hxdefined} applied to $\psi^{\ast}$ the centroaffine metric $\gb_{\psi^{\ast}}$ is $-\fD^{2}\Phi^{\ast}$ restricted to $T\Sigma_{\psi^{\ast}}$. Hence $\gb_{\psi^{\ast}}(\fD^{2}\Phi(e_{i}, \cdot), \fD^{2}\Phi(e_{j}, \cdot)) = \gb_{\psi}^{-1}(\gb_{\psi}(e_{i}, \cdot), \gb_{\psi}(e_{j}, \cdot)) = \gb_{\psi}(e_{i}, e_{j})$, which is \ref{cml3}.

Finally \ref{cml4} is proved. At $x \in C$ choose a basis $(\eul, e_{1}, \dots, e_{n})$ with $e_{j} \in \ker d\psi(x)$. Since $\fD^{2}\psi(\eul, \cdot) = 0$ and $d\psi(\eul) = \psi$, the forms $\fD^{2}\Phi = \psi\,\fD^{2}\psi + d\psi\otimes d\psi$ and $-\fD^{2}\log\psi = \psi^{-2}d\psi\otimes d\psi - \psi^{-1}\fD^{2}\psi$ are block diagonal, with blocks $\psi^{2}$ and $\psi\,\fD^{2}\psi|_{\ker d\psi}$, and $1$ and $-\psi^{-1}\fD^{2}\psi|_{\ker d\psi}$, respectively. Comparing determinants gives the first identity of \eqref{mulogdet}, because $\mu(\psi) = (-1)^{n}\det\fD^{2}\Phi$. At points of $\Sigma_{\psi}$ the second block is $\gb_{\psi}$ by \eqref{hxdefined}, and both sides of the second identity of \eqref{mulogdet} are invariant under the dilations of $C$, so it holds on $C$.
\end{proof}

For a locally strongly convex $\Sigma$ there is a unique equiaffine transversal $\xi$ such that $h_{\xi}$ is positive definite and $|\theta_{\xi}| = \vol_{h_{\xi}}$ \cite[Chapter II, Theorem $3.1$]{Nomizu-Sasaki}. This $\xi$ is the \emph{equiaffine normal} of $\Sigma$, called the affine normal in \cite{Nomizu-Sasaki}, and the corresponding $h = h_{\xi}$, $\nabla = \nabla^{\xi}$, and $S = S_{\xi}$ are the \emph{equiaffine metric}, also called the Blaschke metric, the \emph{equiaffine connection}, and the \emph{equiaffine shape operator}. The \emph{affine mean curvature} is $\tfrac{1}{n}\tr S$. The qualifier equiaffine reflects the dependence of these objects on the volume form $\Psi$. Since $\theta_{\xi}$ is linear in $\Psi$ while $\vol_{h_{\la\xi}} = \la^{-n/2}\vol_{h_{\xi}}$ for constant $\la > 0$, replacing $\Psi$ by $c\Psi$, $c > 0$, replaces $\xi$ by $c^{-2/(n+2)}\xi$, $h$ by $c^{2/(n+2)}h$, and $S$ by $c^{-2/(n+2)}S$, while $\nabla$ is unchanged.

These objects are invariant under the group of affine transformations of $\ste$ preserving $\Psi$ up to sign, but not under $GL(\ste)$. With respect to the equiaffine metric $h$, the equiaffine connection $\nabla$ has a conjugate connection $\bar{\nabla}$, defined by $Xh(Y, Z) = h(\nabla_{X}Y, Z) + h(Y, \bar{\nabla}_{X}Z)$ for vector fields $X$, $Y$, and $Z$ tangent to $\Sigma$, and, when $\Sigma = \Sigma_{\psi}$, $\bar{\nabla}$ represents the flat projective structure induced on $\Sigma_{\psi}$ by its conormal Gauss map \cite{Fox-coupled}.

\begin{definition}\label{affinespheredefinition}
A locally strongly convex hypersurface with equiaffine normal $\xi$ is a \emph{hyperbolic affine sphere centered at the origin} if $\xi = \la\eul$ along the hypersurface for a positive function $\la$.
\end{definition}

\begin{definition}\label{asymptoticdefinition}
A hyperbolic affine sphere $\Sigma \subset C$ centered at the origin is \emph{asymptotic to $\pr C$} if $\Sigma$ meets every ray of $C$ and is the boundary of a closed convex subset of $\ste$ that is contained in $C$. This subset is called the closed convex subset of $C$ \emph{bounded by} $\Sigma$.
\end{definition}
The first of the two conditions of Definition \ref{asymptoticdefinition} cannot be omitted. If $C' \subsetneq C$ is an open properly convex cone, then a hyperbolic affine sphere centered at the origin and asymptotic to $\pr C'$ is the boundary of a closed convex subset of $\ste$ contained in $C$, but it does not meet the rays of $C$ not contained in $C'$.

Let $\cAs(C, \{1\})$ be the set of smooth positive concave functions on $C$, homogeneous of degree $1$, whose Hessian is negative semidefinite with radical spanned by $\eul$, that is, the set defined in the same way as $\cAs(C, \Ga)$ with the trivial group in place of $\Ga$. For $\psi \in \cAs(C, \{1\})$, the sets $\Sigma_{\psi} = \{x \in C: \psi(x) = 1\}$ and $L_{\psi} = \{x \in C: \psi(x) \geq 1\}$, the function $\mu(\psi)$ of \eqref{mudefined}, and the centroaffine metric of $\Sigma_{\psi}$ are defined as for elements of $\cAs(C, \Ga)$, and Lemmas \ref{chartlemma} and \ref{affinespherelemma}, whose proofs are pointwise, apply to $\psi$. Lemma \ref{asymptoticlemma} describes the affine spheres asymptotic to $\pr C$ in terms of such functions.
Asymptotic is used in this sense in \cite{Cheng-Yau-affinespheres, Loftin-affinespheres}. By Lemma \ref{gammaconvexbodylemma}, for $\psi \in \cA(C, \Ga)$ the set $L_{\psi}$ is a closed convex subset of $\ste$ contained in $C$, and its boundary is $\Sigma_{\psi}$. Hence, for $\psi \in \cAs(C, \Ga)$ such that $\Sigma_{\psi}$ is a hyperbolic affine sphere centered at the origin, $\Sigma_{\psi}$ is asymptotic to $\pr C$, and $L_{\psi}$ is the closed convex subset of $C$ bounded by $\Sigma_{\psi}$.

\begin{lemma}\label{affinespherelemma}
For $\psi \in \cAs(C, \Ga)$ and $\Sigma = \Sigma_{\psi}$, with equiaffine normal $\xi$, the following are equivalent.
\begin{enumerate}[label=$(\arabic*)$, ref=(\arabic*)]
\item\label{asl1} $\Sigma$ is a hyperbolic affine sphere centered at the origin.
\item\label{asl2} $\xi = \la\eul$ along $\Sigma$ for a constant $\la > 0$.
\item\label{asl3} $\mu(\psi)$ is constant.
\item\label{asl4} The equiaffine metric of $\Sigma$ is a constant multiple of the centroaffine metric $\gb_{\psi}$.
\end{enumerate}
If \ref{asl1}--\ref{asl4} hold, then $\mu(\psi) = \la^{n+2}$, the equiaffine metric of $\Sigma$ is $h = \la^{-1}\gb_{\psi} = -\la^{-1}\fD^{2}\psi|_{T\Sigma}$, its equiaffine connection is the induced connection $\nabla^{\eul}$ of the transversal $\eul$ and preserves the volume form $\theta = \la\,\iota_{\eul}\Psi$, and its affine mean curvature is $-\la$.
\end{lemma}
\begin{proof}
That \ref{asl2} implies \ref{asl1} is immediate from Definition \ref{affinespheredefinition}.
If \ref{asl1} holds, the equiaffine normal $\xi$ of $\Sigma$ has the form $\xi = \la \eul$ for a positive function $\la$, so $\fD_{X}\xi = X(\la)\eul + \la X$ is tangent because $\xi$ is equiaffine, implying that $X(\la) = 0$ and $\la$ is locally constant. Radial projection identifies $\Sigma$ with the connected set $\Om$, so $\la$ is constant, which is \ref{asl2}. 

Since $\fD_{X}\eul = X$ for every vector $X$, the transversal $\eul$ is equiaffine, with shape operator $S_{\eul} = -\mathrm{id}$. Since $\fD^{2}\psi(\eul, \cdot) = 0$ by homogeneity and $d\psi|_{T\Sigma} = 0$, the bilinear form $\fD^{2}\Phi = \psi\fD^{2}\psi + d\psi\otimes d\psi$ has, along $\Sigma$ and in a basis $X_{1}, \dots, X_{n}, \eul$ with $X_{i}$ tangent, the block form $\mathrm{diag}(-h_{\eul}, 1)$. For a transversal $\eta$ with $h_{\eta}$ positive definite, let
\begin{align}\label{dettheta}
\det\nolimits_{\theta_{\eta}}h_{\eta} = \frac{\det(h_{\eta}(X_{i}, X_{j}))}{\theta_{\eta}(X_{1}, \dots, X_{n})^{2}},
\end{align}
which does not depend on the basis $X_{1}, \dots, X_{n}$ of $T\Sigma$, and equals $1$ exactly when $|\theta_{\eta}| = \vol_{h_{\eta}}$. By the definition \eqref{mudefined} of $\mu$ and the block form of $\fD^{2}\Phi$,
\begin{align}\label{mudettheta}
\mu(\psi) = \det\nolimits_{\theta_{\eul}}h_{\eul}.
\end{align}
Since the equiaffine normal $\xi$ is the unique equiaffine transversal $\eta$ with $h_{\eta}$ positive definite and $|\theta_{\eta}| = \vol_{h_{\eta}}$, an equiaffine transversal $\eta$ with $h_{\eta}$ positive definite equals $\xi$ if and only if $\det\nolimits_{\theta_{\eta}}h_{\eta} = 1$. For a constant $\la > 0$, the transversal $\la\eul$ is equiaffine, with $h_{\la\eul} = \la^{-1}h_{\eul}$ positive definite and $\theta_{\la\eul} = \la\theta_{\eul}$, so by \eqref{dettheta} and \eqref{mudettheta}
\begin{align}\label{blaschkecondition}
\det\nolimits_{\theta_{\la\eul}}h_{\la\eul} = \det\nolimits_{\la\theta_{\eul}}(\la^{-1}h_{\eul}) = \la^{-n-2}\det\nolimits_{\theta_{\eul}}h_{\eul} = \la^{-n-2}\mu(\psi),
\end{align}
and hence $\xi = \la\eul$ if and only if $\mu(\psi) = \la^{n+2}$. If \ref{asl2} holds, this shows that $\mu(\psi) = \la^{n+2}$ is constant, which is \ref{asl3}. If \ref{asl3} holds, then $\la = \mu(\psi)^{1/(n+2)}$ is a positive constant with $\mu(\psi) = \la^{n+2}$, so $\xi = \la\eul$, which is \ref{asl2}. 

If \ref{asl2} holds, then $h_{\xi} = h_{\la\eul} = \la^{-1}h_{\eul} = \la^{-1}\gb_{\psi}$, which is \ref{asl4}. Suppose \ref{asl4} holds, $h_{\xi} = c\,\gb_{\psi}$ for a constant $c > 0$. Since $\eul$ is transverse to $\Sigma$, $\xi = \la\eul + Z$ for a positive function $\la$ and a tangent vector field $Z$, and $h_{\xi} = h_{\la\eul} = \la^{-1}h_{\eul}$, so $\la = c^{-1}$ is constant. Since $\xi$ is equiaffine, $\fD_{X}\xi = \la X + \fD_{X}Z$ is tangent for every tangent $X$, while by \eqref{gaussformula} with the transversal $\eul$, $\fD_{X}Z = \nabla^{\eul}_{X}Z + h_{\eul}(X, Z)\eul$, so $h_{\eul}(X, Z) = 0$ for every tangent $X$, and $Z = 0$ because $h_{\eul}$ is nondegenerate. Hence $\xi = \la\eul$ with $\la$ constant, which is \ref{asl2}.

The final claims follow from $h = h_{\la\eul} = \la^{-1}h_{\eul}$, $\nabla^{\la\eul} = \nabla^{\eul}$ for constant $\la$, $\theta = \theta_{\la\eul} = \la\theta_{\eul}$, and $S = S_{\la\eul} = \la S_{\eul} = -\la\,\mathrm{id}$, because the affine mean curvature is $\tfrac{1}{n}\operatorname{tr}S$.
\end{proof}

For an affine sphere, $\gb_{\psi}$ is the metric denoted $h$ in Theorem A, and when $\mu(\psi) = 1$ it is the equiaffine metric of $\Sigma_{\psi}$, the equiaffine metric of the domain $\Om$.
For every $\psi \in \cAs(C, \Ga)$, by \eqref{mulogdet} the Riemannian Hessian metric $-\fD^{2}\log\psi$ on $C$ is the Riemannian product of a line with $(\Sigma_{\psi}, \gb_{\psi})$ \cite[Theorem $1$]{Loftin-KE}. When $\mu(\psi)$ is constant, \eqref{mulogdet} also shows that $F = -N\log\psi$ satisfies $\det\fD^{2}F = N^{N}\mu(\psi)e^{2F}$, so that, after the addition of a constant, $F$ is the canonical potential of $C$ \cite{Fox-schwarz}, and $-\fD^{2}\log\psi = N^{-1}\fD^{2}F$ determines a K\"ahler--Einstein metric on the tube domain over $C$ \cite[Theorem $3$]{Loftin-KE}; see also \cite[Theorem $4.1$]{Fox-schwarz}. 
In an affine chart, with $u = -\psi(z, 1)$, the equation $\mu(\psi) = 1$ is $\det \fD^{2}u = (-u)^{-n-2}$ by Lemma \ref{chartlemma}, and radial projection identifies $\gb_{\psi}$ with $-u^{-1}\fD^{2}u$. Indeed, at $x = (z, 1)/p(z) \in \Sigma_{\psi}$, with $p(z) = \psi(z, 1)$, the vectors $v_{i} = \tfrac{\pr}{\pr z^{i}}((z,1)/p(z))$ span $T_{x}\Sigma_{\psi}$, $\fD^{2}\psi$ annihilates the radial direction and is homogeneous of degree $-1$, so $-\fD^{2}\psi_{x}(v_{i}, v_{j}) = -p^{-2}\,p\,\fD^{2}p(e_{i}, e_{j}) = -p^{-1}\fD^{2}p(e_{i}, e_{j}) = -u^{-1}\fD^{2}u(e_{i}, e_{j})$.

By Lemma \ref{affinespherelemma}, $\Sigma_{\psi}$ is a hyperbolic affine sphere centered at the origin if and only if $\mu(\psi)$ is constant. By Lemma \ref{chartlemma}, this means that $u = -p$ solves $\det \fD^{2}u = \ka(-u)^{-n-2}$ in every affine chart, with $\ka$ constant, which is the equation of hyperbolic affine spheres of Calabi \cite{Calabi-completeaffine} and Cheng--Yau \cite{Cheng-Yau-affinespheres}. The equivalence of this equation with Definition \ref{affinespheredefinition}, reproved in the present notation by Lemma \ref{affinespherelemma}, is classical; see \cite[Section $5$, (5.12)--(5.13)]{Cheng-Yau-affinespheres}, \cite[Theorem $3$]{Loftin-KE}, and \cite{Loftin-affinespheres}. Lemma \ref{uniquenessproposition} shows that two functions in $\cAs(C, \Ga)$ with this property differ by a constant factor, so that such an affine sphere, if it exists, is unique up to homothety. The proof is an application of the maximum principle on the compact quotient $M$.
\begin{lemma}\label{uniquenessproposition}
If $\psi_{1}, \psi_{2} \in \cAs(C, \Ga)$ satisfy $\mu(\psi_{1}) = \ka_{1}$ and $\mu(\psi_{2}) = \ka_{2}$ with $\ka_{1}, \ka_{2}$ constant, then $\psi_{2}$ is a constant multiple of $\psi_{1}$.
\end{lemma}

\begin{proof}
The quotient $f = \psi_{2}/\psi_{1}$ is homogeneous of degree $0$ and $\Ga$-invariant, so it is identified with a positive smooth function on the compact manifold $M$, on which it attains its maximum $f_{\max}$ and minimum $f_{\min}$. At a point $z_{0}$ of $\Om$ over a maximum point, take an affine chart as in Lemma \ref{chartlemma} and let $p_{i}(z) = \psi_{i}(z, 1)$ for $i = 1,2$, so that $p_{2} = fp_{1}$. At $z_{0}$, $\fD f = 0$ and $\fD^{2}f \leq 0$, so $-\fD^{2}p_{2} = -f\fD^{2}p_{1} - p_{1}\fD^{2}f \geq -f\fD^{2}p_{1} > 0$, and by Lemma \ref{chartlemma}
\begin{align}
\ka_{2} = p_{2}^{n+2}\det(-\fD^{2}p_{2}) \geq f^{n+2}p_{1}^{n+2}f^{n}\det(-\fD^{2}p_{1}) = f_{\max}^{2N}\ka_{1}.
\end{align}
At a minimum point the same argument gives $\ka_{2} \leq f_{\min}^{2N}\ka_{1}$. Hence $f_{\max}^{2N} \leq \ka_{2}/\ka_{1} \leq f_{\min}^{2N}$, so $f$ is constant.
\end{proof}

Lemmas \ref{asymptoticlemma} and \ref{generaluniquenesslemma} show that the uniqueness of Lemma \ref{uniquenessproposition} holds without the group, for affine spheres asymptotic to $\pr C$, and that such an affine sphere is invariant under every unimodular linear transformation preserving $C$.

\begin{lemma}\label{asymptoticlemma}
Let $C \subset \ste$ be an open properly convex cone.
\begin{enumerate}[label=$(\arabic*)$, ref=(\arabic*)]
\item\label{asyl1} For $\psi \in \cAs(C, \{1\})$, the set $L_{\psi}$ is closed in $\ste$ if and only if $\psi(x) \to 0$ as $x \to \pr C\setminus\{0\}$.
\item\label{asyl2} A hypersurface $\Sigma \subset C$ is a hyperbolic affine sphere centered at the origin and asymptotic to $\pr C$ if and only if $\Sigma = \Sigma_{\psi}$ for a $\psi \in \cAs(C, \{1\})$ such that $\mu(\psi)$ is constant and $L_{\psi}$ is closed in $\ste$. In this case $L_{\psi}$ is the closed convex subset of $C$ bounded by $\Sigma$.
\end{enumerate}
\end{lemma}
\begin{proof}
First \ref{asyl1} is proved. Suppose $\psi(x) \to 0$ as $x \to \pr C\setminus\{0\}$. Then $\psi$ extends to a continuous function on $\bar{C}\setminus\{0\}$ vanishing on $\pr C\setminus\{0\}$, and, being homogeneous of degree $1$ and bounded on the compact intersection of $\bar{C}$ with the unit sphere of a Euclidean norm $|\cdot|$ on $\ste$, it satisfies $\psi(x) \leq c|x|$ for a constant $c$ and all $x \in C$. Let $x_{k} \in L_{\psi}$ converge to $x \in \ste$. Then $x \in \bar{C}$, and $|x| \geq c^{-1}$ because $c|x_{k}| \geq \psi(x_{k}) \geq 1$. If $x$ were in $\pr C\setminus\{0\}$, then $\psi(x_{k})$ would tend to $0$. Hence $x \in C$, and $\psi(x) \geq 1$ by continuity, so $x \in L_{\psi}$. Conversely, suppose that $L_{\psi}$ is closed in $\ste$ and that there are $\delta > 0$ and $x_{k} \in C$ converging to $x \in \pr C\setminus\{0\}$ such that $\psi(x_{k}) \geq \delta$. Then $\psi(\delta^{-1}x_{k}) \geq 1$, so the points $\delta^{-1}x_{k}$ of $L_{\psi}$ converge to $\delta^{-1}x$, which is not in $C$, a contradiction.

Next \ref{asyl2} is proved. Let $\psi \in \cAs(C, \{1\})$ be such that $\mu(\psi)$ is constant and $L_{\psi}$ is closed in $\ste$. By Lemma \ref{affinespherelemma}, $\Sigma_{\psi}$ is a hyperbolic affine sphere centered at the origin. It meets every ray of $C$, because $\psi$ is positive and homogeneous of degree $1$, and it is the boundary of $L_{\psi}$, which is convex because $\psi$ is concave, and contained in $C$.

Conversely, let $\Sigma$ be a hyperbolic affine sphere centered at the origin and asymptotic to $\pr C$, and let $K$ be the closed convex subset of $C$ bounded by $\Sigma$. Because the second fundamental form of $\Sigma$ relative to its equiaffine normal $\xi = \la\eul$ is positive definite, $K$ lies, near each point of $\Sigma$, on the side of $\Sigma$ into which $\xi$ points, so, for $y \in \Sigma$ and all sufficiently small $s > 0$, $(1 + s)y \in K$ and $(1 - s)y \notin K$. For $x \in C$ let $R_{x} = \{t > 0: tx \in K\}$. It is nonempty, because the ray through $x$ meets $\Sigma \subset K$, it is an interval, because $K$ is convex, and it is a closed subset of $(0, \infty)$ bounded away from $0$, because $K$ is closed in $\ste$ and does not contain $0$. If $t \in R_{x}$ and $tx \in \Sigma$, then, by the preceding, $t$ is the left endpoint of $R_{x}$ and is not its right endpoint. Since an endpoint $t$ of $R_{x}$ satisfies $tx \in \pr K = \Sigma$, it follows that $R_{x} = [a(x), \infty)$ for some $a(x) > 0$ and that $a(x)x$ is the only point of $\Sigma$ on the ray through $x$.

Let $\psi(x) = a(x)^{-1}$. Then $\psi$ is positive and homogeneous of degree $1$, $K = L_{\psi}$, and $\Sigma = \Sigma_{\psi}$. For $x, y \in C$, the points $\psi(x)^{-1}x$ and $\psi(y)^{-1}y$ lie in $K$, so their convex combination $(\psi(x) + \psi(y))^{-1}(x + y)$ lies in $K$, which means that $\psi(x + y) \geq \psi(x) + \psi(y)$. With the homogeneity of $\psi$, this shows that $\psi$ is concave. Since $\eul$ is transverse to $\Sigma$, the map $(s, y) \mapsto sy$ from $(0, \infty)\times\Sigma$ to $C$ is a diffeomorphism, and $\psi(sy) = s$, so $\psi$ is smooth. The remaining arguments of the proof of Lemma \ref{smoothboundarylemma}, which are pointwise, show that the radical of $\fD^{2}\psi$ is spanned by $\eul$, because $\Sigma$ is locally strongly convex. Hence $\psi \in \cAs(C, \{1\})$, and $\mu(\psi)$ is constant by Lemma \ref{affinespherelemma}.
\end{proof}

\begin{lemma}\label{generaluniquenesslemma}
Let $\psi_{1}, \psi_{2} \in \cAs(C, \{1\})$ satisfy $\mu(\psi_{1}) = \ka_{1}$ and $\mu(\psi_{2}) = \ka_{2}$ with $\ka_{1}, \ka_{2}$ constant, and $\psi_{i}(x) \to 0$ as $x \to \pr C\setminus\{0\}$ for $i = 1, 2$. Then $\psi_{2} = (\ka_{2}/\ka_{1})^{1/(2N)}\psi_{1}$. Consequently:
\begin{enumerate}[label=$(\arabic*)$, ref=(\arabic*)]
\item\label{gul1} A hyperbolic affine sphere in $C$ centered at the origin and asymptotic to $\pr C$, if there is one, is unique up to homothety.
\item\label{gul2} Such an affine sphere is mapped to itself by every $g \in SL^{\pm}(\ste)$ such that $gC = C$. In particular, it is $\Ga$-invariant, and it is the level set $\Sigma_{\psi}$ of a $\psi \in \cAs(C, \Ga)$ with $\mu(\psi)$ constant.
\end{enumerate}
\end{lemma}
\begin{proof}
Since $\mu(c\psi) = c^{2N}\mu(\psi)$ for $c > 0$, after replacing $\psi_{2}$ by $(\ka_{1}/\ka_{2})^{1/(2N)}\psi_{2}$ it can be supposed that $\ka_{1} = \ka_{2} = \ka$, and it is to be shown that $\psi_{1} = \psi_{2}$. Choose linear coordinates $(z, t)$ on $\ste$ as in Lemma \ref{chartlemma} with $t \in C^{\ast}$, so that the slice $\Om_{0} = \{z \in \rea^{n}: (z, 1) \in C\}$ is bounded and meets every ray of $C$, and let $p_{i}(z) = \psi_{i}(z, 1)$. By \eqref{chartformula}, $p_{i}$ is a positive function on $\Om_{0}$ such that $-\fD^{2}p_{i}$ is positive definite and
\begin{align}\label{slicemongeampere}
\det(-\fD^{2}p_{i}) = \ka p_{i}^{-n-2},
\end{align}
and $p_{i}(z) \to 0$ as $z \to \pr\Om_{0}$, because $(z, 1) \to \pr C\setminus\{0\}$. Let $c > 1$ and suppose that $p_{1} - cp_{2}$ is positive somewhere in $\Om_{0}$. Since $p_{1} - cp_{2}$ extends continuously to the compact set $\bar{\Om}_{0}$ and vanishes on $\pr\Om_{0}$, it attains its maximum at a point $z_{0} \in \Om_{0}$, where $-\fD^{2}p_{1} \geq -c\,\fD^{2}p_{2} > 0$ and $p_{1} > cp_{2}$. By \eqref{slicemongeampere}, at $z_{0}$,
\begin{align}
\ka p_{1}^{-n-2} = \det(-\fD^{2}p_{1}) \geq c^{n}\det(-\fD^{2}p_{2}) = c^{n}\ka p_{2}^{-n-2} > c^{2n+2}\ka p_{1}^{-n-2},
\end{align}
which is impossible because $c > 1$. Hence $p_{1} \leq cp_{2}$ on $\Om_{0}$ for every $c > 1$, so $p_{1} \leq p_{2}$. Exchanging the roles of $p_{1}$ and $p_{2}$ gives $p_{1} = p_{2}$, and $\psi_{1} = \psi_{2}$ by homogeneity.

By Lemma \ref{asymptoticlemma}, the hyperbolic affine spheres in $C$ centered at the origin and asymptotic to $\pr C$ are the level sets $\Sigma_{\psi}$ of the $\psi \in \cAs(C, \{1\})$ such that $\mu(\psi)$ is constant and $\psi(x) \to 0$ as $x \to \pr C\setminus\{0\}$. Any two such functions are proportional, by the first part of the lemma, so their level sets are homothetic. This shows \ref{gul1}. Let $\psi$ be such a function and let $g \in SL^{\pm}(\ste)$ satisfy $gC = C$. Then $\psi\circ g \in \cAs(C, \{1\})$, $\psi\circ g$ tends to $0$ at $\pr C\setminus\{0\}$ because $g$ preserves $\pr C\setminus\{0\}$, and $\mu(\psi\circ g) = \mu(\psi)\circ g$ is the same constant as $\mu(\psi)$, because $g^{\ast}\Psi = \pm\Psi$ in \eqref{mudefined}. By the first part of the lemma, $\psi\circ g = \psi$, so $g$ preserves $\Sigma_{\psi}$. Applied to the elements of $\Ga$, this shows that $\psi \in \cAs(C, \Ga)$. This shows \ref{gul2}.
\end{proof}

The following lemma is likewise classical.
\begin{lemma}\label{gausslemma}
Let $\psi \in \cAs(C, \Ga)$ satisfy $\mu(\psi) = 1$, and let $\Sigma = \Sigma_{\psi}$, with equiaffine metric $h = \gb_{\psi}$ and equiaffine connection $\nabla$. Let $D$ be the Levi-Civita connection of $h$, and let $\dift = \nabla - D$ be the difference tensor.
\begin{enumerate}[label=$(\arabic*)$, ref=(\arabic*)]
\item\label{gl1} The \emph{cubic form} $C_{h}(X, Y, Z) = h(\dift_{X}Y, Z)$ is completely symmetric, and $\dift$ is \emph{apolar}, meaning that $\operatorname{tr}\dift_{X} = 0$ for every $X$.
\item\label{gl2} The Levi-Civita Ricci tensor and scalar curvature of $h$ satisfy
\begin{align}\label{gaussricci}
\begin{aligned}
\mathrm{Ric}_{h}(Y, Z) &= -(n-1)h(Y, Z) + \tr(\dift_{Y}\dift_{Z}) \geq -(n-1)h(Y, Z),\\
\mathrm{Scal}_{h} &= -n(n-1) + |\dift|^{2}_{h},
\end{aligned}
\end{align}
with equality in the first inequality for $Y$ at a point if and only if $\dift_{Y} = 0$ there.
\end{enumerate}
\end{lemma}
\begin{proof}
These facts are classical. See \cite[Chapter II]{Nomizu-Sasaki}, \cite{Calabi-completeaffine}, and \cite{Loftin-survey}. The proof is sketched to fix the sign conventions and normalizations. By Lemma \ref{affinespherelemma}, the equiaffine normal is $\xi = \eul$, the shape operator is $S = -\mathrm{id}$, and $\nabla$ preserves $\theta = \theta_{\eul} = \vol_{h}$. Since $\fD$ is flat and torsion free, differentiating the Gauss formula $\fD_{X}Y = \nabla_{X}Y + h(X, Y)\eul$ gives the Gauss and Codazzi equations
\begin{align}\label{gausscodazzi}
\begin{aligned}
R^{\nabla}(X, Y)Z &= h(Y, Z)SX - h(X, Z)SY = -h(Y, Z)X + h(X, Z)Y, \\
(\nabla_{X}h)(Y, Z) &= (\nabla_{Y}h)(X, Z).
\end{aligned}
\end{align}
Because $\nabla$ and $D$ are torsion free, $\dift_{X}Y = \dift_{Y}X$, and because $Dh = 0$, $(\nabla_{X}h)(Y, Z) = -C_{h}(X, Y, Z) - C_{h}(X, Z, Y)$. With the Codazzi equation these give that $C_{h}$ is completely symmetric and $\nabla h = -2C_{h}$. Since $\nabla$ and $D$ both preserve $\vol_{h}$, $\operatorname{tr}\dift_{X} = 0$. This proves \ref{gl1}.

Next \ref{gl2} is proved. By the complete symmetry of $C_{h}$, the connection $\bar{\nabla} = D - \dift$ satisfies $Xh(Z, W) = h(\nabla_{X}Z, W) + h(Z, \bar{\nabla}_{X}W)$, and so $h(R^{\nabla}(X, Y)Z, W) + h(Z, R^{\bar{\nabla}}(X, Y)W) = 0$. By the Gauss equation, $h(R^{\nabla}(X, Y)Z, W)$ is antisymmetric in $Z$ and $W$, so $R^{\bar{\nabla}} = R^{\nabla}$. Averaging the identities $R^{D \pm \dift}(X, Y) = R^{D}(X, Y) \pm \big((D_{X}\dift)_{Y} - (D_{Y}\dift)_{X}\big) + [\dift_{X}, \dift_{Y}]$ gives $R^{D}(X, Y) = R^{\nabla}(X, Y) - [\dift_{X}, \dift_{Y}]$. Taking the trace over $X$, with $R^{\nabla}$ given by \eqref{gausscodazzi} and $\operatorname{tr}\dift_{\dift_{Y}Z} = 0$, yields $\mathrm{Ric}_{h}(Y, Z) = -(n-1)h(Y, Z) + \tr(\dift_{Y}\dift_{Z})$. Since $\dift_{Y}$ is $h$-self-adjoint, $\tr(\dift_{Y}\dift_{Y}) = |\dift_{Y}|^{2}_{h} \geq 0$, with equality if and only if $\dift_{Y} = 0$. Tracing gives the formula for the scalar curvature.
\end{proof}

The standing assumption that $\Ga$ divides $\Om$ is not used in Lemma \ref{productcogaugelemma}, which concerns products of cones. For $i = 1, \dots, k$, let $\ste_{i}$ be a real vector space of dimension $N_{i} \geq 1$ with a parallel volume form $\Psi_{i}$, and let $C_{i} \subset \ste_{i}$ be an open properly convex cone, a half-line when $N_{i} = 1$. Let $N = \sum_{i = 1}^{k}N_{i}$, $a_{i} = N_{i}/N$, $\ste = \ste_{1}\oplus\dots\oplus\ste_{k}$ with the volume form $\Psi = \Psi_{1}\wedge\dots\wedge\Psi_{k}$, and $C = C_{1}\times\dots\times C_{k}$, an open properly convex cone in $\ste$. A function on $C_{i}$ is identified with its pullback to $C$ under the projection. For positive functions $\psi_{i}$ on $C_{i}$ homogeneous of degree $1$, let
\begin{align}\label{productcogauge}
\psi_{1}\boxtimes\dots\boxtimes\psi_{k} = \psi_{1}^{a_{1}}\cdots\psi_{k}^{a_{k}},
\end{align}
a positive function on $C$ homogeneous of degree $\sum_{i = 1}^{k}a_{i} = 1$. Here and below, for an open properly convex cone $C$, $\cAs(C, \{1\})$ is the set introduced after Definition \ref{asymptoticdefinition}. For $\phi \in \cAs(C, \{1\})$, the level set $\Sigma_{\phi}$, the function $\mu(\phi)$ of \eqref{mudefined}, and the centroaffine metric $\gb_{\phi}$ are defined as for elements of $\cAs(C, \Ga)$, and Lemma \ref{polarmaplemma}\ref{pml1}, Lemma \ref{centroaffinemetriclemma}\ref{cml4}, and Lemma \ref{affinespherelemma}, whose proofs do not use the group, apply to $\phi$. When $N_{i} = 1$, $\psi_{i}$ is a positive multiple $cx$ of a linear coordinate $x$, so $\psi_{i} \in \cAs(C_{i}, \{1\})$, its level set is a point, $\gb_{\psi_{i}} = 0$, and, by \eqref{mudefined}, $\mu(\psi_{i}) = c^{2}$ when $\Psi_{i} = dx$; in this case both identities of \eqref{mulogdet} hold on $C_{i}$ by direct computation.

\begin{lemma}\label{productcogaugelemma}
Let $\psi_{i}$ be positive concave functions on the $C_{i}$, homogeneous of degree $1$, and let $\psi = \psi_{1}\boxtimes\dots\boxtimes\psi_{k}$.
\begin{enumerate}[label=$(\arabic*)$, ref=(\arabic*)]
\item\label{pcl1} The function $\psi$ is concave, and if each $\psi_{i}$ is invariant under a group $\Ga_{i} \subset SL^{\pm}(\ste_{i})$, then $\psi$ is invariant under $\Ga_{1}\times\dots\times\Ga_{k}$.
\item\label{pcl3} If $\psi_{i} \in \cAs(C_{i}, \{1\})$ for each $i$, then $\psi \in \cAs(C, \{1\})$, and, with $\vartheta_{i} = \log\psi_{i}$,
\begin{align}\label{productmu}
&\mu(\psi) = \prod_{i=1}^{k}a_{i}^{N_{i}}\,\mu(\psi_{i}),& &\gb_{\psi} = \sum_{i=1}^{k}a_{i}\gb_{\psi_{i}} + \sum_{i=1}^{k}a_{i}\,d\vartheta_{i}^{2} - \Big(\sum_{i=1}^{k}a_{i}\,d\vartheta_{i}\Big)^{2}.&
\end{align}
In the second identity of \eqref{productmu}, $\gb_{\psi}$ is pulled back to $C$ by radial projection onto $\Sigma_{\psi}$, as in Lemma \ref{centroaffinemetriclemma}\ref{cml4}, and each $\gb_{\psi_{i}}$ is pulled back to $C_{i}$ by radial projection onto $\Sigma_{\psi_{i}}$ and then to $C$ by the projection of $C$ onto $C_{i}$. In particular, if each $\mu(\psi_{i})$ is constant, then so is $\mu(\psi)$, and $\Sigma_{\psi}$ is a hyperbolic affine sphere centered at the origin.
\end{enumerate}
\end{lemma}
\begin{proof}
The map $(u_{1}, \dots, u_{k}) \mapsto u_{1}^{a_{1}}\cdots u_{k}^{a_{k}}$ is concave and nondecreasing in each variable on $(0, \infty)^{k}$, so $\psi$ is concave, and the invariance is immediate. This proves \ref{pcl1}.

Next \ref{pcl3} is proved. Suppose $\psi_{i} \in \cAs(C_{i}, \{1\})$. Since $\log\psi = \sum_{i = 1}^{k}a_{i}\vartheta_{i}$, $\psi^{-1}\fD^{2}\psi = \fD^{2}\log\psi + (d\log\psi)^{2}$, and similarly for each $\psi_{i}$,
\begin{align}\label{producthessian}
\psi^{-1}\fD^{2}\psi = \sum_{i=1}^{k}a_{i}\,\psi_{i}^{-1}\fD^{2}\psi_{i} - \tfrac{1}{2}\sum_{i,j=1}^{k}a_{i}a_{j}(d\vartheta_{i} - d\vartheta_{j})^{2},
\end{align}
because $\sum_{i = 1}^{k}a_{i}\,d\vartheta_{i}^{2} - (\sum_{i = 1}^{k}a_{i}\,d\vartheta_{i})^{2} = \tfrac{1}{2}\sum_{i, j = 1}^{k}a_{i}a_{j}(d\vartheta_{i} - d\vartheta_{j})^{2}$, squares meaning symmetric squares. Each term on the right of \eqref{producthessian} is negative semidefinite, so $\fD^{2}\psi$ is negative semidefinite. If $v = (v_{1}, \dots, v_{k})$ is in its radical, then every term vanishes on $v$, so $v_{i} = s_{i}\eul_{i}$ for constants $s_{i}$, where $\eul_{i}$ is the Euler field of $\ste_{i}$, and then $s_{i} - s_{j} = d\vartheta_{i}(v) - d\vartheta_{j}(v) = 0$. Hence the radical of $\fD^{2}\psi$ is spanned by $\eul$, and $\psi \in \cAs(C, \{1\})$. Since $\vartheta_{i}$ depends only on $x_{i}$, $-\fD^{2}\log\psi = \sum_{i = 1}^{k}a_{i}(-\fD^{2}\vartheta_{i})$ is block diagonal with respect to $\ste = \ste_{1}\oplus\dots\oplus\ste_{k}$, so, by the first identity of \eqref{mulogdet} on $C$ and on each $C_{i}$,
\begin{align}\label{productmuproof}
\mu(\psi) = \psi^{2N}\prod_{i=1}^{k}a_{i}^{N_{i}}\det(-\fD^{2}\vartheta_{i}) = \psi^{2N}\prod_{i=1}^{k}a_{i}^{N_{i}}\mu(\psi_{i})\psi_{i}^{-2N_{i}} = \prod_{i=1}^{k}a_{i}^{N_{i}}\mu(\psi_{i}),
\end{align}
because $\psi^{2N} = \prod_{i = 1}^{k}\psi_{i}^{2N_{i}}$. By the second identity of \eqref{mulogdet} on $C$ and on each $C_{i}$, $\gb_{\psi} = -\fD^{2}\log\psi - (d\log\psi)^{2} = \sum_{i = 1}^{k}a_{i}(d\vartheta_{i}^{2} + \gb_{\psi_{i}}) - (\sum_{i = 1}^{k}a_{i}\,d\vartheta_{i})^{2}$, which is the second identity of \eqref{productmu}. The final claim follows from Lemma \ref{affinespherelemma}.
\end{proof}

The final claim of Lemma \ref{productcogaugelemma} recovers the product of hyperbolic affine spheres of Calabi \cite{Calabi-completeaffine}.

\begin{example}\label{hyperboloidexample}
If $\Ga \subset SO^{+}(n, 1)$ preserves the future cone $C$ of a Lorentzian inner product $\lb\,,\,\ra$ of signature $(n, 1)$ with volume form $\pm\Psi$, let $\psi_{0}(x) = (-\lb x, x\ra)^{1/2}$. In coordinates in which $\lb x, x\ra = |z|^{2} - t^{2}$, $p(z) = (1 - |z|^{2})^{1/2}$ on the unit ball, $\det(-\fD^{2}p) = (1 - |z|^{2})^{-(n+2)/2}$, and so $\mu(\psi_{0}) = 1$. Differentiating $\psi_{0}^{2} = -\lb x, x\ra$ twice gives $-\fD^{2}\psi_{0} = \psi_{0}^{-1}(\lb\,,\,\ra + d\psi_{0}\otimes d\psi_{0})$, whose restriction to the tangent spaces of $\Sigma_{\psi_{0}}$ is the hyperbolic metric of curvature $-1$. By Lemma \ref{affinespherelemma}, $\Sigma_{\psi_{0}}$ has equiaffine normal $\eul$ and affine mean curvature $-1$.
\end{example}

\begin{example}\label{orthantexample}
Let $\ell_{1}, \dots, \ell_{N}$ be linearly independent linear forms on $\ste$, let $C = \{x \in \ste: \ell_{i}(x) > 0, 1 \leq i \leq N\}$, a simplicial cone, so that $\Om$ is an open simplex, and let $c_{\ell} > 0$ be the constant such that $\Psi = \pm c_{\ell}\, d\ell_{1}\wedge\dots\wedge d\ell_{N}$. The group comprising $g \in SL(\ste)$ such that $\ell_{i}\circ g = a_{i}\ell_{i}$, $a_{i} > 0$, is abelian and isomorphic to $\rea^{n}$, and a lattice $\Ga$ in it acts cocompactly on $\Om$ with quotient an $n$-torus. The cone $C$ is the product of the half-lines $\{\ell_{i} > 0\}$, and, by Lemma \ref{productcogaugelemma}\ref{pcl3} applied with $\psi_{i} = \ell_{i}$, for which $\mu(\psi_{i}) = 1$ with respect to $d\ell_{i}$ and $\gb_{\psi_{i}} = 0$, the function $(\ell_{1}\cdots\ell_{N})^{1/N}$ has $\mu = N^{-N}$ with respect to $d\ell_{1}\wedge\dots\wedge d\ell_{N}$. Since $\mu(t\psi) = t^{2N}\mu(\psi)$ and replacing $d\ell_{1}\wedge\dots\wedge d\ell_{N}$ by $\Psi$ divides $\mu$ by $c_{\ell}^{2}$, the function $\psi_{0} = N^{1/2}c_{\ell}^{1/N}(\ell_{1}\cdots\ell_{N})^{1/N}$ satisfies $\mu(\psi_{0}) = 1$. By \eqref{productmu}, in the coordinates $y_{i} = \log \ell_{i}$, $\gb_{\psi_{0}}$ is the restriction to $\Sigma_{\psi_{0}}$ of $\tfrac{1}{N}\sum_{i=1}^{N}dy_{i}^{2} - \tfrac{1}{N^{2}}(\sum_{i=1}^{N}dy_{i})^{2}$, a flat metric, so $(M, \gb_{\psi_{0}})$ is a flat torus.
\end{example}

The simplicial cone of Example \ref{orthantexample} is linearly equivalent to the orthant, a homogeneous cone, and for every homogeneous cone $\mu(\rho_{C})$ is constant, so that the level sets of $\rho_{C}$ are the hyperbolic affine spheres asymptotic to $\pr C$. By the transformation rule $\phi_{C}(gx) = |\det g|^{-1}\phi_{C}(x)$ for $g \in GL(\ste)$ with $gC = C$ \cite{Vinberg}, the function $\mu(\rho_{C})$ is invariant under every such $g$, hence constant when these act transitively on $C$, and $\rho_{C} \to 0$ at $\pr C\setminus\{0\}$ by Lemma \ref{vinberglemma}. In Example \ref{orthantexample}, $\phi_{C}$ is a constant multiple of $(\ell_{1}\cdots\ell_{N})^{-1}$, so $\rho_{C} = \phi_{C}^{-1/N}$ is a constant multiple of $\psi_{0}$. That the canonical potential of a homogeneous cone differs by an additive constant from the logarithm of its characteristic function is observed in \cite[Section $4.5$]{Fox-schwarz}.

If $\psi \in \cA(C, \Ga)$, then, by Lemma \ref{polartransformlemma}, $\psi^{\ast} \in \cA(C^{\ast}, \Ga^{\ast})$. Let $K^{\ast}_{\psi} = \{\zeta \in C^{\ast}: \psi^{\ast}(\zeta) \leq 1\}$, and define
\begin{align}\label{pdefined}
\cP(\psi) = |K_{\psi}/\Ga|\,|K^{\ast}_{\psi}/\Ga^{\ast}| = \frac{1}{N^{2}}\int_{M}\psi^{-N}\int_{M^{\ast}}(\psi^{\ast})^{-N},
\end{align}
where the second equality follows from Lemma \ref{conevolumelemma}. Since a concave function on an open set is continuous, $\psi^{-N}$ and $(\psi^{\ast})^{-N}$ are continuous positive functions on the compact manifolds $M$ and $M^{\ast}$, so both integrals in \eqref{pdefined} are finite and positive. In particular $0 < \cP(\psi) < \infty$ for every $\psi \in \cA(C, \Ga)$. The notation $\cP(\psi)$ is consistent with that of \eqref{volprodintro} because, since $C \setminus L_{\psi}$ and $K_{\psi}$, and $C^{\ast}\setminus L^{\circ}_{\psi}$ and $K^{\ast}_{\psi}$, differ by null sets, $\cP(\psi) = \cP(L_{\psi})$. By \eqref{polardefined}, $(c\psi)^{\ast} = c^{-1}\psi^{\ast}$ for $c > 0$, and in \eqref{pdefined} this yields the scaling invariance $\cP(c\psi) = \cP(\psi)$. Moreover, $\cP$ is continuous in the following sense. If $\psi, \psi_{k} \in \cA(C, \Ga)$, then $\psi_{k}/\psi$ is $\Ga$-invariant and homogeneous of degree $0$ and so is identified with a function on $M$. Suppose that $\psi_{k}/\psi$ converges to $1$ uniformly, so that, given $\ep \in (0, 1)$, for $k$ large there holds $(1-\ep)\psi \leq \psi_{k} \leq (1+\ep)\psi$. Because, by \eqref{polardefined}, polarity reverses order, $(1+\ep)^{-1}\psi^{\ast} \leq \psi^{\ast}_{k} \leq (1-\ep)^{-1}\psi^{\ast}$. Substituting these bounds in \eqref{pdefined} gives
\begin{align}\label{pcontinuity}
\left(\frac{1-\ep}{1+\ep}\right)^{N} \leq \frac{\cP(\psi_{k})}{\cP(\psi)} \leq \left(\frac{1+\ep}{1-\ep}\right)^{N},
\end{align}
so $\cP(\psi_{k}) \to \cP(\psi)$.

\begin{remark}\label{hausdorffremark}
For $\psi, \psi_{k} \in \cA(C, \Ga)$, uniform convergence of $\psi_{k}/\psi$ to $1$ is equivalent to local Hausdorff convergence of $L_{\psi_{k}}$ to $L_{\psi}$, meaning that for every $R > 0$ each point of $L_{\psi_{k}}$ of norm at most $R$ lies within $\delta_{k}(R)$ of $L_{\psi}$ and each point of $L_{\psi}$ of norm at most $R$ lies within $\delta_{k}(R)$ of $L_{\psi_{k}}$, with $\delta_{k}(R) \to 0$. If $(1-\ep)\psi \leq \psi_{k} \leq (1+\ep)\psi$, then $y \in L_{\psi}$ implies $y/(1-\ep) \in L_{\psi_{k}}$, and $y \in L_{\psi_{k}}$ implies $y/(1+\ep) \in L_{\psi}$, which gives local Hausdorff convergence.

Conversely, suppose that $L_{\psi_{k}}$ converges locally to $L_{\psi}$, let $K \subset C$ be compact, and let $0 < \ep < \min_{K}\psi$, which suffices, since a smaller $\ep$ gives a stronger conclusion. If $\psi_{k}(x_{k}) \geq \psi(x_{k}) + \ep$ for infinitely many $k$ and some $x_{k} \in K$, then the points $x_{k}/(\psi(x_{k}) + \ep) \in L_{\psi_{k}}$ have a limit point $y$, which lies in $L_{\psi}$ although $\psi(y) < 1$. If $\psi_{k}(x_{k}) \leq \psi(x_{k}) - \ep$ for infinitely many $k$ and some $x_{k} \in K$, then the points $z_{k} = x_{k}/(\psi(x_{k}) - \ep/2)$ satisfy $\psi(z_{k}) \geq 1 + \eta$ for some $\eta > 0$, so the balls of some fixed radius $r$ about them lie in $L_{\psi}$, and a hyperplane through $z_{k}$ separating it from $L_{\psi_{k}}$ would leave a point of such a ball at distance at least $r$ from $L_{\psi_{k}}$; hence $z_{k} \in L_{\psi_{k}}$, that is, $\psi_{k}(x_{k}) \geq \psi(x_{k}) - \ep/2$, a contradiction.

So $\psi_{k} \to \psi$ uniformly on compact subsets of $C$, and, since $\psi_{k}/\psi$ is $\Ga$-invariant and homogeneous of degree $0$, uniform convergence on the compact set $\Fset$ of Lemma \ref{fundamentalsetlemma} is uniform convergence on $C$.

Consequently \eqref{pcontinuity} says that $\cP$ is continuous on $\cK(C, \Ga)$ for local Hausdorff convergence, and the same holds for each covolume separately and for polarity. By Lemma \ref{densitylemma}, $\Ga$-convex bodies with smooth, locally strongly convex boundary are dense for this convergence.
\end{remark}

The covolume product behaves under products of cones as the Mahler product behaves under products of convex bodies, \eqref{mahlerproductintro}. Consider the setting of Lemma \ref{productcogaugelemma} with $k = 2$, write $a = a_{1}$ and $b = a_{2}$, and suppose that $\Ga_{i} \subset SL^{\pm}(\ste_{i})$ preserves $C_{i}$ and acts freely, properly, and cocompactly on $\Om_{i} = \projp(C_{i})$, with quotient $M_{i}$; when $N_{i} = 1$, $\Ga_{i}$ is trivial and $M_{i}$ is a point. For $\tau > 0$ let
\begin{align}\label{relativedilation}
g_{\tau} = \left(e^{N_{2}\tau}\,\mathrm{id}_{\ste_{1}}, e^{-N_{1}\tau}\,\mathrm{id}_{\ste_{2}}\right) \in SL(\ste),
\end{align}
and let $\Ga = \Ga_{1}\times\Ga_{2}\times\lb g_{\tau}\ra$.

\begin{lemma}\label{productdivisiblelemma}
The group $\Ga$ preserves $C$ and acts freely, properly, and cocompactly on $\Om = \projp(C)$. More precisely, the map
\begin{align}\label{productidentification}
\Om \to \Om_{1}\times\Om_{2}\times\rea, \qquad [x_{1}, x_{2}] \mapsto \left([x_{1}], [x_{2}], \log\rho_{C_{1}}(x_{1}) - \log\rho_{C_{2}}(x_{2})\right)
\end{align}
is a homeomorphism that intertwines the action of $\Ga_{1}\times\Ga_{2}$ with its action on the first two factors and the action of $g_{\tau}$ with translation by $N\tau$ on the third factor, so $M = \Om/\Ga$ is homeomorphic to $M_{1}\times M_{2}\times S^{1}$.
\end{lemma}
\begin{proof}
The map \eqref{productidentification} is well defined because $\log\rho_{C_{1}}(tx_{1}) - \log\rho_{C_{2}}(tx_{2})$ does not depend on $t > 0$, and it is continuous because the $\rho_{C_{i}}$ are continuous and positive (Lemma \ref{vinberglemma}; when $N_{i} = 1$, $\rho_{C_{i}}$ is a positive multiple of a linear coordinate on the half-line $C_{i}$). Its inverse sends $([y_{1}], [y_{2}], s)$, with $\rho_{C_{i}}(y_{i}) = 1$, to $[e^{s}y_{1}, y_{2}]$, which is also continuous. Because $\Ga_{i}$ preserves $\rho_{C_{i}}$, the group $\Ga_{1}\times\Ga_{2}$ preserves the third coordinate, and $g_{\tau}$ adds $N_{2}\tau + N_{1}\tau = N\tau$ to it. Hence $\Ga$ acts on $\Om_{1}\times\Om_{2}\times\rea$ as the product of the actions of $\Ga_{1}$, $\Ga_{2}$, and $N\tau\mathbb{Z}$, each free, proper, and cocompact. Finally $\det g_{\tau} = e^{N_{1}N_{2}\tau - N_{2}N_{1}\tau} = 1$.
\end{proof}

By theorems of Vey \cite[Th\'eor\`emes $3$ and $5$ and Lemme $2$]{Vey} (see \cite[Section $5.1$]{Benoist-survey}), a group dividing a properly convex cone in a vector space acts semisimply on the vector space, and to every decomposition of the vector space as a direct sum of two subspaces stable under the group there corresponds a factorization of the cone as a product of properly convex cones in these subspaces, each stable under the group. The factors of a divisible product of cones are divisible \cite[Section $2.4$]{Benoist-survey}, and Lemma \ref{productdivisiblelemma} constructs groups dividing a product from groups dividing its factors.

Not every group dividing a product of cones has a subgroup of finite index of the form $\Ga_{1}\times\Ga_{2}\times\lb g_{\tau}\ra$. For example, when $N = 3$ and the simplicial cone of Example \ref{orthantexample} is regarded as the product of the quadrant $\{\ell_{1} > 0, \ell_{2} > 0\}$ and the half-line $\{\ell_{3} > 0\}$, a group of this form is generated by two maps $g$ satisfying $\ell_{i}\circ g = a_{i}\ell_{i}$, $1 \leq i \leq 3$, for which the vectors $(\log a_{1}, \log a_{2}, \log a_{3})$ are multiples of $(1, -1, 0)$ and $(1, 1, -2)$, respectively, while a lattice $\Ga$ as in Example \ref{orthantexample} need not contain any nontrivial element of the first kind.

For $\Ga_{i}$-convex bodies $L_{i} = L_{\psi_{i}} \subset C_{i}$ let $L_{1}\boxtimes L_{2} = L_{\psi_{1}\boxtimes\psi_{2}}$. By Lemma \ref{productcogaugelemma}\ref{pcl1} and because $(\psi_{1}\boxtimes\psi_{2})\circ g_{\tau} = e^{aN_{2}\tau - bN_{1}\tau}\psi_{1}\boxtimes\psi_{2} = \psi_{1}\boxtimes\psi_{2}$, $\psi_{1}\boxtimes\psi_{2} \in \cA(C, \Ga)$, so $L_{1}\boxtimes L_{2}$ is a $\Ga$-convex body by Lemma \ref{gammaconvexbodylemma}. It is the union of the images of $L_{1}\times L_{2}$ under the relative dilations,
\begin{align}\label{productbody}
L_{1}\boxtimes L_{2} = \{(x_{1}, x_{2}) \in C: \psi_{1}(x_{1})^{N_{1}}\psi_{2}(x_{2})^{N_{2}} \geq 1\} = \bigcup_{s > 0}\left(s^{N_{2}}L_{1}\times s^{-N_{1}}L_{2}\right),
\end{align}
and, because $\phi_{C} = \phi_{C_{1}}\phi_{C_{2}}$, $\rho_{C} = \rho_{C_{1}}\boxtimes\rho_{C_{2}}$.

More generally, let $Q = (0, \infty)^{2}$, let $d_{\tau}(u_{1}, u_{2}) = (e^{N_{2}\tau}u_{1}, e^{-N_{1}\tau}u_{2})$, and let $\mathcal{G}_{\tau}$ be the set of positive concave functions on $Q$ homogeneous of degree $1$ and invariant under $d_{\tau}$, which contains $G_{0}(u_{1}, u_{2}) = u_{1}^{a}u_{2}^{b}$. By \eqref{superadditive} on $Q$, $G(u + v) \geq G(u) + G(v) > G(u)$ for $u, v \in Q$, so, by the continuity of the concave function $G$, each $G \in \mathcal{G}_{\tau}$ is nondecreasing in each variable, so, for $\psi_{i} \in \cA(C_{i}, \Ga_{i})$, $G(\psi_{1}, \psi_{2})$ is concave, and it is positive, homogeneous of degree $1$, and $\Ga$-invariant. Thus $G(\psi_{1}, \psi_{2}) \in \cA(C, \Ga)$, and $G_{0}(\psi_{1}, \psi_{2}) = \psi_{1}\boxtimes\psi_{2}$. For $\Ga_{i}$-convex bodies $L_{i} = L_{\psi_{i}}$, write $G(L_{1}, L_{2})$ for the $\Ga$-convex body $L_{G(\psi_{1}, \psi_{2})}$.

The polar $G^{\ast}(A_{1}, A_{2}) = \inf_{r \in Q}(r_{1}A_{1} + r_{2}A_{2})/G(r_{1}, r_{2})$ of $G \in \mathcal{G}_{\tau}$ again lies in $\mathcal{G}_{\tau}$. Indeed, $G^{\ast}$ is concave and homogeneous of degree $1$, as an infimum of linear functions of $(A_{1}, A_{2}) \in Q$; it is positive, because the positive concave function $G$ is bounded above on the segment $\{r \in Q: r_{1} + r_{2} = 1\}$, by a constant $c_{G}$, so that $G(r) \leq c_{G}(r_{1} + r_{2})$ for $r \in Q$ and $G^{\ast}(A_{1}, A_{2}) \geq c_{G}^{-1}\min\{A_{1}, A_{2}\}$; and it is invariant under $d_{\tau}$, by the substitution $r \mapsto d_{\tau}r$ in the infimum and the invariance of $G$ under $d_{\tau}$. Moreover, $G_{0}^{\ast} = a^{-a}b^{-b}G_{0}$ by the weighted inequality between arithmetic and geometric means. For $G \in \mathcal{G}_{\tau}$ let
\begin{align}\label{planarvolume}
&V(G) = \int_{\{r \in R_{\tau}:\, G(r) \leq 1\}}r_{1}^{N_{1}-1}r_{2}^{N_{2}-1}\,dr_{1}\,dr_{2},& &\cP_{N_{1}, N_{2}}(G) = V(G)\,V(G^{\ast}),&
\end{align}
where $R_{\tau} = \{r \in Q: 0 \leq \log r_{1} - \log r_{2} < N\tau\}$ is a fundamental domain for $\lb d_{\tau}\ra$. In the variables $P = N_{1}\log r_{1} + N_{2}\log r_{2}$ and $D = \log r_{1} - \log r_{2}$, for which $r_{1}^{N_{1}-1}r_{2}^{N_{2}-1}dr_{1}\,dr_{2} = N^{-1}e^{P}dP\,dD$,
\begin{align}\label{productcovolumecomputation}
V(G_{0}) = \frac{1}{N}\int_{0}^{N\tau}dD\int_{-\infty}^{0}e^{P}\,dP = \tau,
\end{align}
and, since $G_{0}^{\ast} = a^{-a}b^{-b}G_{0}$, $V(G_{0}^{\ast}) = (a^{a}b^{b})^{N}\tau$, so $\cP_{N_{1}, N_{2}}(G_{0}) = a^{N_{1}}b^{N_{2}}\tau^{2} = \tau^{2}N_{1}^{N_{1}}N_{2}^{N_{2}}N^{-N}$.

\begin{lemma}\label{productcovolumelemma}
For $\Ga_{i}$-convex bodies $L_{i} \subset C_{i}$ and $G \in \mathcal{G}_{\tau}$,
\begin{align}\label{productcovolume}
\left|(C\setminus G(L_{1}, L_{2}))/\Ga\right| = N_{1}N_{2}\,V(G)\,\left|(C_{1}\setminus L_{1})/\Ga_{1}\right|\,\left|(C_{2}\setminus L_{2})/\Ga_{2}\right|.
\end{align}
\end{lemma}
\begin{proof}
Let $\psi_{i}$ be the cogauge of $L_{i}$, let $\psi = G(\psi_{1}, \psi_{2})$, and let $\Sigma_{i} = \{y \in C_{i}: \psi_{i}(y) = 1\}$. Let $\nu_{i}$ be the measure on $\Sigma_{i}$ that corresponds, under the radial projection of $\Sigma_{i}$ onto $\Om_{i}$, to the measure on $\Om_{i}$ determined by the density $\psi_{i}^{-N_{i}}|\iota_{\eul_{i}}\Psi_{i}|$; when $\Sigma_{i}$ is smooth, $\nu_{i}$ is the restriction of $|\iota_{\eul_{i}}\Psi_{i}|$ to $\Sigma_{i}$. By \eqref{conevolumelocal} applied to $t^{-1}\psi_{i}$, for $t > 0$ and a Borel set $E \subset \Sigma_{i}$ the set $\{ry: 0 < r \leq t,\ y \in E\}$ has volume $t^{N_{i}}\nu_{i}(E)/N_{i}$, so, in the coordinates $x_{i} = r_{i}y_{i}$, $r_{i} > 0$, $y_{i} \in \Sigma_{i}$, Lebesgue measure on $C_{i}$ is $r_{i}^{N_{i}-1}dr_{i}\,d\nu_{i}$.

Since $\psi_{i} = 1$ on $\Sigma_{i}$, if $R_{i} \subset \Sigma_{i}$ is a Borel fundamental domain for the action of $\Ga_{i}$, then $\nu_{i}(R_{i}) = \int_{M_{i}}\psi_{i}^{-N_{i}} = N_{i}|K_{\psi_{i}}/\Ga_{i}|$, by Lemma \ref{conevolumelemma} when $N_{i} \geq 2$ and directly when $N_{i} = 1$. By Lemma \ref{productdivisiblelemma}, applied with $\psi_{i}$ in place of $\rho_{C_{i}}$, which changes nothing in its proof, the set of $(r_{1}y_{1}, r_{2}y_{2})$ with $y_{i} \in R_{i}$ and $(r_{1}, r_{2}) \in R_{\tau}$ is a fundamental domain for the action of $\Ga$ on $C$. Since $\psi(r_{1}y_{1}, r_{2}y_{2}) = G(r_{1}, r_{2})$, the volume of $K_{\psi}$ intersected with this fundamental domain is $\nu_{1}(R_{1})\,\nu_{2}(R_{2})\,V(G)$. Since $K_{\psi}$ and $C\setminus L_{\psi}$ differ by a null set, this is \eqref{productcovolume}.
\end{proof}

\begin{lemma}\label{productpolarlemma}
Let $\psi_{i} \in \cA(C_{i}, \Ga_{i})$ and $G \in \mathcal{G}_{\tau}$. Then $C^{\ast} = C_{1}^{\ast}\times C_{2}^{\ast}$, the group $\Ga^{\ast}$ is the group of the same form for the $\Ga_{i}^{\ast}$, with the same $\tau$, and
\begin{align}\label{productpolar}
&G(\psi_{1}, \psi_{2})^{\ast} = G^{\ast}(\psi_{1}^{\ast}, \psi_{2}^{\ast}),& &G(L_{1}, L_{2})^{\circ} = G^{\ast}(L_{1}^{\circ}, L_{2}^{\circ}).&
\end{align}
In particular $(L_{1}\boxtimes L_{2})^{\circ} = a^{a}b^{b}\left(L_{1}^{\circ}\boxtimes L_{2}^{\circ}\right)$.
\end{lemma}
\begin{proof}
That $C^{\ast} = C_{1}^{\ast}\times C_{2}^{\ast}$ is immediate. The inverse transpose of $g_{\tau}$ is $(e^{-N_{2}\tau}\mathrm{id}, e^{N_{1}\tau}\mathrm{id})$, the inverse of the element \eqref{relativedilation} for $\std_{1}\oplus\std_{2}$, so $\Ga^{\ast}$ is generated by $\Ga_{1}^{\ast}\times\Ga_{2}^{\ast}$ and that element. Let $\zeta = (\zeta_{1}, \zeta_{2}) \in C^{\ast}$ and $A_{i} = \psi_{i}^{\ast}(\zeta_{i}) > 0$. Writing $x_{i} = r_{i}y_{i}$ with $\psi_{i}(y_{i}) = 1$, the infimum of $\zeta_{i}(y_{i})$ over $y_{i} \in \Sigma_{i}$ is $A_{i}$, and the denominator below does not depend on the $y_{i}$, so, by \eqref{polardefined},
\begin{align}\label{productpolarcomputation}
G(\psi_{1}, \psi_{2})^{\ast}(\zeta) = \inf_{x \in C}\frac{\zeta_{1}(x_{1}) + \zeta_{2}(x_{2})}{G(\psi_{1}(x_{1}), \psi_{2}(x_{2}))} = \inf_{r \in Q}\frac{r_{1}A_{1} + r_{2}A_{2}}{G(r_{1}, r_{2})} = G^{\ast}(A_{1}, A_{2}).
\end{align}
This is the first identity of \eqref{productpolar}, and the second follows because $L^{\circ}_{\psi} = L_{\psi^{\ast}}$ by the discussion following \eqref{polardefined}. The last claim follows from $G_{0}^{\ast} = a^{-a}b^{-b}G_{0}$ and $L_{c\phi} = c^{-1}L_{\phi}$ for $c > 0$.
\end{proof}

\begin{theorem}\label{producttheorem}
For $\Ga_{i}$-convex bodies $L_{i} \subset C_{i}$ and $G \in \mathcal{G}_{\tau}$,
\begin{align}\label{compositeidentity}
\cP(G(L_{1}, L_{2})) = N_{1}^{2}N_{2}^{2}\,\cP_{N_{1}, N_{2}}(G)\,\cP(L_{1})\,\cP(L_{2}),
\end{align}
where $\cP(L_{i}) = 1$ when $N_{i} = 1$. In particular,
\begin{align}\label{productidentity}
\cP(L_{1}\boxtimes L_{2}) = \tau^{2}\,\frac{N_{1}^{N_{1}+2}N_{2}^{N_{2}+2}}{N^{N}}\,\cP(L_{1})\,\cP(L_{2}).
\end{align}
\end{theorem}
\begin{proof}
When $N_{i} = 1$ every $\Ga_{i}$-convex body is a half-line $[t, \infty)$, whose polar is $[t^{-1}, \infty)$, so $\cP(L_{i}) = 1$. Applying Lemma \ref{productcovolumelemma} to $G(L_{1}, L_{2})$ and, in the dual setting, to $G(L_{1}, L_{2})^{\circ} = G^{\ast}(L_{1}^{\circ}, L_{2}^{\circ})$ (Lemma \ref{productpolarlemma}), and multiplying, gives \eqref{compositeidentity}, and \eqref{productidentity} is the case $G = G_{0}$.
\end{proof}

The identity \eqref{productidentity} is the product of the identity \eqref{productcovolume} for $G = G_{0}$, in which $V(G_{0}) = \tau$, and the corresponding identity for the polars,
\begin{align}\label{productpolarcovolume}
N^{N}\left|(C^{\ast}\setminus (L_{1}\boxtimes L_{2})^{\circ})/\Ga^{\ast}\right| = \tau N_{1}^{N_{1}+1}N_{2}^{N_{2}+1}\left|(C_{1}^{\ast}\setminus L_{1}^{\circ})/\Ga_{1}^{\ast}\right|\,\left|(C_{2}^{\ast}\setminus L_{2}^{\circ})/\Ga_{2}^{\ast}\right|,
\end{align}
which follows from Lemmas \ref{productcovolumelemma} and \ref{productpolarlemma}. The factor $N_{1}N_{2}$ in each of these two identities comes from the relation $\nu_{i}(R_{i}) = N_{i}|K_{\psi_{i}}/\Ga_{i}|$ between the cone measure and the covolume, used in the proof of Lemma \ref{productcovolumelemma}, and these factors account for the factor $N_{1}^{2}N_{2}^{2}$ in \eqref{productidentity}. The constant $N_{1}^{N_{1}}N_{2}^{N_{2}}N^{-N}$ appears only in \eqref{productpolarcovolume}, and it plays the role of $m_{1}!\,m_{2}!/(m_{1}+m_{2})!$ in \eqref{mahlerproductintro}.

For convex bodies, $|K_{1}\times K_{2}| = |K_{1}||K_{2}|$, and the factorials are the ratio of the volume of the polar of $K_{1}\times K_{2}$, which is the free sum of $K_{1}^{\circ}$ and $K_{2}^{\circ}$, to $|K_{1}^{\circ}||K_{2}^{\circ}|$; after integration over the factors, this ratio is $m_{1}m_{2}$ times the area of the triangle $\{r_{1} + r_{2} \leq 1\}$ in the quadrant for the weight $r_{1}^{m_{1}-1}r_{2}^{m_{2}-1}$. Here the constant is the factor $(a^{a}b^{b})^{N}$ by which the dilation relating $(L_{1}\boxtimes L_{2})^{\circ}$ to $L_{1}^{\circ}\boxtimes L_{2}^{\circ}$ in Lemma \ref{productpolarlemma} multiplies covolumes.

The factor $\tau$ in each identity is a modulus of $M$ with no counterpart for convex bodies. By Lemma \ref{productdivisiblelemma}, $M$ is homeomorphic to $M_{1}\times M_{2}\times S^{1}$, and, writing $a = N_{1}/N$, $b = N_{2}/N$, and $\vartheta_{i} = \log\psi_{i}$, the second identity of \eqref{productmu} gives $\gb_{\psi_{1}\boxtimes\psi_{2}} = a\gb_{\psi_{1}} + b\gb_{\psi_{2}} + ab\,(d\vartheta_{1} - d\vartheta_{2})^{2}$ for a product $L_{1}\boxtimes L_{2}$ of bodies with cogauges $\psi_{i} \in \cAs(C_{i}, \Ga_{i})$. Since $g_{\tau}$ translates $\vartheta_{1} - \vartheta_{2}$ by $N\tau$, the centroaffine metric is a Riemannian product in which the circle factor has length $(ab)^{1/2}N\tau = (N_{1}N_{2})^{1/2}\tau$.

\begin{lemma}\label{volumeproposition}
For $\psi \in \cAs(C, \Ga)$, $G_{\psi}(K_{\psi}) = K^{\ast}_{\psi}$, and
\begin{align}\label{volumeformulas}
|K^{\ast}_{\psi}/\Ga^{\ast}| = \int_{K_{\psi}/\Ga}\mu(\psi)\,dx = \frac{1}{N}\int_{M}\mu(\psi)\psi^{-N}.
\end{align}
\end{lemma}
\begin{proof}
By Lemma \ref{polarmaplemma}\ref{pml3}--\ref{pml4}, $G_{\psi}:C \to C^{\ast}$ is a diffeomorphism such that $\psi^{\ast}\circ G_{\psi} = \psi$, so $G_{\psi}(K_{\psi}) = K^{\ast}_{\psi}$. It is equivariant and its Jacobian has absolute value $\mu(\psi)$, which gives the first equality in \eqref{volumeformulas}, while Lemma \ref{conevolumelemma} gives the second equality in \eqref{volumeformulas}.
\end{proof}

For $\psi \in \cAs(C, \Ga)$, by \eqref{pdefined} and Lemma \ref{volumeproposition},
\begin{align}
&\cS_{0}(\psi) = N|K_{\psi}/\Ga|,&& \cS_{1}(\psi) = N|K^{\ast}_{\psi}/\Ga^{\ast}|, && N^{2}\cP(\psi) = \cS_{0}(\psi)\cS_{1}(\psi). 
\end{align}

\begin{lemma}\label{dualityproposition}
For $\psi \in \cAs(C, \Ga)$, $\mu(\psi^{\ast})\circ G_{\psi} = \mu(\psi)^{-1}$, and, for all $\be \in \rea$,
\begin{align}\label{sbetaduality}
\cS_{\be}(\psi^{\ast}) = \cS_{1-\be}(\psi).
\end{align}
In particular $\cS_{1/2}(\psi^{\ast}) = \cS_{1/2}(\psi)$ and $\cP(\psi^{\ast}) = \cP(\psi)$.
\end{lemma}
\begin{proof}
By Lemma \ref{polarmaplemma}\ref{pml4}, $\fD^{2}\Phi^{\ast}\circ G_{\psi}$ is the bilinear form inverse to $\fD^{2}\Phi$, so $\mu(\psi^{\ast})\circ G_{\psi} = \mu(\psi)^{-1}$. By \eqref{sbetadefined} applied on $M^{\ast}$, the change of variables $\zeta = G_{\psi}(x)$, which maps $K_{\psi}$ onto $K^{\ast}_{\psi} = K_{\psi^{\ast}}$ with $d\zeta = \mu(\psi)dx$,
\begin{align}
\cS_{\be}(\psi^{\ast}) = N\int_{K_{\psi^{\ast}}/\Ga^{\ast}}\mu(\psi^{\ast})^{\be}d\zeta = N\int_{K_{\psi}/\Ga}\mu(\psi)^{1-\be}dx = \cS_{1-\be}(\psi).
\end{align}
\end{proof}
In particular, $\mu(\psi)$ is constant if and only if $\mu(\psi^{\ast})$ is constant, so the conormal map $G_{\psi}$ of an affine sphere $\Sigma_{\psi}$ maps it onto the affine sphere $\Sigma_{\psi^{\ast}}$. This correspondence between the affine spheres foliating $C$ and those foliating $C^{\ast}$ is \cite[Theorem $4.2$]{Fox-schwarz}, where it is derived from the Legendre transform of the canonical potential.

\begin{remark}\label{klartagmahlerremark}
By Lemma \ref{dualityproposition}, polarity maps the cogauge of the affine sphere of $C$ to that of $C^{\ast}$. The characteristic cogauge $\rho_{C} = \phi_{C}^{-1/N}$ has this property only up to a factor given by Mahler products. For a compact convex set $K$ with nonempty interior in an affine space of dimension $n$, let $\bar{\mahler}(K) = \min_{p \in \operatorname{int}K}\mahler(K - p)$, where $\operatorname{int}K$ is the interior of $K$, be its Mahler product about its Santal\'o point; it is an affine invariant. For $\zeta \in C^{\ast}$ let $K_{\zeta} = \{x \in \bar{C}: \zeta(x) = 1\}$, a compact convex set with nonempty interior in the affine hyperplane $\{x \in \ste: \zeta(x) = 1\}$, whose interior in that hyperplane is $K_{\zeta}\cap C$. Since $K_{t\zeta} = t^{-1}K_{\zeta}$ and $K_{\zeta\circ\gamma^{-1}} = \gamma K_{\zeta}$, the function $\zeta \mapsto \bar{\mahler}(K_{\zeta})$ is the lift of a function on $M^{\ast}$. The following identity is a reformulation, in the conventions for dual cones used here, of the identity of Klartag \cite[Corollary $2.4$]{Klartag-isotropic} between the logarithmic Laplace transform of a cone and the Legendre transform of that of its dual:
\begin{align}\label{klartagidentity}
\rho_{C}^{\ast}(\zeta) = (n!)^{2/N}\,\bar{\mahler}(K_{\zeta})^{1/N}\rho_{C^{\ast}}(\zeta), \qquad \zeta \in C^{\ast}.
\end{align}
It can be proved directly as follows. Let $\zeta \in C^{\ast}$. Since $x \mapsto \zeta(x)/\rho_{C}(x)$ is homogeneous of degree $0$ on $C$, by \eqref{polardefined},
\begin{align}\label{klartagstep1}
\rho_{C}^{\ast}(\zeta)^{N} = \inf_{x \in C}\zeta(x)^{N}\phi_{C}(x) = \inf_{x \in K_{\zeta}\cap C}\phi_{C}(x).
\end{align}
Let $x \in K_{\zeta}\cap C$ and let $H^{\ast}_{x} = \{\theta \in \std: \theta(x) = 1\}$. In the coordinates $\eta = t\theta$ on $C^{\ast}$, with $t > 0$ and $\theta \in C^{\ast}\cap H^{\ast}_{x}$, $d\eta = t^{n}\,dt\,d\si_{x}(\theta)$, where $\si_{x}$ is the restriction of $|\iota_{\eul}\Psi^{\ast}|$ to $H^{\ast}_{x}$, as in the proof of Lemma \ref{conevolumelemma}, and $\eta(x) = t$, so
\begin{align}\label{klartagstep2}
\phi_{C}(x) = \int_{0}^{\infty}e^{-t}t^{n}\,dt\,\si_{x}(C^{\ast}\cap H^{\ast}_{x}) = n!\,\si_{x}(C^{\ast}\cap H^{\ast}_{x}).
\end{align}
In the same way, $\phi_{C^{\ast}}(\zeta) = n!\,|K_{\zeta}|$, where $|K_{\zeta}|$ is the measure of $K_{\zeta}$ with respect to the restriction of $|\iota_{\eul}\Psi|$ to the hyperplane $\{y \in \ste: \zeta(y) = 1\}$, which is a Lebesgue measure on this hyperplane. The map $\theta \mapsto (\theta - \zeta)|_{\ker\zeta}$ is an affine isomorphism from $H^{\ast}_{x}$ onto the dual of $\ker\zeta$, because $\ste$ is the direct sum of the line spanned by $x$ and $\ker\zeta$. For $\theta \in H^{\ast}_{x}$ and $y \in K_{\zeta}$, $\theta(y) = 1 + (\theta - \zeta)(y - x)$, because $(\theta - \zeta)(x) = 0$ and $y - x \in \ker\zeta$, so $\theta$ is positive on $K_{\zeta}$, that is, $\theta \in C^{\ast}$, if and only if the negative of $(\theta - \zeta)|_{\ker\zeta}$ lies in the interior of the polar body $(K_{\zeta} - x)^{\circ}$ of the convex body $K_{\zeta} - x \subset \ker\zeta$. If $x, v_{1}, \dots, v_{n}$ is a basis of $\ste$ with $v_{i} \in \ker\zeta$ and $\Psi(x, v_{1}, \dots, v_{n}) = 1$, then its dual basis has the form $\zeta, \theta_{1}, \dots, \theta_{n}$, with $\theta_{i}(x) = 0$, and $\Psi^{\ast}(\zeta, \theta_{1}, \dots, \theta_{n}) = 1$, so the isomorphism carries $\si_{x}$ to the Lebesgue measure on the dual of $\ker\zeta$ dual to the Lebesgue measure on $\ker\zeta$ with respect to which $|K_{\zeta}|$ is computed. Hence $\si_{x}(C^{\ast}\cap H^{\ast}_{x}) = |(K_{\zeta} - x)^{\circ}|$, and, by \eqref{klartagstep2},
\begin{align}\label{klartagstep3}
\phi_{C^{\ast}}(\zeta)\,\phi_{C}(x) = (n!)^{2}\,|K_{\zeta}|\,|(K_{\zeta} - x)^{\circ}| = (n!)^{2}\,\mahler(K_{\zeta} - x).
\end{align}
Taking the infimum of \eqref{klartagstep3} over $x \in K_{\zeta}\cap C$ and using \eqref{klartagstep1} gives $\phi_{C^{\ast}}(\zeta)\rho_{C}^{\ast}(\zeta)^{N} = (n!)^{2}\bar{\mahler}(K_{\zeta})$, which is \eqref{klartagidentity}, because $\phi_{C^{\ast}} = \rho_{C^{\ast}}^{-N}$. In particular, by \eqref{klartagidentity}, the function $\zeta \mapsto \bar{\mahler}(K_{\zeta})$ is continuous on $C^{\ast}$, and hence so is the function it induces on $M^{\ast}$.
Consequently $\rho_{C}^{\ast}$ is a constant multiple of $\rho_{C^{\ast}}$ if and only if $\bar{\mahler}(K_{\zeta})$ does not depend on $\zeta$. This holds for homogeneous cones, because the adjoints of the linear automorphisms of a homogeneous cone act transitively on its dual \cite[Chapter I, Proposition $9$]{Vinberg}, so that the sections $K_{\zeta}$ are linear images of one another. Since $\rho_{C}^{-N} = \phi_{C}$ and $\rho_{C^{\ast}}^{-N} = \phi_{C^{\ast}}$, \eqref{pdefined} and \eqref{klartagidentity} give
\begin{align}\label{klartagvolumeproduct}
N^{2}(n!)^{2}\,\cP(\rho_{C}) = \int_{M}\phi_{C}\int_{M^{\ast}}\frac{\phi_{C^{\ast}}(\zeta)}{\bar{\mahler}(K_{\zeta})}.
\end{align}
When $C$ is homogeneous, $\rho_{C}$ is a constant multiple of the cogauge of the affine sphere, so, by Theorem A, $\cP(\rho_{C})$ is the minimum of $\cP$ on $\cK(C, \Ga)$, and, by \eqref{klartagvolumeproduct}, it equals $\int_{M}\phi_{C}\int_{M^{\ast}}\phi_{C^{\ast}}$ divided by $N^{2}(n!)^{2}\bar{\mahler}(K)$, for any section $K = K_{\zeta}$ of $C$. For a simplicial cone $K$ is a simplex, for which $(n!)^{2}\bar{\mahler}(K) = N^{N}$, so the minimum is $\int_{M}\phi_{C}\int_{M^{\ast}}\phi_{C^{\ast}}/N^{N+2}$, and for the Lorentz cone $K$ is an ellipsoid. Thus the simplex and the ellipsoid, the conjectured minimizer and the maximizer of the Mahler product, appear as the sections of the cones in the two model cases.
\end{remark}

\begin{lemma}\label{cauchyschwarzproposition}
For $\psi \in \cAs(C, \Ga)$ and $\be \in \rea$,
\begin{align}\label{cauchyschwarz}
\cS_{\be}(\psi)\cS_{1-\be}(\psi) \geq \cS_{1/2}(\psi)^{2} = \vol(M, \gb_{\psi})^{2},
\end{align}
with equality for some $\be \neq 1/2$ if and only if $\Sigma_{\psi}$ is an affine sphere. In particular $N^{2}\cP(\psi) \geq \vol(M, \gb_{\psi})^{2}$.
\end{lemma}
\begin{proof}
This is the Cauchy--Schwarz inequality for $\mu^{\be/2}$ and $\mu^{(1-\be)/2}$ in $L^{2}(\psi^{-N})$, together with Lemma \ref{centroaffinemetriclemma}\ref{cml1}.
\end{proof}

The functionals $\cS_{\be}$ have a description in terms of Euclidean quantities, which is the form \eqref{sbetaintro} given in the introduction and which connects them with the $L_{p}$ affine surface areas.

Fix a Euclidean inner product $g$ on $\ste$ whose volume form is $\pm\Psi$, and let $d\si$ be the induced hypersurface measure. For $\psi \in \cAs(C, \Ga)$, let $\nrm_{\psi}$ be the $g$-unit normal to $\Sigma_{\psi}$ pointing into $L_{\psi}$, let $s_{\psi} = g(\eul, \nrm_{\psi})$ be the support function of $\Sigma_{\psi}$, which is positive because $L_{\psi}$ lies on the side of $\Sigma_{\psi}$ away from the origin. For a convex body containing the origin in its interior, the corresponding positive quantity is computed with the outer normal and is the usual support function; for a convex set not containing the origin, such as $L_{\psi}$, the usual support function, computed with the normal pointing away from the set, is $-s_{\psi}$.

Let $\varkappa_{\psi}$ be the Gauss--Kronecker curvature of $\Sigma_{\psi}$ with respect to $\nrm_{\psi}$. The \emph{centroaffine curvature} of $\Sigma_{\psi}$ is $\det\nolimits_{\theta_{\eul}}h_{\eul}$, the determinant of the centroaffine metric with respect to the induced volume form $\theta_{\eul}$, which by \eqref{mudettheta} equals $\mu(\psi)$. For the transversal $\nrm_{\psi}$ the analogous quantity $\det\nolimits_{\theta_{\nrm_{\psi}}}h_{\nrm_{\psi}}$ is $|\varkappa_{\psi}|$, since $h_{\nrm_{\psi}}$ is the second fundamental form and $\theta_{\nrm_{\psi}}$ is $d\si$. Lemma \ref{centroaffineproposition} shows that the centroaffine curvature equals $|\varkappa_{\psi}|s_{\psi}^{-n-2}$, which is the centroaffine curvature of \cite{Hug-polar} and the function written $\kappa_{0}$ in \cite{Ludwig-general}.

\begin{lemma}\label{centroaffineproposition}
Let $\psi \in \cAs(C, \Ga)$. Then, on $\Sigma_{\psi}$, $g(\nrm_{\psi}, \cdot) = G_{\psi}/|G_{\psi}|_{g}$, where $|\cdot|_{g}$ is the dual norm, $s_{\psi} = |G_{\psi}|_{g}^{-1} = \psi^{\ast}(g(\nrm_{\psi}, \cdot))$, and $\mu(\psi) = |\varkappa_{\psi}|s_{\psi}^{-n-2}$. Consequently
\begin{align}\label{sbetacentroaffine}
\cS_{\be}(\psi) = \int_{\Sigma_{\psi}/\Ga}s_{\psi}^{1 - (n+2)\be}|\varkappa_{\psi}|^{\be}\,d\si.
\end{align}
\end{lemma}

\begin{proof}
Write $\omega = G_{\psi} = \fD\Phi$, so that $\fD^{2}\Phi(\eul, \cdot) = \omega$ and $\omega(\eul) = 2\Phi = 1$ on $\Sigma_{\psi} = \{x \in C: \Phi(x) = 1/2\}$. Since $T\Sigma_{\psi} = \ker \omega$ and $\Phi > 1/2$ on the interior of $L_{\psi}$, the covector $g(\nrm_{\psi}, \cdot)$ is a positive multiple of $\omega$ of dual norm $1$, so it equals $\omega/|\omega|_{g}$. Hence $s_{\psi} = \omega(\eul)/|\omega|_{g} = |\omega|_{g}^{-1}$, and $\psi^{\ast}(\omega/|\omega|_{g}) = |\omega|_{g}^{-1}\psi^{\ast}(G_{\psi}) = s_{\psi}$ by Lemma \ref{polarmaplemma}. For vector fields $X$ and $Y$ tangent to $\Sigma_{\psi}$, differentiating $\omega(Y) = 0$ along $X$ gives $\fD^{2}\Phi(X, Y) = -\omega(\fD_{X}Y) = -|\omega|_{g}\,g(\fD_{X}Y, \nrm_{\psi})$, so the second fundamental form $g(\fD_{X}Y, \nrm_{\psi})$ of $\Sigma_{\psi}$ with respect to $\nrm_{\psi}$ is $-|\omega|_{g}^{-1}\fD^{2}\Phi|_{T\Sigma_{\psi}}$, which is positive definite.

Let $e_{1}, \dots, e_{n}$ be a $g$-orthonormal basis of $T_{x}\Sigma_{\psi}$. Then $|\varkappa_{\psi}| = |\omega|_{g}^{-n}|\det(\fD^{2}\Phi(e_{i}, e_{j}))| = s_{\psi}^{n}|\det(\fD^{2}\Phi(e_{i}, e_{j}))|$. Since $\fD^{2}\Phi(e_{i}, \eul) = \omega(e_{i}) = 0$ and $\fD^{2}\Phi(\eul, \eul) = 1$, the matrix of $\fD^{2}\Phi$ in the basis $e_{1}, \dots, e_{n}, \eul_{x}$ is block diagonal with blocks $(\fD^{2}\Phi(e_{i}, e_{j}))$ and $1$, and $|\Psi(e_{1}, \dots, e_{n}, \eul_{x})| = |g(\eul_{x}, \nrm_{\psi})| = s_{\psi}$ because $\Psi = \pm\vol_{g}$. By \eqref{mudefined}, $\mu(\psi) = |\det(\fD^{2}\Phi(e_{i}, e_{j}))|s_{\psi}^{-2}$. Hence $|\varkappa_{\psi}| = s_{\psi}^{n+2}\mu(\psi)$.

For a function $F$ on $M$, writing $x = ry$ with $y \in \Sigma_{\psi}$ and $0 < r \leq 1$ gives $dx = r^{n}s_{\psi}\,dr\,d\si$, so $N\int_{K_{\psi}/\Ga}F\,dx = \int_{\Sigma_{\psi}/\Ga}Fs_{\psi}\,d\si$, and \eqref{sbetacentroaffine} follows from \eqref{sbetadefined}.
\end{proof}

\begin{remark}\label{lpremark}
With $p = N\be/(1-\be)$ for $\be \neq 1$, equivalently $\be = p/(N+p)$, the integrand $s^{1-(n+2)\be}|\varkappa|^{\be}$, where $s$ is the support function in the inward-normal convention fixed above, in which it is positive, and is written $h$ in \cite{Hug-contributions, Lutwak-Lp}, equals
\begin{align}\label{lpintegrand}
|\varkappa|^{p/(N+p)}s^{-N(p-1)/(N+p)},
\end{align}
which, for a convex body in $\rea^{N}$ containing the origin, is the integrand of the $L_{p}$ affine surface area $\Omega_{p}$ of Lutwak \cite{Lutwak-Lp}, in the representation of Hug \cite{Hug-contributions} valid for all $p > 0$. See \cite[equation (2)]{Ludwig-general}.

For $p < 0$, $p \neq -N$, the same integrand defines the $L_{p}$ affine surface area of \cite[Section $3$]{Schutt-Werner}, which is the definition used in \cite{Werner-Ye}. Ludwig \cite{Ludwig-general} introduces two families of general affine surface areas that contain $\Omega_{p}$ for $-N < p < 0$ \cite[Corollary $7$]{Ludwig-general} and, for $p < -N$, a functional defined by an integral over the sphere in place of the integral of \eqref{lpintegrand} over the boundary \cite[equation (5) and Corollary $11$]{Ludwig-general}. The values $\be = 0, 1/(n+2), 1/2, 1$ correspond to $p = 0, 1, N, \pm\infty$ which respectively yield the volume, the classical affine surface area, the self-dual centroaffine surface area, and the polar volume. The ranges $0 < \be < 1$, $\be < 0$, and $\be > 1$ correspond to $p > 0$, $-N < p < 0$, and $p < -N$, and no $\be$ corresponds to $p = -N$. The duality of Lemma \ref{dualityproposition} corresponds, under $\be \leftrightarrow 1 - \be$, that is, for nonzero finite $p$, $p \leftrightarrow N^{2}/p$, which exchanges the two ranges of negative $p$, to the identity $\Omega_{p}(K^{\circ}) = \Omega_{N^{2}/p}(K)$ for convex bodies, due to Hug for $p > 0$ \cite{Hug-polar}, \cite[Theorem $4$]{Ludwig-general}, and to Werner and Ye for all $p \neq -N$ and bodies with boundary of class $C^{2}$ and positive curvature \cite[Corollary $3.1$]{Werner-Ye}.

Since $(N-p)/(N+p) = 1 - 2\be$ and $\Omega_{0}(K) = N|K|$, the functional $\mathcal{T}_{\be} = \cS_{\be}\cS_{0}^{2\be - 1}$ of Section \ref{familysection} is the counterpart of the quantity $\Omega_{p}(K)(N|K|)^{-(N-p)/(N+p)}$ compared with its value for the ball in the $L_{p}$ affine isoperimetric inequalities. These are due to Lutwak \cite{Lutwak-Lp} for $p \geq 1$ and to Werner and Ye \cite[Theorem $4.2$]{Werner-Ye} for the remaining $p \neq -N$; see also \cite[Theorems $3$ and $8$]{Ludwig-general}, where the inequalities for $p > 0$ and for $-N < p < 0$ are obtained for all convex bodies with centroid at the origin from the Jensen and Blaschke--Santal\'o inequalities. For a convex body with centroid at the origin, this quantity is at most its value for the ball when $p \geq 0$, and at least that value when $-N < p < 0$, with equality in both cases exactly for ellipsoids when $p \neq 0$. When $p < -N$ and the boundary is of class $C^{2}$ with positive curvature, it is bounded below by its value for the ball multiplied by a power of the constant of the inequality of Bourgain and Milman \cite{Bourgain-Milman}, and this bound is not sharp.

The proofs in \cite{Werner-Ye} combine H\"older inequalities relating $\Omega_{p}(K)$, $|K|$, and $|K^{\circ}|$ with the Blaschke--Santal\'o inequality, when $p > -N$, and with the inequality of Bourgain and Milman, when $p < -N$. Here Theorem A takes the place of both, and it has the direction of the inequality of Bourgain and Milman. Accordingly, for $\be \geq 1$ the same H\"older inequality with Theorem A gives a sharp inequality with the affine sphere as the only extremal, \eqref{holder4} of Theorem \ref{reductionproposition}, whereas for $\be < 1$ the bound for $\mathcal{T}_{\be}$ by a power of the covolume product obtained in this way does not yield a bound by the value at the affine sphere, and what can be proved is described in Theorem \ref{reductionproposition}.

The comparison made in this remark is between integrands on hypersurfaces that are smooth and locally strongly convex. For general convex bodies the functionals $\Omega_{p}$ with $p > 0$ are finite, vanish on polytopes, and are upper semicontinuous, whereas those with $p < 0$ are infinite on polytopes and lower semicontinuous \cite[Theorem $1$ and Corollaries $7$ and $11$]{Ludwig-general}; the corresponding extensions of $\cS_{\be}$ from $\cAs(C, \Ga)$ to $\cA(C, \Ga)$ are not considered here.
\end{remark}

\section{Existence of minimizers and the proofs of Theorems A and B}\label{existencesection}
This section proves Theorems A and B. A minimizer of $\cP$ exists by the Harnack inequality (Lemma \ref{minimizerexistsproposition}); the variational lemma of Aleksandrov (Lemma \ref{aleksandrovlemma}) and the identification of the cone volume measure with a Monge--Amp\`ere measure (Lemma \ref{conevolumemeasurelemma}) show that its polar solves the affine sphere equation weakly (Lemma \ref{weakeulerlagrangeproposition}); the regularity theory of Caffarelli makes this solution smooth (Lemma \ref{regularityproposition}); and Theorem A then follows from Theorem B and the uniqueness of affine spheres (Lemma \ref{uniquenessproposition}).

\begin{lemma}\label{minimizerexistsproposition}
On $\cA(C, \Ga)$ the functional $\cP$ attains its infimum.
\end{lemma}

\begin{proof}
By \eqref{pdefined}, $\cP$ is positive, so $\inf_{\cA(C, \Ga)}\cP$ is finite and there is a minimizing sequence $\psi_{k} \in \cA(C, \Ga)$. Since a concave function on an open set is continuous and $\rho_{C}$ is smooth and positive, the ratio $\psi_{k}/\rho_{C}$ is positive, continuous, homogeneous of degree $0$, and $\Ga$-invariant, so it descends to a continuous function on the compact manifold $M$. Consequently $s_{k} = \sup_{x \in C}(\psi_{k}(x)/\rho_{C}(x))$ is finite and positive, and, because $\cP$ is scale invariant, replacing $\psi_{k}$ by $s_{k}^{-1}\psi_{k}$ yields a minimizing sequence satisfying $\sup_{x \in C}\left(\psi_{k}(x)/\rho_{C}(x)\right) = 1$. By the Harnack inequality \eqref{harnack}, $c\rho_{C} \leq \psi_{k} \leq \rho_{C}$ for all $k$.

Let $\Fset$ be a compact set as in Lemma \ref{fundamentalsetlemma}. Because the convex hull of $\Fset$ is a compact subset of the open convex set $C$, there is a bounded open convex set $U$ with $\Fset \subset U$ and $\bar{U} \subset C$. On $U$ the functions $-\psi_{k}$ are finite and convex, and they are pointwise bounded, because $-\rho_{C} \leq -\psi_{k} \leq 0$ and $\rho_{C}$ is bounded on $\bar{U}$. By \cite[Theorem 10.6]{Rockafellar}, the $-\psi_{k}$ are uniformly bounded and equi-Lipschitz on the closed bounded subset $\Fset$ of $U$, so, by the Arzel\`a--Ascoli theorem, after passing to a subsequence, which is again written $\psi_{k}$, the sequence $\psi_{k}$ converges uniformly on $\Fset$. Since $\rho_{C}$ is bounded below by a positive constant on the compact set $\Fset$, the ratios $\psi_{k}/\rho_{C}$ converge uniformly on $\Fset$, and since these ratios are $\Ga$-invariant and homogeneous of degree $0$, by \eqref{fundamentalset} they converge uniformly on $C$.

Hence $\psi_{k}$ converges locally uniformly on $C$ to a function $\psi$, which is concave, homogeneous of degree $1$, and $\Ga$-invariant, as a pointwise limit of functions with these properties, and satisfies $c\rho_{C} \leq \psi \leq \rho_{C}$, so that $\psi \in \cA(C, \Ga)$. Finally, $\psi_{k}/\psi = (\psi_{k}/\rho_{C})/(\psi/\rho_{C})$, where the numerator converges uniformly on $C$ to the denominator and the denominator is bounded below by $c > 0$, so $\psi_{k}/\psi \to 1$ uniformly on $C$, and \eqref{pcontinuity} gives $\cP(\psi) = \lim_{k \to \infty}\cP(\psi_{k}) = \inf_{\cA(C, \Ga)}\cP$.
\end{proof}
By Remark \ref{hausdorffremark}, Lemma \ref{minimizerexistsproposition} is an equivariant Blaschke selection theorem: the $\Ga$-convex bodies normalized by $\sup_{C}\psi/\rho_{C} = 1$ form a compact set for local Hausdorff convergence, on which $\cP$ is continuous. That a selection theorem of this kind should hold for Fuchsian convex bodies is remarked at the end of \cite[Section $5.2$]{Fillastre-Fuchsian}.

Lemma \ref{aleksandrovlemma} is an equivariant form of a variational lemma of Aleksandrov from his solution of the Minkowski problem, for which see \cite[Lemma $7.5.3$ and Section $8.2$]{Schneider-book}. The proof here follows the radial form of \cite[Lemma $4.3$ and Theorem $4.5$]{HLYZ}, and \cite{Li-Ye-Zhu} treats unbounded closed convex sets in cones.

\begin{lemma}\label{aleksandrovlemma}
Let $\varphi$ be a positive continuous $\Ga^{\ast}$-invariant function on $C^{\ast}$, homogeneous of degree $1$. Let $[\varphi] = \{x \in C: \zeta(x) \geq \varphi(\zeta) \text{ for all } \zeta \in C^{\ast}\}$ 
and let $\psi_{[\varphi]}(x) = \inf_{\zeta \in C^{\ast}}(\zeta(x)/\varphi(\zeta))$ for $x \in C$. Then:
\begin{enumerate}[label=$(\arabic*)$, ref=(\arabic*)]
\item\label{al1} $\psi_{[\varphi]} \in \cA(C, \Ga)$ and $[\varphi] = L_{\psi_{[\varphi]}}$. In particular $[\varphi]$ is a $\Ga$-convex body with cogauge $\psi_{[\varphi]}$.
\item\label{al2} $\psi^{\ast}_{[\varphi]} \geq \varphi$ on $C^{\ast}$, with equality if $\varphi = \psi^{\ast}$ for some $\psi \in \cA(C, \Ga)$, in which case $[\varphi] = L_{\psi}$ and $\psi_{[\varphi]} = \psi$.
\item\label{al3} $\psi_{[\varphi]}$ is differentiable at almost every $x \in C$. At such $x$, the set $[\varphi]$ has a unique supporting hyperplane at $p = x/\psi_{[\varphi]}(x)$. It has the form $\{y \in \ste: \zeta(y) = \zeta(p)\}$ for some $\zeta \in C^{\ast}$, which is determined up to a positive multiple, so the point $\al_{\varphi}(x) = \{t\zeta: t > 0\} \in \Om^{\ast}$ is well defined, and it satisfies $\al_{\varphi}(\gamma x) = \gamma^{-\ast}\al_{\varphi}(x)$ for $\gamma \in \Ga$.
\item\label{al4} Let $g$ be a continuous function on $M^{\ast}$ and let $\varphi_{t} = \varphi e^{tg}$ for $t \in \rea$. For every $x \in C$ and $|t| \leq 1$, $|\log\psi_{[\varphi_{t}]}(x) - \log\psi_{[\varphi]}(x)| \leq |t|\max_{M^{\ast}}|g|$, and at every $x$ at which $\psi_{[\varphi]}$ is differentiable,
\begin{align}\label{logradialderivative}
\tfrac{d}{dt}_{|t=0}\log\psi_{[\varphi_{t}]}(x) = -g(\al_{\varphi}(x)).
\end{align}
Consequently
\begin{align}\label{aleksandrovvariation}
\tfrac{d}{dt}_{|t=0}\int_{M}\psi_{[\varphi_{t}]}^{-N} = N\int_{M}(g\circ\al_{\varphi})\,\psi_{[\varphi]}^{-N}.
\end{align}
\end{enumerate}
\end{lemma}
\begin{proof}
First \ref{al1} is proved. The function $\psi_{[\varphi]}$ is an infimum of linear functions, so it is concave and homogeneous of degree $1$. It is $\Ga$-invariant because, for $\gamma \in \Ga$, the substitution $\zeta \mapsto \gamma^{\ast}\zeta$, which preserves $C^{\ast}$ and $\varphi$, shows $\psi_{[\varphi]}(\gamma x) = \psi_{[\varphi]}(x)$. It is finite and nonnegative because it is an infimum of positive numbers. Since $\varphi/\rho_{C^{\ast}}$ is the lift of a continuous function on the compact manifold $M^{\ast}$, there is $a > 0$ with $\varphi \leq a\,\rho_{C^{\ast}}$, so $\psi_{[\varphi]} \geq a^{-1}\rho^{\ast}_{C^{\ast}}$, where $\rho^{\ast}_{C^{\ast}}$ is the polar of $\rho_{C^{\ast}} \in \cAs(C^{\ast}, \Ga^{\ast})$, which is positive by Lemma \ref{polarmaplemma}\ref{pml4}. Hence $\psi_{[\varphi]} \in \cA(C, \Ga)$. By definition, $x \in [\varphi]$ if and only if $\zeta(x)/\varphi(\zeta) \geq 1$ for all $\zeta \in C^{\ast}$, that is, if and only if $\psi_{[\varphi]}(x) \geq 1$, so $[\varphi] = L_{\psi_{[\varphi]}}$. That $[\varphi]$ is a $\Ga$-convex body with cogauge $\psi_{[\varphi]}$ follows from Lemma \ref{gammaconvexbodylemma}.

Next \ref{al2} is proved. By the paragraph following the definition of polarity, $\psi^{\ast}_{[\varphi]}$ is the support function $s_{[\varphi]}(\zeta) = \inf_{x \in [\varphi]}\zeta(x)$, and $s_{[\varphi]} \geq \varphi$ because $\zeta(x) \geq \varphi(\zeta)$ for every $x \in [\varphi]$, by the definition of $[\varphi]$. If $\varphi = \psi^{\ast}$ with $\psi \in \cA(C, \Ga)$, then $[\varphi] = \{x \in C: \zeta(x) \geq \psi^{\ast}(\zeta) \text{ for all } \zeta \in C^{\ast}\} = L_{\psi^{\ast\ast}} = L_{\psi}$ by Lemma \ref{polartransformlemma}, so $\psi_{[\varphi]} = \psi$ and $\psi^{\ast}_{[\varphi]} = \varphi$.

Claims \ref{al3} and \ref{al4} are proved together. Throughout, $g$ is identified with its lift to $C^{\ast}$, which is continuous, $\Ga^{\ast}$-invariant, and homogeneous of degree $0$. Fix $x \in C$. Fix a norm on $\std$, let $\mathsf{Z}$ be the set of unit vectors in $\bar{C}^{\ast}$, and let $\mathsf{Z}^{\circ} = \mathsf{Z}\cap C^{\ast}$. Because $x \in C$ and $\mathsf{Z}$ is compact, the function $\mathsf{Z} \ni \zeta \mapsto \zeta(x)$ is bounded below by a positive constant. For $|t| \leq 1$ and $\zeta \in \mathsf{Z}^{\circ}$ let $F_{t}(\zeta) = \varphi_{t}(\zeta)/\zeta(x)$. Since $e^{-\max_{M^{\ast}}|g|}\varphi \leq \varphi_{t} \leq e^{\max_{M^{\ast}}|g|}a\,\rho_{C^{\ast}}$ for $|t| \leq 1$, and $\rho_{C^{\ast}} \to 0$ at $\pr C^{\ast}\setminus\{0\}$ by Lemma \ref{vinberglemma}, the positive continuous function $F_{t}$ extends continuously by $0$ to $\mathsf{Z}\setminus \mathsf{Z}^{\circ}$, with $F_{t}(\zeta) \leq e^{\max_{M^{\ast}}|g|}a\,\rho_{C^{\ast}}(\zeta)/\zeta(x)$ for all $\zeta \in \mathsf{Z}^{\circ}$ and $|t| \leq 1$, while $\max_{\mathsf{Z}} F_{t} \geq F_{t}(\zeta_{1}) \geq m := e^{-\max_{M^{\ast}}|g|}\varphi(\zeta_{1})/\zeta_{1}(x)$ for any fixed $\zeta_{1} \in \mathsf{Z}^{\circ}$. Hence $F_{t}$ attains its maximum on the compact set $\mathsf{Z}$, necessarily at a point of $\mathsf{Z}^{\circ}$, and every maximizer of $F_{t}$, for every $|t| \leq 1$, lies in
\begin{align}\label{zzerodefined}
\mathsf{Z}_{0} = \{\zeta \in \mathsf{Z}^{\circ}: e^{\max_{M^{\ast}}|g|}a\,\rho_{C^{\ast}}(\zeta)/\zeta(x) \geq m\},
\end{align}
which is a compact subset of $\mathsf{Z}^{\circ}$, depending on $x$ but not on $t$, because it is closed in $\mathsf{Z}$ and disjoint from $\mathsf{Z}\setminus\mathsf{Z}^{\circ}$. Since $\zeta \mapsto \zeta(x)/\varphi_{t}(\zeta)$ is homogeneous of degree $0$, its infimum over $C^{\ast}$ equals its infimum over $\mathsf{Z}^{\circ}$, which is $1/\max_{\mathsf{Z}_{0}}F_{t}$, so
\begin{align}\label{thetadefined}
\Theta(t, x) = -\log\psi_{[\varphi_{t}]}(x) = \max_{\zeta \in \mathsf{Z}_{0}}\left(\log \varphi(\zeta) + tg(\zeta) - \log\zeta(x)\right).
\end{align}
As a maximum of affine functions of $t$ with slopes $g(\zeta)$, $\Theta(\cdot, x)$ is convex and $(\max_{M^{\ast}}|g|)$-Lipschitz on $[-1, 1]$, which gives the first claim of \ref{al4}.

A \emph{supporting hyperplane} of $[\varphi]$ at $p$ is a set $\{y \in \ste: \eta(y) = c\}$ with $\eta \in \std\setminus\{0\}$ and $c \in \rea$ such that $\eta(p) = c$ and $\eta(y) \geq c$ for all $y \in [\varphi]$; the pair $(\eta, c)$ is determined by the hyperplane up to a common positive factor. Since $\eta$ is then bounded below on $[\varphi]$, Lemma \ref{gammaconvexbodylemma} gives $\eta \in \bar{C}^{\ast}$, and since $\eta \neq 0$ and $p \in C$, $c = \eta(p) > 0$, so every supporting hyperplane of $[\varphi]$ at $p$ can be written as $\{y \in \ste: \eta(y) = 1\}$ with $\eta \in \bar{C}^{\ast}$, $\eta(p) = 1$, and $\eta(y) \geq 1$ for all $y \in [\varphi]$.

Let $p = x/\psi_{[\varphi]}(x)$. If $\zeta \in \mathsf{Z}_{0}$ is a maximizer at $t = 0$ of the expression on the right-hand side of \eqref{thetadefined}, then $\zeta(x)/\varphi(\zeta) = \psi_{[\varphi]}(x)$, so $\zeta(p) = \varphi(\zeta)$, and $\zeta(y) \geq \varphi(\zeta)$ for all $y \in [\varphi]$ by the definition of $[\varphi]$, so $\{y \in \ste: \zeta(y) = \varphi(\zeta)\}$ is a supporting hyperplane of $[\varphi]$ at $p$.

On the other hand, let $\eta \in \std$ be a supergradient of $\psi_{[\varphi]}$ at $x$, meaning $\psi_{[\varphi]}(y) \leq \psi_{[\varphi]}(x) + \eta(y - x)$ for all $y \in C$. Taking $y = sx$ with $s > 0$ and using homogeneity gives $(s - 1)\psi_{[\varphi]}(x) \leq (s-1)\eta(x)$ for all $s > 0$, hence $\eta(x) = \psi_{[\varphi]}(x)$, hence $\eta(y) \geq \psi_{[\varphi]}(y)$ for all $y \in C$, hence $\eta(y) \geq 1 = \eta(p)$ for all $y \in [\varphi] = \{y' \in C: \psi_{[\varphi]}(y') \geq 1\}$, so $\{y \in \ste: \eta(y) = 1\}$ is a supporting hyperplane of $[\varphi]$ at $p$. Conversely, if $\{y \in \ste: \eta(y) = 1\}$ is a supporting hyperplane of $[\varphi]$ at $p$, normalized as above, then for $y \in C$ the point $y/\psi_{[\varphi]}(y)$ lies in $[\varphi]$, so $\eta(y) \geq \psi_{[\varphi]}(y)$, while $\eta(x) = \psi_{[\varphi]}(x)\eta(p) = \psi_{[\varphi]}(x)$, so $\eta$ is a supergradient of $\psi_{[\varphi]}$ at $x$. Thus $\eta \mapsto \{y \in \ste: \eta(y) = 1\}$ is a bijection from the supergradients of $\psi_{[\varphi]}$ at $x$ to the supporting hyperplanes of $[\varphi]$ at $p$.

The concave function $\psi_{[\varphi]}$ is differentiable almost everywhere, by \cite[Theorem 25.5]{Rockafellar} applied to $-\psi_{[\varphi]}$, and it is differentiable at $x$ if and only if it has a unique supergradient there \cite[Theorem 25.1]{Rockafellar}, that is, if and only if $[\varphi]$ has a unique supporting hyperplane at $p$. At such $x$ the unique supergradient is $d\psi_{[\varphi]}(x)$.

Distinct elements of $\mathsf{Z}_{0}$ are unit vectors in $C^{\ast}$, so no one is a positive multiple of another, and they determine distinct supporting hyperplanes $\{y \in \ste: \zeta(y) = \varphi(\zeta)\}$ at $p$; hence at such $x$ the maximizer $\zeta$ at $t = 0$ is unique, and $\{y \in \ste: \zeta(y) = \varphi(\zeta)\} = \{y \in \ste: d\psi_{[\varphi]}(x)(y) = 1\}$, so $\zeta = \varphi(\zeta)\,d\psi_{[\varphi]}(x)$ with $\varphi(\zeta) > 0$. Therefore $d\psi_{[\varphi]}(x) \in C^{\ast}$, and $\al_{\varphi}(x)$, the ray of $d\psi_{[\varphi]}(x)$, lies in $\Om^{\ast}$.

The map $\al_{\varphi}$ is measurable, being defined almost everywhere as the composition of the measurable map $d\psi_{[\varphi]}$ with the projection $C^{\ast} \to \Om^{\ast}$, and differentiating $\psi_{[\varphi]}(\gamma y) = \psi_{[\varphi]}(y)$ at $y = x$ gives $d\psi_{[\varphi]}(\gamma x)\circ\gamma = d\psi_{[\varphi]}(x)$, that is, $d\psi_{[\varphi]}(\gamma x) = \gamma^{-\ast}d\psi_{[\varphi]}(x)$, which is the equivariance in \ref{al3}. This proves \ref{al3}.

For \eqref{logradialderivative}, let $x$ be a point at which $\psi_{[\varphi]}$ is differentiable and let $\zeta_{0} \in \mathsf{Z}_{0}$ be the unique maximizer at $t = 0$, which spans $\al_{\varphi}(x)$. For $|t| \leq 1$ let $\zeta_{t} \in \mathsf{Z}_{0}$ be a maximizer at time $t$. Evaluating the expression in \eqref{thetadefined} at $\zeta_{0}$ gives $\Theta(t, x) \geq \Theta(0, x) + tg(\zeta_{0})$, and evaluating it at $\zeta_{t}$ and using that $\zeta_{t}$ is not a better choice than $\zeta_{0}$ at $t = 0$ gives $\Theta(t, x) \leq \Theta(0, x) + tg(\zeta_{t})$. Hence, for $0 < |t| \leq 1$, the difference quotient $(\Theta(t, x) - \Theta(0, x))/t$ lies between $g(\zeta_{0})$ and $g(\zeta_{t})$. Every limit point of $\zeta_{t}$ as $t \to 0$ lies in the compact set $\mathsf{Z}_{0}$ and, by the continuity of $\Theta(\cdot, x)$ and of the expression in \eqref{thetadefined} in $(t, \zeta)$, is a maximizer at $t = 0$, so $\zeta_{t} \to \zeta_{0}$, and the continuity of $g$ gives $\tfrac{d}{dt}_{|t=0}\Theta(t, x) = g(\zeta_{0}) = g(\al_{\varphi}(x))$, which is \eqref{logradialderivative}.

Finally, $g\circ\al_{\varphi}$ is a $\Ga$-invariant function homogeneous of degree $0$ defined almost everywhere on $C$, so it is a measurable function on $M$, and $\psi_{[\varphi_{t}]}^{-N} = e^{N(\Theta(t, \cdot) - \Theta(0, \cdot))}\psi_{[\varphi]}^{-N}$. Since $|e^{u} - 1| \leq |u|e^{|u|}$ and $|\Theta(t, \cdot) - \Theta(0, \cdot)| \leq |t|\max_{M^{\ast}}|g|$, the difference quotients in $t$ of the exponential factor are bounded by $N\max_{M^{\ast}}|g|\,e^{N\max_{M^{\ast}}|g|}$ for $0 < |t| \leq 1$, and they converge almost everywhere to $Ng\circ\al_{\varphi}$ by \eqref{logradialderivative}; since $\psi_{[\varphi]}^{-N}$ is continuous on the compact manifold $M$, \eqref{aleksandrovvariation} follows by dominated convergence.
\end{proof}

\begin{lemma}\label{conevolumemeasurelemma}
Choose linear coordinates $x = (y, t)$ on $\ste$ with $\Psi = dy^{1}\wedge\dots\wedge dy^{n}\wedge dt$, such that $(0, 1) \in C$, and dual coordinates $\zeta = (\eta, s)$ with $\zeta(x) = \lb\eta, y\ra + st$. Then the slice $\{(\eta, s) \in \std: s = 1\}$ meets every ray of $C^{\ast}$, and $\Om^{\ast}_{0} = \{\eta \in \rea^{n}: (\eta, 1) \in C^{\ast}\}$ is a bounded convex domain. For $\psi \in \cA(C, \Ga)$ let $L = L_{\psi}$ and $w(\eta) = -\psi^{\ast}(\eta, 1)$. Then $w$ is a negative convex function on $\Om^{\ast}_{0}$. Let $\MA(w)(Z) = |\pr w(Z)|$ be its Monge--Amp\`ere measure in the sense of Aleksandrov. The map $\al_{\psi^{\ast}}$ of Lemma \ref{aleksandrovlemma}\ref{al3} is defined almost everywhere on $\Om$ and takes values in $\Om^{\ast}$. Let $(\al_{\psi^{\ast}})_{\ast}(\psi^{-N})$ be the pushforward under $\al_{\psi^{\ast}}$ of the $\Ga$-invariant measure on $\Om$ determined by the density $\psi^{-N}|\iota_{\eul}\Psi|$, regarded as a measure on $\Om^{\ast}_{0}$ by identifying directions of $C^{\ast}$ with points of the slice. Then
\begin{align}\label{conevolumemeasure}
(\al_{\psi^{\ast}})_{\ast}(\psi^{-N}) = (-w)\MA(w).
\end{align}
\end{lemma}
\begin{proof}
If $\zeta \in C^{\ast}$ then $s = \zeta((0, 1)) > 0$, so every ray of $C^{\ast}$ meets the slice, and $\Om^{\ast}_{0}$ is bounded because $(0,1)$ lies in the interior of $C$. Since $L + C \subset L$, $(0,1) \in C$, and $L$ is closed by Lemma \ref{gammaconvexbodylemma}, $L$ is the epigraph $\{(y, t): t \geq \tau(y)\}$ of a lower semicontinuous convex function $\tau$ on the set of those $y \in \rea^{n}$ for which $(y, t) \in L$ for some $t$. This set is all of $\rea^{n}$. Indeed, let $(y_{0}, t_{0}) \in L$ and $y \in \rea^{n}$. Since $(0, 1)$ is an interior point of the open cone $C$, for $T$ sufficiently large $(y - y_{0}, T - t_{0}) = (T - t_{0})\big((y - y_{0})/(T - t_{0}), 1\big) \in C$, so $(y, T) \in L + C \subset L$. Moreover $\tau(y) > -\infty$ for every $y$, because $L \subset C$ and $\bar{C}$ contains no line. Thus $\tau$ is a finite convex function on $\rea^{n}$, and in particular it is locally Lipschitz. Then
\begin{align}
w(\eta) = -\inf_{(y,t) \in L}\left(\lb\eta, y\ra + t\right) = \sup_{y \in \rea^{n}}\left(\lb -\eta, y\ra - \tau(y)\right) = \tau^{\ast}(-\eta),
\end{align}
where $\tau^{\ast}$ is the Legendre transform, so $w$ is convex, and $w < 0$ because $\psi^{\ast} > 0$. For almost every direction $x$, the boundary point $x/\psi(x)$ has a unique supporting hyperplane, whose normal lies in $C^{\ast}$ by Lemma \ref{aleksandrovlemma}\ref{al3} applied with $\varphi = \psi^{\ast}$. Writing $x/\psi(x) = (y, \tau(y))$, it follows that $\tau$ is differentiable at $y$, because the supporting hyperplane is unique.

Since $\psi$ is concave, it is locally Lipschitz on $C$, so the map $x \mapsto x/\psi(x)$ from $C$ to $\pr L$ is locally Lipschitz, and its composition with the projection $(y, t) \mapsto y$ is a locally Lipschitz map from directions to $\rea^{n}$. In the other direction, $y \mapsto (y, \tau(y))$ is locally Lipschitz, and its composition with the projection of $C$ onto its directions, which is locally Lipschitz, is a locally Lipschitz map from $\rea^{n}$ to directions. These two maps are inverse to each other, so the correspondence between a direction $x$ and the point $y \in \rea^{n}$ such that $x/\psi(x) = (y, \tau(y))$ is locally bi-Lipschitz, and it identifies null sets of directions with null sets of $y$. At such a point $p = (y, \tau(y))$ the inner normal direction is $\eta(y) = -\fD\tau(y)$, and $-w(\eta(y)) = \lb\eta(y), y\ra + \tau(y) = \tau(y) - \lb y, \fD\tau(y)\ra$. Since $w(\eta) = \tau^{\ast}(-\eta)$ and, $\tau$ being lower semicontinuous, $y \in \pr\tau^{\ast}(\xi)$ if and only if $\xi \in \pr\tau(y)$ \cite[Theorem $23.5$]{Rockafellar}, the subdifferential image $\pr w(Z)$ of a Borel set $Z \subset \Om^{\ast}_{0}$ is $-\{y : \eta(y) \in Z\}$ up to the null set of $y$ at which $\tau$ is not differentiable, and Lebesgue measure is invariant under $y \mapsto -y$, so $\MA(w) = \eta_{\ast}(dy)$.

On the other hand, let $B$ be a Borel subset of $\rea^{n}$ consisting of points at which $\tau$ is differentiable. The map $(s, y) \mapsto (sy, s\tau(y))$ from $(0, 1)\times B$ to $\ste$ is injective, because $\pr L$ meets each ray of $C$ exactly once, and it is locally Lipschitz, because $\tau$ is locally Lipschitz. Its Jacobian determinant at $(s, y)$ is $s^{n}(\tau(y) - \lb y, \fD\tau(y)\ra)$, and its image is the cone from the origin over the part of $\pr L$ above $B$. By the change of variables formula for injective locally Lipschitz maps, the volume of this cone is 
\begin{align}\label{conevolumegraph}
\int_{B}\int_{0}^{1}s^{n}\left(\tau(y) - \lb y, \fD\tau(y)\ra\right)\,ds\,dy = \frac{1}{N}\int_{B}\left(\tau(y) - \lb y, \fD\tau(y)\ra\right)\,dy.
\end{align}
By \eqref{conevolumelocal} with $f = 1$, the volume of this cone is also $\tfrac{1}{N}$ times the measure assigned by the density $\psi^{-N}|\iota_{\eul}\Psi|$ on $\Om$ to the set $E \subset \Om$ of the directions of the points of $\pr L$ above $B$. For a Borel set $Z \subset \Om^{\ast}_{0}$, taking $B = \{y: \eta(y) \in Z\}$, for which $E$ equals $\al_{\psi^{\ast}}^{-1}(Z)$ up to a null set, and using $\tau(y) - \lb y, \fD\tau(y)\ra = -w(\eta(y))$ and $\MA(w) = \eta_{\ast}(dy)$, gives $(\al_{\psi^{\ast}})_{\ast}(\psi^{-N})(Z) = \int_{B}(-w(\eta(y)))\,dy = \int_{Z}(-w)\,d\MA(w)$.
\end{proof}

\begin{lemma}\label{weakeulerlagrangeproposition}
Let $\psi \in \cA(C, \Ga)$ minimize $\cP$, and, in the coordinates of Lemma \ref{conevolumemeasurelemma}, let $w(\eta) = -\psi^{\ast}(\eta, 1)$. Then $w$ is a negative convex Aleksandrov solution of
\begin{align}\label{weakaffinesphere}
\det \fD^{2}w = c_{\psi}(-w)^{-n-2} \,\,\text{in}\,\, \Om^{\ast}_{0}, \qquad c_{\psi} = \frac{|K_{\psi}/\Ga|}{|K^{\ast}_{\psi}/\Ga^{\ast}|},
\end{align}
meaning that $\MA(w) = c_{\psi}(-w)^{-n-2}d\eta$ as Borel measures on $\Om^{\ast}_{0}$, and $w(\eta) \to 0$ as $\eta \to \pr\Om^{\ast}_{0}$.
\end{lemma}
\begin{proof}
For $\varphi$ as in Lemma \ref{aleksandrovlemma} let $\widetilde{\cP}(\varphi) = \tfrac{1}{N^{2}}\int_{M}\psi_{[\varphi]}^{-N}\int_{M^{\ast}}\varphi^{-N}$. Since $\psi^{\ast}_{[\varphi]} \geq \varphi$ by Lemma \ref{aleksandrovlemma}\ref{al2}, \eqref{pdefined} gives $\widetilde{\cP}(\varphi) \geq \cP(\psi_{[\varphi]}) \geq \cP(\psi)$, and $\widetilde{\cP}(\psi^{\ast}) = \cP(\psi)$ by the equality case of Lemma \ref{aleksandrovlemma}\ref{al2}. Hence, for every continuous $g$ on $M^{\ast}$, $t \mapsto \widetilde{\cP}(\psi^{\ast}e^{tg})$ has a minimum at $t = 0$. It is differentiable at $t = 0$, because, by Lemma \ref{aleksandrovlemma}\ref{al2}, $\psi_{[\psi^{\ast}]} = \psi$, so that the first factor $\int_{M}\psi_{[\psi^{\ast}e^{tg}]}^{-N}$ is differentiable at $t = 0$ by \eqref{aleksandrovvariation}, while the second factor $\int_{M^{\ast}}(\psi^{\ast})^{-N}e^{-Ntg}$ is a smooth function of $t \in \rea$. Hence
\begin{align}
0 = \tfrac{d}{dt}_{|t=0}\log \widetilde{\cP}(\psi^{\ast}e^{tg}) = N\frac{\int_{M}(g\circ\al_{\psi^{\ast}})\psi^{-N}}{\int_{M}\psi^{-N}} - N\frac{\int_{M^{\ast}}g(\psi^{\ast})^{-N}}{\int_{M^{\ast}}(\psi^{\ast})^{-N}}.
\end{align}
Here, by Lemma \ref{aleksandrovlemma}\ref{al3}, $\al_{\psi^{\ast}}$ is defined almost everywhere on $\Om$ and is equivariant, so it induces a map, also written $\al_{\psi^{\ast}}$, defined almost everywhere on $M$ and taking values in $M^{\ast}$. Since $g$ is an arbitrary continuous function on the compact manifold $M^{\ast}$, and a finite Borel measure on $M^{\ast}$ is determined by the integrals against it of the continuous functions, there holds $(\al_{\psi^{\ast}})_{\ast}(\psi^{-N}) = c_{\psi} (\psi^{\ast})^{-N}$ as measures on $M^{\ast}$, where $c_{\psi} = \int_{M}\psi^{-N}/\int_{M^{\ast}}(\psi^{\ast})^{-N}$ equals the constant $c_{\psi}$ of \eqref{weakaffinesphere} by \eqref{pdefined}.

Because $\al_{\psi^{\ast}}$ is $\Ga$-equivariant, the pushforward under $\al_{\psi^{\ast}}$ of the $\Ga$-invariant measure $\psi^{-N}$ on $\Om$ is a $\Ga^{\ast}$-invariant measure on $\Om^{\ast}$ whose quotient is $(\al_{\psi^{\ast}})_{\ast}(\psi^{-N})$ on $M^{\ast}$, so the identity holds also as an identity of $\Ga^{\ast}$-invariant measures on $\Om^{\ast}$ (if a Borel set $Z \subset \Om^{\ast}$ meets each $\Ga^{\ast}$-orbit at most once, then, by the equivariance of $\al_{\psi^{\ast}}$ and the freeness of the action of $\Ga$, $\al_{\psi^{\ast}}^{-1}(Z)$ meets each $\Ga$-orbit at most once, so both sides assign to $Z$ the mass they assign to its image in $M^{\ast}$), and it can be verified on the slice of Lemma \ref{conevolumemeasurelemma}, where neither group appears.

By Lemma \ref{conevolumemeasurelemma} the left-hand side is $(-w)\MA(w)$, and, by the slice description of densities used in Lemma \ref{chartlemma}, $(\psi^{\ast})^{-N}$ is $(-w)^{-N}d\eta$. Hence $(-w)\MA(w) = c_{\psi}(-w)^{-N}d\eta$, which is \eqref{weakaffinesphere} because $N + 1 = n + 2$. The boundary behavior follows from Lemma \ref{harnacklemma} applied to $\psi^{\ast} \in \cA(C^{\ast}, \Ga^{\ast})$.
\end{proof}

Lemma \ref{weakeulerlagrangeproposition} says that a minimizer of $\cP$ solves, in the $\Ga$-equivariant setting, the $L_{-N}$ Minkowski problem for $C$-close sets with constant data; the $L_{p}$ Minkowski problem for $C$-close sets without an equivariance hypothesis is studied in \cite{Ai-Ye-Zhu, Yang-Ye-Zhu}. The measure $(\al_{\psi^{\ast}})_{\ast}(\psi^{-N})$ of Lemma \ref{conevolumemeasurelemma} is the equivariant counterpart of the cone-volume measure of a pseudo-cone, and, without a group, the measures that arise as weighted cone-volume measures of $C$-pseudo-cones are characterized in \cite{Schneider-conevolume}, while a Minkowski existence theorem for weighted surface area measures of $C$-pseudo-cones is proved in \cite{Schneider-weightedminkowski}. In the Fuchsian case of Example \ref{hyperboloidexample}, the classical equivariant Minkowski problem is solved in \cite{Oliker-Simon}, and in \cite{Bonsante-Fillastre} for convex sets in Minkowski space invariant under affine deformations of cocompact lattices of $SO^{+}(n, 1)$.

\begin{lemma}\label{regularityproposition}
Let $\psi \in \cA(C, \Ga)$ minimize $\cP$, and let $w$ and $c_{\psi}$ be as in Lemma \ref{weakeulerlagrangeproposition}. Then $w$ is smooth with $\fD^{2}w > 0$. Consequently $\psi^{\ast} \in \cAs(C^{\ast}, \Ga^{\ast})$ and $\mu(\psi^{\ast}) = c_{\psi}$.
\end{lemma}
\begin{proof}
For $n = 1$, the Monge--Amp\`ere measure of the convex function $w$ on the interval $\Om^{\ast}_{0}$ is its second derivative in the sense of distributions, so \eqref{weakaffinesphere} says that $w^{\prime\prime}$ equals the positive continuous function $c_{\psi}(-w)^{-3}$. Hence $w$ is of class $C^{2}$ with $w^{\prime\prime} > 0$, and differentiating the equation $w^{\prime\prime} = c_{\psi}(-w)^{-3}$ repeatedly shows that $w$ is smooth.

Henceforth it is assumed that $n \geq 2$, as is needed to apply the results of Caffarelli used below. By Lemma \ref{harnacklemma}, on each compact subset of $\Om^{\ast}_{0}$ the function $-w$ is bounded above and below by positive constants, so the right-hand side of \eqref{weakaffinesphere} is bounded between positive constants on compact subsets. Suppose that for some supporting affine function $\ell$ of $w$ the contact set $\{\eta \in \Om^{\ast}_{0}: w(\eta) = \ell(\eta)\}$ contains more than one point. The contact set $S$ is convex and relatively closed in $\Om^{\ast}_{0}$.

It is shown first that $S$ has no extreme points in $\Om^{\ast}_{0}$. Suppose that $p \in \Om^{\ast}_{0}$ is an extreme point of $S$, and choose a convex domain $\Om'$ with $p \in \Om'$ and $\bar{\Om'} \subset \Om^{\ast}_{0}$. The Monge--Amp\`ere measure of $w$ is bounded above and below by positive multiples of Lebesgue measure on $\Om'$. The set $S \cap \Om'$ is the contact set of $\ell$ and $w|_{\Om'}$. It contains more than one point, because $S$ is convex with more than one point and $\Om'$ contains a neighborhood of $p$ in $S$, and $p$ is an extreme point of $S \cap \Om'$. This contradicts \cite[Theorem $1$]{Caffarelli-strictconvexity}, applied to $w - \ell$ on $\Om'$, which says that if the Monge--Amp\`ere measure of a convex function is bounded above and below by positive multiples of Lebesgue measure, then its contact set with a supporting affine function, if not a single point, has no extreme points in the interior of the domain. (That theorem is stated for viscosity solutions. Aleksandrov solutions of the two-sided inequality are viscosity solutions, as noted in the introduction of \cite{Caffarelli-strictconvexity}.)

Extend $w$ continuously to $\bar{\Om^{\ast}_{0}}$ by $0$, which is possible by Lemma \ref{weakeulerlagrangeproposition}. The closure $\bar{S}$ is a compact convex set, equal to the convex hull of its extreme points. Since $S$ is relatively closed, $\bar{S}\cap\Om^{\ast}_{0} = S$, so an extreme point of $\bar{S}$ in $\Om^{\ast}_{0}$ would be an extreme point of $S$; hence all extreme points of $\bar{S}$ lie on $\pr\Om^{\ast}_{0}$, and $\ell = w = 0$ at these points by continuity. Since $\ell$ is affine, $\ell = 0$ on $\bar{S}$, so $w = 0$ on $S$, contradicting $w < 0$. Hence $w$ is strictly convex, meaning that each of its supporting affine functions touches it at exactly one point (compare \cite[Corollary $2$]{Caffarelli-strictconvexity}, which gives the strict convexity of a solution with boundary values $0$ when the right-hand side is bounded above and below by positive constants on the whole domain; here the right-hand side of \eqref{weakaffinesphere} is not bounded above near $\pr\Om^{\ast}_{0}$, which is why the argument has been given).

Since $w$ is convex and finite on $\Om^{\ast}_{0}$, it is locally Lipschitz. Consequently, the right-hand side $c_{\psi}(-w)^{-n-2}$ of \eqref{weakaffinesphere} is locally Lipschitz on $\Om^{\ast}_{0}$, and it is positive. Let $\eta_{0} \in \Om^{\ast}_{0}$ and let $\ell$ be a supporting affine function of $w$ at $\eta_{0}$. By the continuity of $w$ on $\bar{\Om^{\ast}_{0}}$, $w \geq \ell$ on $\bar{\Om^{\ast}_{0}}$. If $w(q) = \ell(q)$ for some $q \in \bar{\Om^{\ast}_{0}}$ different from $\eta_{0}$, then, by the convexity of $w$, $w = \ell$ on the segment joining $\eta_{0}$ to $q$, which contains points of $\Om^{\ast}_{0}$ different from $\eta_{0}$, contradicting the strict convexity of $w$. Hence $w - \ell$ is positive on $\bar{\Om^{\ast}_{0}}\setminus\{\eta_{0}\}$, so $m = \min_{\pr\Om^{\ast}_{0}}(w - \ell) > 0$, by the compactness of $\pr\Om^{\ast}_{0}$, and, for $0 < t < m$, the closure of the section $S_{t} = \{\eta \in \Om^{\ast}_{0}: w(\eta) < \ell(\eta) + t\}$, which is an open convex set containing $\eta_{0}$, is contained in $\Om^{\ast}_{0}$, so that $w - \ell - t$ vanishes on $\pr S_{t}$. 

The interior estimate \cite[Theorem $2$]{Caffarelli-interior} is stated for a convex viscosity solution of $\det\fD^{2}u = f$, with $f$ of class $C^{\al}$, on a convex domain normalized by an affine change of variables, that vanishes on the boundary of the domain. Applied to $u = w - \ell - t$ on $S_{t}$, it shows that $w$ is locally of class $C^{2,\al}$ on $S_{t}$, which is a neighborhood of $\eta_{0}$. Since $\eta_{0} \in \Om^{\ast}_{0}$ is arbitrary, $w \in C^{2, \al}_{\mathrm{loc}}(\Om^{\ast}_{0})$, and $\fD^{2}w > 0$ because $w$ is convex and $\det\fD^{2}w > 0$. Schauder theory applied to the linearized equation gives $w \in \cinf(\Om^{\ast}_{0})$. 

The function $-\psi^{\ast}$ is smooth with $-\fD^{2}\psi^{\ast}$ positive semidefinite with radical spanned by $\eul$, since in the basis of Lemma \ref{chartlemma} it is block diagonal with blocks $\fD^{2}w$ and $0$. So $\psi^{\ast} \in \cAs(C^{\ast}, \Ga^{\ast})$, and Lemma \ref{chartlemma} gives $\mu(\psi^{\ast}) = (-w)^{n+2}\det \fD^{2}w = c_{\psi}$.
\end{proof}

\begin{proof}[Proof of Theorem B]
By Lemma \ref{gammaconvexbodylemma}, the map $L \mapsto \psi_{L}$ is a bijection from the set $\cK(C, \Ga)$ of $\Ga$-convex bodies in $C$ onto $\cA(C, \Ga)$, with inverse $\psi \mapsto L_{\psi}$, and $\cP(L) = \cP(\psi_{L})$ by the remark following \eqref{pdefined}. Hence a $\Ga$-convex body $L$ minimizes $\cP$ on $\cK(C, \Ga)$ if and only if its cogauge $\psi_{L}$ minimizes $\cP$ on $\cA(C, \Ga)$. By Lemma \ref{minimizerexistsproposition}, $\cP$ attains its infimum on $\cA(C, \Ga)$.

Let $\psi \in \cA(C, \Ga)$ be any minimizer. By Lemma \ref{regularityproposition}, $\psi^{\ast} \in \cAs(C^{\ast}, \Ga^{\ast})$ with $\mu(\psi^{\ast})$ constant. By Lemma \ref{polarmaplemma}\ref{pml4} and Lemma \ref{dualityproposition}, applied on $M^{\ast}$, $\psi^{\ast\ast} \in \cAs(C, \Ga)$ and $\mu(\psi^{\ast\ast})\circ G_{\psi^{\ast}} = \mu(\psi^{\ast})^{-1}$. Since $G_{\psi^{\ast}}$ maps $C^{\ast}$ onto $C$, by Lemma \ref{polarmaplemma}\ref{pml3} applied on $M^{\ast}$, it follows that $\mu(\psi^{\ast\ast})$ is constant on $C$, and $\psi^{\ast\ast} = \psi$ by Lemma \ref{polartransformlemma}. Thus every minimizer $\psi$ lies in $\cAs(C, \Ga)$ and has constant $\mu(\psi)$. By Lemma \ref{smoothboundarylemma} the boundary $\Sigma_{\psi}$ of $L_{\psi}$ is a smooth, locally strongly convex hypersurface, and by Lemma \ref{affinespherelemma} it is a hyperbolic affine sphere centered at the origin. It is $\Ga$-invariant because $\psi$ is, and it is asymptotic to $\pr C$ by the remark following Definition \ref{asymptoticdefinition}.

Finally, if $\psi_{1}, \psi_{2} \in \cAs(C, \Ga)$ are minimizers, then, by Lemma \ref{uniquenessproposition}, $\psi_{2} = c\psi_{1}$ for a constant $c > 0$, so that $L_{\psi_{2}} = c^{-1}L_{\psi_{1}}$, which shows that any two minimizers are homothetic.
\end{proof}

\begin{proof}[Proof of Theorem A]
Let $L_{0}$ be the closed convex subset of $C$ bounded by $\Sigma_{0}$, in the sense of Definition \ref{asymptoticdefinition}. By Lemma \ref{generaluniquenesslemma}\ref{gul2}, $\Sigma_{0}$ is $\Ga$-invariant, so $L_{0}$ is $\Ga$-invariant, and $L_{0}$ is a $\Ga$-convex body by Lemma \ref{gammaconvexbodylemma}. Let $\psi_{1} \in \cA(C, \Ga)$ be the cogauge of $L_{0}$. Since $\pr L_{0} = \Sigma_{0}$ is a smooth, locally strongly convex hypersurface, Lemma \ref{smoothboundarylemma} shows that $\psi_{1} \in \cAs(C, \Ga)$ and $\Sigma_{\psi_{1}} = \Sigma_{0}$, so that the metric $h$ of Theorem A is $\gb_{\psi_{1}}$. Since $\Sigma_{0}$ is a hyperbolic affine sphere centered at the origin, by Lemma \ref{affinespherelemma}, $\mu(\psi_{1})$ is equal to a positive constant $\ka$. By \eqref{pdefined}, \eqref{volumeformulas}, and Lemma \ref{centroaffinemetriclemma}\ref{cml1}, 
\begin{align}\label{theoremAminimum}
\begin{aligned}
\cP(\psi_{1}) &= \frac{1}{N^{2}}\int_{M}\psi_{1}^{-N}\int_{M}\mu(\psi_{1})\psi_{1}^{-N} = \frac{\ka}{N^{2}}\Big(\int_{M}\psi_{1}^{-N}\Big)^{2}\\
& = \frac{1}{N^{2}}\Big(\int_{M}\mu(\psi_{1})^{1/2}\psi_{1}^{-N}\Big)^{2} = \frac{\vol(M, h)^{2}}{N^{2}}.
\end{aligned}
\end{align}
By Theorem B there is $\psi_{0} \in \cAs(C, \Ga)$ that minimizes $\cP$ on $\cA(C, \Ga)$ and for which $\mu(\psi_{0})$ is constant. By Lemma \ref{uniquenessproposition}, $\psi_{1}$ is a positive constant multiple of $\psi_{0}$. Since $\cP(c\psi) = \cP(\psi)$ for $c > 0$ and $\psi \in \cA(C, \Ga)$, it follows that $\psi_{1}$ also minimizes $\cP$ on $\cA(C, \Ga)$. Hence, for every $\Ga$-convex body $L \subset C$, with cogauge $\psi_{L} \in \cA(C, \Ga)$, by \eqref{theoremAminimum},
\begin{align}\label{theoremAproofinequality}
\cP(L) = \cP(\psi_{L}) \geq \cP(\psi_{1}) = \frac{\vol(M, h)^{2}}{N^{2}}.
\end{align}
If $\pr L$ is a homothetic image of $\Sigma_{0}$, then $\psi_{L}$ is a positive constant multiple of $\psi_{1}$, and equality holds in \eqref{theoremAproofinequality} by the invariance of $\cP$ under rescaling. Conversely, if equality holds in \eqref{theoremAproofinequality}, then $\psi_{L}$ minimizes $\cP$ on $\cA(C, \Ga)$, so, by Theorem B, $\psi_{L}$ is a positive constant multiple of the minimizer $\psi_{1}$, and $\pr L = \Sigma_{\psi_{L}}$ is a homothetic image of $\Sigma_{\psi_{1}} = \Sigma_{0}$.
\end{proof}

\begin{corollary}\label{fuchsiancorollary}
Suppose that $\Ga$ preserves a Lorentzian quadratic form $\lb\,,\,\ra$ on $\ste$ of signature $(n, 1)$ and that $C$ is a component of $\{x \in \ste: \lb x, x\ra < 0\}$, so that, when $\Ga \subset SO^{+}(n, 1)$, the $\Ga$-convex bodies are the Fuchsian convex bodies of \cite{Fillastre-Fuchsian}. Let $\vol_{\mathrm{hyp}}(M)$ be the volume of $M$ with the hyperbolic metric of curvature $-1$. For every $\Ga$-convex body $L \subset C$,
\begin{align}\label{fuchsianinequality}
\cP(L) \geq \frac{\vol_{\mathrm{hyp}}(M)^{2}}{N^{2}},
\end{align}
with equality if and only if $\pr L$ is a homothetic image of the hyperboloid $\{x \in C: \lb x, x\ra = -1\}$.
\end{corollary}
\begin{proof}
The computation of Example \ref{hyperboloidexample} uses only that $\Ga$ preserves $\lb\,,\,\ra$ and $C$, so the hyperboloid is a $\Ga$-invariant hyperbolic affine sphere asymptotic to $\pr C$, and its centroaffine metric is the hyperbolic metric of curvature $-1$. Theorem A with this $\Sigma_{0}$ gives the claim.
\end{proof}

\begin{corollary}\label{existencecorollary}
Let $C$ and $\Ga$ be as in Theorem A. There is a $\Ga$-invariant hyperbolic affine sphere in $C$ asymptotic to $\pr C$, namely the level set $\Sigma_{\psi_{0}}$ of a $\psi_{0} \in \cAs(C, \Ga)$ with $\mu(\psi_{0})$ constant, and it is unique up to homothety.
\end{corollary}
\begin{proof}
By Theorem B there is $\psi_{0} \in \cAs(C, \Ga)$ with $\mu(\psi_{0})$ constant, and its level set $\Sigma_{\psi_{0}}$ is a hyperbolic affine sphere centered at the origin, which is $\Ga$-invariant because $\psi_{0}$ is. Any two such $\psi_{0}$ are proportional by Lemma \ref{uniquenessproposition}, which is the uniqueness up to homothety. By Lemma \ref{gammaconvexbodylemma}, $L_{\psi_{0}}$ is a closed convex subset of $\ste$ contained in $C$, and its boundary is $\Sigma_{\psi_{0}} = \{x \in C: \psi_{0}(x) = 1\}$, so $\Sigma_{\psi_{0}}$ is asymptotic to $\pr C$ in the sense of Definition \ref{asymptoticdefinition}.
\end{proof}

\begin{definition}\label{maxmindefinition}
Let $\psi_{0}$ be as in Corollary \ref{existencecorollary}. The \emph{sharp Harnack constant} of $M$ is $\mathfrak{c}(M) = \mathfrak{c}_{\psi_{0}}(M)$ (Definition \ref{sharpharnackdefinition}), and the \emph{minimal} and \emph{maximal covolume products} of $M$ are
\begin{align}\label{maxmindefined}
&\mathfrak{m}(M) = \min_{L \in \cK(C, \Ga)}\cP(L) = \cP(\psi_{0}),& &\mathfrak{M}(M) = \sup_{L \in \cK(C, \Ga)}\cP(L).&
\end{align}
\end{definition}
The invariant $\mathfrak{m}(M)$ is studied further in Section \ref{invariantsection}.

\begin{corollary}\label{maximumcorollary}
Let $M$, $\mathfrak{c}(M)$, $\mathfrak{m}(M)$, and $\mathfrak{M}(M)$ be as in Definition \ref{maxmindefinition}.
\begin{enumerate}[label=$(\arabic*)$, ref=(\arabic*)]
\item\label{mc1} $\mathfrak{c}(M)$ depends only on the flat projective structure of $M$, and $\mathfrak{c}(M^{\ast}) = \mathfrak{c}(M)$.
\item\label{mc2} The supremum $\mathfrak{M}(M)$ is attained, $\mathfrak{M}(M^{\ast}) = \mathfrak{M}(M)$, and
\begin{align}\label{maxminbound}
\mathfrak{m}(M) < \mathfrak{M}(M) \leq \mathfrak{c}(M)^{-N}\mathfrak{m}(M).
\end{align}
\end{enumerate}
\end{corollary}
\begin{proof}
Since $\psi_{0}$ is unique up to a positive factor, and the ratio in \eqref{sharpharnackdefined} is unchanged when $\phi$ is multiplied by a positive constant, $\mathfrak{c}(M)$ depends only on $C$ and $\Ga$, and it is unchanged by conjugation in $GL(\ste)$. By Lemma \ref{dualityproposition}, $\mu(\psi_{0}^{\ast})$ is constant, so $\mathfrak{c}(M^{\ast}) = \mathfrak{c}_{\psi_{0}^{\ast}}(M^{\ast})$. By Lemma \ref{polartransformlemma} and \eqref{polardefined}, $\psi \mapsto \psi^{\ast}$ is an order reversing bijection from $\cA(C, \Ga)$ to $\cA(C^{\ast}, \Ga^{\ast})$ satisfying $(c\psi)^{\ast} = c^{-1}\psi^{\ast}$, so $\al\psi_{0} \leq \psi \leq \be\psi_{0}$ if and only if $\be^{-1}\psi_{0}^{\ast} \leq \psi^{\ast} \leq \al^{-1}\psi_{0}^{\ast}$, and the ratios in \eqref{sharpharnackdefined} for $\psi$ and $\psi^{\ast}$ coincide. This proves \ref{mc1}.

Let $\psi \in \cA(C, \Ga)$ be normalized by $\sup_{C}\psi/\psi_{0} = 1$, so that $\mathfrak{c}(M)\psi_{0} \leq \psi \leq \psi_{0}$. Then $\psi^{-N} \leq \mathfrak{c}(M)^{-N}\psi_{0}^{-N}$ and $\psi^{\ast} \geq \psi_{0}^{\ast}$, so, by \eqref{pdefined}, $\cP(\psi) \leq \mathfrak{c}(M)^{-N}\cP(\psi_{0})$. The proof of Lemma \ref{minimizerexistsproposition} uses only the normalization, Lemma \ref{harnacklemma}, and \eqref{pcontinuity}, so it applies to a maximizing sequence, and the supremum is attained. Since polarity is a bijection from $\cK(C, \Ga)$ to $\cK(C^{\ast}, \Ga^{\ast})$ with $\cP(L^{\circ}) = \cP(L)$, $\mathfrak{M}(M^{\ast}) = \mathfrak{M}(M)$. 

It remains to show that $\mathfrak{m}(M) < \mathfrak{M}(M)$. Let $x \in C$ and let $\psi_{x}$ be the cogauge of the body $L_{x}$ of Theorem \ref{sharpharnacktheorem}. It is claimed that $\psi_{x} \notin \cAs(C, \Ga)$. Suppose, to the contrary, that $\psi_{x} \in \cAs(C, \Ga)$, and let $y \in \Sigma_{\psi_{x}} = \pr L_{x}$. By Lemma \ref{polarmaplemma}\ref{pml3} and its proof, $\zeta = G_{\psi_{x}}(y)$ lies in $C^{\ast}$, and the hyperplane $\{z \in \ste: \zeta(z) = 1\}$ supports $L_{x}$ exactly at $y$, so that $\zeta(z) > 1$ for $z \in L_{x}\setminus\{y\}$. Since $\zeta$ is positive on $C$, the infimum of $\zeta$ over $L_{x}$, which equals $1$, is equal to $\inf_{\gamma \in \Ga}\zeta(\gamma x)$.

Fix a norm $|\cdot|$ on $\ste$. Since $\zeta \in C^{\ast}$, there is $c_{1} > 0$ such that $\zeta(z) \geq c_{1}|z|$ for $z \in \bar{C}$, and, since $\rho_{C}$ extends continuously by $0$ to $\bar{C}$, as shown in the proof of Lemma \ref{gammaconvexbodylemma}, for every $R > 0$ the set $\{z \in C: \rho_{C}(z) = \rho_{C}(x),\ \zeta(z) \leq R\}$ is a compact subset of $C$. It contains every point $\gamma x$, $\gamma \in \Ga$, such that $\zeta(\gamma x) \leq R$, because $\rho_{C}$ is $\Ga$-invariant. Because $\Ga$ acts properly on $\Om$ and preserves $\rho_{C}$, it acts properly on $C$, so only finitely many points of the orbit $\Ga x$ lie in a given compact subset of $C$. Hence $\inf_{\gamma \in \Ga}\zeta(\gamma x)$ is attained, at a point $\gamma_{0}x \in L_{x}$ such that $\zeta(\gamma_{0}x) = 1$, and so $y = \gamma_{0}x$. This shows that $\Sigma_{\psi_{x}} \subset \Ga x$, which is impossible, because $\Ga$ is countable, being a discrete subgroup of $GL(\ste)$, while $\Sigma_{\psi_{x}}$, which meets every ray of $C$, is uncountable. Hence $\psi_{x} \notin \cAs(C, \Ga)$, so, by Theorem B, $L_{x}$ is not a minimizer of $\cP$, and $\mathfrak{M}(M) \geq \cP(L_{x}) > \mathfrak{m}(M)$. This completes the proof of \ref{mc2}.
\end{proof}

Combined with Remark \ref{klartagmahlerremark}, \eqref{maxminbound} bounds $\mathfrak{m}(M)$ in terms of the characteristic functions of $C$ and $C^{\ast}$. Let $|B^{n}|$ be the volume of the Euclidean unit ball and let $c_{0} > 0$ be an absolute constant such that $c_{0}^{n}N^{N} \leq (n!)^{2}\bar{\mahler}(K)$ for every $n \geq 1$ and every compact convex set $K$ with nonempty interior in an $n$-dimensional affine space, where $N = n + 1$. Such a constant exists by \cite[Corollary $6.1$]{Bourgain-Milman}, which gives $\bar{\mahler}(K) \geq c_{2}^{n}|B^{n}|^{2}$ for an absolute constant $c_{2} > 0$ and all such $K$, the inequality proved there being valid for every convex body containing the origin in its interior, because, by Stirling's formula, there is an absolute constant $c_{3} > 0$ such that $(n!)^{2}|B^{n}|^{2} \geq c_{3}^{n}N^{N}$ for all $n \geq 1$. See also \cite[Theorem $2$]{Mastrantonis-Rubinstein}, where it is shown that $n!\,\bar{\mahler}(K) \geq (\pi/4)^{n}$.

The constant $c_{0}$ is not the same as the constant $c$ of the introduction, because the comparison here is with the Mahler product $N^{N}/(n!)^{2}$ of the simplex about its centroid and not with that of the ball. With the Blaschke--Santal\'o inequality, $c_{0}^{n}N^{N} \leq (n!)^{2}\bar{\mahler}(K_{\zeta}) \leq (n!)^{2}|B^{n}|^{2}$ for $\zeta \in C^{\ast}$. Since $\mathfrak{m}(M) \leq \cP(\rho_{C}) \leq \mathfrak{c}(M)^{-N}\mathfrak{m}(M)$, \eqref{klartagvolumeproduct} gives
\begin{align}\label{characteristicbounds}
\frac{\mathfrak{c}(M)^{N}}{N^{2}(n!)^{2}|B^{n}|^{2}} \leq \frac{\mathfrak{m}(M)}{\int_{M}\phi_{C}\int_{M^{\ast}}\phi_{C^{\ast}}} \leq \frac{1}{c_{0}^{n}N^{N+2}}.
\end{align}

\begin{example}\label{simplicialharnackexample}
Let $C$ and $\Ga$ be as in Example \ref{orthantexample}, and let $\Lambda \subset \Pi_{0} = \{y \in \rea^{N}: \sum_{i = 1}^{N}y_{i} = 0\}$ be the lattice of the vectors $(\log a_{1}, \dots, \log a_{N})$ for $g \in \Ga$ with $\ell_{i}\circ g = a_{i}\ell_{i}$. The linear maps multiplying each $\ell_{i}$ by a positive factor commute with $\Ga$, act transitively on $C$, and multiply $\psi_{0}$ by constants, so in Theorem \ref{sharpharnacktheorem} the point $x$ can be taken to satisfy $\ell_{i}(x) = 1$ for all $i$. By the argument proving \ref{mc1}, the ratio in \eqref{sharpharnack} for $\psi_{x}$ equals the corresponding ratio for $\psi_{x}^{\ast}$ relative to $\psi_{0}^{\ast}$.

For $\zeta = \sum_{i = 1}^{N}\zeta_{i}\ell_{i} \in C^{\ast}$, $\psi_{x}^{\ast}(\zeta)$ is the minimum of $\zeta$ over the orbit $\Ga x$, that is, $\min_{\lambda \in \Lambda}\sum_{i = 1}^{N}\zeta_{i}e^{\lambda_{i}}$, while, for $\psi_{0} = (\ell_{1}\cdots\ell_{N})^{1/N}$, $\psi_{0}^{\ast}(\zeta) = N(\zeta_{1}\cdots\zeta_{N})^{1/N}$ by the inequality between arithmetic and geometric means. On the slice where $\eta_{i} = \log\zeta_{i}$ satisfies $\eta \in \Pi_{0}$, $\psi_{0}^{\ast}(\zeta) = N$ and $\psi_{x}^{\ast}(\zeta)/\psi_{0}^{\ast}(\zeta) = \min_{\lambda \in \Lambda}\tfrac{1}{N}\sum_{i=1}^{N}e^{\eta_{i} + \lambda_{i}}$. By the inequality between arithmetic and geometric means, this quantity is at least $1$ for every $\eta \in \Pi_{0}$, because $\sum_{i = 1}^{N}(\eta_{i} + \lambda_{i}) = 0$, and it equals $1$ when $\eta = 0$, so its infimum over $\eta \in \Pi_{0}$ is $1$. It is a continuous function of $\eta \in \Pi_{0}$ invariant under translation by elements of the lattice $\Lambda$, so it attains its maximum. Since $\mathfrak{c}(M)$ is the quotient of the infimum by the supremum, this gives
\begin{align}\label{simplicialharnack}
\mathfrak{c}(M)^{-1} = \max_{\eta \in \Pi_{0}}\min_{\lambda \in \Lambda}\frac{1}{N}\sum_{i=1}^{N}e^{\eta_{i} + \lambda_{i}}.
\end{align}
For $N = 2$, with $\Ga$ generated by $g$ satisfying $\ell_{1}\circ g = e^{s}\ell_{1}$ and $\ell_{2}\circ g = e^{-s}\ell_{2}$, $s > 0$, \eqref{simplicialharnack} gives $\mathfrak{c}(M) = 1/\cosh(s/2)$. Areas in $\ste$ are computed in the coordinates $(\ell_{1}, \ell_{2})$ and areas in $\std$ in the dual coordinates; this does not affect the covolume product, which is unchanged when $\Psi$ is rescaled. In these coordinates, with $x$ the point such that $\ell_{1}(x) = \ell_{2}(x) = 1$, a period of $C \setminus L_{x}$ is the triangle with vertices $0$, $x$, $gx$, of area $\tfrac{1}{2}(e^{s} - e^{-s})$, and a period of $C^{\ast}\setminus L_{x}^{\circ}$ is the triangle with vertices $0$, $\zeta_{0}$, $\zeta_{0}\circ g^{-1}$, where $\zeta_{0} = (e^{s}+1)^{-1}(\ell_{1} + e^{s}\ell_{2})$ equals $1$ at $x$ and $gx$, of area $\tfrac{1}{2}(e^{s} - 1)/(e^{s} + 1)$, so $\cP(L_{x}) = \sinh^{2}(s/2)$. The affine sphere is the hyperbola $\ell_{1}\ell_{2} = 1$, which bounds a region of area $s$ per period, while its polar curve $\zeta_{1}\zeta_{2} = 1/4$ bounds one of area $s/4$, so $\mathfrak{m}(M) = s^{2}/4$. With \eqref{maxminbound},
\begin{align}\label{nonebounds}
\frac{\sinh^{2}(s/2)}{(s/2)^{2}} \leq \frac{\mathfrak{M}(M)}{\mathfrak{m}(M)} \leq \cosh^{2}(s/2),
\end{align}
so $\mathfrak{M}(M)/\mathfrak{m}(M)$ is not bounded by any function of the dimension alone.
\end{example}

\begin{corollary}\label{productinvariantcorollary}
In the setting of Theorem \ref{producttheorem},
\begin{align}\label{productinvariant}
\mathfrak{m}(M) = \tau^{2}\,\frac{N_{1}^{N_{1}+2}N_{2}^{N_{2}+2}}{N^{N}}\,\mathfrak{m}(M_{1})\,\mathfrak{m}(M_{2}),
\end{align}
with $\mathfrak{m}(M_{i}) = 1$ when $N_{i} = 1$, and the minimizers of $\cP$ on $\cK(C, \Ga)$ are the products $L_{1}\boxtimes L_{2}$ of minimizers of $\cP$ on the $\cK(C_{i}, \Ga_{i})$.
\end{corollary}
\begin{proof}
By Corollary \ref{existencecorollary} there are $\psi_{0, i} \in \cAs(C_{i}, \Ga_{i})$ with $\mu(\psi_{0,i})$ constant; when $N_{i} = 1$, any positive linear function on $C_{i}$ will do. By Lemma \ref{productcogaugelemma}\ref{pcl3} and the $\Ga$-invariance shown before Lemma \ref{productcovolumelemma}, $\psi_{0} = \psi_{0,1}\boxtimes\psi_{0,2} \in \cAs(C, \Ga)$ and $\mu(\psi_{0})$ is constant. By Theorem B and Lemma \ref{uniquenessproposition}, the minimizers of $\cP$ on $\cK(C, \Ga)$ are exactly the homothetic images of $L_{\psi_{0}}$, and those on $\cK(C_{i}, \Ga_{i})$ are exactly the homothetic images of $L_{\psi_{0,i}}$. Since $(t_{1}^{-1}\psi_{1})\boxtimes(t_{2}^{-1}\psi_{2}) = t_{1}^{-a}t_{2}^{-b}\,\psi_{1}\boxtimes\psi_{2}$, the products of minimizers are exactly the minimizers, and \eqref{productinvariant} follows from \eqref{productidentity}.
\end{proof}
The same proofs give \eqref{productidentity} for $h_{i}$-invariant $L_{i}$, and \eqref{productinvariant} in general, when $g_{\tau}$ is replaced by $(h_{1}e^{N_{2}\tau}, h_{2}e^{-N_{1}\tau})$ with $h_{i} \in SL^{\pm}(\ste_{i})$ preserving $C_{i}$ and normalizing $\Ga_{i}$, because such $h_{i}$ preserves the $\Ga_{i}$-invariant affine sphere by Corollary \ref{existencecorollary}.

By Lemma \ref{densitylemma}\ref{dl1} and the continuity of $\cP$, $\mathfrak{M}(M)$ is the supremum of $\cP$ over $\Ga$-convex bodies cut out by finitely many $\Ga$-orbits of half-spaces, and, applying this to $M^{\ast}$ and taking polars, also over the closures in $C$ of the sets $\operatorname{conv}(\Ga X) + C$ with $X \subset C$ finite.

\begin{question}\label{maximizerquestion}
Determine $\mathfrak{M}(M)$ and the maximizers of $\cP$ on $\cK(C, \Ga)$. Is there always a maximizer of the form $L_{x}$, a polytope with one orbit of vertices, or, equivalently by polarity, one with one orbit of facets, $\{x \in C: \zeta(\gamma x) \geq 1 \text{ for all } \gamma \in \Ga\}$ for some $\zeta \in C^{\ast}$?
\end{question}

Products give evidence for an affirmative answer and suggest a route to it. In the setting of Theorem \ref{producttheorem}, let $\mathfrak{M}_{N_{1}, N_{2}}(\tau) = \sup_{G \in \mathcal{G}_{\tau}}\cP_{N_{1}, N_{2}}(G)$. Applying \eqref{compositeidentity} to maximizers $L_{i}$ gives
\begin{align}\label{maxproduct}
\mathfrak{M}(M) \geq N_{1}^{2}N_{2}^{2}\,\mathfrak{M}_{N_{1}, N_{2}}(\tau)\,\mathfrak{M}(M_{1})\,\mathfrak{M}(M_{2}),
\end{align}
while, by Corollary \ref{productinvariantcorollary}, the minimum factors in the same way with $\cP_{N_{1}, N_{2}}(G_{0})$ in place of the planar supremum. Moreover $\mathfrak{M}_{N_{1}, N_{2}}(\tau) > \cP_{N_{1}, N_{2}}(G_{0})$. Indeed, when the $L_{i}$ are minimizers, $G(L_{1}, L_{2})$ is homothetic to the affine sphere only if $G$ is a constant multiple of $G_{0}$, because the pairs $(\psi_{1}(x_{1}), \psi_{2}(x_{2}))$ fill $Q$, so, by Theorem A and \eqref{compositeidentity}, $\cP_{N_{1}, N_{2}}(G) > \cP_{N_{1}, N_{2}}(G_{0})$ for every other $G$. Consequently $\mathfrak{M}/\mathfrak{m}$ is supermultiplicative, with an additional factor $\mathfrak{M}_{N_{1}, N_{2}}(\tau)/\cP_{N_{1}, N_{2}}(G_{0}) > 1$ coming from the circle factor of Lemma \ref{productdivisiblelemma}, and the product $L_{1}\boxtimes L_{2}$ of maximizers is never a maximizer.

The class of polytopes in Question \ref{maximizerquestion} is stable under this construction. If $L_{i} = \{\zeta_{i}(\gamma_{i}x_{i}) \geq 1 \text{ for all } \gamma_{i} \in \Ga_{i}\}$ and $G(u) = \min_{k \in \mathbb{Z}}(\al e^{kN_{2}\tau}u_{1} + \be e^{-kN_{1}\tau}u_{2})$ is a polygon with one orbit of edges, with $\al, \be > 0$, then, because $\al$ and $\be$ are positive and $\psi_{i} = \min_{\gamma_{i} \in \Ga_{i}}\zeta_{i}\circ\gamma_{i}$, $G(\psi_{1}, \psi_{2}) = \min_{\gamma \in \Ga}(\al\zeta_{1} + \be\zeta_{2})\circ\gamma$, so $G(L_{1}, L_{2})$ is a polytope with one orbit of facets. This suggests that equality holds in \eqref{maxproduct}, with $\mathfrak{M}_{N_{1}, N_{2}}(\tau)$ attained by a polygon with one orbit of edges, and, for a simplicial cone whose lattice is a product of lattices of lower rank, that the maximizers are polytopes with one orbit of facets. For $N = 2$ the planar problem is the whole problem, and whether the polygon $L_{x}$ of Example \ref{simplicialharnackexample} is a maximizer is open already in that case.

\section{First and second variations of covolumes and the covolume product}\label{variationsection}

\begin{lemma}\label{coneopenlemma}
The set $\cAs(C, \Ga)$ is a convex cone, open in the $C^{2}$ topology of the space of smooth $\Ga$-invariant functions on $C$ homogeneous of degree one. In particular, for $\psi \in \cAs(C, \Ga)$ and $f \in \cinf(M)$, $(1 + tf)\psi \in \cAs(C, \Ga)$ for $|t|$ small.
\end{lemma}

\begin{proof}
Positivity, concavity, homogeneity of degree one, $\Ga$-invariance, and the negative semidefiniteness of $\fD^{2}\psi$ are preserved by sums and positive multiples, and the radical of $\fD^{2}(\psi_{1} + \psi_{2}) = \fD^{2}\psi_{1} + \fD^{2}\psi_{2}$, a sum of negative semidefinite forms, is the intersection of their radicals, which is spanned by $\eul$. Hence $\cAs(C, \Ga)$ is a convex cone. Its defining conditions, positivity and negative definiteness of $\fD^{2}\psi$ transversally to $\eul$, are open conditions on the $2$-jet of $\psi$ that, by $\Ga$-invariance and the compactness of $M$, need only be checked on the compact set $\Fset$ of Lemma \ref{fundamentalsetlemma}, by homogeneity, $\Ga$-invariance, and \eqref{fundamentalset}, so they are open in the $C^{2}$ topology. The last claim follows because $f\psi$ is a smooth $\Ga$-invariant function homogeneous of degree one.
\end{proof}

Let $\psi \in \cAs(C, \Ga)$ and $\mu = \mu(\psi)$, and let $\Phi$, $\fD^{2}\Phi$, and $G_{\psi}$ be as in \eqref{phiqdefined}. By Lemma \ref{polarmaplemma}\ref{pml1}, $\sign\det \fD^{2}\Phi = (-1)^{n}$, so $\mu = (-1)^{n}\det \fD^{2}\Phi$. For $f \in \cinf(M)$ define
\begin{align}\label{lopdefined}
\Lop_{\psi}f = (-1)^{n}\mu^{-1}\psi^{2}\lb\adj(\fD^{2}\Phi), \fD^{2}f\ra = \psi^{2}(\det \fD^{2}\Phi)^{-1}\lb\adj(\fD^{2}\Phi), \fD^{2}f\ra,
\end{align}
the trace of $\fD^{2}f$ with respect to $\fD^{2}\Phi$, multiplied by $\psi^{2}$. It is homogeneous of degree $0$ and $\Ga$-invariant, hence a function on $M$. By Lemma \ref{laplacianlemma}, $\Lop_{\psi} = -\Delta_{\gb_{\psi}}$ when $\mu(\psi)$ is constant. 

\begin{lemma}\label{ibplemma}
Let $\psi \in \cAs(C, \Ga)$ and $\mu = \mu(\psi)$. For $f, f_{1} \in \cinf(M)$,
\begin{align}\label{ibp}
&\int_{M}\mu\psi^{-N}\Lop_{\psi}f = 0,& &\int_{M}\mu\psi^{-N}f_{1}\Lop_{\psi}f = \int_{M}\mu\psi^{-N}\lb df_{1}, df\ra_{\gb_{\psi}}.&
\end{align}
\end{lemma}
\begin{proof}
Let $Y = \psi^{2-N}f_{1}\,\adj(\fD^{2}\Phi)(df, \cdot)$, the vector field obtained by contracting $df$ into one slot of $\adj(\fD^{2}\Phi)$. It is homogeneous of degree $1-N$, and it is $\Ga$-invariant, because $\fD^{2}\Phi$ is $\Ga$-invariant and the tensor $\Psi^{\ast}\otimes\Psi^{\ast}$, by means of which $\adj(\fD^{2}\Phi)$ is defined, is invariant under $SL^{\pm}(\ste)$. Since $\adj(\fD^{2}\Phi)$ is divergence free,
\begin{align}
\begin{aligned}
\operatorname{div}Y &= \psi^{2-N}f_{1}\lb\adj(\fD^{2}\Phi), \fD^{2}f\ra + \psi^{2-N}\adj(\fD^{2}\Phi)(df_{1}, df) \\
&\quad + (2-N)\psi^{1-N}f_{1}\adj(\fD^{2}\Phi)(d\psi, df).
\end{aligned}
\end{align}
Since $\fD^{2}\Phi(\eul, \cdot) = \fD\Phi = \psi\,d\psi$, there holds $\adj(\fD^{2}\Phi)(d\psi, \cdot) = (\det \fD^{2}\Phi)\psi^{-1}\eul$. Hence the last term of $\operatorname{div}Y$ vanishes, because $\adj(\fD^{2}\Phi)(d\psi, df) = (\det \fD^{2}\Phi)\psi^{-1}df(\eul)$ and $f$ is homogeneous of degree $0$. By \eqref{lopdefined}, $\lb\adj(\fD^{2}\Phi), \fD^{2}f\ra = (-1)^{n}\mu\psi^{-2}\Lop_{\psi}f$, and by \eqref{dualmetricadjugate} with polarization, $\adj(\fD^{2}\Phi)(df_{1}, df) = (-1)^{n+1}\mu\psi^{-2}\lb df_{1}, df\ra_{\gb_{\psi}}$. Hence
\begin{align}
\operatorname{div}Y = (-1)^{n}\mu\psi^{-N}\left(f_{1}\Lop_{\psi}f - \lb df_{1}, df\ra_{\gb_{\psi}}\right).
\end{align}
By Lemma \ref{divergencelemma}, $\int_{M}\operatorname{div}Y = 0$, which is the second identity of \eqref{ibp}, and taking $f_{1} = 1$ gives the first.
\end{proof}
Equivalently, $\Lop_{\psi}f = -\Delta_{\gb_{\psi}}f - \tfrac{1}{2}\lb d\log\mu(\psi), df\ra_{\gb_{\psi}}$, where $\Delta_{\gb_{\psi}}$ is the Laplacian of $\gb_{\psi}$, so that $\Lop_{\psi}$ is the Laplacian of $\gb_{\psi}$ weighted by the measure $\mu(\psi)\psi^{-N} = \mu(\psi)^{1/2}\vol_{\gb_{\psi}}$. It reduces to $-\Delta_{\gb_{\psi}}$ when $\mu(\psi)$ is constant, which is Lemma \ref{laplacianlemma}.

\begin{lemma}\label{laplacianlemma}
Let $\psi_{0} \in \cAs(C, \Ga)$ satisfy $\mu(\psi_{0}) = \ka$ constant, and let $\Delta$ be the Laplacian of $\gb_{\psi_{0}}$. Then $\Lop_{\psi_{0}}f = -\Delta f$ for every $f \in \cinf(M)$.
\end{lemma}
\begin{proof}
Since $\mu(\psi_{0})$ is constant, $\psi_{0}^{-N}$ is a constant multiple of $\vol_{\gb_{\psi_{0}}}$ by Lemma \ref{centroaffinemetriclemma}\ref{cml1}, so the second identity of \eqref{ibp} becomes $\int_{M}f_{1}\Lop_{\psi_{0}}f\,\psi_{0}^{-N} = \int_{M}\lb df_{1}, df\ra_{\gb_{\psi_{0}}}\psi_{0}^{-N}$ for all $f_{1} \in \cinf(M)$. Integrating by parts on the closed Riemannian manifold $(M, \gb_{\psi_{0}})$ gives $\int_{M}\lb df_{1}, df\ra_{\gb_{\psi_{0}}}\psi_{0}^{-N} = -\int_{M}f_{1}\Delta f\,\psi_{0}^{-N}$, so $\int_{M}f_{1}(\Lop_{\psi_{0}}f + \Delta f)\psi_{0}^{-N} = 0$ for all $f_{1} \in \cinf(M)$, and $\Lop_{\psi_{0}}f = -\Delta f$.
\end{proof}

\begin{lemma}
\label{firstvariationproposition}
Let $\psi \in \cAs(C, \Ga)$, $\mu = \mu(\psi)$, $f \in \cinf(M)$, and $\psi_{t} = (1 + tf)\psi$, which lies in $\cAs(C, \Ga)$ for $|t|$ small. Then
\begin{align}\label{firstvariation}
&\tfrac{d}{dt}_{|t=0}|K_{\psi_{t}}/\Ga| = -\int_{M}f\psi^{-N},& &\tfrac{d}{dt}_{|t=0}|K^{\ast}_{\psi_{t}}/\Ga^{\ast}| = \int_{M}f\mu\psi^{-N}.&
\end{align}
Consequently $\psi$ is critical for $\cP$ on $\cAs(C, \Ga)$, meaning that $\tfrac{d}{dt}_{|t=0}\cP((1 + tf)\psi) = 0$ for every $f \in \cinf(M)$, if and only if $\mu(\psi) = |K^{\ast}_{\psi}/\Ga^{\ast}|/|K_{\psi}/\Ga|$ is constant, that is, if and only if $\Sigma_{\psi}$ is a hyperbolic affine sphere.
\end{lemma}
\begin{proof}
Write $\Phi_{t} = \Phi_{\psi_{t}} = \tfrac{1}{2}\psi_{t}^{2}$ and $\varphi = f\psi$. The first formula follows from \eqref{conevolume}, since $\psi_{t}^{-N} = (1 + tf)^{-N}\psi^{-N}$, by differentiation under the integral sign, which is justified because the integrand depends smoothly on $t$ and $M$ is compact. The same justification applies to the other derivatives of integrals over $M$ computed in this section. For the second, $\tfrac{d}{dt}_{|t=0}\Phi_{t} = \psi\varphi = 2f\Phi$, and
\begin{align}
\fD^{2}(2f\Phi) = 2f\,\fD^{2}\Phi + 2(df\otimes \fD\Phi + \fD\Phi\otimes df) + 2\Phi \fD^{2}f.
\end{align}
By the Jacobi formula, since $\lb\adj(\fD^{2}\Phi), \fD^{2}\Phi\ra = N\det \fD^{2}\Phi$, $\adj(\fD^{2}\Phi)(\fD\Phi, \cdot) = (\det \fD^{2}\Phi)\eul$, and $df(\eul) = 0$,
\begin{align}\label{deltamu}
\tfrac{d}{dt}_{|t=0}\mu(\psi_{t}) = (-1)^{n}\lb\adj(\fD^{2}\Phi), \fD^{2}(2f\Phi)\ra = \mu\left(2Nf + \Lop_{\psi}f\right).
\end{align}
Hence, by Lemma \ref{ibplemma},
\begin{align}
\tfrac{d}{dt}_{|t=0}\int_{M}\mu(\psi_{t})\psi_{t}^{-N} = \int_{M}\mu\psi^{-N}(2Nf - Nf) + \int_{M}\mu\psi^{-N}\Lop_{\psi}f = N\int_{M}f\mu\psi^{-N},
\end{align}
and \eqref{volumeformulas} gives the second formula. Then
\begin{align}
\tfrac{d}{dt}_{|t=0}\log\cP(\psi_{t}) = \int_{M}f\psi^{-N}\left(\frac{\mu}{|K^{\ast}_{\psi}/\Ga^{\ast}|} - \frac{1}{|K_{\psi}/\Ga|}\right),
\end{align}
which vanishes for all $f$ if and only if $\mu = |K^{\ast}_{\psi}/\Ga^{\ast}|/|K_{\psi}/\Ga|$.
\end{proof}
By Lemma \ref{firstvariationproposition} and Corollary \ref{existencecorollary}, the positive multiples of the cogauge $\psi_{0}$ of the affine sphere are the only critical points of $\cP$ on $\cAs(C, \Ga)$. Theorem \ref{secondvariationtheorem} shows that the critical point is nondegenerate transverse to rescalings.
\begin{corollary}\label{nosmoothmaximizercorollary}
No maximizer of $\cP$ on $\cK(C, \Ga)$ has smooth, locally strongly convex boundary.
\end{corollary}
\begin{proof}
By Lemma \ref{smoothboundarylemma}, the cogauge of a $\Ga$-convex body with smooth, locally strongly convex boundary lies in $\cAs(C, \Ga)$. If such a body maximized $\cP$, its cogauge would be critical for $\cP$ on $\cAs(C, \Ga)$, hence a positive multiple of $\psi_{0}$ by Lemma \ref{firstvariationproposition} and Corollary \ref{existencecorollary}, hence a minimizer, contradicting $\mathfrak{m}(M) < \mathfrak{M}(M)$ in \eqref{maxminbound}.
\end{proof}

\begin{lemma}\label{covolumeconvexitylemma}
The function $\psi \mapsto |K^{\ast}_{\psi}/\Ga^{\ast}|$ is strictly convex on the convex cone $\cAs(C, \Ga)$. More precisely, for $\psi \in \cAs(C, \Ga)$, $\mu = \mu(\psi)$, $f \in \cinf(M)$, and $\varphi = f\psi$,
\begin{align}\label{secondcovolume}
\tfrac{d^{2}}{dt^{2}}_{|t=0}|K^{\ast}_{\psi + t\varphi}/\Ga^{\ast}| = \int_{M}\mu\psi^{-N}\left((N-1)f^{2} + |df|^{2}_{\gb_{\psi}}\right) \geq (N-1)\int_{M}f^{2}\mu\psi^{-N}.
\end{align}
\end{lemma}
\begin{proof}
Write $\psi_{t} = \psi + t\varphi$, which lies in $\cAs(C, \Ga)$ for $|t|$ sufficiently small, by Lemma \ref{coneopenlemma}. For such $t$ and for $|s|$ sufficiently small, $\psi_{t + s} = (1 + s\varphi/\psi_{t})\psi_{t}$, so \eqref{firstvariation}, applied with $\psi_{t}$ in place of $\psi$ and $\varphi/\psi_{t} \in \cinf(M)$ in place of $f$, gives $\tfrac{d}{dt}|K^{\ast}_{\psi_{t}}/\Ga^{\ast}| = \int_{M}\varphi\,\mu(\psi_{t})\psi_{t}^{-N-1}$, the direction $\varphi$ being fixed. Differentiating again with \eqref{deltamu},
\begin{align}
\tfrac{d^{2}}{dt^{2}}_{|t=0}|K^{\ast}_{\psi_{t}}/\Ga^{\ast}| = \int_{M}\mu\psi^{-N}\left(2Nf^{2} + f\Lop_{\psi}f - (N+1)f^{2}\right),
\end{align}
and Lemma \ref{ibplemma} with $f_{1} = f$ gives \eqref{secondcovolume}. 
Applying \eqref{secondcovolume} at each point of a segment $\psi + t\varphi$ contained in $\cAs(C, \Ga)$, in place of $\psi$, shows that along the segment the second derivative of $t \mapsto |K^{\ast}_{\psi + t\varphi}/\Ga^{\ast}|$ is positive unless $\varphi = 0$, because $N - 1 > 0$, which is the strict convexity.
\end{proof}
In the Fuchsian case of Example \ref{hyperboloidexample}, Lemma \ref{covolumeconvexitylemma} is the convexity of the covolume of Fuchsian convex bodies due to Fillastre \cite{Fillastre-Fuchsian}, applied to the polar sets $L^{\circ}_{\psi}$, whose support functions are the $\psi$. Without a group action, the corresponding convexity of the covolume, in the form of Brunn--Minkowski and Aleksandrov--Fenchel inequalities, is proved for coconvex bodies in \cite{Khovanskii-Timorin} and for coconvex sets of finite volume in \cite{Schneider-coconvex}.

The segments $(1-t)\psi_{0} + t\psi_{1}$, $0 \leq t \leq 1$, in $\cAs(C, \Ga)$ along which Lemma \ref{covolumeconvexitylemma} asserts convexity are Minkowski combinations of polar bodies. By Lemma \ref{polartransformlemma} and the paragraph following the definition \eqref{polardefined} of the polar, $\psi$ is the support function $x \mapsto \inf_{\zeta \in L^{\circ}_{\psi}}\zeta(x)$ of the polar body $L^{\circ}_{\psi}$ on $C$, defined with an infimum rather than the usual supremum, and a $\Ga^{\ast}$-convex body is determined by its support function. The Minkowski combination $(1-t)L^{\circ}_{\psi_{0}} + tL^{\circ}_{\psi_{1}}$ is a $\Ga^{\ast}$-convex body, being convex, $\Ga^{\ast}$-invariant, and closed, the latter because both summands have recession cone $\bar{C}^{\ast}$, which contains no line \cite[Corollary $9.1.2$]{Rockafellar}, and its support function is $(1-t)\psi_{0} + t\psi_{1}$. Hence $L^{\circ}_{(1-t)\psi_{0} + t\psi_{1}} = (1-t)L^{\circ}_{\psi_{0}} + tL^{\circ}_{\psi_{1}}$. In the Fuchsian case this is the Minkowski addition of Fuchsian convex bodies of \cite{Fillastre-Fuchsian}.

\begin{theorem}\label{secondvariationtheorem}
Let $\psi_{0} \in \cAs(C, \Ga)$ satisfy $\mu(\psi_{0}) = \ka$ for a constant $\ka > 0$, let $\nu = \psi_{0}^{-N} = \ka^{-1/2}\vol_{\gb_{\psi_{0}}}$, let $V = \int_{M}\nu$, and, for $f \in \cinf(M)$, let $\bar{f} = V^{-1}\int_{M}f\nu$. Then
\begin{align}\label{secondvariation}
\tfrac{d^{2}}{dt^{2}}_{|t=0}\cP((1 + tf)\psi_{0}) = \frac{\ka V}{N}\int_{M}\left(|df|^{2}_{\gb_{\psi_{0}}} + 2N(f - \bar{f})^{2}\right)\nu.
\end{align}
In particular the second variation is a positive semidefinite quadratic form whose radical consists of the constant functions, which correspond to rescalings.
\end{theorem}
\begin{proof}
By Lemma \ref{centroaffinemetriclemma}\ref{cml1}, $\ka^{1/2}\psi_{0}^{-N}$ is the Riemannian density of $\gb_{\psi_{0}}$, which is the equality $\nu = \ka^{-1/2}\vol_{\gb_{\psi_{0}}}$ of the statement. By \eqref{conevolume} with $f = 1$, $|K_{\psi_{0}}/\Ga| = V/N$, and by \eqref{volumeformulas}, $|K^{\ast}_{\psi_{0}}/\Ga^{\ast}| = \ka V/N$. Along $\psi_{t} = (1+tf)\psi_{0}$, $|K_{\psi_{t}}/\Ga| = \tfrac{1}{N}\int_{M}(1+tf)^{-N}\nu$ has first and second derivatives $-\int_{M} f\nu$ and $(N+1)\int_{M} f^{2}\nu$ at $t = 0$, and, by \eqref{firstvariation} and \eqref{secondcovolume}, $|K^{\ast}_{\psi_{t}}/\Ga^{\ast}|$ has first and second derivatives $\ka\int_{M} f\nu$ and $\ka\int_{M}((N-1)f^{2} + |df|^{2}_{\gb_{\psi_{0}}})\nu$. Hence the second derivative of the product is
\begin{align}
\begin{aligned}
&\frac{\ka V}{N}(N+1)\int_{M} f^{2}\nu - 2\ka\Big(\int_{M} f\nu\Big)^{2} + \frac{\ka V}{N}\int_{M}\left((N-1)f^{2} + |df|^{2}\right)\nu\\
&\quad = \frac{\ka V}{N}\left(\int_{M}|df|^{2}\nu + 2N\int_{M} f^{2}\nu - \frac{2N}{V}\Big(\int_{M} f\nu\Big)^{2}\right),
\end{aligned}
\end{align}
which is \eqref{secondvariation}, because $\int_{M}(f - \bar{f})^{2}\nu = \int_{M}f^{2}\nu - V^{-1}(\int_{M}f\nu)^{2}$. The right-hand side of \eqref{secondvariation} is nonnegative, and it vanishes if and only if $df = 0$ and $f = \bar{f}$, that is, since $M$ is connected, if and only if $f$ is constant. For a constant $c$, the function $(1 + tc)\psi_{0}$ is a rescaling of $\psi_{0}$.
\end{proof}

\begin{proof}[Proof of Theorem C]
Let $\Sigma_{0}$, $L_{0}$, $\psi_{0}$, and $h$ be as in the statement of Theorem C. By the first part of the proof of Theorem A, $L_{0}$ is a $\Ga$-convex body, its cogauge $\psi_{0}$ lies in $\cAs(C, \Ga)$ and satisfies $\Sigma_{\psi_{0}} = \Sigma_{0}$ and $\gb_{\psi_{0}} = h$, and $\mu(\psi_{0})$ is equal to a positive constant $\ka$.

Let $f \in \cinf(M)$. By Lemma \ref{coneopenlemma}, $\psi_{t} = (1 + tf)\psi_{0}$ lies in $\cAs(C, \Ga)$ for $|t|$ sufficiently small, so that $\{x \in C: \psi_{t}(x) \geq 1\} = L_{\psi_{t}}$ is a $\Ga$-convex body, and $\cP(L_{\psi_{t}}) = \cP(\psi_{t})$. Hence Theorem \ref{secondvariationtheorem} applies. In the notation of that theorem, $\nu = \ka^{-1/2}\vol_{h}$ and $\int_{M}\nu = \ka^{-1/2}\vol(M, h)$, so that $\ka(\int_{M}\nu)\,\nu = \vol(M, h)\vol_{h}$, and the mean of $f$ with respect to $\nu$ is its $h$-mean $\bar{f}$. Therefore \eqref{secondvariation} becomes \eqref{theoremCformula}, in which $V$ denotes $\vol(M, h)$ and not the quantity $\int_{M}\nu$ that is denoted by $V$ in Theorem \ref{secondvariationtheorem}. In particular \eqref{theoremCformula} holds for every choice of $\Sigma_{0}$ among its homothetic images. The last assertion of Theorem C follows from the last assertion of Theorem \ref{secondvariationtheorem}, because, for a constant $c$ and $|t|$ sufficiently small, $(1 + tc)\psi_{0}$ is the cogauge of the image of $L_{0}$ under the homothety $x \mapsto (1 + tc)^{-1}x$, which preserves $C$.
\end{proof}
\begin{remark}\label{ellipticremark}
For a convex body $K \subset \ste$ with smooth boundary of positive curvature and the origin in its interior, the gauge function $\gamma$ is positive, convex, and homogeneous of degree $1$ on $\ste\setminus\{0\}$, and $|K| = \tfrac{1}{N}\int_{S^{n}}\gamma^{-N}$ and $|K^{\circ}| = \tfrac{1}{N}\int_{S^{n}}\mu(\gamma)\gamma^{-N}$, with $\mu(\gamma) = \det \fD^{2}(\tfrac{1}{2}\gamma^{2}) > 0$, where integrals over $S^{n} = \projp(\ste)$ are defined as in Section \ref{settingsection}.

The computations of this section apply with three changes of sign. Since $\fD^{2}(\tfrac{1}{2}\gamma^{2})$ is positive definite, the centroaffine metric $\gb_{\gamma}$ is the restriction of $+\fD^{2}\gamma$ to $\pr K$. Accordingly Lemma \ref{centroaffinemetriclemma}\ref{cml2} reads $|df|^{2}_{\gb_{\gamma}} = +\gamma^{2}\mu(\gamma)^{-1}\adj(\fD^{2}\tfrac{1}{2}\gamma^{2})(df, df)$, the integration by parts of Lemma \ref{ibplemma} reads $\int_{S^{n}}\mu\gamma^{-N}f_{1}\Lop_{\gamma}f = -\int_{S^{n}}\mu\gamma^{-N}\langle df_{1}, df\rangle_{\gb_{\gamma}}$, and consequently the gradient term in \eqref{secondcovolume}, and hence in Theorem \ref{secondvariationtheorem}, changes sign, while the terms without derivatives, which come from \eqref{deltamu} and from the expansion of $(1 + tf)^{-N}$, are unchanged.

The result is that the second variation of the Mahler product $|K||K^{\circ}|$ at a centered ellipsoid, in the direction $(1 + tf)\gamma$, is $\tfrac{\ka V}{N}\int_{S^{n}}(-|df|^{2} + 2N(f - \bar{f})^{2})\nu$, where, for the unit ball, $\ka = 1$ and $\nu$ is the round measure. On spherical harmonics $f_{k}$ of degree $k \geq 1$ this equals $\tfrac{\ka V}{N}(2N - k(k+n-1))\int_{S^{n}}f_{k}^{2}\nu$, which is positive for $k = 1$, zero for $k = 2$, and negative for $k \geq 3$. The directions $k = 1$ are the translations of $K$ relative to the origin, since translating by $a$ changes $\gamma$ to first order by $-\gamma\,d\gamma(a)$. Those with $k = 2$ are the infinitesimal actions of $\mathfrak{sl}(\ste)$, which preserve the volume product.

The centroaffine volume has the same sign pattern at the ellipsoid: by the computation of Wang, who concludes that the ellipsoid is an unstable centroaffine minimal hypersurface \cite[Corollary $3.2$]{Wang-centroaffine}, its second variation on $f_{k}$ is a positive multiple of $-k(k+n-1)\left(k(k+n-1) - 2N\right)\int_{S^{n}}f_{k}^{2}\nu$.

This agrees with the Blaschke--Santal\'o inequality and with the minimizing property of the Santal\'o point \cite{Santalo-affineinvariant}, the point $x_{0}$ at which $|K - x_{0}||(K - x_{0})^{\circ}|$ is minimal, which is the origin exactly when the barycenter of $K^{\circ}$ is at the origin \cite[Section $1$]{AKM-santalo}. The same barycenter condition is the one in Remark \ref{klartagremark}. The contrast with Theorem \ref{secondvariationtheorem}, where all directions transverse to rescalings are positive, is the infinitesimal form of Theorem A.
\end{remark}

\section{A linear inequality for covolumes and a second proof of the Blaschke--Santal\'o inequality for the covolume product}\label{reversesantalosection}

This section gives a second proof of Theorem A, including the characterization of the equality case, together with a stability estimate (Theorem \ref{stabilitytheorem}). The inequality is proved first for $\Ga$-convex bodies with smooth, locally strongly convex boundary, that is, by Lemma \ref{smoothboundarylemma}, for cogauges $\psi \in \cAs(C, \Ga)$ (Theorem \ref{smoothreversesantalotheorem}), and it is then extended to all $\Ga$-convex bodies by approximation (Lemma \ref{densitylemma}). The proof does not use the direct method of Section \ref{existencesection}. It requires only the existence of the affine sphere, given by Corollary \ref{existencecorollary} or by the theorems of Cheng and Yau, and it yields another proof of the uniqueness of the affine sphere up to homothety.

The key inequality, Theorem \ref{suminequalityproposition}, is linear rather than multiplicative. By the inequality between arithmetic and geometric means, it is a consequence of Theorem A, while the family of such inequalities obtained by rescaling $\psi_{0}$ is equivalent to Theorem A for such bodies (Corollary \ref{weightedcorollary} and Theorem \ref{smoothreversesantalotheorem}). Its proof uses only the convexity of the polar covolume and the convexity of $x \mapsto x^{-N}$. 

Theorem \ref{relativetheorem} in Section \ref{relativesection} is the analogue of Theorem \ref{suminequalityproposition} for bounded convex domains.

\begin{lemma}\label{balancedlemma}
Let $\psi_{0} \in \cAs(C, \Ga)$ satisfy $\mu(\psi_{0}) = \ka$ for a constant $\ka > 0$. Then the polar $\psi_{0}^{\ast}$ satisfies
\begin{align}\label{dualconstant}
\mu(\psi_{0}^{\ast}) = \ka^{\ast} := \ka^{-1},
\end{align}
the level sets $\Sigma_{\psi_{0}}$ and $\Sigma_{\psi_{0}^{\ast}}$ are hyperbolic affine spheres with affine mean curvatures $-\ka^{1/(n+2)}$ and $-\ka^{-1/(n+2)}$, respectively, and
\begin{align}\label{balanced}
N\ka^{1/2}|K_{\psi_{0}}/\Ga| = N(\ka^{\ast})^{1/2}|K^{\ast}_{\psi_{0}}/\Ga^{\ast}| = \vol(M, \gb_{\psi_{0}}) = \vol(M^{\ast}, \gb_{\psi_{0}^{\ast}}).
\end{align}
The common value in \eqref{balanced} is the same for all positive multiples of $\psi_{0}$.
\end{lemma}
\begin{proof}
The identity \eqref{dualconstant} is Lemma \ref{dualityproposition}, and the affine mean curvatures are given by Lemma \ref{affinespherelemma}. By \eqref{volumeformulas} and Lemma \ref{centroaffinemetriclemma}\ref{cml1}, $|K^{\ast}_{\psi_{0}}/\Ga^{\ast}| = \ka|K_{\psi_{0}}/\Ga|$ and $N\ka^{1/2}|K_{\psi_{0}}/\Ga| = \vol(M, \gb_{\psi_{0}})$, and the last equality in \eqref{balanced} is \eqref{sbetaduality} with $\be = 1/2$. The last claim holds because $\vol(M, \gb_{\psi}) = \cS_{1/2}(\psi)$ is invariant under rescaling.
\end{proof}

\begin{definition}\label{blaschkevolumedefinition}
The \emph{equiaffine volume} of $M = \Om/\Ga$ is $\cV(M) = \vol(M, \gb_{\psi_{0}})$, where $\psi_{0} \in \cAs(C, \Ga)$ is any element with $\mu(\psi_{0})$ constant.
\end{definition}

By Corollary \ref{existencecorollary} and Lemma \ref{balancedlemma}, $\cV(M)$ is well defined, it equals $\vol(M, h)$ for the metric $h$ of Theorem A, and it is preserved by polar duality, $\cV(M^{\ast}) = \cV(M)$. Replacing $\Psi$ by $c\Psi$ multiplies $\mu$ by $c^{-2}$, so does not change the set of $\psi_{0}$ with $\mu(\psi_{0})$ constant, and $\gb_{\psi_{0}}$ does not depend on $\Psi$. Hence $\cV(M)$ does not depend on $\Psi$. Moreover, $\cV(M)$ depends only on the flat projective structure of $M$, which determines the pair $(C, \Ga)$ up to the action of $GL(\ste)$. Indeed, for $g \in GL(\ste)$ the function $\psi_{0}\circ g^{-1}$ is the cogauge of the image under $g$ of $L_{\psi_{0}}$, its centroaffine curvature is again constant, because $g$ multiplies $\Psi$ by a constant, and, since the centroaffine metric is defined using only the vector space structure of $\ste$, $g$ induces an isometry between the corresponding centroaffine metrics on the quotients. 

\begin{theorem}\label{suminequalityproposition}
Let $\psi_{0} \in \cAs(C, \Ga)$ satisfy $\mu(\psi_{0}) = \ka$ for a constant $\ka > 0$ and let $\ka^{\ast} = \ka^{-1}$. For every $\psi \in \cAs(C, \Ga)$,
\begin{align}\label{suminequality}
\ka^{1/2}|K_{\psi}/\Ga| + (\ka^{\ast})^{1/2}|K^{\ast}_{\psi}/\Ga^{\ast}| \geq \ka^{1/2}|K_{\psi_{0}}/\Ga| + (\ka^{\ast})^{1/2}|K^{\ast}_{\psi_{0}}/\Ga^{\ast}| = \frac{2}{N}\cV(M),
\end{align}
where $\cV(M)$ is the equiaffine volume, or, equivalently,
\begin{align}\label{ratioinequality}
\frac{|K_{\psi}/\Ga|}{|K_{\psi_{0}}/\Ga|} + \frac{|K^{\ast}_{\psi}/\Ga^{\ast}|}{|K^{\ast}_{\psi_{0}}/\Ga^{\ast}|} \geq 2.
\end{align}
Equality holds in \eqref{suminequality} and \eqref{ratioinequality} if and only if $\psi = \psi_{0}$.
\end{theorem}

\begin{proof}
Write $\psi = f\psi_{0}$ and $\nu = \psi_{0}^{-N}$. By Lemma \ref{covolumeconvexitylemma}, $t \mapsto |K^{\ast}_{(1-t)\psi_{0} + t\psi}/\Ga^{\ast}|$ is convex on $[0,1]$, and, by \eqref{firstvariation}, its derivative at $t = 0$ is $\int_{M}(f - 1)\mu(\psi_{0})\psi_{0}^{-N} = \ka\int_{M}(f-1)\nu$. Hence
\begin{align}\label{rin1}
|K^{\ast}_{\psi}/\Ga^{\ast}| \geq |K^{\ast}_{\psi_{0}}/\Ga^{\ast}| + \ka\int_{M}(f-1)\nu.
\end{align}
Since $x^{-N} - 1 \geq -N(x - 1)$ for $x > 0$, with equality only at $x = 1$,
\begin{align}\label{rin2}
|K_{\psi}/\Ga| - |K_{\psi_{0}}/\Ga| = \frac{1}{N}\int_{M}(f^{-N} - 1)\nu \geq -\int_{M}(f - 1)\nu,
\end{align}
and, because $f$ is continuous and $\nu$ is positive, equality holds in \eqref{rin2} only if $f = 1$. Adding $\ka$ times \eqref{rin2} to \eqref{rin1} and dividing by $\ka^{1/2}$ gives the inequality in \eqref{suminequality}; the equality in \eqref{suminequality} is \eqref{balanced}. If equality holds in \eqref{suminequality}, then it holds in \eqref{rin2}, so $f = 1$, that is, $\psi = \psi_{0}$, and conversely. Dividing \eqref{suminequality} by $\tfrac{1}{N}\cV(M)$, which by \eqref{balanced} equals both $\ka^{1/2}|K_{\psi_{0}}/\Ga|$ and $(\ka^{\ast})^{1/2}|K^{\ast}_{\psi_{0}}/\Ga^{\ast}|$, gives \eqref{ratioinequality}, with the same equality case.
\end{proof}

Since $\psi^{\ast\ast} = \psi$ and $\mu(\psi_{0}^{\ast}) = \ka^{\ast}$ by Lemma \ref{dualityproposition}, the inequalities \eqref{suminequality} and \eqref{ratioinequality} are unchanged when $(\psi, \psi_{0}, \ka)$ is replaced by $(\psi^{\ast}, \psi_{0}^{\ast}, \ka^{\ast})$.

Since every positive multiple of $\psi_{0}$ is also the cogauge of an affine sphere, \eqref{ratioinequality} is a one-parameter family of inequalities, one for each level set of the affine sphere. This observation yields Corollary \ref{weightedcorollary}.
\begin{corollary}\label{weightedcorollary}
Let $\cV(M)$ be the equiaffine volume. For $a, b > 0$ and every $\psi \in \cAs(C, \Ga)$,
\begin{align}\label{weightedinequality}
a|K_{\psi}/\Ga| + b|K^{\ast}_{\psi}/\Ga^{\ast}| \geq \frac{2}{N}(ab)^{1/2}\cV(M),
\end{align}
with equality if and only if $\Sigma_{\psi}$ is the affine sphere with $\mu(\psi) = a/b$, whose polar satisfies $\mu(\psi^{\ast}) = b/a$.
\end{corollary}
\begin{proof}
Apply Theorem \ref{suminequalityproposition} with $\psi_{0}$ the multiple of the affine sphere for which $\ka = a/b$, and multiply by $(ab)^{1/2}$.
\end{proof}

\begin{remark}\label{klartagremark}
Corollary \ref{weightedcorollary} is the hyperbolic counterpart of the construction of affine hemispheres of elliptic type by Klartag \cite{Klartag-hemispheres}. For $p > 0$ and a probability measure $\mathsf{m}$ on $\rea^{n}$ with barycenter at the origin whose support spans $\rea^{n}$, Klartag minimizes
\begin{align}\label{klartagfunctional}
\int_{\rea^{n}}\psi\,d\mathsf{m} + \left(\int_{\rea^{n}}(\psi^{\ast})^{-(n+p)}\,dy\right)^{-2/p},
\end{align}
over proper convex $\psi$ with $\psi(0) < 0$ \cite[Theorem $2.3$]{Klartag-hemispheres}. When $p = 1$ and $\mathsf{m}$ is the uniform probability measure on a bounded open convex set $L \subset \rea^{n}$ with barycenter at the origin, the graph of the minimizer over $L$ is an elliptic affinely spherical hypersurface centered at the origin, and its polar hypersurface is an affine hemisphere \cite[Theorems $1.2$, $3.10$, and $5.10$]{Klartag-hemispheres}.

In both settings an affine sphere is an extremal of a sum of a primal and a dual volume. The dual term is convex, by the Borell--Brascamp--Lieb inequality \cite[Lemma $2.1$]{Klartag-hemispheres} in the elliptic setting and by Lemma \ref{covolumeconvexitylemma} here, and the extremal is characterized by polarity carrying the cone measure of the hypersurface to a constant multiple of the cone measure of its polar \cite[Proposition $5.8$]{Klartag-hemispheres}, which here is the equation $\mu(\psi) = \ka$. The barycenter condition of Klartag, which is necessary for solvability \cite[Proposition $3.3$]{Klartag-hemispheres}, reflects the freedom to add linear functions to $\psi$. Here there is no such freedom. A $\Ga$-invariant linear function $\ell$ vanishes, by the last claim of Lemma \ref{vinberglemma}, because $\ell/\rho_{C}$ is a continuous function on the compact manifold $M$, hence bounded. So the only symmetry is scaling, whose first variation condition is the balance relation \eqref{balanced}. After optimization over scale, the inequality of Klartag bounds a volume product from above, in the direction of the Blaschke--Santal\'o inequality, whereas Corollary \ref{weightedcorollary} is equivalent to the lower bound of Theorem \ref{smoothreversesantalotheorem}.
\end{remark}

\begin{theorem}\label{smoothreversesantalotheorem}
Let $\cV(M)$ be the equiaffine volume. If $\psi_{0} \in \cAs(C, \Ga)$ satisfies $\mu(\psi_{0}) = 1$, then any $\psi \in \cAs(C, \Ga)$ satisfies
\begin{align}\label{smoothreversesantalo}
\cP(\psi) \geq \cP(\psi_{0}) = \frac{\cV(M)^{2}}{N^{2}},
\end{align}
with equality if and only if $\psi$ is a constant multiple of $\psi_{0}$.
\end{theorem}
\begin{proof}
Minimizing the left-hand side of \eqref{weightedinequality} over $a, b$ with $ab = 1$ gives $2\cP(\psi)^{1/2} \geq \tfrac{2}{N}\cV(M)$. Equality forces equality in \eqref{weightedinequality} for the minimizing $a, b$, so $\psi$ is a multiple of $\psi_{0}$. The value at $\psi_{0}$ is \eqref{balanced}.
\end{proof}

\begin{theorem}\label{stabilitytheorem}
Let $\psi_{0} \in \cAs(C, \Ga)$ satisfy $\mu(\psi_{0}) = \ka$ for a constant $\ka > 0$, let $\ka^{\ast} = \ka^{-1}$, let $\nu = \psi_{0}^{-N}$, and, for $t > 0$, let
\begin{align}\label{qdefined}
\mathsf{q}(t) = t^{-N} - 1 + N(t - 1),
\end{align}
which is nonnegative, strictly convex, and vanishes only at $t = 1$.
\begin{enumerate}[label=$(\arabic*)$, ref=(\arabic*)]
\item\label{stt1} For every $\psi \in \cA(C, \Ga)$,
\begin{align}\label{sumdeficit}
\ka^{1/2}|K_{\psi}/\Ga| + (\ka^{\ast})^{1/2}|K^{\ast}_{\psi}/\Ga^{\ast}| - \frac{2}{N}\cV(M) \geq \frac{\ka^{1/2}}{N}\int_{M}\mathsf{q}\left(\psi/\psi_{0}\right)\nu.
\end{align}
In particular \eqref{suminequality} and \eqref{weightedinequality} hold for every $\psi \in \cA(C, \Ga)$, with equality in \eqref{suminequality} if and only if $\psi = \psi_{0}$.
\item\label{stt2} Suppose $\ka = 1$, so that $\nu = \vol_{\gb_{\psi_{0}}}$. For $\psi \in \cA(C, \Ga)$ let $s = \left(|K_{\psi}/\Ga|/|K^{\ast}_{\psi}/\Ga^{\ast}|\right)^{1/(2N)}$. Then
\begin{align}\label{productdeficit}
\cP(\psi)^{1/2} - \frac{\cV(M)}{N} \geq \frac{1}{2N}\int_{M}\mathsf{q}\left(s\psi/\psi_{0}\right)\nu.
\end{align}
In particular \eqref{smoothreversesantalo} holds for every $\psi \in \cA(C, \Ga)$, with equality if and only if $\psi$ is a constant multiple of $\psi_{0}$.
\item\label{stt3} In the setting of \ref{stt2}, $\mathfrak{c}(M) \leq s\psi/\psi_{0} \leq \mathfrak{c}(M)^{-1}$, and
\begin{align}\label{stability}
\cP(\psi)^{1/2} - \frac{\cV(M)}{N} \geq \frac{N+1}{4}\,\mathfrak{c}(M)^{N+2}\int_{M}\left(\frac{s\psi}{\psi_{0}} - 1\right)^{2}\nu.
\end{align}
\end{enumerate}
\end{theorem}
\begin{proof}
Write $f = \psi/\psi_{0}$. Suppose first that $\psi \in \cAs(C, \Ga)$. Adding to \eqref{rin1} the product of $\ka$ with the equality $|K_{\psi}/\Ga| - |K_{\psi_{0}}/\Ga| = \tfrac{1}{N}\int_{M}(f^{-N} - 1)\nu$ of \eqref{rin2} gives
\begin{align}
\ka|K_{\psi}/\Ga| + |K^{\ast}_{\psi}/\Ga^{\ast}| \geq \ka|K_{\psi_{0}}/\Ga| + |K^{\ast}_{\psi_{0}}/\Ga^{\ast}| + \frac{\ka}{N}\int_{M}\mathsf{q}(f)\nu,
\end{align}
and dividing by $\ka^{1/2}$ and using \eqref{balanced} gives \eqref{sumdeficit}. For $\psi \in \cA(C, \Ga)$, by Lemma \ref{densitylemma}\ref{dl2} there are $\psi_{k} \in \cAs(C, \Ga)$ such that $\psi_{k}/\psi \to 1$ uniformly. Both covolumes of $\psi_{k}$ converge to those of $\psi$, as shown in the derivation of \eqref{pcontinuity}, and, since $f$ is a positive continuous function on the compact manifold $M$, $\psi_{k}/\psi_{0} = (\psi_{k}/\psi)f$ converges uniformly to $f$ and takes values in a compact subinterval of $(0, \infty)$, so $\int_{M}\mathsf{q}(\psi_{k}/\psi_{0})\nu \to \int_{M}\mathsf{q}(f)\nu$. Hence \eqref{sumdeficit} holds for $\psi$. Since $\mathsf{q}(f)$ is continuous and nonnegative, the right-hand side of \eqref{sumdeficit} is nonnegative, and it vanishes if and only if $f = 1$. This gives the claims about \eqref{suminequality}, and those about \eqref{weightedinequality} follow as in the proof of Corollary \ref{weightedcorollary}. This proves \ref{stt1}.

Since $(s\psi)^{\ast} = s^{-1}\psi^{\ast}$, \eqref{pdefined} gives $|K_{s\psi}/\Ga| = s^{-N}|K_{\psi}/\Ga|$ and $|K^{\ast}_{s\psi}/\Ga^{\ast}| = s^{N}|K^{\ast}_{\psi}/\Ga^{\ast}|$, and by the choice of $s$ both equal $\cP(\psi)^{1/2}$. Applying \eqref{sumdeficit} with $\ka = 1$ to $s\psi$ and dividing by $2$ gives \eqref{productdeficit}. If equality holds in \eqref{smoothreversesantalo}, the left-hand side of \eqref{productdeficit} vanishes, so $s\psi = \psi_{0}$. This proves \ref{stt2}.

Let $a$ and $b$ be the minimum and the maximum of $sf$ on $M$. By \eqref{balanced} with $\ka = 1$ and by \ref{stt2}, $|K_{s\psi}/\Ga| = |K^{\ast}_{s\psi}/\Ga^{\ast}| = \cP(\psi)^{1/2} \geq \tfrac{1}{N}\cV(M) = |K_{\psi_{0}}/\Ga| = |K^{\ast}_{\psi_{0}}/\Ga^{\ast}|$. If $a > 1$, then $|K_{s\psi}/\Ga| = \tfrac{1}{N}\int_{M}(sf)^{-N}\nu < |K_{\psi_{0}}/\Ga|$, a contradiction. If $b < 1$, then $s\psi \leq b\psi_{0}$, so $(s\psi)^{\ast} \geq b^{-1}\psi_{0}^{\ast}$, because polarity reverses order, and $|K^{\ast}_{s\psi}/\Ga^{\ast}| \leq b^{N}|K^{\ast}_{\psi_{0}}/\Ga^{\ast}| < |K^{\ast}_{\psi_{0}}/\Ga^{\ast}|$, a contradiction.

Hence $a \leq 1 \leq b$, and, by Definitions \ref{sharpharnackdefinition} and \ref{maxmindefinition}, $a \geq \mathfrak{c}(M)b \geq \mathfrak{c}(M)$ and $b \leq \mathfrak{c}(M)^{-1}a \leq \mathfrak{c}(M)^{-1}$. Since $\mathsf{q}(1) = 0 = \mathsf{q}'(1)$ and $\mathsf{q}''(t) = N(N+1)t^{-N-2} \geq N(N+1)\mathfrak{c}(M)^{N+2}$ for $t \leq \mathfrak{c}(M)^{-1}$, the Taylor formula gives $\mathsf{q}(t) \geq \tfrac{1}{2}N(N+1)\mathfrak{c}(M)^{N+2}(t-1)^{2}$ for $\mathfrak{c}(M) \leq t \leq \mathfrak{c}(M)^{-1}$, and \eqref{stability} follows from \eqref{productdeficit}. This proves \ref{stt3}.
\end{proof}

Theorem \ref{stabilitytheorem} also gives the uniqueness of the affine sphere up to homothety. If $\psi_{1} \in \cAs(C, \Ga)$ satisfies $\mu(\psi_{1}) = \ka$, with the same constant $\ka$ as $\psi_{0}$, which can always be arranged by multiplying $\psi_{1}$ by a positive constant, then exchanging the roles of $\psi_{0}$ and $\psi_{1}$ in \eqref{suminequality} and using \eqref{balanced} shows that $\vol(M, \gb_{\psi_{1}}) = \vol(M, \gb_{\psi_{0}})$, so that equality holds in \eqref{suminequality} for $\psi = \psi_{1}$, and $\psi_{1} = \psi_{0}$.

The integral on the right-hand side of \eqref{stability} is unchanged when $\psi$ is rescaled, because $s\psi$ is, and it is bounded below by $\inf_{a > 0}\int_{M}(a\psi/\psi_{0} - 1)^{2}\nu$, the square of the distance in $L^{2}(\nu)$ from $1$ to the ray of positive multiples of $\psi/\psi_{0}$. By Theorem C, as $t \to 0$ the left-hand side of \eqref{productdeficit} for $\psi = (1 + tf)\psi_{0}$ is $\tfrac{t^{2}}{4}\int_{M}(|df|^{2}_{\gb_{\psi_{0}}} + 2N(f - \bar{f})^{2})\nu + O(t^{3})$, whereas \eqref{stability} retains only the remainder of the convexity of $x \mapsto x^{-N}$, which accounts for the constant $\tfrac{N+1}{4}$ in place of $\tfrac{N}{2}$ and for the absence of a gradient term.

\begin{remark}\label{deficitremark}
Let $\psi_{0}$, $\ka$, and $\ka^{\ast}$ be as in Theorem \ref{suminequalityproposition}. The pointwise identity $\ka^{1/2} + \ka^{-1/2}\mu - 2\mu^{1/2} = (\ka^{1/4} - (\ka^{\ast})^{1/4}\mu^{1/2})^{2}$ and Lemma \ref{centroaffinemetriclemma}\ref{cml1} give, for every $\psi \in \cAs(C, \Ga)$,
\begin{align}\label{deficitidentity}
\begin{aligned}
 - 2\vol(M, \gb_{\psi}) & + N\left(\ka^{1/2}|K_{\psi}/\Ga| + (\ka^{\ast})^{1/2}|K^{\ast}_{\psi}/\Ga^{\ast}|\right) 
\\ &
= \int_{M}\left(\ka^{1/4} - (\ka^{\ast})^{1/4}\mu(\psi)^{1/2}\right)^{2}\psi^{-N},
\end{aligned}
\end{align}
the right-hand side of which measures the failure of $\Sigma_{\psi}$ to be the affine sphere $\Sigma_{\psi_{0}} = \{x \in C: \psi_{0}(x) = 1\}$. Thus \eqref{suminequality} is equivalent to
\begin{align}\label{deficitform}
\cV(M) - \vol(M, \gb_{\psi}) \leq \tfrac{1}{2}\int_{M}\left(\ka^{1/4} - (\ka^{\ast})^{1/4}\mu(\psi)^{1/2}\right)^{2}\psi^{-N},
\end{align}
an upper bound for the defect of the centroaffine volume in terms of the deficit in the affine sphere equation.
\end{remark}

\begin{remark}\label{fillastreremark}
In the setting of Corollary \ref{fuchsiancorollary}, with $\Ga \subset SO^{+}(n, 1)$, the inequality \eqref{fuchsianinequality} and its equality case also follow from results of Fillastre \cite{Fillastre-Fuchsian}. Let $B$ be the $\Ga$-convex body bounded by the hyperboloid $\{x \in C: \lb x, x\ra = -1\}$, let $V = \vol_{\mathrm{hyp}}(M)$, and let $\mathsf{h}$ be the support function of a $\Ga$-convex body $L$ in the convention of \cite{Fillastre-Fuchsian}, in which it is negative, regarded as a function on $M$ by restriction to the hyperboloid as in \cite[Sections $2.3$ and $3.1$]{Fillastre-Fuchsian}. Integrals over $M$ are taken with respect to the hyperbolic volume $d\vol_{\mathrm{hyp}}$.

In the notation of \cite{Fillastre-Fuchsian} for covolumes and mixed covolumes, $\operatorname{covol}(B) = V/N$ and $\operatorname{covol}(L, B, \dots, B) = \tfrac{1}{N}\int_{M}|\mathsf{h}|\,d\vol_{\mathrm{hyp}}$, by \cite[Lemma $3.9$(v)]{Fillastre-Fuchsian} when $\pr L$ is smooth and locally strongly convex, and by the continuity of the mixed covolume in general. The first reversed Minkowski inequality \cite[Theorem $5.5$]{Fillastre-Fuchsian}, applied with $K_{1} = L$ and $K_{2} = B$, therefore gives
\begin{align}\label{fillastrecovolume}
\operatorname{covol}(L) \geq \frac{1}{NV^{N-1}}\left(\int_{M}|\mathsf{h}|\,d\vol_{\mathrm{hyp}}\right)^{N}.
\end{align}
Under the identification of $\std$ with $\ste$ given by $-\lb\,,\,\ra$, the polar $L^{\circ}$ is the polar dual of \cite{Fillastre-Fuchsian}, whose boundary is the radial graph of $1/|\mathsf{h}|$ over the hyperboloid. Hence the covolume of $L^{\circ}$ is $\tfrac{1}{N}\int_{M}|\mathsf{h}|^{-N}\,d\vol_{\mathrm{hyp}}$, and, by Jensen's inequality,
\begin{align}\label{fillastrepolarcovolume}
\frac{1}{N}\int_{M}|\mathsf{h}|^{-N}\,d\vol_{\mathrm{hyp}} \geq \frac{V^{N+1}}{N}\left(\int_{M}|\mathsf{h}|\,d\vol_{\mathrm{hyp}}\right)^{-N},
\end{align}
with equality if and only if $\mathsf{h}$ is constant. Multiplying \eqref{fillastrecovolume} and \eqref{fillastrepolarcovolume} gives \eqref{fuchsianinequality}, and equality forces $\mathsf{h}$ to be constant, so that $L$ is a homothetic image of $B$.

The product of the covolumes of a body and its polar is not considered in \cite{Fillastre-Fuchsian}. The proof of Theorem \ref{suminequalityproposition} extends this argument to divisible cones, for which no explicit reference body is available, with the convexity of the polar covolume (Lemma \ref{covolumeconvexitylemma}) in place of the reversed Minkowski inequality.
\end{remark}

\section{Critical points and second variation of the scale-invariant functionals \texorpdfstring{$\mathcal{T}_{\be}$}{T beta}}\label{familysection}

For $\be \in \rea$ let
\begin{align}\label{tbetadefined}
&\mathcal{T}_{\be}(\psi) = \cS_{\be}(\psi)\cS_{0}(\psi)^{2\be - 1},& &\mathcal{T}'_{\be}(\psi) = \cS_{\be}(\psi)\cS_{1}(\psi)^{1 - 2\be}.&
\end{align}
By \eqref{csscaling}, these functionals are scale invariant, and, by Lemma \ref{dualityproposition}, $\mathcal{T}'_{\be}(\psi) = \mathcal{T}_{1-\be}(\psi^{\ast})$, so polarity exchanges the two families. In particular, $\mathcal{T}_{1} = \mathcal{T}'_{0} = N^{2}\cP$ and $\mathcal{T}_{1/2} = \mathcal{T}'_{1/2} = \cS_{1/2}$.

This section treats the critical points and the second variations of the functionals \eqref{tbetadefined}.
Their versions on bounded domains are treated in Section \ref{relativesection}. 

\begin{lemma}\label{sbetafirstvariationlemma}
For $\psi \in \cAs(C, \Ga)$, $f \in \cinf(M)$, $\mu = \mu(\psi)$, $\nu = \psi^{-N}$, and $F = \mu^{\be - 1}$,
\begin{align}\label{sbetafirstvariation}
\tfrac{d}{dt}_{|t=0}\cS_{\be}((1+tf)\psi) = N(2\be - 1)\int_{M}f\mu^{\be}\nu + \be\int_{M}\lb dF, df\ra_{\gb_{\psi}}\mu\nu.
\end{align}
\end{lemma}
\begin{proof}
By \eqref{deltamu}, the derivative is $\int_{M}\mu^{\be}\nu\left(\be(2Nf + \Lop_{\psi}f) - Nf\right)$, and applying Lemma \ref{ibplemma} with $f_{1} = \mu^{\be-1}$ converts $\be\int_{M}\mu^{\be}\psi^{-N}\Lop_{\psi}f$ into the last term.
\end{proof}

\begin{theorem}\label{criticalpointstheorem}
Let $\psi \in \cAs(C, \Ga)$. If $\Sigma_{\psi}$ is a hyperbolic affine sphere, then $\psi$ is critical for $\mathcal{T}_{\be}$ and $\mathcal{T}'_{\be}$ for every $\be \in \rea$. Conversely:
\begin{enumerate}[label=$(\arabic*)$, ref=(\arabic*)]
\item\label{cpt1} If $\be \in (-\infty, 0)\cup(0, \tfrac{1}{2}]\cup[1, \infty)$ and $\psi$ is critical for $\mathcal{T}_{\be}$, then $\Sigma_{\psi}$ is an affine sphere.
\item\label{cpt2} If $\be \in (-\infty, 0]\cup[\tfrac{1}{2}, 1)\cup(1, \infty)$ and $\psi$ is critical for $\mathcal{T}'_{\be}$, then $\Sigma_{\psi}$ is an affine sphere.
\end{enumerate}
\end{theorem}

\begin{proof}
For $\mu = \mu(\psi)$, let $\delta_{\psi}$ denote the formal adjoint of $d$ with respect to $\gb_{\psi}$ and the measure $\mu\nu = \mu^{1/2}\vol_{\gb_{\psi}}$, so that the second term of \eqref{sbetafirstvariation} takes the form $\be\int_{M}f(\delta_{\psi}dF)\mu\nu$.
Write $c_{0} = \cS_{\be}(\psi)/\cS_{0}(\psi)$ and $c_{1} = \cS_{\be}(\psi)/\cS_{1}(\psi)$. By Lemma \ref{sbetafirstvariationlemma} and \eqref{firstvariation}, $\psi$ is critical for $\mathcal{T}_{\be}$, respectively $\mathcal{T}'_{\be}$, if and only if
\begin{align}
\label{eltbeta}&N(2\be - 1)(\mu^{\be} - c_{0}) + \be\mu\,\delta_{\psi}dF = 0,\\
\label{eltprimebeta}&N(2\be - 1)(F - c_{1}) + \be\,\delta_{\psi}dF = 0,
\end{align}
the second after division by $\mu$, using $\mu^{\be} = \mu F$. If $\mu$ is constant, then $F$ is constant, $c_{0} = \mu^{\be}$, and $c_{1} = F$, so both \eqref{eltbeta} and \eqref{eltprimebeta} hold. This proves that if $\Sigma_{\psi}$ is a hyperbolic affine sphere, then $\psi$ is critical for $\mathcal{T}_{\be}$ and $\mathcal{T}'_{\be}$ for every $\be \in \rea$. 

Suppose \eqref{eltprimebeta} holds. Multiplying by $(F - c_{1})\mu\nu$ and integrating gives
\begin{align}\label{betadirichlet}
\be\int_{M}|dF|^{2}_{\gb_{\psi}}\mu\nu = -N(2\be - 1)\int_{M}(F - c_{1})^{2}\mu\nu.
\end{align}
If $\be > 1/2$, the left-hand side of \eqref{betadirichlet} is nonnegative while the right-hand side is nonpositive. If $\be \leq 0$, the left-hand side is nonpositive while the right-hand side is nonnegative. If $\be = 1/2$, the right-hand side vanishes. In all three cases both sides of \eqref{betadirichlet} vanish, so $F$ is constant: for $\be \neq 0$ because the left-hand side vanishes, and for $\be = 0$ because the right-hand side does. Since $F = \mu^{\be - 1}$ and $\be \neq 1$ in the range of \ref{cpt2}, $\mu$ is constant as well. This proves \ref{cpt2}.

Since $\mathcal{T}_{\be}(\psi) = \mathcal{T}'_{1-\be}(\psi^{\ast})$ and polarity maps $\cAs(C, \Ga)$ bijectively onto $\cAs(C^{\ast}, \Ga^{\ast})$, carrying smooth variations to smooth variations, $\psi$ is critical for $\mathcal{T}_{\be}$ if and only if $\psi^{\ast}$ is critical for $\mathcal{T}'_{1-\be}$. By Lemma \ref{dualityproposition}, $\mu(\psi)$ is constant if and only if $\mu(\psi^{\ast})$ is, and $\be \mapsto 1 - \be$ maps the range of \ref{cpt1} onto that of \ref{cpt2}, so \ref{cpt1} follows from \ref{cpt2} applied to $(C^{\ast}, \Ga^{\ast})$.
\end{proof}

The cases left open by Theorem \ref{criticalpointstheorem} reduce to one: $\mathcal{T}_{0} = \mathcal{T}'_{1} = 1$ identically, and, since $\mathcal{T}_{\be}(\psi) = \mathcal{T}'_{1-\be}(\psi^{\ast})$, the functionals $\mathcal{T}_{\be}$ with $\tfrac{1}{2} < \be < 1$ correspond to the functionals $\mathcal{T}'_{\be}$ with $0 < \be < \tfrac{1}{2}$, for which the two sides of \eqref{betadirichlet} have the same sign and the argument does not decide whether every critical point is an affine sphere.

By Theorem \ref{criticalpointstheorem}, for every $\be \in \rea$ there is a scale invariant functional, $\mathcal{T}_{\be}$ or $\mathcal{T}'_{\be}$, whose only critical points are the affine spheres. For $\be = 0$ this functional is $\mathcal{T}'_{0} = N^{2}\cP$, so this is Lemma \ref{firstvariationproposition}. For the classical equiaffine surface area, $\be = 1/(n+2)$, this applies to $\mathcal{T}_{1/(n+2)}$. For the affine maximal case $\be = (n+1)/(n+2)$ of Remark \ref{affinemaximalremark} it applies to $\mathcal{T}'_{(n+1)/(n+2)}$, which by duality equals $\mathcal{T}_{1/(n+2)}$ evaluated on the polar.

For $\be = 1/2$ the critical points of the centroaffine volume $\cS_{1/2}$ are the equivariant counterparts of the centroaffine minimal hypersurfaces of Wang \cite{Wang-centroaffine}, the critical points of the area of the centroaffine metric under compactly supported variations, now often called centroaffine extremal hypersurfaces \cite{Li-Li-Simon}. Wang showed that proper affine spheres centered at the origin are centroaffine minimal \cite[Example $1$]{Wang-centroaffine}, and that so are the hypersurfaces $\ell_{1}^{\al_{1}}\cdots\ell_{N}^{\al_{N}} = 1$ in the cone $\{\ell_{1} > 0, \dots, \ell_{N} > 0\}$, where $\ell_{1}, \dots, \ell_{N}$ is a basis of $\std$ and $\al_{i} > 0$; Wang states this for coordinate functions, and the general case follows because centroaffine minimality is invariant under $GL(\ste)$. For $\al_{i}$ not all equal these hypersurfaces are not affine spheres \cite[Example $2$]{Wang-centroaffine}, and, consistent with Theorem \ref{criticalpointstheorem}, they are invariant under no lattice $\Ga$ as in Example \ref{orthantexample}. Remark \ref{relativecriticalremark} gives the analogous statement for compactly supported perturbations of the level sets of the Cheng--Yau solution of a bounded domain.

Complete centroaffine extremal hypersurfaces in $\rea^{N}$ that are not affine spheres are plentiful; see \cite{Xing-Sun-Xu} and the references there.

\begin{theorem}\label{tbetasecondvariationtheorem}
Let $\psi_{0} \in \cAs(C, \Ga)$ satisfy $\mu(\psi_{0}) = \ka$ for a constant $\ka > 0$, let $\nu = \psi_{0}^{-N}$ and $V = \int_{M}\nu$, and let $\Delta$ be the Laplacian of $\gb_{\psi_{0}}$. For $f \in \cinf(M)$ and $\psi_{t} = (1+tf)\psi_{0}$,
\begin{align}\label{sbetasecond}
\tfrac{d^{2}}{dt^{2}}_{|t=0}\cS_{\be}(\psi_{t}) = 2\ka^{\be}\int_{M}\left(b_{0}f^{2} + b_{1}|df|^{2} + b_{2}(\Delta f)^{2}\right)\nu,
\end{align}
with $b_{0} = 2N^{2}\be^{2} - N(2N+1)\be + \tfrac{N(N+1)}{2}$, $b_{1} = 2N\be^{2} - \tfrac{3N}{2}\be$, and $b_{2} = \tfrac{\be(\be - 1)}{2}$, and
\begin{align}\label{tbetasecond}
\begin{aligned}
\tfrac{d^{2}}{dt^{2}}_{|t=0}&\log\mathcal{T}_{\be}(\psi_{t}) \\
&= \frac{1}{V}\int_{M}\left(2N^{2}\be(2\be - 1)(f - \bar{f})^{2} + N\be(4\be - 3)|df|^{2} + \be(\be - 1)(\Delta f)^{2}\right)\nu,
\end{aligned}
\end{align}
where norms are taken with respect to $\gb_{\psi_{0}}$.
\end{theorem}
\begin{proof}
By Lemma \ref{laplacianlemma}, $\Lop_{\psi_{0}} = -\Delta$. Write $\mu(\psi_{t}) = \ka(1 + ta_{1} + t^{2}a_{2} + O(t^{3}))$. By \eqref{deltamu}, $a_{1} = 2Nf - \Delta f$. Write $I_{0} = \int_{M} f^{2}\nu$, $I_{1} = \int_{M}|df|^{2}\nu$, and $I_{2} = \int_{M}(\Delta f)^{2}\nu$, so that $\int_{M} fa_{1}\nu = 2NI_{0} + I_{1}$ and $\int_{M} a_{1}^{2}\nu = 4N^{2}I_{0} + 4NI_{1} + I_{2}$. Expanding $\cS_{1}(\psi_{t}) = \ka\int_{M}(1 + ta_{1} + t^{2}a_{2})(1+tf)^{-N}\nu$ to second order and comparing with \eqref{secondcovolume}, which gives $\tfrac{d^{2}}{dt^{2}}\cS_{1} = N\ka((N-1)I_{0} + I_{1})$, yields $\int_{M} a_{2}\nu = N(2N-1)I_{0} + \tfrac{3N}{2}I_{1}$. Expanding $\cS_{\be}(\psi_{t}) = \ka^{\be}\int_{M}(1 + ta_{1} + t^{2}a_{2})^{\be}(1+tf)^{-N}\nu$, the coefficient of $t^{2}$ is
\begin{align}
\ka^{\be}\int_{M}\left(\be a_{2} + \tfrac{\be(\be-1)}{2}a_{1}^{2} - N\be fa_{1} + \tfrac{N(N+1)}{2}f^{2}\right)\nu,
\end{align}
and substituting gives \eqref{sbetasecond}. Similarly $\tfrac{d}{dt}_{|t=0}\cS_{\be}(\psi_{t}) = \ka^{\be}N(2\be - 1)\int_{M} f\nu$. With $\cS_{0}(\psi_{0}) = V$, $\cS_{\be}(\psi_{0}) = \ka^{\be}V$, $\tfrac{d}{dt}\cS_{0} = -N\int_{M} f\nu$, and $\tfrac{d^{2}}{dt^{2}}\cS_{0} = N(N+1)I_{0}$, the second derivative of $\log\mathcal{T}_{\be} = \log\cS_{\be} + (2\be-1)\log\cS_{0}$ is
\begin{align}
\frac{2}{V}(b_{0}I_{0} + b_{1}I_{1} + b_{2}I_{2}) - N^{2}(2\be-1)^{2}\bar{f}^{2} + (2\be-1)\left(\frac{N(N+1)}{V}I_{0} - N^{2}\bar{f}^{2}\right),
\end{align}
in which the coefficient of $I_{0}/V$ is $2b_{0} + (2\be-1)N(N+1) = 2N^{2}\be(2\be-1)$ and that of $\bar{f}^{2}$ is $-2N^{2}\be(2\be-1)$. This gives \eqref{tbetasecond}.
\end{proof}

\begin{corollary}\label{tbetaspectralcorollary}
Let $\psi_{0} \in \cAs(C, \Ga)$ satisfy $\mu(\psi_{0}) = \ka$ for a constant $\ka > 0$, let $\nu = \psi_{0}^{-N}$ and $V = \int_{M}\nu$, let $\Delta$ be the Laplacian of $\gb_{\psi_{0}}$, let $\be \in \rea$, and let $\psi_{t} = (1 + tf)\psi_{0}$. If $f \in \cinf(M)$ satisfies $-\Delta f = \la f$ with $\la > 0$, then $\tfrac{d^{2}}{dt^{2}}_{|t=0}\log\mathcal{T}_{\be}(\psi_{t}) = \tfrac{1}{V}Q_{\be}(\la)\int_{M}f^{2}\nu$, where
\begin{align}\label{hbeta}
Q_{\be}(\la) = \be\left(2N^{2}(2\be - 1) + N(4\be - 3)\la + (\be - 1)\la^{2}\right).
\end{align}
Consequently, on variations transverse to rescalings, the second variation of $\mathcal{T}_{\be}$ at $\psi_{0}$ is positive definite for $\be < 0$ and for $\be \geq 1$, it is negative definite for $0 < \be \leq 1/2$, and it is not positive semidefinite for $1/2 < \be < 1$. In particular,
\begin{align}\label{centroaffinesecond}
\tfrac{d^{2}}{dt^{2}}_{|t=0}\cS_{1/2}(\psi_{t}) = -\ka^{1/2}\int_{M}\left(\tfrac{N}{2}|df|^{2} + \tfrac{1}{4}(\Delta f)^{2}\right)\nu
\end{align}
for every $f \in \cinf(M)$, so the second variation of the centroaffine volume $\cS_{1/2}$ at $\psi_{0}$ is negative definite on variations transverse to rescalings. By duality (Lemmas \ref{dualityproposition} and \ref{centroaffinemetriclemma}\ref{cml3}), the same conclusions hold for $\mathcal{T}'_{1-\be}$ in place of $\mathcal{T}_{\be}$.
\end{corollary}

\begin{proof}
For $\be \geq 1$ all three coefficients inside the parentheses of \eqref{hbeta} are nonnegative and the first is positive; for $\be < 0$ all three are negative, and so is $\be$; for $0 < \be \leq 1/2$ the first two are nonpositive, not both zero, and the third is negative. For $1/2 < \be < 1$, $Q_{\be}(\la) \to -\infty$ as $\la \to \infty$, and the eigenvalues of $-\Delta$ are unbounded. Formula \eqref{centroaffinesecond} is \eqref{sbetasecond} with $\be = 1/2$, for which $b_{0} = 0$.
\end{proof}

Formula \eqref{centroaffinesecond} is the equivariant form of the second variation formula of Wang at a hyperbolic affine sphere \cite[Corollary $3.1$]{Wang-centroaffine}. Wang varies the hypersurface radially, by $x \mapsto (1 + t\phi(x))x$ for a function $\phi$ homogeneous of degree $0$. Since $\psi_{t}(x/(1 + tf(x))) = 1$ when $\psi_{0}(x) = 1$, the hypersurface $\Sigma_{\psi_{t}}$ is the image of $\Sigma_{\psi_{0}}$ under $x \mapsto x/(1 + tf(x))$, which agrees with the radial variation with $\phi = -f$ up to terms of order $t^{2}$. At a critical point the second derivative of a functional along a curve depends only on the first order term of the curve, so the two second variation formulas agree. 

In the terminology of Wang, in which a centroaffine minimal hypersurface is stable if the second variation of its centroaffine area is nonpositive, \eqref{centroaffinesecond} says that the affine sphere is stable, indeed strictly stable transverse to rescalings, with respect to $\Ga$-invariant variations, as it is with respect to compactly supported ones \cite[Corollary $3.1$]{Wang-centroaffine}. This is the Hamiltonian stability of the corresponding Lagrangian submanifold in Remark \ref{lagrangianremark}.

\begin{remark}\label{affinemaximalremark}
By Lemma \ref{dualityproposition}, $\cS_{(n+1)/(n+2)}(\psi) = \cS_{1/(n+2)}(\psi^{\ast})$ is the equivariant equiaffine surface area of the polar hypersurface $\Sigma_{\psi^{\ast}}$. Since $(n+1)/(n+2) \in [1/2, 1)$, Theorem \ref{criticalpointstheorem}\ref{cpt2} shows that the affine spheres are the only critical points of $\mathcal{T}'_{(n+1)/(n+2)}$, and, since $\mathcal{T}'_{(n+1)/(n+2)}(\psi) = \mathcal{T}_{1/(n+2)}(\psi^{\ast})$ with $1/(n+2) \in (0, 1/2]$, Corollary \ref{tbetaspectralcorollary} shows that their second variations are negative definite transverse to rescalings. Remark \ref{relativecriticalremark} gives the analogous uniqueness for compactly supported perturbations of the level sets of the Cheng--Yau solution of a bounded domain.
\end{remark}

\begin{theorem}\label{reductionproposition}
Let $\psi_{0} \in \cAs(C, \Ga)$ satisfy $\mu(\psi_{0}) = 1$, and let $\cV(M)$ be the equiaffine volume. For every $\psi \in \cAs(C, \Ga)$,
\begin{align}
\label{holder1}&\mathcal{T}_{\be}(\psi) \leq \cS_{1/2}(\psi)^{2\be} && (0 \leq \be \leq \tfrac{1}{2}),\\
\label{holder2}&\mathcal{T}_{\be}(\psi) \geq \cS_{1/2}(\psi)^{2\be} && (\be < 0),\\
\label{holder3}&\mathcal{T}_{\be}(\psi) \geq \left(N^{2}\cP(\psi)\right)^{2\be - 1}\cS_{1/2}(\psi)^{2 - 2\be} \geq \cV(M)^{2(2\be-1)}\cS_{1/2}(\psi)^{2-2\be} && (\be \geq 1),
\end{align}
and $\mathcal{T}_{\be}(\psi_{0}) = \cV(M)^{2\be}$. Moreover,
\begin{align}\label{holder4}
&\mathcal{T}_{\be}(\psi) \geq \left(N^{2}\cP(\psi)\right)^{\be} \geq \cV(M)^{2\be} = \mathcal{T}_{\be}(\psi_{0})&& (\be \geq 1),
\end{align}
with equality in the second inequality if and only if $\psi$ is a constant multiple of $\psi_{0}$, so $\psi_{0}$ minimizes $\mathcal{T}_{\be}$ for $\be \geq 1$, and its positive multiples are the only minimizers. If the centroaffine volume satisfies $\cS_{1/2}(\psi) \leq \cV(M)$ for all $\psi \in \cAs(C, \Ga)$, then $\psi_{0}$ maximizes $\mathcal{T}_{\be}$ for $0 \leq \be \leq 1/2$ and minimizes it for $\be < 0$.
\end{theorem}

\begin{proof}
The inequalities \eqref{holder1}, \eqref{holder2}, and \eqref{holder3} are the H\"older inequality for the log-convex function $\be \mapsto \cS_{\be}(\psi)$. The inequality \eqref{holder1} interpolates $\be$ between $0$ and $1/2$, \eqref{holder2} interpolates $0$ between $\be$ and $1/2$, and \eqref{holder3} interpolates $1$ between $1/2$ and $\be$, giving $\cS_{1} \leq \cS_{1/2}^{(2\be-2)/(2\be-1)}\cS_{\be}^{1/(2\be-1)}$. The second inequality in \eqref{holder3} is Theorem A. At $\psi_{0}$, $\mu = 1$, so $\cS_{\be}(\psi_{0}) = \cV(M)$ for all $\be$. For \eqref{holder4}, the first inequality is an identity when $\be = 1$, and, when $\be > 1$, interpolating $1$ between $0$ and $\be$ gives $\cS_{1} \leq \cS_{0}^{(\be - 1)/\be}\cS_{\be}^{1/\be}$, so that $\cS_{\be} \geq \cS_{1}^{\be}\cS_{0}^{1 - \be}$ and $\mathcal{T}_{\be} \geq (\cS_{1}\cS_{0})^{\be} = (N^{2}\cP)^{\be}$, and the second inequality in \eqref{holder4} and its equality case are Theorem A. The last statement follows by comparing exponents, since $2\be \geq 0$ in \eqref{holder1} and $2\be < 0$ in \eqref{holder2}.
\end{proof}
The inequalities \eqref{holder1}, \eqref{holder2}, and \eqref{holder4} cover all $\be$ outside the interval $(1/2, 1)$, on which, by Corollary \ref{tbetaspectralcorollary}, the affine sphere is not a local minimum of $\mathcal{T}_{\be}$. By \eqref{holder4}, for $\be \geq 1$ the affine sphere is the minimizer of $\mathcal{T}_{\be}$, and, by Lemma \ref{dualityproposition}, it is the minimizer of $\mathcal{T}'_{\be}$ for $\be \leq 0$. This is the sharp equivariant counterpart of the inequality of \cite{Werner-Ye} for $p < -N$ described in Remark \ref{lpremark}. In \eqref{holder3} the first inequality is stronger than the first inequality of \eqref{holder4}, by Lemma \ref{cauchyschwarzproposition}, but the conclusion drawn from it requires the bound $\cS_{1/2}(\psi) \leq \cV(M)$. For $\be \leq 1/2$, Theorem \ref{reductionproposition} reduces the question whether the affine sphere is the extremum of $\mathcal{T}_{\be}$, which is the equivariant form of the $L_{p}$ affine isoperimetric inequalities of Remark \ref{lpremark} for $-N < p \leq N$, to the single inequality $\cS_{1/2}(\psi) \leq \cV(M)$, which is Question \ref{centroaffineconjecture} below. The second variation \eqref{centroaffinesecond} is consistent with that inequality.

\begin{question}\label{centroaffineconjecture}
For a closed properly convex flat projective manifold $M = \Om/\Ga$ of dimension $n \geq 3$, is it true that every $\psi \in \cAs(C, \Ga)$ satisfies $\vol(M, \gb_{\psi}) \leq \cV(M)$, with equality only if $\mu(\psi)$ is constant?
\end{question}
The examples of \cite{Blayac-Viaggi}, for which $\Om$ contains properly embedded cones, so that $\pr\Om$ is not $C^{1}$, are test cases for Question \ref{centroaffineconjecture} in which the affine sphere has flat directions at infinity. In the case of convex bodies the analogous inequality, $\Omega_{N}(K) \leq \Omega_{N}(B)$ for the self-dual affine surface area of Remark \ref{lpremark}, $K$ having Santal\'o point the origin and $B$ the unit ball, follows from the Blaschke--Santal\'o inequality via \eqref{cauchyschwarz}. The inequality of Theorem A goes the other way, so a different argument is needed. In fact every splitting of the geometric mean $\mu^{1/2}\psi^{-N}$ into the integrands of $|K_{\psi}|$ and $|K^{\ast}_{\psi}|$ by the arithmetic--geometric mean inequality yields a bound on $\vol(M, \gb_{\psi})$ by a quantity that is minimized, not maximized, at the affine sphere. In dimension $n = 2$ a different splitting works, as Theorem \ref{lowdimcentroaffinetheorem} below shows.

\begin{lemma}\label{centroaffineformulalemma}
Let $\psi_{0} \in \cAs(C, \Ga)$ satisfy $\mu(\psi_{0}) = 1$, so that $\gb_{\psi_{0}}$ is the equiaffine metric of $\Sigma_{\psi_{0}}$, let $\nu = \psi_{0}^{-N} = \vol_{\gb_{\psi_{0}}}$, let $D$ and $\Delta$ be the Levi-Civita connection and the Laplacian of $\gb_{\psi_{0}}$, let $\dift$ be the difference tensor as in Lemma \ref{gausslemma}, and for a positive function $f \in \cinf(M)$ let $H_{ij} = D_{i}D_{j}f - \dift_{ij}\,^{p}D_{p}f$ be the Hessian $\nabla^{2}f$ of $f$ for the equiaffine connection $\nabla = D + \dift$ of $\Sigma_{\psi_{0}}$. Let $\mathsf{T}$ be the endomorphism of $TM$ determined by the bilinear form $\gb_{\psi_{0}} - f^{-1}H$ with respect to $\gb_{\psi_{0}}$, that is, $\gb_{\psi_{0}}(\mathsf{T}Y, Z) = \gb_{\psi_{0}}(Y, Z) - f^{-1}H(Y, Z)$ for all $Y, Z$. Then $\psi = f\psi_{0}$ is in $\cAs(C, \Ga)$ if and only if $\gb_{\psi_{0}} - f^{-1}H$ is positive definite, and in that case
\begin{align}\label{centroaffineformula}
\mu(f\psi_{0}) = f^{2N}\det \mathsf{T}, \qquad \vol(M, \gb_{f\psi_{0}}) = \int_{M}(\det \mathsf{T})^{1/2}\,\nu,
\end{align}
and $\tr \mathsf{T} = n - f^{-1}\Delta f$. Moreover, as metrics on $M$, $\gb_{f\psi_{0}} = \gb_{\psi_{0}}(\mathsf{T}\,\cdot\,, \,\cdot\,)$.
\end{lemma}
\begin{proof}
In the notation of \eqref{phiqdefined}, $\Phi_{\psi} = \tfrac{1}{2}\psi^{2} = f^{2}\Phi_{\psi_{0}}$. Compute $\fD^{2}\Phi_{\psi}$ at a point of $\Sigma_{\psi_{0}} = \{x \in C: \psi_{0}(x) = 1\}$ in a basis $\eul, e_{1}, \dots, e_{n}$ with $e_{i}$ tangent to $\Sigma_{\psi_{0}}$ and $\gb_{\psi_{0}}$-orthonormal.
\begin{itemize}
\item By homogeneity, $\fD^{2}\Phi_{\psi}(\eul, \eul) = 2\Phi_{\psi} = f^{2}$ and $\fD^{2}\Phi_{\psi}(\eul, e_{i}) = d\Phi_{\psi}(e_{i}) = f\,df(e_{i})$.
\item On vectors tangent to $\Sigma_{\psi_{0}}$, $\fD^{2}\Phi_{\psi} = f^{2}\fD^{2}\Phi_{\psi_{0}} + 2f(df\otimes \fD\Phi_{\psi_{0}} + \fD\Phi_{\psi_{0}}\otimes df) + 2\Phi_{\psi_{0}}(df\otimes df + f\fD^{2}f) = -f^{2}\gb_{\psi_{0}} + df\otimes df + fH$. Here $\fD\Phi_{\psi_{0}} = \psi_{0}d\psi_{0}$ annihilates $T\Sigma_{\psi_{0}}$, $\fD^{2}\Phi_{\psi_{0}} = -\gb_{\psi_{0}}$ there, and $\fD^{2}f$ restricted to $T\Sigma_{\psi_{0}}$ is the Hessian of $f$ for the induced connection with transversal $\eul$, because $f$ is homogeneous of degree $0$. By Lemma \ref{affinespherelemma}, that connection is the equiaffine connection $D + \dift$, because $\eul$ is the equiaffine normal.
\end{itemize}
Let $A_{\psi}$ and $A_{\psi_{0}}$ be the matrices of $\fD^{2}\Phi_{\psi}$ and $\fD^{2}\Phi_{\psi_{0}}$ in the basis $\eul, e_{1}, \dots, e_{n}$, and write $H_{ij}$ and $f_{i}$ for $H(e_{i}, e_{j})$ and $df(e_{i})$. By the two computations above, the first row and column of $A_{\psi}$ are $(f^{2}, ff_{1}, \dots, ff_{n})$, and its remaining entries are $-f^{2}\delta_{ij} + f_{i}f_{j} + fH_{ij}$, while $A_{\psi_{0}} = \operatorname{diag}(1, -I)$. Subtracting $f^{-1}f_{i}$ times the first row of $A_{\psi}$ from its $(i+1)$st row for $i = 1, \dots, n$, which does not change the determinant, makes the first column $(f^{2}, 0, \dots, 0)$ and replaces the remaining entries by $-f^{2}\delta_{ij} + fH_{ij}$. Equivalently, the vectors $e'_{i} = e_{i} - f^{-1}f_{i}\eul$ satisfy $\fD^{2}\Phi_{\psi}(e'_{i}, \eul) = 0$, so that they span $\ker d\psi$, because $\fD^{2}\Phi_{\psi}(\eul, \cdot) = \psi\,d\psi$, and $\fD^{2}\Phi_{\psi}(e'_{i}, e'_{j}) = -f^{2}\delta_{ij} + fH_{ij}$. Hence
\begin{align}
\det A_{\psi} = f^{2}\det\left(-f^{2}I + fH\right) = (-1)^{n}f^{n+2}\det(fI - H),
\end{align}
where $H$ denotes the matrix $(H_{ij})$. By \eqref{mudefined}, for any basis $e_{0}, \dots, e_{n}$ of $\ste$ and any $\psi \in \cAs(C, \Ga)$, $\mu(\psi) = |\det(\fD^{2}\Phi_{\psi}(e_{i}, e_{j}))|\,\Psi(e_{0}, \dots, e_{n})^{-2}$. Since $\det A_{\psi_{0}} = (-1)^{n}$ and $\mu(\psi_{0}) = 1$, this gives $|\Psi(\eul, e_{1}, \dots, e_{n})| = 1$, and then $\mu(f\psi_{0}) = |\det A_{\psi}| = f^{n+2}|\det(fI - H)| = f^{2N}|\det \mathsf{T}|$, because $(I - f^{-1}H)$ is the matrix of $\mathsf{T}$ in the $\gb_{\psi_{0}}$-orthonormal basis $e_{1}, \dots, e_{n}$. By Lemma \ref{polarmaplemma}\ref{pml1}, $\psi \in \cAs(C, \Ga)$ if and only if $\fD^{2}\Phi_{\psi}$ has signature $(1, n)$, and, since $\fD^{2}\Phi_{\psi}(\eul, \eul) = f^{2} > 0$, this holds if and only if the form with matrix $-f^{2}I + fH$ is negative definite, that is, if and only if $\gb_{\psi_{0}} - f^{-1}H$ is positive definite. Then $\mu(\psi)^{1/2}\psi^{-N} = f^{N}(\det \mathsf{T})^{1/2}f^{-N}\psi_{0}^{-N}$, which gives the volume formula by Lemma \ref{centroaffinemetriclemma}\ref{cml1}.

Next, $\tr H = \Delta f$, because $\dift$ is trace-free. Finally, by the second identity of \eqref{mulogdet}, the pullback of $\gb_{\psi}$ to $C$ equals $-\psi^{-1}\fD^{2}\psi = -\psi^{-2}\fD^{2}\Phi_{\psi}$ on $\ker d\psi$, and $\psi = f$ at the point considered, so $\gb_{\psi}(e'_{i}, e'_{j}) = \delta_{ij} - f^{-1}H_{ij}$. Since $e'_{i}$ and $e_{i}$ project to the same tangent vector of $M$ and the $e_{i}$ are $\gb_{\psi_{0}}$-orthonormal, this shows that $\gb_{f\psi_{0}} = \gb_{\psi_{0}}(\mathsf{T}\,\cdot\,, \,\cdot\,)$ on $M$.
\end{proof}

\begin{lemma}\label{blaschkehessianlemma}
Let $n \geq 2$, let $C \subset \ste$ be an open properly convex cone, not necessarily divisible, let $\psi_{0} \in \cAs(C, \{1\})$ satisfy $\mu(\psi_{0}) = \ka$ constant, and let $\Sigma = \Sigma_{\psi_{0}}$, with equiaffine metric $h$, equiaffine connection $\nabla$, and equiaffine normal $\xi = \ka^{1/(n+2)}\eul$. A smooth function $f$ on $\Sigma$ satisfies $\nabla^{2}f = q h$ for some function $q$ if and only if $f = a + \ell|_{\Sigma}$ for a constant $a$ and a linear function $\ell \in \std$, and $a$ and $\ell$ are determined by $f$. Consequently, if $f$ is constant outside a compact subset of $\Sigma$, then $f$ is constant, and if $\Ga$ divides $\Om = \projp(C)$, $\psi_{0}$ is $\Ga$-invariant, and $f$ is $\Ga$-invariant, then $f$ is constant.
\end{lemma}
\begin{proof}
By Lemma \ref{affinespherelemma} the equiaffine normal is $\xi = \ka^{1/(n+2)}\eul$, so the Gauss formula reads $\fD_{X}Y = \nabla_{X}Y + h(X, Y)\xi$ for vector fields $X, Y$ tangent to $\Sigma$. For $\ell \in \std$ this gives $X(Y\ell) = \ell(\fD_{X}Y) = \ell(\nabla_{X}Y) + h(X, Y)\ell(\xi)$, that is,
\begin{align}\label{linearhessian}
\nabla^{2}\ell = \ka^{1/(n+2)}\ell\,h \quad\text{on } \Sigma,
\end{align}
and constants satisfy $\nabla^{2}a = 0$. Conversely, suppose that $\nabla^{2}f = q h$, and write $\la = \ka^{1/(n+2)}$. Differentiating gives $\nabla_{k}\nabla_{i}\nabla_{j}f = (\nabla_{k}q)h_{ij} + q\nabla_{k}h_{ij}$. Antisymmetrize in $k$ and $i$. The terms $q\nabla h$ cancel by the Codazzi equation in \eqref{gausscodazzi}, and, by the Ricci identity and the Gauss equation in \eqref{gausscodazzi}, which with $S = -\la\,\mathrm{id}$ reads $R^{\nabla}(X, Y)Z = \la(h(X, Z)Y - h(Y, Z)X)$,
\begin{align}\label{qgradient}
(\nabla_{k}q)h_{ij} - (\nabla_{i}q)h_{kj} = 2\nabla_{[k}\nabla_{i]}\nabla_{j}f = \la\left(h_{ij}\nabla_{k}f - h_{kj}\nabla_{i}f\right).
\end{align}
Contracting \eqref{qgradient} with $h^{ij}$ gives $(n-1)(dq - \la\,df) = 0$, so $dq = \la\,df$, because $n \geq 2$. Hence the triple $(f, df, q)$ satisfies the first order linear system $\nabla f = df$, $\nabla df = qh$, $\nabla q = \la\,df$, and so is determined along any curve by its value at one point. Since $\Sigma$ is connected, being identified with $\Om$ by radial projection, the solutions of $\nabla^{2}f = q h$ form a vector space of dimension at most $n + 2$. The functions $a + \ell|_{\Sigma}$ form a subspace of dimension $n+2$. If $a + \ell = 0$ on $\Sigma$ with $(a, \ell) \neq (0, 0)$, then $\Sigma$ lies in the hyperplane $\{y \in \ste: \ell(y) = -a\}$, so $\fD_{X}Y$ is tangent to $\Sigma$ for tangent $X$ and $Y$, contradicting the positive definiteness of $h$. Hence every solution has the form $a + \ell|_{\Sigma}$, with $a$ and $\ell$ unique.

If $f = a + \ell|_{\Sigma}$ is $\Ga$-invariant, then for $\gamma \in \Ga$ the linear function $\ell\circ\gamma - \ell$ vanishes on $\Sigma$, hence is zero because $\Sigma$ is not contained in a hyperplane. So $\ell$ is $\Ga$-invariant, $\ell/\rho_{C}$ is the lift of a continuous function on the compact manifold $M$, hence bounded, and $\ell = 0$ by the last claim of Lemma \ref{vinberglemma}. If $f = b$ outside a compact set $Q \subset \Sigma$, then $\ell = b - a$ on the open set $\Sigma\setminus Q$, which is nonempty because $\Sigma$ is not compact. The preceding argument, which is local, shows that no nonempty open subset of $\Sigma$ lies in a hyperplane, so $\ell = 0$ and $f = a$.
\end{proof}

\begin{theorem}\label{lowdimcentroaffinetheorem}
Let $n \leq 2$, let $\psi_{0} \in \cAs(C, \Ga)$ satisfy $\mu(\psi_{0}) = 1$, and let $\nu = \psi_{0}^{-N}$. Every $\psi = f\psi_{0} \in \cAs(C, \Ga)$ satisfies
\begin{align}\label{lowdimcentroaffine}
\vol(M, \gb_{\psi}) \leq \vol(M, \gb_{\psi_{0}}) - \tfrac{1}{2}\int_{M}|d\log f|^{2}_{\gb_{\psi_{0}}}\,\nu,
\end{align}
with equality if and only if $f$ is constant.
\end{theorem}

\begin{proof}
Let $H = \nabla^{2}f$ and $\mathsf{T}$ be as in Lemma \ref{centroaffineformulalemma}, so that $\mathsf{T}$ is the endomorphism of $TM$ defined by $\gb_{\psi_{0}}(\mathsf{T}Y, Z) = \gb_{\psi_{0}}(Y, Z) - f^{-1}H(Y, Z)$, which is positive definite because $\psi \in \cAs(C, \Ga)$. By that lemma, $\vol(M, \gb_{\psi}) = \int_{M}(\det\mathsf{T})^{1/2}\nu$ and $\tr \mathsf{T} - n = -f^{-1}\Delta f$, and, by integration by parts on the closed manifold $M$, $\int_{M}f^{-1}\Delta f\,\nu = \int_{M}f^{-2}|df|^{2}_{\gb_{\psi_{0}}}\nu$. Hence, in every dimension, with $R = 1 + \tfrac{1}{2}(\tr\mathsf{T} - n) - (\det\mathsf{T})^{1/2}$,
\begin{align}\label{lowdimidentity}
\vol(M, \gb_{\psi_{0}}) - \vol(M, \gb_{\psi}) = \frac{1}{2}\int_{M}|d\log f|^{2}_{\gb_{\psi_{0}}}\,\nu + \int_{M}R\,\nu.
\end{align}
Let $\la_{1}, \dots, \la_{n}$ be the eigenvalues of the positive definite endomorphism $\mathsf{T}$. Then $R = \tfrac{1}{2}(\la_{1}^{1/2} - \la_{2}^{1/2})^{2}$ when $n = 2$, and $R = \tfrac{1}{2}(\la_{1}^{1/2} - 1)^{2}$ when $n = 1$. So, when $n \leq 2$, $R$ is nonnegative, which gives \eqref{lowdimcentroaffine}. Moreover, $R$ vanishes at a point if and only if $\mathsf{T}$ is a multiple of the identity at that point, when $n = 2$, and if and only if $\mathsf{T} = 1$ at that point, when $n = 1$. If equality holds in \eqref{lowdimcentroaffine}, then, since $R$ is continuous and nonnegative, $R$ vanishes at every point. If $n = 2$, then $\mathsf{T}$ is a multiple of the identity at every point, that is, the equiaffine Hessian $H = \nabla^{2}f$ is pointwise a multiple of $\gb_{\psi_{0}}$, and, since $f$ is $\Ga$-invariant, Lemma \ref{blaschkehessianlemma} shows that $f$ is constant. If $n = 1$, then $\mathsf{T} = 1$ at every point, so $\Delta f = \tr H = 0$, and $f$ is constant because $M$ is closed.
\end{proof}

In particular the inequality of Question \ref{centroaffineconjecture} holds for $n \leq 2$, so, by Theorem \ref{reductionproposition}, for convex $\mathbb{RP}^{2}$-structures on closed surfaces the affine sphere is extremal for $\mathcal{T}_{\be}$ for every $\be \notin (\tfrac{1}{2}, 1)$. Theorem \ref{lowdimcentroaffinetheorem} bounds the defect $\cV(M) - \vol(M, \gb_{\psi})$ of the centroaffine volume below by $0$, while \eqref{deficitform}, which is equivalent to Theorem \ref{suminequalityproposition}, bounds it above by one half of the deficit $\int_{M}(\ka^{1/4} - (\ka^{\ast})^{1/4}\mu(\psi)^{1/2})^{2}\psi^{-N}$ in the affine sphere equation. The upper bound holds in every dimension, and for $n \geq 3$ the lower bound is the content of Question \ref{centroaffineconjecture}.

For $n \geq 3$ the last integrand in \eqref{lowdimidentity} need not be nonnegative. If all eigenvalues of $\mathsf{T}$ equal $1 + a$ with $a > 0$, then $(\det \mathsf{T})^{1/2} = (1 + a)^{n/2} > 1 + \tfrac{n}{2}a$. So a proof for $n \geq 3$ must use that $\mathsf{T} - I = -f^{-1}H$ comes from a Hessian, together with the constraint $\mathsf{T} > 0$ globally on $M$.

\begin{remark}\label{energyremark}
By Lemma \ref{centroaffineformulalemma}, $\tr\mathsf{T} = \tr_{\gb_{\psi_{0}}}\gb_{\psi}$, so the integration by parts in the proof of Theorem \ref{lowdimcentroaffinetheorem} gives, in every dimension,
\begin{align}\label{energyidentity}
\frac{1}{2}\int_{M}\tr_{\gb_{\psi_{0}}}(\gb_{f\psi_{0}})\,\nu = \frac{n}{2}\cV(M) - \frac{1}{2}\int_{M}|d\log f|^{2}_{\gb_{\psi_{0}}}\,\nu.
\end{align}
The left-hand side is the Dirichlet energy of the identity map from $(M, \gb_{\psi_{0}})$ to $(M, \gb_{f\psi_{0}})$, and the last integrand in \eqref{lowdimidentity} vanishes exactly at the points where this map is conformal, if $n = 2$. When $n = 2$ the area of the target is at most this energy, which is the inequality of Theorem \ref{lowdimcentroaffinetheorem}. For $n \geq 3$ the last integrand in \eqref{lowdimidentity} can be negative for a positive definite $\mathsf{T}$, as the example preceding this remark shows, so the pointwise comparison of volume and energy used when $n \leq 2$ does not extend. This does not settle the global question for the centroaffine metrics arising from cogauges in $\cAs(C, \Ga)$, which is Question \ref{centroaffineconjecture}.
\end{remark}

\begin{remark}\label{lagrangianremark}
For $\psi \in \cAs(C, \Ga)$ the graph $\Lambda_{\psi} = \{(x, G_{\psi}(x)): x \in C\}$ of $\fD\Phi$ is a conical Lagrangian submanifold of $\ste\times\std$ with its canonical symplectic form, and, by Lemma \ref{polarmaplemma}\ref{pml4}, it is also the graph of $\fD\Phi^{\ast}$ over $C^{\ast}$. The neutral metric $\mathsf{G}((u, \al), (v, \be)) = \al(v) + \be(u)$ restricts to $\Lambda_{\psi}$ as $2\fD^{2}\Phi$, which has signature $(1, n)$ by Lemma \ref{polarmaplemma}\ref{pml1}. By Lemma \ref{volumeproposition}, the part of $\Lambda_{\psi}$ lying over $K_{\psi}$ projects onto $K_{\psi}$ and onto $K^{\ast}_{\psi}$, so the two covolumes are the volumes of these projections modulo $\Ga$. The volume modulo $\Ga$ of this part of $\Lambda_{\psi}$ for the induced metric is $2^{N/2}N^{-1}\cS_{1/2}(\psi)$, by \eqref{sbetadefined}, and Lemma \ref{cauchyschwarzproposition} compares the product of the two projected volumes with the square of this volume. When $\mu(\psi)$ is constant, $\Lambda_{\psi}$ has mean curvature zero \cite[Theorem $4.3$]{Fox-schwarz}.

Symplectic reduction of the level set $\{(x, \zeta): \zeta(x) = 1\}$ by the flow $(x, \zeta) \mapsto (e^{t}x, e^{-t}\zeta)$ yields a para-K\"ahler space form \cite[Section $4.6$]{Fox-schwarz}, \cite{Hildebrand-crossratio, Hildebrand-parakahler}, of which $\Om\times\Om^{\ast}$ is an open subset on which $\Ga$ acts freely and properly. The metric induced on this space form by $\mathsf{G}$ has constant para-holomorphic sectional curvature $2$, so its Ricci tensor is $N$ times the metric, and consequently, with the metric induced by $-\mathsf{G}$, the quotient $(\Om\times\Om^{\ast})/\Ga$ is a para-K\"ahler Einstein manifold whose Ricci tensor is $-N$ times the metric. The value of the constant can be checked on the hyperboloid of Example \ref{hyperboloidexample}, for which the Lagrangian submanifold described next is totally geodesic and has induced metric, for $\mathsf{G}$, equal to $-2$ times a hyperbolic metric of curvature $-1$, hence of curvature $\tfrac{1}{2}$.

The intersection of $\Lambda_{\psi}$ with this level set is the graph of the conormal map of $\Sigma_{\psi}$. Its tangent vectors $(v, \fD^{2}\psi(v, \cdot))$, $v \in T\Sigma_{\psi}$, are horizontal, being $\mathsf{G}$-orthogonal to $(x, \zeta)$ and $(x, -\zeta)$, so it descends to a Lagrangian submanifold of $(\Om\times\Om^{\ast})/\Ga$, diffeomorphic to $M$, whose induced metric is $2\gb_{\psi}$, whose volume is $2^{n/2}\cS_{1/2}(\psi)$, and which has mean curvature zero exactly when $\Sigma_{\psi}$ is an affine sphere \cite{Hildebrand-crossratio, Rungi-equiaffine}.

Replacing $\psi$ by $(1 + tf)\psi$ with $f \in \cinf(M)$ deforms this Lagrangian submanifold by a Hamiltonian isotopy, with Hamiltonian potential $f$ at $t = 0$, through Lagrangian submanifolds that project diffeomorphically onto $M$. In these terms, Question \ref{centroaffineconjecture} asks whether the Lagrangian submanifold of mean curvature zero maximizes the volume among the Lagrangian submanifolds obtained from it in this way, a para-K\"ahler counterpart of the problem studied by Oh \cite{Oh-stability} for minimal Lagrangian submanifolds of K\"ahler--Einstein manifolds. Since the Einstein constant $-N$ is negative, the criterion of \cite[Theorem $1$]{Anciaux-Georgiou}, that the first positive eigenvalue of the Laplacian of the induced metric be at least the Einstein constant, is satisfied, so this Lagrangian submanifold is Hamiltonian stable, its second variation being negative for every nonconstant Hamiltonian potential, and the second variation formula of \cite{Anciaux-Georgiou} is $2^{n/2}$ times \eqref{centroaffinesecond}.

By Lemma \ref{centroaffineformulalemma}, when $\mu(\psi_{0})$ is constant the induced metric of the Lagrangian submanifold corresponding to $f\psi_{0}$ is $2\gb_{\psi_{0}}(\mathsf{T}\,\cdot\,, \,\cdot\,)$, and Theorem \ref{lowdimcentroaffinetheorem} shows that, when $n \leq 2$, the Lagrangian submanifold of mean curvature zero has the largest volume among all the Lagrangian submanifolds that arise in this way from cogauges in $\cAs(C, \Ga)$, and not only among those near it. This does not show that it has the largest volume among all the Lagrangian submanifolds Hamiltonian isotopic to it.
\end{remark}

\section{Bounded convex domains and the relative setting}\label{relativesection}
The linear inequality \eqref{suminequality} and the results of Section \ref{familysection} have relative versions on bounded domains, in which the infinite volumes in \eqref{suminequality} are replaced by the relative covolume and the relative polar covolume as defined in Definition \ref{relativesettingdefinition}. This section fixes that setting, proves the relative version of the linear inequality (Theorem \ref{relativetheorem}), and states the relative versions of the results of Section \ref{familysection}, with the modifications needed in their proofs (Remarks \ref{relativecriticalremark} and \ref{relativelowdimremark}). 

In this section the standing assumption that $\Ga$ divides $\Om$ is suspended. Instead, it is supposed that $\Om = \projp(C)$ is a bounded properly convex domain in an affine chart, and $\cAs(C)$ denotes the set $\cAs(C, \{1\})$ introduced after Definition \ref{asymptoticdefinition}.
The notation \eqref{phiqdefined} and the function $\mu(\psi)$ of \eqref{mudefined} are defined for $\psi \in \cAs(C)$ as before, and integrals over $\Om$ of compactly supported functions homogeneous of degree $-N$ are defined as in Section \ref{settingsection}. 

\begin{definition}\label{chengyaudefinition}
A \emph{Cheng--Yau solution} on $C$ is a $\psi_{0} \in \cAs(C)$ such that $\mu(\psi_{0})$ is constant and $\psi_{0}(x) \to 0$ as $x \to \pr C\setminus\{0\}$.
\end{definition}
By Lemma \ref{asymptoticlemma}, the level sets of a Cheng--Yau solution are the hyperbolic affine spheres centered at the origin and asymptotic to $\pr C$. By theorems of Cheng and Yau \cite{Cheng-Yau-dirichlet, Cheng-Yau-affinespheres}, a Cheng--Yau solution exists, and by Lemma \ref{generaluniquenesslemma} it is unique up to a positive constant factor. 

Choose linear coordinates $(z, t)$ on $\ste$ as in Lemma \ref{chartlemma} with $t \in C^{\ast}$, so that the slice $\Om_{0} = \{z: (z, 1) \in C\}$ is bounded and $z \mapsto [(z, 1)]$ identifies it with $\Om$. For a Cheng--Yau solution $\psi_{0}$ with $\mu(\psi_{0}) = \ka$, let $p(z) = \psi_{0}(z, 1)$. By Lemma \ref{chartlemma}, $p$ is a positive concave function on $\Om_{0}$ solving $\det(-\fD^{2}p) = \ka p^{-n-2}$, and it tends to $0$ at $\pr\Om_{0}$ because $\psi_{0}$ tends to $0$ at $\pr C\setminus\{0\}$. In the form used by Cheng and Yau, $u = -p$ is a negative convex solution of $\det\fD^{2}u = \ka(-u)^{-n-2}$ on $\Om_{0}$ vanishing on $\pr\Om_{0}$. The uniqueness of such a solution is what is proved in Lemma \ref{generaluniquenesslemma}; see also \cite{Cheng-Yau-dirichlet}.

The properties of the Cheng--Yau solution are collected in \cite[Theorem $4.1$]{Fox-schwarz}, where, in the notation used there, $\psi_{0} = (n+1)^{1/2}e^{-F/(n+1)}$, with $F$ the canonical potential of $C$, satisfies $\mu(\psi_{0}) = 1$. This uniqueness answers, for solutions homogeneous of degree $2$, the second question at the end of \cite[Section $5.8$]{Fox-schwarz}, on the uniqueness of a negative solution $u$ of $\det\fD^{2}u = -1$ with Hessian of signature $(n, 1)$ vanishing on $\pr C$. For such $u$, the function $\psi = (-2u)^{1/2}$ is a Cheng--Yau solution, because $\fD^{2}u(\eul, \cdot) = -\psi\,d\psi$, so the signature forces $\fD^{2}\psi$ to be negative definite on $\ker d\psi$.

The boundary condition cannot be omitted from Definition \ref{chengyaudefinition}, for the reason explained after Definition \ref{asymptoticdefinition}. The restriction to $C$ of the Cheng--Yau solution of a larger cone has constant $\mu$ but does not tend to $0$ at $\pr C\setminus\{0\}$. In the $\Ga$-invariant setting the boundary condition is automatic, by Lemma \ref{harnacklemma}, which is why Lemma \ref{uniquenessproposition} requires no boundary condition.

\begin{definition}\label{relativesettingdefinition}
Let $\psi_{0}$ be a Cheng--Yau solution on $C$ and let $\ka = \mu(\psi_{0})$. Let $\cAs_{c}(\psi_{0})$ be the set of $\psi \in \cAs(C)$ that coincide with $\psi_{0}$ outside the cone over a compact subset of $\Om$. For $\psi \in \cAs_{c}(\psi_{0})$ and $\be \in \rea$, the \emph{relative functionals} are
\begin{align}\label{relativevolumes}
\begin{aligned}
&\Srel{\be}(\psi) = \int_{\Om}\left(\mu(\psi)^{\be}\psi^{-N} - \ka^{\be}\psi_{0}^{-N}\right),&&\Vrel(\psi) = \tfrac{1}{N}\Srel{0}(\psi),\\
&\Vrelpolar(\psi) = \tfrac{1}{N}\Srel{1}(\psi), &&\cVrel(\psi) = \Srel{1/2}(\psi);
\end{aligned}
\end{align}
$\Vrel$, $\Vrelpolar$, and $\cVrel$ are the \emph{relative covolume}, the \emph{relative polar covolume}, and the \emph{relative centroaffine volume}. 

The data $(\Om, \psi_{0}, \ka, \cAs_{c}(\psi_{0}))$ with these functionals constitute the \emph{relative setting}.
\end{definition}

The covolume product is not finite in this setting. By Lemma \ref{centroaffinemetriclemma}\ref{cml1}, $|K_{\psi_{0}}|$ is a constant multiple of the volume of the equiaffine metric of $\Om$. This metric is bi-Lipschitz to the Hilbert metric of $\Om$ \cite[Proposition $3.4$]{Benoist-Hulin-finite} (see also \cite[Section $4.1$]{Benoist-Hulin}), and, as for every properly convex domain, $\Om$ has infinite volume for the Hilbert metric. Hence $|K_{\psi_{0}}| = \infty$, and therefore $|K_{\psi}| = \infty$ for every $\psi \in \cAs_{c}(\psi_{0})$. The relative functionals replace these infinite covolumes.

The set $\cAs_{c}(\psi_{0})$ is convex, because the proof of Lemma \ref{coneopenlemma} that $\cAs(C, \Ga)$ is convex does not use $\Ga$, and a convex combination of functions that coincide with $\psi_{0}$ outside the cones over compact sets $E_{1}$ and $E_{2}$ coincides with $\psi_{0}$ outside the cone over $E_{1} \cup E_{2}$. The integrands in \eqref{relativevolumes} have compact support, so the relative functionals are finite. By Lemma \ref{centroaffinemetriclemma}\ref{cml1}, $\cVrel(\psi) = \int_{\Om}(\vol_{\gb_{\psi}} - \vol_{\gb_{\psi_{0}}})$, and, by \eqref{conevolumelocal} with $f = 1$, applied to $\psi$ and to $\psi_{0}$ with $E$ a compact set of directions outside which $\psi = \psi_{0}$, $\Vrel(\psi) = |K_{\psi}\setminus K_{\psi_{0}}| - |K_{\psi_{0}}\setminus K_{\psi}|$. Similarly,
\begin{align}\label{relativepolardescription}
\Vrelpolar(\psi) = |G_{\psi}(K_{\psi})\setminus G_{\psi_{0}}(K_{\psi_{0}})| - |G_{\psi_{0}}(K_{\psi_{0}})\setminus G_{\psi}(K_{\psi})|,
\end{align}
where $G_{\psi}(K_{\psi})$ plays the role played by the polar set $K^{\ast}_{\psi}$ in Lemma \ref{volumeproposition}. Indeed, the proofs of Lemma \ref{polarmaplemma}\ref{pml1}--\ref{pml2} and of the injectivity of $G_{\psi}$ in Lemma \ref{polarmaplemma}\ref{pml3} do not use $\Ga$, so $G_{\psi}$ and $G_{\psi_{0}}$ are injective local diffeomorphisms of $C$ into $\std$, with Jacobian determinants of absolute values $\mu(\psi)$ and $\ka$.

Let $R_{E}$ be the cone over a compact set $E$ of directions outside which $\psi = \psi_{0}$. By \eqref{conevolumelocal} with $f = \mu(\psi)$, $|G_{\psi}(K_{\psi}\cap R_{E})| = \int_{K_{\psi}\cap R_{E}}\mu(\psi)\,dx = \tfrac{1}{N}\int_{E}\mu(\psi)\psi^{-N}$, and likewise for $\psi_{0}$, so $\Vrelpolar(\psi) = |G_{\psi}(K_{\psi}\cap R_{E})| - |G_{\psi_{0}}(K_{\psi_{0}}\cap R_{E})|$. On the open set $C\setminus R_{E}$ the maps $G_{\psi}$ and $G_{\psi_{0}}$ coincide, as do the sets $K_{\psi}$ and $K_{\psi_{0}}$, so $G_{\psi}(K_{\psi})$ and $G_{\psi_{0}}(K_{\psi_{0}})$ are the unions of the common set $G_{\psi}(K_{\psi}\setminus R_{E})$ with $G_{\psi}(K_{\psi}\cap R_{E})$ and with $G_{\psi_{0}}(K_{\psi_{0}}\cap R_{E})$, respectively, and both unions are disjoint, because $G_{\psi}$ and $G_{\psi_{0}}$ are injective. This gives \eqref{relativepolardescription}.

In the combination $\ka^{1/2}\Vrel + (\ka^{\ast})^{1/2}\Vrelpolar$ of Theorem \ref{relativetheorem} below, the two subtracted densities, $\ka^{1/2}\psi_{0}^{-N}$ and $(\ka^{\ast})^{1/2}\ka\psi_{0}^{-N}$, both equal the Riemannian density of $\gb_{\psi_{0}}$, by Lemma \ref{centroaffinemetriclemma}\ref{cml1}, which is why the covolumes are weighted in this way.

\begin{lemma}\label{relativelemma}
In the relative setting of Definition \ref{relativesettingdefinition}, let $\psi \in \cAs_{c}(\psi_{0})$ and $\mu = \mu(\psi)$.
\begin{enumerate}[label=$(\arabic*)$, ref=(\arabic*)]
\item\label{rl1} For smooth functions $f, f_{1}$ on $\Om$ such that $f_{1}\,df$ has compact support, the second identity of \eqref{ibp} holds with $\int_{M}$ replaced by $\int_{\Om}$. In particular, taking $f_{1} = 1$, the first identity of \eqref{ibp} holds with $\int_{M}$ replaced by $\int_{\Om}$ for smooth $f$ such that $df$ has compact support.
\item\label{rl2} For every smooth function $f$ on $\Om$, $\Lop_{\psi_{0}}f = -\Delta f$, where $\Delta$ is the Laplacian of $\gb_{\psi_{0}}$.
\item\label{rl3} For a compactly supported smooth $f$ on $\Om$, $\psi_{t} = (1 + tf)\psi$ lies in $\cAs_{c}(\psi_{0})$ for $|t|$ small, and
\begin{align}\label{relativefirstvariation}
&\tfrac{d}{dt}_{|t=0}\Vrel(\psi_{t}) = -\int_{\Om}f\psi^{-N},& &\tfrac{d}{dt}_{|t=0}\Vrelpolar(\psi_{t}) = \int_{\Om}f\mu\psi^{-N},&
\end{align}
and
\begin{align}\label{relativesecondcovolume}
\tfrac{d^{2}}{dt^{2}}_{|t=0}\Vrelpolar(\psi_{t}) = \int_{\Om}\mu\psi^{-N}\left((N-1)f^{2} + |df|^{2}_{\gb_{\psi}}\right).
\end{align}
In particular the functional $\Vrelpolar$ is strictly convex on $\cAs_{c}(\psi_{0})$.
\end{enumerate}
\end{lemma}

\begin{proof}
The vector field $Y$ in the proof of Lemma \ref{ibplemma} has compact support when $f_{1}\,df$ does, so the Stokes theorem gives $\int_{\Om}\operatorname{div}Y = 0$ in place of Lemma \ref{divergencelemma}, and the rest of the proof of Lemma \ref{ibplemma} applies verbatim to the second identity of \eqref{ibp}, of which the first is the case $f_{1} = 1$. This proves \ref{rl1}. Next \ref{rl2} is proved. Let $f$ be a smooth function on $\Om$ and let $f_{1}$ be a compactly supported smooth function on $\Om$, so that $f_{1}\,df$ has compact support. Since $\psi_{0}^{-N}$ is a constant multiple of $\vol_{\gb_{\psi_{0}}}$, by Lemma \ref{centroaffinemetriclemma}\ref{cml1}, \ref{rl1} applied with $\psi = \psi_{0}$ and integration by parts on $(\Om, \gb_{\psi_{0}})$, which is justified because $f_{1}$ has compact support, give
\begin{align}\label{relativelaplacian}
\int_{\Om}f_{1}\Lop_{\psi_{0}}f\,\psi_{0}^{-N} = \int_{\Om}\lb df_{1}, df\ra_{\gb_{\psi_{0}}}\psi_{0}^{-N} = -\int_{\Om}f_{1}\Delta f\,\psi_{0}^{-N}.
\end{align}
Since $f_{1}$ is arbitrary, $\Lop_{\psi_{0}}f = -\Delta f$. Finally \ref{rl3} is proved. Let $f$ be a compactly supported smooth function on $\Om$. The function $\psi_{t} = (1 + tf)\psi$ is smooth and homogeneous of degree $1$, and it coincides with $\psi$, hence with $\psi_{0}$ outside the cone over a compact subset of $\Om$, outside the cone over the support of $f$. As in the proof of Lemma \ref{coneopenlemma}, the positivity of $\psi_{t}$ and the negative definiteness of $\fD^{2}\psi_{t}$ transversally to $\eul$ are open conditions on the $2$-jet of $\psi_{t}$ that, by homogeneity, need only be checked on a compact subset of $C$ meeting every ray of the cone over the support of $f$, so they hold for $|t|$ small, and $\psi_{t} \in \cAs_{c}(\psi_{0})$ for such $t$. The identity \eqref{deltamu} is pointwise, and the differences of integrands defining $\Vrel$ and $\Vrelpolar$ along $\psi_{t}$ have compact support, so the proofs of Lemmas \ref{firstvariationproposition} and \ref{covolumeconvexitylemma}, using \ref{rl1} in place of Lemma \ref{ibplemma}, give \ref{rl3}.
\end{proof}

\begin{theorem}\label{relativetheorem}
Let $\Om = \projp(C)$ be a bounded properly convex domain, let $\psi_{0}$ be a Cheng--Yau solution, let $\ka$, $\cAs_{c}(\psi_{0})$, $\Vrel$, and $\Vrelpolar$ be as in Definition \ref{relativesettingdefinition}, and let $\ka^{\ast} = \ka^{-1}$. For every $\psi \in \cAs_{c}(\psi_{0})$,
\begin{align}\label{relativeinequality}
\ka^{1/2}\Vrel(\psi) + (\ka^{\ast})^{1/2}\Vrelpolar(\psi) \geq 0,
\end{align}
with equality if and only if $\psi = \psi_{0}$. The functional $\ka^{1/2}\Vrel + (\ka^{\ast})^{1/2}\Vrelpolar$ is strictly convex along segments in $\cAs_{c}(\psi_{0})$, $\psi_{0}$ is its unique minimizer and its unique critical point on $\cAs_{c}(\psi_{0})$, and the second variation at $\psi_{0}$ in the direction $(1 + tf)\psi_{0}$, with $f$ compactly supported, is
\begin{align}\label{relativesecondvariation}
\tfrac{d^{2}}{dt^{2}}_{|t=0}\left(\ka^{1/2}\Vrel + (\ka^{\ast})^{1/2}\Vrelpolar\right)((1+tf)\psi_{0}) = \int_{\Om}\left(2Nf^{2} + |df|^{2}_{\gb_{\psi_{0}}}\right)d\vol_{\gb_{\psi_{0}}}.
\end{align}
\end{theorem}

\begin{proof}
The proof of Theorem \ref{suminequalityproposition} applies verbatim, with $\int_{M}$ replaced by $\int_{\Om}$, $|K_{\psi}/\Ga| - |K_{\psi_{0}}/\Ga|$ replaced by $\Vrel(\psi)$, and $|K^{\ast}_{\psi}/\Ga^{\ast}| - |K^{\ast}_{\psi_{0}}/\Ga^{\ast}|$ replaced by $\Vrelpolar(\psi)$, using \eqref{relativefirstvariation} and \eqref{relativesecondcovolume} of Lemma \ref{relativelemma} in place of \eqref{firstvariation} and Lemma \ref{covolumeconvexitylemma}. In that proof the ratio $\psi/\psi_{0}$ plays the role of the function there called $f$, and here $\psi/\psi_{0} - 1$ has compact support because $\psi \in \cAs_{c}(\psi_{0})$, so all the integrands do. The critical point is unique because a strictly convex function on a convex set has at most one critical point. Strict convexity of the functional follows from that of $\Vrelpolar$ and of $\Vrel$, the latter because $x \mapsto x^{-N}$ is strictly convex and, along a segment $\psi_{t} = (1-t)\psi_{a} + t\psi_{b}$ in $\cAs_{c}(\psi_{0})$, $\psi_{t}^{-N} = ((1-t)\psi_{a} + t\psi_{b})^{-N}$. For the second variation, with $f$ as in the statement, $\tfrac{d^{2}}{dt^{2}}_{|t=0}\Vrel((1+tf)\psi_{0}) = (N+1)\int_{\Om}f^{2}\psi_{0}^{-N}$, and combining this with \eqref{relativesecondcovolume} and $\ka^{1/2}\psi_{0}^{-N} = \vol_{\gb_{\psi_{0}}}$ gives \eqref{relativesecondvariation}.
\end{proof}

Theorem \ref{relativetheorem} provides a variational characterization of the Cheng--Yau affine sphere of a bounded convex domain. Among all smooth, locally strongly convex hypersurfaces asymptotic to the cone that coincide with the affine sphere near infinity, the affine sphere minimizes the relative, polarity symmetric sum of the covolume and the polar covolume, weighted by $\ka^{1/2}$ and $(\ka^{\ast})^{1/2}$. The proofs of Lemma \ref{relativelemma} and Theorem \ref{relativetheorem}, and those of the statements of Remarks \ref{relativecriticalremark} and \ref{relativelowdimremark} below, use only that $\mu(\psi_{0})$ is constant and that the perturbations have compact support. They do not use the boundary condition of Definition \ref{chengyaudefinition}, which serves to single out $\psi_{0}$. Hence these results hold, with $\cAs_{c}(\psi_{0})$ and the relative functionals defined in the same way, for every $\psi_{0} \in \cAs(C)$ such that $\mu(\psi_{0})$ is constant, for example for the restriction to $C$ of the Cheng--Yau solution of a larger cone. In this sense the minimizing property is shared by every hyperbolic affine sphere centered at the origin, whether or not it is asymptotic to $\pr C$.

For the relative centroaffine volume $\cVrel$ of \eqref{relativevolumes}, the pointwise identity of Remark \ref{deficitremark} gives
\begin{align}\label{relativedeficit}
N\left(\ka^{1/2}\Vrel(\psi) + (\ka^{\ast})^{1/2}\Vrelpolar(\psi)\right) = \int_{\Om}\left(\ka^{1/4} - (\ka^{\ast})^{1/4}\mu(\psi)^{1/2}\right)^{2}\psi^{-N} + 2\cVrel(\psi),
\end{align}
so that \eqref{relativeinequality} reads $-\cVrel(\psi) \leq \tfrac{1}{2}\int_{\Om}(\ka^{1/4} - (\ka^{\ast})^{1/4}\mu(\psi)^{1/2})^{2}\psi^{-N}$. Perturbations $\psi = f\psi_{0}$ with $f - 1$ decaying at $\pr\Om$, for which the integrals \eqref{relativevolumes} converge absolutely, require control of the boundary terms in the integrations by parts and are not treated here.

\begin{remark}\label{relativecriticalremark}
The results of Section \ref{familysection} have relative versions, obtained by the same computations with Lemma \ref{relativelemma}\ref{rl1} in place of Lemma \ref{ibplemma}. In the relative setting of Definition \ref{relativesettingdefinition}, for $\be \in \rea$, the functionals
\begin{align}\label{rbetadefined}
&\mathcal{R}_{\be} = \Srel{\be} + N(2\be - 1)\ka^{\be}\Vrel,& &\mathcal{R}'_{\be} = \Srel{\be} - (2\be - 1)\ka^{\be-1}\Srel{1},&
\end{align}
on $\cAs_{c}(\psi_{0})$ take the place of $\mathcal{T}_{\be}$ and $\mathcal{T}'_{\be}$. The Cheng--Yau solution $\psi_{0}$ is critical for both for every $\be$, and, for $\be$ in the range of Theorem \ref{criticalpointstheorem}\ref{cpt1}, respectively \ref{cpt2}, it is the only critical point of $\mathcal{R}_{\be}$, respectively $\mathcal{R}'_{\be}$, in $\cAs_{c}(\psi_{0})$. For $\psi \in \cAs_{c}(\psi_{0})$ and compactly supported smooth $f$ on $\Om$, the identity \eqref{sbetafirstvariation} holds with $\cS_{\be}$ replaced by $\Srel{\be}$ and $\int_{M}$ by $\int_{\Om}$, and with \eqref{relativefirstvariation} it shows that the derivatives at $t = 0$ of $\mathcal{R}_{\be}$, respectively $\mathcal{R}'_{\be}$, along $(1 + tf)\psi$ vanish for all such $f$ if and only if \eqref{eltbeta}, respectively \eqref{eltprimebeta}, holds with $c_{0} = \ka^{\be}$ and $c_{1} = \ka^{\be - 1}$. For $\mathcal{R}'_{\be}$ the proof of Theorem \ref{criticalpointstheorem}\ref{cpt2} applies, because $F - \ka^{\be - 1}$ has compact support, and shows that $\mu = \mu(\psi)$ is constant. For $\mathcal{R}_{\be}$ there is no duality reducing the claim to the case of $\mathcal{R}'_{\be}$, and a direct argument is needed. Multiplying \eqref{eltbeta} by $(F - \ka^{\be-1})\nu$, which has compact support, and integrating gives
\begin{align}\label{relativechebyshev}
\be\int_{\Om}|dF|^{2}_{\gb_{\psi}}\mu\nu = -N(2\be-1)\int_{\Om}(\mu^{\be} - \ka^{\be})(\mu^{\be-1} - \ka^{\be-1})\nu,
\end{align}
in which the integrand on the right-hand side is pointwise nonpositive for $0 < \be < 1$ and nonnegative for $\be < 0$ or $\be > 1$. Hence, for $\be$ in the range of Theorem \ref{criticalpointstheorem}\ref{cpt1} other than $\be = 1$, the two sides of \eqref{relativechebyshev} have opposite signs, so both vanish, and, as $\be \neq 0$, $dF = 0$ and $\mu$ is constant. For $\be = 1$, $F = 1$ and \eqref{eltbeta} gives $\mu = \ka$ directly. In both cases the constant $\mu$ equals $\ka$, its value outside the cone over a compact set, and the maximum principle argument of the proof of Lemma \ref{uniquenessproposition}, applied to $\psi/\psi_{0}$, which equals $1$ outside the cone over a compact set, gives $\psi = \psi_{0}$. The second variation of $\mathcal{R}_{\be}$ at $\psi_{0}$ along $(1 + tf)\psi_{0}$, with $f$ compactly supported, is obtained from the right-hand side of \eqref{tbetasecond} by replacing $\bar{f}$ by $0$ and $V^{-1}\int_{M}$ by $\ka^{\be}\int_{\Om}$. In particular the second variation of $\cVrel = \mathcal{R}_{1/2}$ at $\psi_{0}$ is $-\int_{\Om}(\tfrac{N}{2}|df|^{2} + \tfrac{1}{4}(\Delta f)^{2})\,d\vol_{\gb_{\psi_{0}}}$, which is negative definite on compactly supported $f \neq 0$.
\end{remark}

By Remark \ref{relativecriticalremark}, the second variation of the relative centroaffine volume $\cVrel$ at $\psi_{0}$ is negative definite on compactly supported variations. This motivates Question \ref{relativecentroaffinequestion}, the relative counterpart of Question \ref{centroaffineconjecture}.

\begin{question}\label{relativecentroaffinequestion}
Is the relative centroaffine volume maximized by the Cheng--Yau affine sphere among compactly supported perturbations? Precisely, for a bounded properly convex domain $\Om$ of dimension $n \geq 3$, with a Cheng--Yau solution $\psi_{0}$ (Definition \ref{chengyaudefinition}), is it true that every $\psi \in \cAs_{c}(\psi_{0})$ satisfies $\cVrel(\psi) \leq 0$, with equality only if $\psi = \psi_{0}$?
\end{question}

\begin{remark}\label{relativelowdimremark}
For $n \leq 2$, the proof of Theorem \ref{lowdimcentroaffinetheorem} applies in the relative setting. With $\psi_{0}$ the Cheng--Yau solution normalized by $\ka = 1$ and $\nu = \psi_{0}^{-N}$, every $\psi = f\psi_{0} \in \cAs_{c}(\psi_{0})$ satisfies
\begin{align}\label{relativelowdim}
\cVrel(\psi) \leq -\tfrac{1}{2}\int_{\Om}|d\log f|^{2}_{\gb_{\psi_{0}}}\,\nu,
\end{align}
with equality if and only if $\psi = \psi_{0}$. Here Lemma \ref{centroaffineformulalemma}, whose proof is pointwise, gives $\cVrel(f\psi_{0}) = \int_{\Om}((\det \mathsf{T})^{1/2} - 1)\nu$, and the integration by parts $\int_{\Om}f^{-1}\Delta f\,\nu = \int_{\Om}f^{-2}|df|^{2}_{\gb_{\psi_{0}}}\nu$ is justified because $df$ has compact support. In the equality case, Lemma \ref{blaschkehessianlemma} is used for $n = 2$ with the hypothesis that $f = 1$ outside a compact set, and for $n = 1$, $f$ is affine in an arc length parameter of $\gb_{\psi_{0}}$ and equal to $1$ near the ends of $\Om$. As in Section \ref{familysection}, \eqref{relativelowdim} together with \eqref{relativedeficit} and Theorem \ref{relativetheorem} bounds $-\cVrel(\psi)$ below by $0$ and above by half of the relative deficit, and Questions \ref{centroaffineconjecture} and \ref{relativecentroaffinequestion} are closely related versions of the same extremal question.
\end{remark}

\section{The minimal covolume product and the equiaffine volume}\label{invariantsection}
By Theorem A and Definitions \ref{maxmindefinition} and \ref{blaschkevolumedefinition},
\begin{align}\label{minvariant}
\mathfrak{m}(M) = \min_{L \in \cK(C, \Ga)}\cP(L) = \frac{\cV(M)^{2}}{N^{2}},
\end{align}
where the equiaffine volume $\cV(M)$ depends only on the flat projective structure of $M$ and satisfies $\cV(M^{\ast}) = \cV(M)$, as shown after Definition \ref{blaschkevolumedefinition}. This section draws consequences of this identity. 

\begin{example}
In the Fuchsian case of Example \ref{hyperboloidexample}, $\gb_{\psi_{0}}$ is the hyperbolic metric of curvature $-1$, and $\mathfrak{m}(M) = \vol_{\mathrm{hyp}}(M)^{2}/N^{2}$. 
\end{example}

\begin{example}\label{ex:orthantmin}
In the case of Example \ref{orthantexample}, in which $c_{\ell} = 1$ can be assumed because $\mathfrak{m}(M)$ does not depend on $\Psi$, identify $\Sigma_{\psi_{0}}$ with the hyperplane $\Pi = \{y \in \rea^{N}: \sum_{i=1}^{N}y_{i} = \tfrac{N}{2}\log\tfrac{1}{N}\}$ in the coordinates $y_{i} = \log \ell_{i}$, on which $\Ga$ acts by a lattice $\Lambda$ of translations parallel to $\Pi_{0} = \{y \in \rea^{N}: \sum_{i=1}^{N}y_{i} = 0\}$. Since $\gb_{\psi_{0}} = \tfrac{1}{N}\sum_{i=1}^{N}dy_{i}^{2}$ on $\Pi$, $\mathfrak{m}(M) = N^{-n-2}|\Pi_{0}/\Lambda|^{2}$, where $|\Pi_{0}/\Lambda|$ is the Euclidean volume of a fundamental domain for $\Lambda$ in $\Pi_{0}$.

For a simplicial cone with $N_{1} = 1$, the formula of Corollary \ref{productinvariantcorollary} agrees with the formula for $\mathfrak{m}(M)$ just given, because the generator of the additional factor of the lattice corresponds to $w = (N_{2}\tau, -\tau, \dots, -\tau)$, which is orthogonal to the hyperplane of $\Pi_{0}$ consisting of the $y$ with $y_{1} = 0$, and satisfies $|w|^{2} = N_{2}N\tau^{2}$.
\end{example}

Lemma \ref{lorentzconelemma} is an elementary rigidity statement related to the results of Vey and Benoist on groups dividing properly convex cones \cite{Benoist-survey, Vey}. For irreducible cones, its conclusion also follows from the theorem of Benoist \cite{Benoist-divisibles2} that, when $C$ is not symmetric, $\Ga \cap SL(\ste)$ is Zariski dense in $SL(\ste)$, and so preserves no nonzero quadratic form, together with the classification of symmetric cones. The direct proof is included because it uses only the characteristic function of $C$ and requires neither irreducibility nor that classification.

\begin{lemma}\label{lorentzconelemma}
Suppose that $n \geq 2$ and that $\Ga$ preserves a quadratic form $Q$ on $\ste$ of signature $(n, 1)$. Then $C$ is a connected component of $\{x \in \ste: Q(x) < 0\}$.
\end{lemma}

\begin{proof}
Let $B$ be the symmetric bilinear form such that $B(v, v) = Q(v)$. The function $Q/\rho_{C}^{2}$ is continuous on $C$, $\Ga$-invariant, and homogeneous of degree $0$, so it is identified with a continuous function on the compact manifold $M$, and there is a constant $c$ such that $|Q| \leq c\rho_{C}^{2}$ on $C$. Since $\rho_{C}(x) \to 0$ as $x \to \pr C\setminus\{0\}$, by Lemma \ref{vinberglemma}, and $Q$ is continuous, $Q$ vanishes on $\pr C$. Let $x \in C$ and let $P \subset \ste$ be a two-dimensional subspace containing $x$. As in the proof of Lemma \ref{vinberglemma}, $\bar{C}\cap P$ is the sector bounded by the rays spanned by two linearly independent vectors $y_{1}, y_{2} \in P\cap\pr C$, so $P$ contains nonzero vectors on which $Q$ vanishes, and $x = ay_{1} + by_{2}$ with $a, b > 0$, so that
\begin{align}\label{lorentzconeidentity}
Q(x) = 2ab\,B(y_{1}, y_{2}).
\end{align}
If $Q(x) = 0$, then, by \eqref{lorentzconeidentity}, $B$ vanishes identically on $P$, which is impossible, because a form of signature $(n, 1)$ has no totally isotropic subspace of dimension two. If $Q(x) > 0$, then the restriction of $B$ to the $B$-orthogonal complement of $x$ has signature $(n-1, 1)$, so, because $n \geq 2$, this complement contains a vector $v$ such that $Q(v) > 0$, and $Q$ is positive definite on the subspace $P$ spanned by $x$ and $v$, which contradicts that $P$ contains nonzero vectors on which $Q$ vanishes. Hence $Q < 0$ on $C$. Since $C$ is connected, it is contained in a connected component $C_{Q}$ of $\{x \in \ste: Q(x) < 0\}$. The set $C$ is open in $C_{Q}$, and it is closed in $C_{Q}$, because $Q$ vanishes on $\pr C$, so that $\bar{C}\cap C_{Q} = C$. Since $C_{Q}$ is connected, $C = C_{Q}$.
\end{proof}

\begin{theorem}\label{surfacetheorem}
Let $M = \Om/\Ga$ be a closed properly convex flat projective surface with $\chi(M) < 0$. Let $\psi_{0} \in \cAs(C, \Ga)$ be the unique element with $\mu(\psi_{0}) = 1$ (Corollary \ref{existencecorollary}), let $h = \gb_{\psi_{0}}$, let $\dift$ be the difference tensor of Lemma \ref{gausslemma}, and let $C_{h}(X, Y, Z) = h(\dift_{X}Y, Z)$ be the cubic form.
\begin{enumerate}[label=$(\arabic*)$, ref=(\arabic*)]
\item\label{st1} $3\mathfrak{m}(M)^{1/2} = \vol(M, h) = 2\pi|\chi(M)| + \tfrac{1}{2}\int_{M}|\dift|^{2}_{h}\,d\vol_{h}$.
\item\label{st2} $\mathfrak{m}(M) \geq \tfrac{4\pi^{2}}{9}\chi(M)^{2}$, with equality if and only if $\Ga$ preserves a Lorentzian quadratic form, that is, if and only if the projective structure is Fuchsian.
\item\label{st3} If $M$ is orientable, then $C_{h}$ is the real part of a cubic differential $U$ that is holomorphic for the conformal structure of $h$. If $|U|^{2/3}$ denotes the metric equal to $|u|^{2/3}|dz|^{2}$ when $U = u\,dz^{3}$ in a local complex coordinate $z = x + iy$, then $h \geq 2^{1/3}|U|^{2/3}$ pointwise, and
\begin{align}\label{cubicbound}
2^{1/3}\int_{M}|U|^{2/3} \leq 3\mathfrak{m}(M)^{1/2} \leq 2\pi|\chi(M)| + 2^{1/3}\int_{M}|U|^{2/3}.
\end{align}
\end{enumerate}
\end{theorem}
\begin{proof}
By \eqref{minvariant} with $N = 3$, $3\mathfrak{m}(M)^{1/2} = \vol(M, h)$. By Lemma \ref{gausslemma} with $n = 2$, $\mathrm{Scal}_{h} = -2 + |\dift|^{2}_{h}$, so the Gauss--Bonnet theorem gives $4\pi\chi(M) = \int_{M}\mathrm{Scal}_{h}\,d\vol_{h} = -2\vol(M, h) + \int_{M}|\dift|^{2}_{h}\,d\vol_{h}$ (compare \cite[Theorem $5.4$ and (5.15)]{Fox-pisa}). This proves \ref{st1}.

By \ref{st1}, $3\mathfrak{m}(M)^{1/2} \geq 2\pi|\chi(M)|$, with equality if and only if $\dift = 0$. If $\dift = 0$, then $\Sigma_{\psi_{0}}$ is a hyperbolic affine sphere with vanishing cubic form, hence a hyperboloid \cite[Chapter II, Theorem $4.5$]{Nomizu-Sasaki}, that is, a connected component of $\{x \in \ste: Q(x) = -1\}$ for a quadratic form $Q$ of signature $(2, 1)$. For $\gamma \in \Ga$, the quadratic forms $Q\circ\gamma$ and $Q$ agree on $\Sigma_{\psi_{0}}$, because $\gamma$ preserves $\Sigma_{\psi_{0}}$, hence agree on $C$, by homogeneity, and hence on $\ste$, because $C$ is open. So $\Ga$ preserves the Lorentzian quadratic form $Q$. This proves \ref{st2}.

Conversely, suppose that $\Ga$ preserves a quadratic form $Q$ of signature $(n, 1)$, which, after multiplication by a positive constant, can be supposed to have volume form $\pm\Psi$. By Lemma \ref{lorentzconelemma}, $C$ is a connected component of $\{x \in \ste: Q(x) < 0\}$. The computation of Example \ref{hyperboloidexample} uses only that $\Ga$ preserves $Q$ and $C$, and it gives that $(-Q)^{1/2} \in \cAs(C, \Ga)$, that $\mu((-Q)^{1/2}) = 1$, and that the centroaffine metric of the hyperboloid $\{x \in C: Q(x) = -1\}$ is the hyperbolic metric of curvature $-1$. By Lemma \ref{uniquenessproposition}, $\psi_{0} = (-Q)^{1/2}$, so $h$ is the hyperbolic metric of curvature $-1$, and, by the Gauss--Bonnet theorem, $\vol(M, h) = 2\pi|\chi(M)|$, so that equality holds, by \ref{st1}.

Finally \ref{st3} is proved. Since $S = -\mathrm{id}$, the Gauss equations of $\nabla = D + \dift$ and of its conjugate $D - \dift$ coincide (proof of Lemma \ref{gausslemma}). Their difference is $2((D_{X}\dift)_{Y} - (D_{Y}\dift)_{X})$, so $DC_{h}$ is completely symmetric, that is, $C_{h}$ is a Codazzi tensor. On a surface, a trace-free symmetric cubic form is locally the real part of $u\,dz^{3}$ in a conformal coordinate $z$, and it is a Codazzi tensor if and only if $u$ is holomorphic \cite[Lemma $3.5$(1)]{Fox-pisa}. If $h = e^{v}|dz|^{2}$, then $|\dift|^{2}_{h} = |C_{h}|^{2}_{h} = |\mathrm{Re}(u\,dz^{3})|^{2}_{h} = 4|u|^{2}e^{-3v}$.

By Lemma \ref{gausslemma}, $\mathrm{Scal}_{h} = -n(n-1) + |\dift|^{2}_{h}$. Because $h$ is the lift of a Riemannian metric on the closed manifold $M$, it is complete on $\Sigma_{\psi_{0}}$, so, by a theorem of Calabi \cite{Calabi-completeaffine}, $\mathrm{Ric}_{h} \leq 0$, and hence $\mathrm{Scal}_{h} \leq 0$. Together these give $|\dift|^{2}_{h} \leq n(n-1)$ for every closed properly convex flat projective manifold. For $n = 2$ this is $4|u|^{2}e^{-3v} \leq 2$, that is, $e^{v} \geq 2^{1/3}|u|^{2/3}$, which is the pointwise bound $h \geq 2^{1/3}|U|^{2/3}$. Comparing the area densities of the two metrics in this bound and integrating over $M$ gives the lower bound in \eqref{cubicbound}. For the upper bound, the density $\tfrac{1}{2}|\dift|^{2}_{h}\,d\vol_{h}$ on $M$ equals $2|u|^{2}e^{-2v}\,dx\,dy$ in a local complex coordinate, and, by the pointwise bound, $2|u|^{2}e^{-2v} \leq 2\cdot 2^{-2/3}|u|^{2/3} = 2^{1/3}|u|^{2/3}$. It follows that $\tfrac{1}{2}|\dift|^{2}_{h}\,d\vol_{h}$ is at most $2^{1/3}$ times the area density of $|U|^{2/3}$, and integrating over $M$ and using \ref{st1} gives
\begin{align}
3\mathfrak{m}(M)^{1/2} - 2\pi|\chi(M)| = \frac{1}{2}\int_{M}|\dift|^{2}_{h}\,d\vol_{h} \leq 2^{1/3}\int_{M}|U|^{2/3}.
\end{align}
\end{proof}

The pointwise bound $|\dift|^{2}_{h} \leq 2$, equivalently $h \geq 2^{1/3}|U|^{2/3}$, follows, as above, from the theorem of Calabi, and in the form $h \geq 2^{1/3}|U|^{2/3}$ it appears in \cite{Dumas-Wolf} and \cite[(2.5)]{Nie-entropy}. A direct proof on compact surfaces, in a setting more general than that of affine spheres and in terms of a tensor equal to $2C_{h}$ up to sign, is \cite[(7.15)]{Fox-pisa}. The upper bound in \eqref{cubicbound} is due to Nie \cite[Section $2$, Lemma]{Nie-entropy}, and both bounds appear in \cite[Lemma $3.4$]{Ouyang-Tamburelli}, whose normalization of the cubic differential agrees with the one here, as their form $1 = -\kappa_{h} + 2|q|^{2}h^{-3}$ of the equation of Wang \cite{Wang-examples} shows.

\begin{remark}\label{surfaceremark}
By Labourie \cite{Labourie-cubic} and Loftin \cite{Loftin-affinespheres}, convex $\mathbb{RP}^{2}$-structures on a closed oriented surface $S$ of genus at least $2$ correspond to pairs $(J, U)$ consisting of a conformal structure $J$ and a $J$-holomorphic cubic differential $U$, normalized here as in Theorem \ref{surfacetheorem}\ref{st3}. Accordingly, $\mathfrak{m}$ is regarded as a function of $(J, U)$.

Along each ray $s \mapsto (J, sU)$ with $U \neq 0$, the equiaffine metric is pointwise nondecreasing in $s$, while $s^{-2/3}$ times its area is nonincreasing \cite[Theorem $9.9$ and Corollary $9.10$]{Fox-pisa}. Strict monotonicity of the metric, in the language of cyclic Higgs bundles and for all ranks, is \cite[Proposition $4.1$]{Dai-Li-cyclic}. Since the area of the equiaffine metric is $3\mathfrak{m}^{1/2}$ by Theorem \ref{surfacetheorem}\ref{st1}, $\mathfrak{m}(J, sU)$ is nondecreasing and $s^{-4/3}\mathfrak{m}(J, sU)$ is nonincreasing in $s$. By \eqref{cubicbound},
\begin{align}\label{raybound}
2^{1/3}s^{2/3}\int_{S}|U|^{2/3} \leq 3\mathfrak{m}(J, sU)^{1/2} \leq 2\pi|\chi(S)| + 2^{1/3}s^{2/3}\int_{S}|U|^{2/3},
\end{align}
so $s^{-4/3}\mathfrak{m}(J, sU)$ converges to $2^{2/3}\big(\int_{S}|U|^{2/3}\big)^{2}/9$ as $s \to \infty$, which identifies the limit whose existence is shown in \cite[Corollary $9.10$]{Fox-pisa}. Consistently, $s^{-2/3}$ times the equiaffine metric converges to $2^{1/3}|U|^{2/3}$ uniformly on compact subsets of the complement of the zeros of $U$ \cite{Loftin-flat}.

By the lower bound in \eqref{raybound}, for fixed $J$ the sublevel sets of $\mathfrak{m}$ are bounded in the space of cubic differentials. On the other hand, $\mathfrak{m}$ is constant, equal to $4\pi^{2}\chi(S)^{2}/9$, on the Fuchsian locus $U = 0$, a copy of Teichm\"uller space. Hence $\mathfrak{m}$ is not a proper function on the space of convex $\mathbb{RP}^{2}$-structures on $S$, nor on the moduli space, its quotient by the mapping class group.

Two other functionals of convex $\mathbb{RP}^{2}$-structures are affine functions of $\mathfrak{m}^{1/2}$. The first is Mettler's functional $\Efun(g) = \int_{S}|(\nabla - D^{g})_{0}|^{2}_{g}\,d\vol_{g}$ on Riemannian metrics $g$, where $\nabla$ represents the projective structure, $D^{g}$ is the Levi-Civita connection of $g$, and $(\cdot)_{0}$ denotes the trace-free part. The equiaffine metric $h$ is a critical point of $\Efun$ \cite[Theorem A]{Mettler-prescribing}, and, since $(\nabla - D^{h})_{0} = \dift$, Theorem \ref{surfacetheorem}\ref{st1} gives $\Efun(h) = 6\mathfrak{m}^{1/2} - 4\pi|\chi(S)|$. The second is the shifted Blaschke area $3\mathfrak{m}^{1/2} - 2\pi|\chi(S)|$, the weight used by Cao \cite{Cao-goldman}, whose results imply that the region of the moduli space where $\mathfrak{m} \leq c$ has finite Goldman volume, polynomial in $c^{1/2}$.
\end{remark}

Theorem \ref{hyperbolicvolumetheorem} extends \ref{st2} of Theorem \ref{surfacetheorem} to convex projective structures on closed hyperbolic manifolds of dimension at least $3$. In its proof the role of the Gauss--Bonnet theorem in Theorem \ref{surfacetheorem} is played by the main theorem of Besson, Courtois, and Gallot \cite{BCG}. Compare \cite{ABC}, where a version of that theorem for Finsler metrics is applied to the Hilbert metric, with the equiaffine metric serving as an auxiliary Riemannian metric, to bound the Hilbert volume of a closed strictly convex projective manifold of dimension at least $3$ that admits a hyperbolic structure from below by the hyperbolic volume times a constant depending only on the dimension. Here the theorem of \cite{BCG} is applied directly to the equiaffine metric, which is Riemannian, strict convexity is not needed, and the resulting bound for the equiaffine volume is sharp.

\begin{theorem}\label{hyperbolicvolumetheorem}
Let $n \geq 3$, let $M_{0}$ be a closed hyperbolic manifold of dimension $n$, with hyperbolic volume $\vol_{\mathrm{hyp}}(M_{0})$, and let $M = \Om/\Ga$ be a closed properly convex flat projective manifold homeomorphic to $M_{0}$. Then
\begin{align}\label{hyperbolicvolumebound}
\mathfrak{m}(M) = \frac{\cV(M)^{2}}{N^{2}} \geq \frac{\vol_{\mathrm{hyp}}(M_{0})^{2}}{N^{2}},
\end{align}
with equality if and only if $\Ga$ preserves a Lorentzian quadratic form, that is, if and only if the flat projective structure of $M$ is that of a hyperbolic metric.
\end{theorem}
\begin{proof}
Let $\psi_{0} \in \cAs(C, \Ga)$ be the unique element with $\mu(\psi_{0}) = 1$ (Corollary \ref{existencecorollary}), and let $h = \gb_{\psi_{0}}$, so that $\cV(M) = \vol(M, h)$. If $M$ is not orientable, then $M$ and $M_{0}$ can be replaced by their orientation double covers which are homeomorphic because a homeomorphism preserves the first Stiefel-Whitney class so lifts to the orientation double covers. Making this replacement multiplies both $\vol(M, h)$ and $\vol_{\mathrm{hyp}}(M_{0})$ by $2$. This replaces $\Ga$ by a subgroup of index $2$ but does not change $\Om$, $\psi_{0}$, or $\Sigma_{\psi_{0}}$, which remain invariant under the original group $\Ga$. So it can and will be supposed that $M$ and $M_{0}$ are oriented.

The \emph{volume entropy} of a Riemannian metric $g$ on a closed manifold is defined by $\operatorname{ent}(g) = \lim_{r \to \infty}r^{-1}\log\vol(B(x, r))$, where $B(x, r)$ is the ball in the universal cover having radius $r$ with respect to the lifted metric and centered on the point $x$; the limit exists and does not depend on $x$. By Lemma \ref{gausslemma}\ref{gl2}, $\mathrm{Ric}_{h} \geq -(n-1)h$, so, by the Bishop volume comparison theorem, the volume of a ball of radius $r$ in the universal cover $(\Sigma_{\psi_{0}}, h)$ of $(M, h)$ is at most that of a ball of radius $r$ in the hyperbolic space of dimension $n$ and curvature $-1$. Hence $\operatorname{ent}(h) \leq n-1$, and $n - 1$ is the volume entropy of a hyperbolic metric $g_{0}$ on $M_{0}$. Because a homeomorphism of closed oriented manifolds has degree $\pm 1$, by \cite[Théorème Principal]{BCG},
\begin{align}\label{bcginequality}
(n-1)^{n}\vol(M, h) \geq \operatorname{ent}(h)^{n}\vol(M, h) \geq \operatorname{ent}(g_{0})^{n}\vol(M_{0}, g_{0}) = (n-1)^{n}\vol_{\mathrm{hyp}}(M_{0}),
\end{align}
and equality in the second inequality of \eqref{bcginequality} holds only if $(M, h)$ is homothetic to $(M_{0}, g_{0})$. This shows \eqref{hyperbolicvolumebound}.

If equality holds in \eqref{hyperbolicvolumebound}, then both inequalities in \eqref{bcginequality} are equalities, so $h$ is a constant multiple of a hyperbolic metric and $\operatorname{ent}(h) = n-1$. Since multiplying a metric by a constant $c^{2}$ multiplies its volume entropy by $c^{-1}$, $h$ is hyperbolic, so $\mathrm{Ric}_{h} = -(n-1)h$. By the equality case of Lemma \ref{gausslemma}\ref{gl2}, $\dift = 0$, and, as in the proof of \ref{st2} of Theorem \ref{surfacetheorem}, $\Sigma_{\psi_{0}}$ is a hyperboloid. Since $\Sigma_{\psi_{0}}$ is invariant under the original group $\Ga$, the argument there shows that $\Ga$ preserves a Lorentzian quadratic form. Conversely, if $\Ga$ preserves a Lorentzian quadratic form, then, by Lemma \ref{lorentzconelemma}, Example \ref{hyperboloidexample}, and Lemma \ref{uniquenessproposition}, as in the proof of \ref{st2} of Theorem \ref{surfacetheorem}, $h$ is hyperbolic, and $\vol(M, h) = \vol_{\mathrm{hyp}}(M_{0})$ by the Mostow rigidity theorem.
\end{proof}

By the equality case of Theorem \ref{hyperbolicvolumetheorem}, $\mathfrak{m}$ is not constant on a family of convex projective structures on $M_{0}$ that contains its hyperbolic structure and a structure not projectively equivalent to it. The proof of the inequality \eqref{hyperbolicvolumebound} applies without change if, instead of a homeomorphism, there is a continuous map from $M$ to $M_{0}$ of nonzero degree $d$, both manifolds being oriented, and yields $\cV(M) \geq |d|\vol_{\mathrm{hyp}}(M_{0})$.

Theorem \ref{hyperbolicvolumetheorem} gives no lower bound for manifolds carrying a properly convex projective structure but no hyperbolic metric, such as those of Kapovich \cite{Kapovich-GT} and the closed orientable properly convex projective $4$-manifolds of every positive even Euler characteristic constructed by Riolo, Seppi, and Slavich \cite[Theorem $1.3$]{Riolo-Seppi-Slavich}. In dimension $4$, Tamburelli \cite{Tamburelli-euler} proves $\chi(M) \geq 0$ for closed properly convex flat projective $M$ by transferring an identity on the cone over the affine sphere $\Sigma_{0}$ to a Bochner identity on $M$, which suggests looking for a lower bound for $\mathfrak{m}(M)$ in terms of $\chi(M)$ in dimension $4$.

\section{Hyperbolic affine spheres as projective K\"ahler--Einstein metrics}\label{kahlereinsteinsection}
The affine sphere equation is the Einstein equation in a real analogue of K\"ahler geometry in which a flat projective structure replaces the complex structure and a Bernstein--Gelfand--Gelfand complex replaces the Dolbeault complex. In this analogy the hyperbolic affine spheres correspond to the K\"ahler--Einstein metrics with negative Einstein constant, and Theorems A, B, and C are the counterparts of the variational theory of the corresponding energy functional. The correspondence between the functionals is formal, at the level of normalized first variations; the correspondence between the equations becomes exact on the tube domain over the cone (Lemma \ref{tubelemma}).

This point of view has its origin in the work of Cheng and Yau \cite{Cheng-Yau-real} and becomes explicit in the work of Loftin \cite{Loftin-affinespheres, Loftin-KE, Loftin-met}, in which a section $u$ of $\mathcal{E}$ is regarded as a potential analogous to the local potential of a K\"ahler metric \cite{Loftin-met}. It is described in this section.

In this section $M = \Om/\Ga$ is a closed properly convex flat projective manifold, as in the rest of the paper, $\en$ is its flat projective structure, and $\mathcal{E} = |\Det TM|^{1/N}$. Sections of the line bundle $\mathcal{E}$ correspond bijectively to $\Ga$-invariant functions on $C$ homogeneous of degree $1$, the section $u$ corresponding to the function $\tilde{u}$, and for $\psi \in \cAs(C, \Ga)$ the section $u$ corresponding to $\tilde{u} = -\psi$ is negative.
The restriction of $\fD^{2}\tilde{u}$ to a subbundle of $TC$ complementary to the span of $\eul$ does not depend on the choice of the complement, because $\fD^{2}\tilde{u}(\eul, \cdot) = 0$ by homogeneity. It is the invariant lift of the projective Hessian $\Bp(u)_{ij} = \nabla_{i}\nabla_{j}u - P_{ij}u$, a section of $S^{2}(T^{\ast}M)\otimes\mathcal{E}$, in the sense explained in \cite{Fox-conelike}. This follows by a routine computation using the identification of the flat affine connection on $C$ with the Thomas connection of the projective structure of $\Om$ \cite{BEG, Thomas-invariants}, which in the tractor formalism \cite{BEG, Cap-Gover-tractor} amounts to the identification of $D_{A}D_{B}u$ with the Hessian of $\tilde{u}$; see also \cite{Eastwood-notes, Gover-projective}. The projectively invariant Monge--Amp\`ere operator $\monge(u) = u^{n+2}\det\Bp(u)$ satisfies
\begin{align}\label{mumonge}
\mu(\psi) = |\monge(u)| = (-1)^{n}\monge(u).
\end{align}
Here $\det\Bp(u)$ is a section of $\mathcal{E}^{-n-2}$, a bundle of densities, so its sign is defined, and it is positive for $\psi \in \cAs(C, \Ga)$ because $\Bp(u)$ is the horizontal restriction of $-\fD^{2}\psi$, which is positive definite. Hence the sign of $\monge(u)$ is that of $u^{n+2}$, which is $(-1)^{n}$ because $u < 0$. Moreover,
\begin{align}\label{sbetaprojective}
& \cS_{\be}(\psi) = \int_{M}|u|^{\al}|\det\Bp(u)|^{(\al+n+1)/(n+2)}, && \be = \tfrac{\al + n+1}{n+2},
\end{align}
the integrand being a section of $\mathcal{E}^{\al - (n+2)\be} = \mathcal{E}^{-n-1}$, so a volume density. In particular, $\cS_{0}$, $\cS_{1}$, and $\cS_{1/2}$ correspond to $\al = -n-1$, $\al = 1$, and $\al = -n/2$. Reformulated in this way, these expressions make sense for any projective structure, flat or not.

A \emph{projective metric} is a nondegenerate section $h_{ij}$ of $S^{2}(T^{\ast}M)\otimes\mathcal{E}$. The operator $h_{ij} \mapsto \nabla_{[i}h_{j]k}$ on sections of $S^{2}(T^{\ast}M)\otimes\mathcal{E}$ does not depend on the choice of $\nabla \in \en$, and a projective metric $h_{ij}$ is \emph{projective K\"ahler} if $\nabla_{[i}h_{j]k} = 0$ for some, equivalently every, $\nabla \in \en$. Two projective K\"ahler metrics are in the same \emph{projective K\"ahler class} if their difference is the projective Hessian of a section of $\mathcal{E}$, and the \emph{trivial class} consists of the projective K\"ahler metrics that are projective Hessians of sections of $\mathcal{E}$. The projective Hessian is the first operator of the Bernstein--Gelfand--Gelfand sequence associated with the representation of $SL(\ste)$ on $\std$, the second being $h \mapsto \nabla_{[i}h_{j]k}$, and, because $\en$ is flat, this sequence is a complex, which is a resolution of the sheaf of local solutions of $\Bp(u) = 0$, that is, of the locally constant sheaf of linear functions on $\ste$ \cite{Cap-Slovak-Soucek, Eastwood-notes}. Consequently every nondegenerate projective Hessian $h = \Bp(u)$ is projective K\"ahler, and the projective K\"ahler classes determine elements of the first cohomology $H^{1}(\Ga; \std)$ of $\Ga$ with coefficients in the representation of $\Ga$ on $\std$.

\begin{lemma}\label{projkernellemma}
If a section $u$ of $\mathcal{E}$ over $M$ satisfies $\Bp(u) = 0$, then $u = 0$.
\end{lemma}
\begin{proof}
The section $u$ corresponds to a $\Ga$-invariant function $\tilde{u}$ on $C$, homogeneous of degree $1$. The condition $\Bp(u) = 0$ means that $\fD^{2}\tilde{u}$ vanishes on horizontal vectors. Since $\tilde{u}$ is homogeneous of degree $1$, also $\fD^{2}\tilde{u}(\eul, \cdot) = 0$. Hence $\fD^{2}\tilde{u} = 0$ on the connected open set $C$, so $\tilde{u}$ is the restriction to $C$ of an affine function, which is linear because $\tilde{u}$ is homogeneous of degree $1$. This linear function is $\Ga$-invariant, because it is determined by its restriction to the open set $C$, which $\Ga$ preserves. Hence $\tilde{u}/\rho_{C}$ is a continuous function on the compact manifold $M$, so it is bounded, and $\tilde{u} = 0$ by the last claim of Lemma \ref{vinberglemma}.
\end{proof}
In terms of the complex described before Lemma \ref{projkernellemma}, the lemma says that the cohomology in degree $0$ vanishes.

Since $\det h$ is a section of $\mathcal{E}^{-n-2}$, the section $|\det h|^{-1/(n+2)}$ of $\mathcal{E}$ is defined, and it is positive. The \emph{projective Ricci form} of a projective metric $h$ is defined to be $\rho(h) = \Bp(|\det h|^{-1/(n+2)})$. This is the analogue of the Ricci form of a K\"ahler metric, which is $-i\pr\bar{\pr}$ of the logarithm of its volume density. Here the logarithm is replaced by the power of the volume density that is a section of $\mathcal{E}$, on which $\Bp$ is defined. A projective metric $h$ is \emph{projective K\"ahler--Einstein} if $\rho(h) = c\,h$ for a constant $c$, the \emph{Einstein constant} of $h$. Because $\rho(h)$ is by definition a projective Hessian, its class is trivial. In this sense the analogue of the first Chern class vanishes as a class, and what remains of it is the sign of the Einstein constant. Moreover, if $\rho(h) = c\,h$ with $c \neq 0$, then $h = \Bp(c^{-1}|\det h|^{-1/(n+2)})$, so $h$ lies in the trivial class and is the projective Hessian of a section of $\mathcal{E}$ having the sign of $c$, while $c = 0$ is impossible, by Lemma \ref{projkernellemma}, because $|\det h|^{-1/(n+2)}$ vanishes nowhere.

\begin{lemma}\label{kelemma}
Let $\psi \in \cAs(C, \Ga)$ and let $u$ be the corresponding section of $\mathcal{E}$. Then $h = \Bp(u)$ is a positive definite projective metric, and $h$ is projective K\"ahler--Einstein if and only if $\mu(\psi)$ is constant, in which case $\rho(h) = -\mu(\psi)^{-1/(n+2)}h$.
\end{lemma}

\begin{proof}
That $h$ is positive definite was observed after \eqref{mumonge}. By \eqref{mumonge}, $\mu(\psi) = |u|^{n+2}|\det\Bp(u)|$, so $|\det h|^{-1/(n+2)} = \mu(\psi)^{-1/(n+2)}|u| = -\mu(\psi)^{-1/(n+2)}u$, because $u < 0$. Hence $\rho(h) = -\Bp(\mu(\psi)^{-1/(n+2)}u)$. If $\mu(\psi) = \ka$ is constant, then $\rho(h) = -\ka^{-1/(n+2)}h$. Conversely, if $\rho(h) = ch$, then $\Bp\big((\mu(\psi)^{-1/(n+2)} + c)u\big) = 0$. By Lemma \ref{projkernellemma}, $(\mu(\psi)^{-1/(n+2)} + c)u = 0$, and since $u$ vanishes nowhere, $\mu(\psi) = (-c)^{-n-2}$ is constant.
\end{proof}

A positive definite projective Hessian $\Bp(u)$ with $u < 0$ corresponds to a smooth positive $\Ga$-invariant function $\psi = -\tilde{u}$ on $C$, homogeneous of degree $1$, such that $\fD^{2}\psi$ is negative definite on horizontal vectors, and such a $\psi$ lies in $\cAs(C, \Ga)$ by the converse in Lemma \ref{polarmaplemma}\ref{pml1}. Hence, by Lemma \ref{kelemma}, Theorem B and Corollary \ref{existencecorollary} say that $M$ carries a positive definite projective K\"ahler--Einstein metric with negative Einstein constant, that such a metric is unique up to multiplication by a positive constant, and that it equals $\Bp(u_{0})$, where $u_{0}$ is the section of $\mathcal{E}$ such that $-\tilde{u}_{0}$ is a minimizer $\psi_{0}$ of the covolume product $\cP$ on $\cA(C, \Ga)$, so that the Riemannian metric $|u_{0}|^{-1}\Bp(u_{0})$ is the centroaffine metric $\gb_{\psi_{0}}$. Theorem C computes the second variation of $\cP$ at the cogauge corresponding to such a metric.

Conversely, by a theorem of Loftin \cite[Theorem $1$]{Loftin-met}, a closed manifold with a flat projective structure that carries a negative section of $\mathcal{E}$ with positive definite projective Hessian is properly convex. For sections satisfying the affine sphere equation this was obtained earlier in \cite{Loftin-affinespheres} from the theorems of Cheng and Yau \cite{Cheng-Yau-affinespheres} and Li \cite{Li-calabi2} on complete affine spheres, and \cite[Theorem $2$]{Loftin-KE} treats manifolds, not necessarily compact, on which $|u|^{-1}\Bp(u)$ is a complete Riemannian metric. Consequently the closed manifolds with a flat projective structure that carry a positive definite projective K\"ahler--Einstein metric with negative Einstein constant are exactly the properly convex ones.

For a Riemannian Hessian metric $g = \fD^{2}F$ on an open subset of $\ste$, the \emph{second Koszul form} is $\be_{\mathrm{K}} = \tfrac{1}{2}\fD^{2}\log\det g$, the determinant being computed with respect to $\Psi$, and $g$ is \emph{Einstein--Hessian} if $\be_{\mathrm{K}} = \la g$ for a function $\la$ \cite[Definition $3.2$ and Proposition $3.4$]{Shima}. In the terminology of \cite[Section $1$]{Cheng-Yau-real}, in which the Ricci tensor of the Hessian metric $g$ is $-\fD^{2}\log\det g = -2\be_{\mathrm{K}}$, which is not the Ricci tensor of the Riemannian metric $g$, a Hessian metric whose second Koszul form equals $g$ is an Einstein K\"ahler affine metric with Einstein constant $-2$.

\begin{lemma}\label{tubelemma}
Let $C \subset \ste$ be an open properly convex cone, let $\psi \in \cAs(C, \{1\})$, and let $g = \fD^{2}F$ for $F = -N\log\psi$. The following are equivalent.
\begin{enumerate}[label=$(\arabic*)$, ref=(\arabic*)]
\item\label{tl1} $\mu(\psi)$ is constant, that is, the level sets of $\psi$ are hyperbolic affine spheres centered at the origin.
\item\label{tl2} The Riemannian Hessian metric $g$ on $C$ is Einstein--Hessian. In this case its second Koszul form equals $g$.
\item\label{tl3} The K\"ahler metric on the tube domain $C + \sqrt{-1}\,\ste$ over $C$, identified with the tangent bundle of $C$, that is induced by $g$ in the sense of \cite[Theorem $3.1$]{Shima}, and for which $F$ is a potential, is K\"ahler--Einstein. In this case its Einstein constant is negative.
\item\label{tl4} The pseudo-K\"ahler metric on the tube domain $C + \sqrt{-1}\,\ste$ for which $\Phi = \tfrac{1}{2}\psi^{2}$ is a potential, which has complex signature $(1, n)$ by Lemma \ref{polarmaplemma}\ref{pml1}, is Ricci flat.
\end{enumerate}
\end{lemma}
\begin{proof}
By \eqref{mulogdet}, $g = -N\fD^{2}\log\psi$ is positive definite and $\det g = N^{N}\psi^{-2N}\mu(\psi)$, so $\be_{\mathrm{K}} = g + \tfrac{1}{2}\fD^{2}\log\mu(\psi)$. If $\mu(\psi)$ is constant, then $\be_{\mathrm{K}} = g$. Conversely, suppose that $\be_{\mathrm{K}} = \la g$ for a function $\la$ on $C$, so that $\tfrac{1}{2}\fD^{2}\log\mu(\psi) = (\la - 1)g$. Since $\log\mu(\psi)$ is homogeneous of degree $0$ and $\eul\log\psi = 1$, there hold $\fD^{2}\log\mu(\psi)(\eul, \cdot) = -d\log\mu(\psi)$ and $g(\eul, \cdot) = N\,d\log\psi$, so $-\tfrac{1}{2}d\log\mu(\psi) = (\la - 1)N\,d\log\psi$. Evaluating on $\eul$ gives $\la = 1$, and then $d\log\mu(\psi) = 0$, so $\mu(\psi)$ is constant on the connected set $C$. This shows the equivalence of \ref{tl1} and \ref{tl2}. Compare \cite[Theorem $3$]{Loftin-KE}. The equivalence of \ref{tl2} and \ref{tl3} is \cite[Theorem $3.1$]{Shima}, and the K\"ahler--Einstein metric has negative Ricci curvature; see \cite[Theorem $3$ and the remark following Theorem $2$]{Loftin-KE} and \cite[Theorem $4.1$]{Fox-schwarz}.

Finally, the Ricci form of the pseudo-K\"ahler metric of \ref{tl4} is $-\sqrt{-1}\,\pr\bar{\pr}\log|\det\fD^{2}\Phi| = -\sqrt{-1}\,\pr\bar{\pr}\log\mu(\psi)$, which, because $\mu(\psi)$ depends only on the real part of the coordinate in the tube domain, vanishes if and only if $\fD^{2}\log\mu(\psi) = 0$, that is, if and only if $\log\mu(\psi)$ is the restriction to $C$ of an affine function. Since $\log\mu(\psi)$ is homogeneous of degree $0$, this is the case if and only if $\mu(\psi)$ is constant, so \ref{tl4} is equivalent to \ref{tl1}.
\end{proof}

\begin{remark}\label{barrierremark}
Let $\psi_{0} \in \cAs(C, \Ga)$ satisfy $\mu(\psi_{0}) = 1$ and let $F_{0} = -N\log\psi_{0} + \tfrac{N}{2}\log N$. By \eqref{mulogdet}, $\det\fD^{2}F_{0} = N^{N}\psi_{0}^{-2N} = e^{2F_{0}}$. The metric $\fD^{2}F_{0}$ is invariant under $\Ga$ and under the dilations of $C$, which together act cocompactly on $C$, so it is complete, and $F_{0}(x) \to \infty$ as $x \to \pr C$ by Lemma \ref{harnacklemma}. Hence $F_{0}$ is the canonical potential of $C$ in the sense of \cite{Fox-schwarz}. The existence of a convex solution of $\det\fD^{2}F = e^{2F}$ that tends to $+\infty$ at the boundary and has complete Hessian metric, on any convex domain containing no complete line, is \cite[Theorem $4.4$]{Cheng-Yau-real}, where it answers a question of Calabi. By a theorem proved independently in \cite[Theorem $1.4$]{Fox-schwarz} and by Hildebrand in \cite{Hildebrand-canonical}, $F_{0}$ is a self-concordant barrier for $C$ with parameter $N$; see \cite{Hildebrand-threeviews} for the relations between such barriers, centroaffine hypersurfaces, and Lagrangian submanifolds. Thus Theorems A and B characterize the canonical potential of a divisible cone variationally: up to an additive constant it is $-N\log\psi$ for the cogauge $\psi$ of a $\Ga$-convex body of least covolume product.
\end{remark}

\begin{remark}\label{compactquotientremark}
The theorem of Cheng and Yau on compact affine manifolds with negative first Chern class \cite[Section $3$, Theorem $1$]{Cheng-Yau-real} gives another proof of the existence of the $\Ga$-invariant affine sphere, with determinants computed with respect to $\Psi$.

Fix $q > 1$ and $\psi_{1} \in \cAs(C, \Ga)$. The map $x \mapsto ([x], \log\psi_{1}(x))$ is a diffeomorphism from $C$ onto $\Om\times\rea$ that intertwines the action of $\Ga$ with its action on the first factor and the dilation by $q$ with the translation by $\log q$ in the second factor, so the group generated by $\Ga$ and $qI$ acts freely and properly discontinuously on $C$ with compact quotient $Q$, diffeomorphic to $M \times S^{1}$, an affine manifold. The density $\psi_{1}^{-N}|\Psi|$ is invariant under $\Ga$ and under all dilations, so it descends to $Q$, and $\fD^{2}\log\psi_{1}^{-N} = -N\fD^{2}\log\psi_{1}$ is positive definite, because $\fD^{2}\psi_{1}$ is negative definite on $\ker d\psi_{1}$, vanishes on $\eul$, and $d\psi_{1}(\eul) = \psi_{1}$.

Hence $Q$ is a compact K\"ahler affine manifold with negative first Chern class in the sense of \cite[Sections $1$ and $3$]{Cheng-Yau-real}, and by \cite[Section $3$, Theorem $1$]{Cheng-Yau-real} there is a density on $Q$ whose lift to $C$ is $e^{F}|\Psi|$ with $\fD^{2}F$ positive definite and $\det\fD^{2}F = e^{2F}$. Such an $F$ is unique: if $F'$ is another, then $F - F'$ descends to $Q$, and at a point where it attains its maximum, $\fD^{2}F \leq \fD^{2}F'$, so $e^{2F} \leq e^{2F'}$ there, so $F \leq F'$ everywhere, and, by symmetry, $F = F'$. Each dilation $x \mapsto tx$, $t > 0$, commutes with $\Ga$ and with $qI$, and $F(tx) + N\log t$ solves the same equation, so, by uniqueness, $F(tx) = F(x) - N\log t$.

Let $\psi_{0} = N^{1/2}e^{-F/N}$. It is positive, $\Ga$-invariant, and homogeneous of degree $1$, and $-N\fD^{2}\log\psi_{0} = \fD^{2}F$ is positive definite. Since $\fD^{2}\psi_{0} = \psi_{0}(\fD^{2}\log\psi_{0} + (d\log\psi_{0})^{2})$, $\fD^{2}\psi_{0}$ is negative definite on $\ker d\psi_{0}$, and it vanishes on $\eul$ by homogeneity, so $\psi_{0} \in \cAs(C, \Ga)$. By \eqref{mulogdet}, $\mu(\psi_{0}) = \psi_{0}^{2N}N^{-N}\det\fD^{2}F = \psi_{0}^{2N}N^{-N}e^{2F} = 1$, so, by Lemma \ref{affinespherelemma}, the level sets of $\psi_{0}$ are $\Ga$-invariant hyperbolic affine spheres, which are asymptotic to $\pr C$ by Lemmas \ref{harnacklemma} and \ref{asymptoticlemma}.

This combines two observations of Loftin: that the level sets of the volume of the complete Hessian metric of Cheng and Yau on a cone are hyperbolic affine spheres \cite[Theorem $3$]{Loftin-KE}, and that the results of \cite{Cheng-Yau-real} apply to the quotient of the cone by a discrete group of homotheties \cite[remark following Theorem $2$]{Loftin-KE}.
\end{remark}

In the coordinates $x + \sqrt{-1}\,y$ of the tube domain, the K\"ahler form of the pseudo-K\"ahler metric of Lemma \ref{tubelemma}\ref{tl4} is a constant multiple of $\sum_{j, k = 1}^{N}\Phi_{jk}\,dx^{j}\wedge dy^{k}$, where $\Phi_{jk}$ are the components of $\fD^{2}\Phi$ in linear coordinates $x^{j}$ on $\ste$ and $y^{j}$ are the corresponding coordinates of the imaginary part, the translations in the imaginary directions are Hamiltonian with moment map a constant multiple of $G_{\psi}$, the symplectic volume density is a constant multiple of $\mu(\psi)\,dx\,dy$, and the volume density of the parallel holomorphic volume form is a constant multiple of $dx\,dy$. Hence, when $\Ga$ divides $C$, by Lemma \ref{volumeproposition}, per unit volume in the imaginary directions and up to constant factors, the covolumes $|K_{\psi}/\Ga|$ and $|K^{\ast}_{\psi}/\Ga^{\ast}|$ are the holomorphic and the symplectic volumes of the part of the tube domain over $K_{\psi}$ modulo $\Ga$, and Theorem A says that their product is least exactly when the metric of Lemma \ref{tubelemma}\ref{tl4} is Ricci flat.

By Lemma \ref{tubelemma}, the K\"ahler potential on the tube domain corresponding to a cogauge $\psi$ is $\phi = -N\log\psi$, and with this correspondence the variational results of this paper correspond formally to the variational description of the K\"ahler--Einstein metrics on a compact K\"ahler manifold $(X, \omega)$ of complex dimension $n$ with negative first Chern class, in terms of K\"ahler potentials $\phi$ such that $\omega_{\phi} = \omega + i\pr\bar{\pr}\phi > 0$. This correspondence, unlike Lemma \ref{tubelemma}, is formal, and it is not used elsewhere in the paper.

Rescaling $\psi$ by a positive constant corresponds to adding a constant to $\phi$. The functional $\cS_{0}(\psi) = \int_{M}\psi^{-N}$ corresponds to $\int_{X}e^{\phi}\omega^{n}$, and the density $\mu(\psi)\psi^{-N}$, whose integral is $\cS_{1}(\psi)$, corresponds to $\omega_{\phi}^{n}$. By \eqref{firstvariation}, the first variation of $\cS_{1}$ in the direction $\vr\psi$ is $N\int_{M}(\vr\psi/\psi)\,\mu(\psi)\psi^{-N}$, which corresponds to $-\int_{X}\vr\phi\,\omega_{\phi}^{n}$, the negative of the first variation of the Aubin--Yau energy $E$. The correspondence is between normalized measures. With $\vr\phi = -N\vr\psi/\psi$ and $V = \int_{X}\omega^{n}$,
\begin{align}\label{normalizedcorrespondence}
&\vr\log\cS_{1}(\psi) = -\int_{M}\vr\phi\,\frac{\mu(\psi)\psi^{-N}}{\cS_{1}(\psi)},& &\vr\left(-\frac{E(\phi)}{V}\right) = -\int_{X}\vr\phi\,\frac{\omega_{\phi}^{n}}{V},&
\end{align}
so the probability measure $\mu(\psi)\psi^{-N}/\cS_{1}(\psi)$ corresponds to $\omega_{\phi}^{n}/V$. The total mass $\cS_{1}(\psi)$ varies with $\psi$, whereas $\int_{X}\omega_{\phi}^{n} = V$ does not depend on $\phi$, and it is $\cS_{1}$ itself, not a normalizing constant, that plays the role of $e^{-E/V}$.

Thus $\log\cP$ corresponds to the functional $\phi \mapsto -E(\phi)/\int_{X}\omega^{n} + \log\int_{X}e^{\phi}\omega^{n}$, which is invariant under the addition of constants to $\phi$ and whose critical points are the solutions of $\omega_{\phi}^{n} = c\,e^{\phi}\omega^{n}$ for a constant $c > 0$. When $\omega$ represents the negative of the first Chern class and $\omega^{n}$ is replaced by its multiple by the exponential of the Ricci potential of $\omega$, this is the K\"ahler--Einstein equation of the theorem of Aubin and Yau, with negative Einstein constant. This functional is convex along segments $t \mapsto (1-t)\phi_{0} + t\phi_{1}$, $t \in [0, 1]$, because $E$ is concave along them and by H\"older's inequality, as the covolume and the polar covolume are convex along segments in $\cAs(C, \Ga)$ (Lemma \ref{covolumeconvexitylemma} and Section \ref{reversesantalosection}), and its minimizers are the K\"ahler--Einstein potentials, as the minimizers of $\cP$ are the cogauges of the bodies bounded by hyperbolic affine spheres (Theorem A).

The linearized operator $-\Delta + 2N$ of the affine sphere equation, which, for $\psi_{0}$ normalized by $\mu(\psi_{0}) = 1$, is by \eqref{deltamu} and Lemma \ref{laplacianlemma} the derivative of $\mu((1+tf)\psi_{0})$ at $t = 0$, and which by Theorem C determines the second variation of $\cP$, corresponds to the linearization $-\Delta_{\omega_{\phi}} + 1$ of the K\"ahler--Einstein equation with negative Einstein constant. Both are positive, and in both cases the positivity accounts for the uniqueness of the solution and the absence of obstructions.

The elliptic counterpart, in which convex bodies and elliptic affine spheres replace $\Ga$-convex bodies and hyperbolic affine spheres, corresponds to the case of positive first Chern class. There the Ding functional \cite{Ding-existence, Ding-Tian} of a toric Fano manifold is a functional of convex functions on $\rea^{n}$, and its boundedness follows from the functional form \cite{AKM-santalo} of the Blaschke--Santal\'o inequality \cite{Berman-Berndtsson}; see also \cite{Cordero-Klartag}, to which the variational problem of Klartag in Remark \ref{klartagremark} is related.

Several aspects of this point of view are not pursued here. The Einstein constant of a positive definite projective K\"ahler--Einstein metric on a closed manifold with a flat projective structure that is not properly convex can be positive or zero, and these cases correspond to the elliptic and the improper affine spheres. The projective K\"ahler classes other than the trivial class are not considered. The vanishing of the class of the projective Ricci form reflects the unimodularity of the holonomy: on the compact quotient of the cone by a group of homotheties used in the remark following \cite[Theorem $2$]{Loftin-KE}, the Hessian metric $-N\fD^{2}\log\psi$ has a first Koszul form whose cohomology class is determined by the affine structure and is not zero, and its second Koszul form represents a multiple of the first affine Chern class \cite[Section $7.6$]{Shima}.

Suppose that $n \geq 2$. The equation $\monge(u) = \ka$ makes sense for any projective structure, flat or not, and its solutions persist under small deformations of a properly convex flat projective structure. For a projective structure $\en$ that need not be flat, there is always a representative $\nabla \in \en$ that preserves a volume density. For such a representative the Ricci tensor, and so the projective Schouten tensor $P_{ij}$, is symmetric, and the projective Hessian $\Bp(u)_{ij} = \nabla_{i}\nabla_{j}u - P_{ij}u$ of a section $u$ of $\mathcal{E}$, computed with any such representative, is a section of $S^{2}(T^{\ast}M)\otimes\mathcal{E}$ that does not depend on the choice of the representative \cite{BEG, Eastwood-notes}. The operator $\monge(u) = u^{n+2}\det\Bp(u)$, projective metrics, and the projective Ricci form $\rho(h) = \Bp(|\det h|^{-1/(n+2)})$ are defined as in the flat case. A one-parameter family $\en_{\ep}$ of projective structures is \emph{smooth} if there are representatives $\nabla^{\ep} \in \en_{\ep}$ depending smoothly on $\ep$. These can be chosen to preserve a fixed volume density, and then the coefficients of the operator $\Bp_{\ep}$ depend smoothly on $\ep$.

\begin{lemma}\label{nonflatlemma}
Let $n \geq 2$, let $M = \Om/\Ga$ be closed and properly convex, with flat projective structure $\en_{0}$, let $\en_{\ep}$, $\ep \in (-1, 1)$, be a smooth one-parameter family of projective structures on $M$, and let $\Bp_{\ep}$, $\monge_{\ep}$, and $\rho_{\ep}$ be the projective Hessian, the operator $\monge$, and the projective Ricci form of $\en_{\ep}$. Let $u_{0}$ be the section of $\mathcal{E}$ corresponding to $-\psi_{0}$ for $\psi_{0} \in \cAs(C, \Ga)$ with $\mu(\psi_{0})$ constant, and let $\ka = \monge_{0}(u_{0})$. Let $\si \in (0, 1)$. There are $\ep_{0} > 0$ and a neighborhood $\mathcal{U}$ of $u_{0}$ in the space $C^{2,\si}(M; \mathcal{E})$ of sections of $\mathcal{E}$ of class $C^{2,\si}$ such that for $|\ep| < \ep_{0}$ the equation $\monge_{\ep}(u) = \ka$ has a unique solution $u_{\ep} \in \mathcal{U}$. The section $u_{\ep}$ is smooth, negative, and depends smoothly on $\ep$, and $h_{\ep} = \Bp_{\ep}(u_{\ep})$ is a positive definite projective metric satisfying $\rho_{\ep}(h_{\ep}) = -|\ka|^{-1/(n+2)}h_{\ep}$.
\end{lemma}
\begin{proof}
The map $(\ep, u) \mapsto \monge_{\ep}(u) - \ka$ is smooth from a neighborhood of $(0, u_{0})$ in $(-1, 1)\times C^{2,\si}(M; \mathcal{E})$ to the space of functions on $M$ of class $C^{0,\si}$, because $\monge_{\ep}(u)$ is a polynomial in $u$ and its first and second derivatives with coefficients depending smoothly on $\ep$. By \eqref{mumonge}, \eqref{deltamu}, and Lemma \ref{laplacianlemma}, its derivative in $u$ at $(0, u_{0})$, in the direction $fu_{0}$ for a function $f$ on $M$ of class $C^{2,\si}$, is $\ka(2Nf - \Delta f)$, where $\Delta$ is the Laplacian of $\gb_{\psi_{0}}$. Since $-\Delta + 2N$ is a positive elliptic operator on the closed manifold $M$, this derivative is an isomorphism from $C^{2,\si}$ to $C^{0,\si}$, and the implicit function theorem gives $\ep_{0}$, $\mathcal{U}$, and the solutions $u_{\ep}$, depending smoothly on $\ep$. After shrinking $\ep_{0}$ and $\mathcal{U}$, $u_{\ep} < 0$ and $\Bp_{\ep}(u_{\ep})$ is positive definite, because these conditions hold for $u_{0}$ and are open in $C^{2}$, so the equation is elliptic at $u_{\ep}$, and $u_{\ep}$ is smooth by elliptic regularity. Since $|\det h_{\ep}| = |\ka||u_{\ep}|^{-n-2}$, $|\det h_{\ep}|^{-1/(n+2)} = -|\ka|^{-1/(n+2)}u_{\ep}$, which gives the last claim.
\end{proof}
In the affine formulation of \cite[Section $3$]{Cheng-Yau-real}, the linearized operator $g^{ij}\partial_{i}\partial_{j}$ of a Hessian metric $g$ is not, in general, formally self-adjoint with respect to the Riemannian volume of $g$, because $g^{ij}\partial_{i}\partial_{j}v = \Delta_{g}v + \tfrac{1}{2}\lb d\log\det g, dv\ra_{g}$, and Cheng and Yau pass to the tube domain to apply the continuity method. In the projective formulation the linearization at an affine sphere is $\ka(-\Delta + 2N)$, which is self-adjoint, so the implicit function theorem applies directly on $M$.

In the setting of Lemma \ref{nonflatlemma}, let $\nabla \in \en_{\ep}$ be the unique representative such that $\nabla u_{\ep} = 0$. Its Ricci tensor $R_{ij}$ is symmetric, its projective Schouten tensor is $P_{ij} = \tfrac{1}{1-n}R_{ij}$, and $h_{\ep} = \Bp_{\ep}(u_{\ep}) = -P_{ij}u_{\ep}$. Hence the Riemannian metric $g_{\ep} = |u_{\ep}|^{-1}h_{\ep}$ equals $P_{ij}$, the Ricci tensor of $\nabla$ is $(1-n)g_{\ep}$, and $\det g_{\ep}$ is $\nabla$-parallel, because $\det h_{\ep} = \ka u_{\ep}^{-n-2}$. Since the operator $h \mapsto \nabla_{[i}h_{j]k}$ does not depend on the choice of the representative of $\en_{\ep}$, the projective metric $h_{\ep}$ is projective K\"ahler if and only if
\begin{align}\label{statisticalcondition}
\nabla_{[i}P_{j]k} = 0,
\end{align}
that is, if and only if the projective Cotton tensor of $\nabla$ vanishes. For a projective structure that is not flat this need not be the case. When \eqref{statisticalcondition} holds, the pair $(\nabla, g_{\ep})$ is a \emph{special statistical structure}, meaning that $\nabla_{[i}(g_{\ep})_{j]k} = 0$ and $\nabla\det g_{\ep} = 0$, and its Ricci tensor is a constant multiple of $g_{\ep}$. By definition \cite[Section $4.2$]{Fox-coupled}, \cite[Definition $9.18$]{Fox-conelike}, and, for surfaces, \cite[Definition $6.2$]{Fox-pisa}, a special statistical structure is \emph{Einstein} if the Ricci tensors of $\nabla$ and of the conjugate connection $\bar{\nabla} = \nabla + g_{\ep}^{kp}\nabla_{i}(g_{\ep})_{jp}$ are multiples of the metric and the scalar curvature is constant, and the special statistical structure induced on a hypersurface in flat affine space by its equiaffine metric and equiaffine connection is Einstein exactly when the hypersurface is an affine sphere. Here the scalar curvature is the constant $-n(n-1)$, and, since the curvature of $\nabla$ is $R_{ijk}\,^{l} = B_{ijk}\,^{l} - 2\delta_{[i}\,^{l}(g_{\ep})_{j]k}$, where $B_{ijk}\,^{l}$ is the projective Weyl tensor of $\en_{\ep}$, the Ricci tensor of $\bar{\nabla}$, which is the contraction of the curvature of $\nabla$ with $g_{\ep}$ on its second and third indices \cite[Lemma $8.15$]{Fox-conelike}, is a multiple of $g_{\ep}$ if and only if
\begin{align}\label{einsteincondition}
g_{\ep}^{jk}B_{ijk}\,^{l} = 0.
\end{align}
Thus $(\nabla, g_{\ep})$ is an Einstein special statistical structure if and only if there hold \eqref{statisticalcondition} and \eqref{einsteincondition}. That these conditions, together with the parallelism of $\det g_{\ep}$, suffice is \cite[Example $9.22$]{Fox-conelike}, where, owing to a misprint, the curvature tensor of $\nabla$ appears in place of the projective Weyl tensor in the last condition; as printed, that condition contradicts the normalization of the metric as $\tfrac{1}{1-n}$ times the Ricci tensor made there, since the trace of $g_{\ep}^{jk}R_{ijk}\,^{l}$ is $g_{\ep}^{jk}R_{jk} = -n(n-1)$.

Both are automatic when $\en_{\ep}$ is flat, and, when $n = 2$, \eqref{einsteincondition} is automatic, because the projective Weyl tensor of a surface vanishes. For $\ep = 0$ the result is the Einstein special statistical structure induced on $M$ by the equiaffine metric and the equiaffine connection of the affine sphere. It is in this sense that the equation $\rho_{\ep}(h_{\ep}) = -|\ka|^{-1/(n+2)}h_{\ep}$ of Lemma \ref{nonflatlemma} is an Einstein equation: it is the Einstein equation for special statistical structures when \eqref{statisticalcondition} and \eqref{einsteincondition} hold, and without them it is only the equation $\monge_{\ep}(u_{\ep}) = \ka$.

Whether $\monge(u) = \ka$ can be solved for projective structures far from the properly convex flat ones appears to be an open question.

\medskip
\noindent\textbf{Acknowledgment of the use of AI tools.}
In preparing this paper I used Claude (Anthropic; Opus 5.5 and Fable 5.2) and ChatGPT (OpenAI; GPT-6) as tools. Both were used for literature searches and comparisons, for error checking, and for auditing drafts of the manuscript. Claude was also used to propose proof strategies, to draft and revise text, including statements and proofs of some lemmas and remarks, to check computations, to test arguments for gaps, and to supply or correct intermediate steps. Some of the results first arose as suggestions made by these tools. 

\medskip

\noindent\textbf{Funding Acknowledgment.} This work was supported in part by Project PID2024-156578NB-I00 funded by MICIU/AEI/10.13039/501100011033/FEDER, EU.

\end{document}